\documentclass{article}
\newtheorem{theorem}{Theorem}
\usepackage{graphicx}
\DeclareSymbolFont{AMSa}{U}{msa}{m}{n}
\DeclareMathDelimiter\ulcorner{\mathopen} {AMSa}{"70}{AMSa}{"70}
\DeclareMathDelimiter\urcorner{\mathclose}{AMSa}{"71}{AMSa}{"71}
\makeatletter
\def\uufill{$\m@th\mathopen\ulcorner\mkern-7mu%
  \cleaders\hbox{\rule[6pt]{1dd}{1dd}}\hfill
  \mkern-7mu\mathclose\urcorner$}
\def\overbrack#1{\vbox{\m@th\ialign{##\crcr
      \uufill\crcr\noalign{\kern-\p@\nointerlineskip}%
      $\hfil\displaystyle{#1}\hfil$\crcr}}}
\makeatother

\title{Quantum Knots and New Quantum Field Theory }
\author{Sze Kui Ng
\\ {\small Institute of Mathematics and Information Tecnology, Hanshan Normal University, China}
}
\begin{document}
\date{}
\maketitle

\begin{abstract}
We propose a new quantum gauge field of 
electrodynamics (QED) and  chromodynamics (QCD)
from which we construct quantum knots and derive 
knot invariants such as the Jones polynomial. Our approach is inspired by the work of Witten who derived knot invariants
from quantum field theory based on the Chern-Simon Lagrangian.
From our approach we can derive new knot invariants which extend
the Jones polynomial and give a complete classification of knots.
Then we model elementary particles by  quantum knots. From the construction of quantum knots we construct quantum  photon propagator and quantum gloun propagator for the  nuclear force beween elementary particles. Then from the propagators  a renormalization group equation is derived for the critical phenomena of QED and QCD  of electrons and quarks including superconductivity and color superconductivity.
\end{abstract}

\section{Introduction}\label{sec00}

In 1989 Witten derived knot invariants such as
the Jones polynomial from quantum field theory based on the
Chern-Simon Lagrangian \cite{Wit}.
Inspired by Witten's work
in this paper we shall derive quantum knots and  knot invariants from a new
quantum gauge field of  electrodynamics (QED) and chromodynamics (QCD). 
In our approach we shall first construct quantum knots and their knot invaraints and then derive the Jones
polynomial.

Then we model elementary particles by  quantum knots. From the construction of quantum knots we construct quantum  photon propagator and quantum gloun propagator for the nuclear force beween elementary particles. Then from the propagators a renormalization group equation is derived for the critical phenomena of QED and QCD.
Then from this  renormalization group equation for the critical phenomena we derive a model of superconductivty for electrons and a model of superfluid and color superconductivity for quarks in neutron stars.

This paper is organized as follows. From section \ref{sec2} 
to section \ref{exampleofknot}
we give a brief description of a new quantum gauge field theory
of QED and QCD. 
Fromt this gauge theory 
 we derive 
 two quantum  Knizhnik-Zamolodchikov equations.
From these two KZ equations
we derive a quantum group structure for the  product of Wilson lines. When 
the  products of Wilson lines form closed curves they are as quantum knots.
From product of  Wilson lines and quantum knots
we derive new knot invariants which extend the Jones
polynomial and gives a classification of knots.

Then   from  section \ref{glounpropagator} to   section \ref{nuclearpotential} the quantum knot based on the $U(1)$ group is as the model of photon and quantum knots based on $SU(2)$ and $SU(3)$ are as glouns. When the quantum knots act on complex vectors representing spinors (for electrons or quarks) this gives models of electron, mesons and baryons.
Then 
 the quantum  Wilson line 
 is as a gloun propagator  connecting  elementary particles represented by quantum knots at $z_0$ to elementary particle represented by quantum knots at $z$. When the gauge group is $U(1)$ it is as the photon propagator. We show that the quantum  gloun propagator  
 gives properties of
nuclear force (or strong interaction) that it is strong in the short distance and is weak in the long distance between the  elementary particles and is with the property of asymptotic freedom.

Then from section  \ref{reQED}  to section \ref{color}  with the gloun propagator and the photon propagator we derive critical phenomena of electrons and quarks including superconductivity and color superconductivity.

\section{New gauge model of QED and QCD}\label{sec2}

Let us 
first describe a new quantum field theory.
Similar to the Wiener measure for the Brownian motion
which is constructed from the
integral $\int_{t_0}^{t_1}\left(\frac{dx}{dt}\right)^2dt$
we construct a measure for QED from the 
following energy integral:
\begin{equation}
\begin{array}{rl}
&-\int_{s_0}^{s_1}D ds:= -\int_{s_0}^{s_1}[\frac12\left(\frac{\partial A_1}{\partial
x^2}-\frac{\partial A_2}{\partial x^1}\right)^*
\left(\frac{\partial A_1}{\partial x^2}-\frac{\partial A_2}{\partial x^1}\right)+ 
 \left(
\frac{dZ^*}{ds} +ie_0 A Z^*\right)
\left( \frac{dZ}{ds}
-ie_0A Z\right)]ds \\
\end{array}
\label{1.1}
\end{equation}
where the complex variable $Z=Z(z(s))$ and the real variables
$A_1=A_1(z(s))$ and $A_2=A_2(z(s))$ are continuous functions in a
form that they are in terms of a (continuously differentiable)
curve $z(s)=C(s)=(x^1(s),x^2(s)), s_0\leq s\leq s_1,
z(s_0)=z(s_1)$ in the complex plane where $s$ is a parameter
representing the proper time in relativity (We shall also write
$z(s)$ in the complex variable form
$C(s)=z(s)=x^1(s)+ix^2(s),s_0\leq s\leq s_1$). The complex
variable $Z=Z(z(s))$ represents a field of matter (such as the
electron) ($Z^*$ denotes its complex conjugate) and the real
variables $A_1=A_1(z(s))$ and $A_2=A_2(z(s))$ represent a
connection (or the gauge field of the photon), and $A=\sum_{j=1}^2A_j\frac{dx^j}{ds}$,  and $e_0$ denotes
the (bare) electric charge.

The integral (\ref{1.1}) has the following gauge symmetry:
\begin{equation}
\begin{array}{rl}
Z^{\prime}(z(s))  := Z(z(s))e^{ie_0a(z(s))},  \quad  A'_j (z(s))  := A_j(z(s))+
\frac{\partial a}{\partial x^j} \quad j=1,2
\end{array}
\label{1.2}
\end{equation}
where $a=a(z)$ is a continuously differentiable real-valued
function of $z$.

Similar to the usual Yang-Mills gauge theory we can generalize
this gauge theory with $U(1)$ gauge symmetry to nonabelian gauge
theories. As an illustration let us consider the gauge symmetry $SU(2)\otimes U(1)$
for the model of QCD where $SU(2)\otimes U(1)$ denotes the  direct
product of the groups $SU(2)$ and $U(1)$.
Similar to (\ref{1.1}) we consider the following energy integral:
\begin{equation}
L := 
\int_{s_0}^{s_1} [\frac12 Tr (D_1A_2-D_2A_1)^{*}(D_1A_2-D_2A_1) +
(D_0^*Z^*)(D_0Z)]ds \label{n1}
\end{equation}
where $Z= (z_1, z_2)^{T}$ is a two dimensional complex vector;
$A_j =\sum_{k=0}^{3}A_j^k e_0 t^k $ $(j=1,2)$ where $A_j^k$ denotes a
component of a gauge field $A^k$; $t^k=i T^k$ denotes a
generator of $SU(2)\otimes U(1)$ where $T^k$ denotes a
self-adjoint generator of $SU(2)\otimes U(1)$ 
 (here for simplicity we choose a
convention that the complex $i$ is absorbed by $t^k$ and $t^k$ is
absorbed by $A_j$; and  the notation $A_j$ is with a little
confusion with the notation $A_j$ in the above formulation of
(\ref{1.1}) where $A_j, j=1,2$ are real valued); and
$D_i=\frac{\partial}{\partial x_i}-A_j $ for $i=1,2$; and
$D_0=\frac{d}{ds}- A$ where $A=\sum_{j=1}^2A_j\frac{dx^j}{ds}$.
We remark that the interaction charge is the bare electric charge $e_0$ which is for general interactions
including the strong and weak interactions.

We have that
(\ref{n1}) is invariant under the following gauge transformation:
\begin{equation}
\begin{array}{rl}
Z^{\prime}(z(s)) & :=U(a(z(s)))Z(z(s)) \\
A_j^{\prime}(z(s)) & := U(a(z(s)))A_j(z(s))U^{-1}(a(z(s)))+
 U(a(z(s)))\frac{\partial U^{-1}}{\partial x^j}(a(z(s))),
j =1,2
\end{array}
\label{n2}
\end{equation}
where $U(a(z(s)))=e^{a(z(s))}$; $a(z(s))=\sum_k e_0 a^k (z(s))t^k$  for some functions $a^k$.
We shall mainly consider the case that $a$ is a function of the form $a(z(s))
=\sum_k \mbox{Re}\, \omega^k(z(s))t^k$ where $\omega^k$ are
analytic functions of $z$ (We let
$\omega(z(s)):=\sum_k\omega^k(z(s))t^k$ and we write
$a(z)=\mbox{Re}\,\omega(z)$).

Similarly we have a model of QCD based on the gauge symmetry $SU(3)\times U(1)$.

The above gauge theory is based on the Banach space $X$ of
continuous functions $Z(z(s))$, $A_j(z(s))$, $j=1,2, s_0\leq s\leq
s_1$ on the one dimensional interval $\lbrack s_0, s_1 \rbrack$.
Since $L$ is positive and the theory is one dimensional (and thus
is simpler than the usual two dimensional Yang-Mills gauge theory)
we have that this gauge theory is similar to the Wiener measure
except that this gauge theory has a gauge symmetry.

\section{ Dirac-Wilson loop } \label{sec4}

Similar to the Wilson loop in quantum field theory \cite{Wit},  from
our quantum theory we introduce an analogue of Wilson loop, as
follows (We shall also call a Wilson loop as a Dirac-Wilson loop).

{\bf Definition}.
 A classical Wilson loop $W_R(C)$ is defined by :
\begin{equation}
W_R(C):= Pe^{\int_C A_jdx^j} \label{n4}
\end{equation}
where $R$ denotes a representation of $SU(2)$; $C(\cdot)=z(\cdot)$
is a fixed closed curve where the quantum gauge theories are based
on it as specific in the above section. As usual the notation $P$
in the definition of $W_R(C)$ denotes a path-ordered product
\cite{Wit}\cite{Kau}\cite{Baez}.

When the curve  $C(\cdot)=z(\cdot)$ is not a closed curve, we call  $W_R(C)$ as a classical Wilson line.

 Let  $z_0$ and $z_1$ denotes the  end points of  $C(\cdot)=z(\cdot)$.   Then,  in \cite{Ng1}, we show that,  based on the above quantum gauge model with the gauge invariance property, the quantum version of  $W_R(C)$, denoted by $W(z_0, z_1)$, satisfies the following
 quantum Knizhnik-Zamolodchikov equation \cite{Ng1,Fra,Fuc, Kni}:
\begin{equation}
\begin{array}{ll}
\partial_{z_i}
 W(z_1, z_1^{\prime})\cdots
W(z_n, z_n^{\prime})
 = \frac{e_0^2}{k_0+g_0} \sum_{j\neq i}^{n}
\frac{\sum_a t_i^a \otimes t_j^a}{z_i-z_j}
 W(z_1, z_1^{\prime})\cdots
W(z_n, z_n^{\prime})\,,
\end{array}
\label{n9}
\end{equation}
for $i\,{=}\,1, \dots, n$ where $g_0$ denotes the dual Coxeter
number of the gauge group.
We remark that in (\ref{n9}) we
have defined $t_i^a\,{:=}\, t^a$ and:
\begin{equation}
\begin{array}{ll}
 t_i^a \otimes t_j^a W(z_1, z_1^{\prime})\cdots
W(z_n, z_n^{\prime})  
 := 
 W(z_1, z_1^{\prime})
\cdots[t^aW(z_i, z_i^{\prime})]\cdots [t^aW(z_j,
z_j^{\prime})]\cdots W(z_n, z_n^{\prime})\,.
\end{array}
\label{n9a}
\end{equation}

The above quantum Knizhnik-Zamolodchikov equation is with respect to the left side variables $z_i$.
It is interesting  that we also have the following
quantum Knizhnik-Zamolodchikov equation with respect to the right side
$z_i^{\prime}$ variables which is dual to (\ref{n9}):
\begin{equation}
\begin{array}{ll}
\partial_{z_i^{\prime}}
 W(z_1,z_1^{\prime})\cdots W(z_n,z_n^{\prime})
 = \frac{e_0^2}{k_0+g_0}\sum_{j\neq i}^{n}
 W(z_1, z_1^{\prime})\cdots
W(z_n, z_n^{\prime})
\frac{\sum_a t_i^a\otimes t_j^a}{z_j^{\prime}-z_i^{\prime}}
\end{array}
\label{d8}
\end{equation}
for $i\,{=}\,1, \dots, n$ where we have defined:
\begin{equation}
\begin{array}{ll}
 W(z_1, z_1^{\prime})\cdots
W(z_n, z_n^{\prime})t_i^a \otimes t_j^a 
 := 
W(z_1, z_1^{\prime})
\cdots [W(z_i, z_i^{\prime})t^a]\cdots [W(z_j,
z_j^{\prime})t^a]\cdots W(z_n, z_n^{\prime})\,.
\end{array}
\label{d8a}
\end{equation}

\section{Solving quantum KZ equation}\label{sec8a}

Let us consider the following product of two quantum
Wilson lines:
\begin{equation}
G(z_1,z_2, z_3, z_4):=
 W(z_1, z_2)W(z_3, z_4)
\label{m1}
\end{equation}
where the two quantum Wilson lines $W(z_1, z_2)$ and
$W(z_3, z_4)$ represent two pieces
of curves starting at $z_1$ and $z_3$ and ending at
$z_2$ and $z_4$ respectively.
We have that this product $G$ satisfies the KZ equation for the
variables $z_1$, $z_3$ and satisfies the dual KZ equation
for the variables $z_2$ and $z_4$.
Then
by solving the two-variables-KZ equation in (\ref{n9}) we have that a form of $G$ is
given by \cite{Chari}\cite{Koh}\cite{Dri}:
\begin{equation}
e^{-\hat{t}\log [\pm (z_1-z_3)]}C_1
\label{m2}
\end{equation}
where $\hat{t}:=\frac{-e_0^2}{k_0+g_0}\sum_a t^a \otimes t^a$
and $C_1$ is a constant matrix which is independent
of the variable $z_1-z_3$.
We see that $G$ is a multi-valued analytic function where the
determination of the $\pm$ sign depended on the choice of the
branch.
Then by solving the dual two-variable-KZ equation
 in (\ref{d8}) we have that
$G$ is of the form:
\begin{equation}
C_2e^{\hat{t}\log [\pm (z_4-z_2)]}
\label{m3}
\end{equation}
where $C_2$ is a constant matrix which is independent
of the variable $z_4-z_2$.
From (\ref{m2}), (\ref{m3}) and letting:
\begin{equation}
C_1=Ae^{\hat{t}\log[\pm (z_4-z_2)]}, \quad
C_2=e^{-\hat{t}\log[\pm(z_1-z_3)]}A
 \label{am3}
\end{equation}
 where $A$ is an initial operator  we have that
$G$ is given by
\begin{equation}
G(z_1, z_2, z_3, z_4)=
e^{-\hat{t}\log [\pm (z_1-z_3)]}Ae^{t\log [\pm (z_4-z_2)]}
\label{m4}
\end{equation}
where at the singular case that $z_1=z_3$ we simply define $\log [\pm (z_1-z_3)]=0$. Similarly
for $z_2=z_4$.

\section{Computation of quantum Wilson lines and loops}\label{sec 8aa}

Let us consider the following product of two quantum
Wilson lines:
\begin{equation}
G(z_1,z_2, z_3, z_4):=
W(z_1, z_2)W(z_3, z_4)
\label{h1}
\end{equation}
where the two quantum Wilson lines $W(z_1, z_2)$ and
$W(z_3, z_4)$ represent two pieces
of curves starting at $z_1$ and $z_3$ and ending at
$z_2$ and $z_4$ respectively.
As shown in the above section we have that $G$
is given by the following formula:
\begin{equation}
G(z_1, z_2, z_3, z_4)=
e^{-\hat{t}\log [\pm (z_1-z_3)]}Ae^{\hat{t}\log [\pm (z_4-z_2)]}
\label{m4a}
\end{equation}
where the product is
a 4-tensor.
Let us set $z_2=z_3$. Then
the 4-tensor $W(z_1, z_2)W(z_3, z_4)$ is reduced to the 2-tensor
$W(z_1, z_2)W(z_2, z_4)$.
By using (\ref{m4a}) the 2-tensor
$W(z_1, z_2)W(z_2, z_4)$
is given by:
\begin{equation}
W(z_1, z_2)W(z_2, z_4)
=e^{-\hat{t}\log [\pm (z_1-z_2)]}A_{14}e^{\hat{t}\log [\pm (z_4-z_2)]}
\label{closed1}
\end{equation}
where $A_{14}=A_1\otimes A_4$ is a 2-tensor reduced from the 4-tensor
$A=A_1\otimes A_2\otimes A_3\otimes A_4$ in (\ref{m4a}). In this reduction the $\hat{t}$ operator
of $\Phi=e^{-\hat{t}\log [\pm (z_1-z_2)]}$ acting on the left side of $A_1$ and $A_3$ in $A$ is reduced
to acting on the left side of $A_1$ and $A_4$ in $A_{14}$. Similarly
the $\hat{t}$ operator of $\Psi=e^{-\hat{t}\log [\pm (z_4-z_2)]}$ acting on the right side of $A_2$
and $A_4$ in $A$ is reduced to acting on the right side of $A_1$ and $A_4$ in $A_{14}$.

Then since $\hat{t}$ is a 2-tensor operator we have that $\hat{t}$ is as a matrix acting on the two sides of
the 2-tensor $A_{14}$ which is also as a matrix with the same dimension as $\hat{t}$.
Thus $\Phi$ and $\Psi$ are as matrices of the same dimension as the matrix
$A_{14}$  acting on $A_{14}$ by the usual matrix operation.
Then since $\hat{t}$ is a Casimir operator for the 2-tensor group representation of $SU(2)$ we have that
$\Phi $  and $\Psi $ commute
with $A_{14}$ since  $\Phi $  and $\Psi$ are exponentials
of $t$ (We remark that $\Phi $  and $\Psi $ are in general not commute with
the 4-tensor initial operator $A$).
Thus we have
\begin{equation}
e^{-\hat{t}\log [\pm (z_1-z_2)]}A_{14}e^{\hat{t}\log[\pm (z_4-z_2)]}
=e^{-\hat{t}\log [\pm (z_1-z_2)]}e^{\hat{t}\log[\pm (z_4-z_2)]}A_{14}
\label{closed1a}
\end{equation}

 Let us set $z_1=z_4$. In this case the quantum
Wilson line forms a closed loop. Now in (\ref{closed1a}) with
$z_1=z_4$ we have that $e^{-\hat{t}\log  \pm (z_1-z_2)}$ and $e^{\hat{t}\log
\pm (z_1-z_2)}$ which come from the two-side KZ equations cancel
each other and from the multi-valued property of the $\log$
function we have
\begin{equation}
W(z_1, z_1) =R^{N}A_{14} \quad\quad N=0, \pm 1, \pm 2, ...
\label{closed2}
\end{equation}
where $R=e^{-i\pi \hat{t}}$ is the monodromy of the KZ equation \cite{Chari}.

\section{Representing braiding of curves by quantum Wilson Lines}\label{sec 9aa}

Consider again the product $G(z_1, z_2, z_3, z_4)=W(z_1,z_2)W(z_3,z_4)$.
We have that $G$ is a multivalued analytic function
where the determination of the $\pm$ sign depended on the choice of the
branch.

Let the two pieces of curves be crossing at $w$. Then we have $W(z_1,z_2)=W(z_1,w)W(w,z_2)$ and
 $W(z_3,z_4)=W(z_3,w)W(w,z_4)$. Thus we have
\begin{equation}
W(z_1,z_2)W(z_3,z_4)=
W(z_1,w)W(w,z_2)W(z_3,w)W(w,z_4)
\label{h2}
\end{equation}

If we interchange $z_1$ and $z_3$, then from
(\ref{h2}) we have the following ordering:
\begin{equation}
 W(z_3,w)W(w, z_2)W(z_1,w)W(w,z_4)
\label{h3}
\end{equation}

Now let us choose a  branch. Suppose that
these two curves are cut from a knot and that
following the orientation of a knot the
curve represented by  $W(z_1,z_2)$ is before the
curve represented by  $W(z_3,z_4)$. Then we fix a branch such that the  product in (\ref{m4a}) is
with two positive signs :
\begin{equation}
W(z_1,z_2)W(z_3,z_4)=
e^{-t\log(z_1-z_3)}Ae^{t\log(z_4-z_2)}
\label{h4}
\end{equation}

Then if we interchange $z_1$ and $z_3$ we have
\begin{equation}
W(z_3,w)W(w, z_2)W(z_1,w)W(w,z_4) =
e^{-t\log[-(z_1-z_3)]}Ae^{t\log(z_4-z_2)}
\label{h5}
\end{equation}
From (\ref{h4}) and (\ref{h5}) as a choice of branch we have
\begin{equation}
W(z_3,w)W(w, z_2)W(z_1,w)W(w,z_4) =
R W(z_1,w)W(w,z_2)W(z_3,w)W(w,z_4)
\label{m7a}
\end{equation}
where $R=e^{-i\pi t}$ is the monodromy of the KZ equation.
In (\ref{m7a}) $z_1$ and $z_3$ denote two points on a closed curve
such that along the direction of the curve the point
$z_1$ is before the point $z_3$ and in this case we choose
a branch such that the angle of $z_3-z_1$ minus the angle
of $z_1-z_3$ is equal to $\pi$.

Now from (\ref{m7a}) we can take a convention that the ordering (\ref{h3}) represents that
the curve represented by  $W(z_1,z_2)$ is upcrossing
the curve represented by  $W(z_3,z_4)$ while
(\ref{h2}) represents zero crossing of these two
curves.

Similarly from the dual KZ equation as a choice of branch which
is consistent with the above formula we have
\begin{equation}
W(z_1,w)W(w,z_4)W(z_3,w)W(w,z_2)=
W(z_1,w)W(w,z_2)W(z_3,w)W(w,z_4)R^{-1}
\label{m8a}
\end{equation}
where $z_2$ is before $z_4$. We take a convention that the ordering (\ref{m8a}) represents that
the curve represented by $W(z_1,z_2)$ is undercrossing the curve represented by $W(z_3,z_4)$.
Here along the orientation of a closed curve the piece of curve
represented by $W(z_1,z_2)$ is before the piece of curve represented by
$W(z_3,z_4)$. In this case since the angle of $z_3-z_1$ minus the angle
of $z_1-z_3$ is equal to $\pi$ we have that the
angle of $z_4-z_2$ minus the angle of $z_2-z_4$ is
also equal to $\pi$ and this gives the $R^{-1}$ in this formula
(\ref{m8a}).

From (\ref{m7a}) and (\ref{m8a}) we have
\begin{equation}
 W(z_3,z_4)W(z_1,z_2)=
RW(z_1,z_2)W(z_3,z_4)R^{-1} \label{m9}
\end{equation}
where $z_1$ and $z_2$ denote the end points of a curve which is before a curve with end points $z_3$ and $z_4$.
From (\ref{m9}) we see that the algebraic structure of these
quantum Wilson lines $W(z,z')$
is analogous to the quasi-triangular quantum
group \cite{Fuc}\cite{Chari}.

\section{Skein relation for the HOMFLY polynomial }\label{sec9a}

In this section let us apply the above result which is from the KZ equation in dual form
to derive the skein relation for the HOMFLY polynomial which includes the Jones polynomial. 
From this relation we then have the skein relation for the Jones polynomial which
is a special case of the HOMFLY polynomial \cite{Kau}\cite{Mur}\cite{Lic}.  

It is well known that from the one-side KZ equation
we can derive a braid group representation which is
related to the the derivation of the
skein relation for the Jones polynomial \cite{Koh}\cite{Dri}. 
We shall see
here that by applying the two KZ equations of the KZ equation in dual form
we also have a  way to derive the
skein relation for the HOMFLY polynomial.

Let us first consider the following theorem
of Kohno and Drinfield \cite{Chari}\cite{Koh}\cite{Dri}:
\begin{theorem}[Kohno-Drinfield]
Let $R$ be the monodromy of the KZ equation for the group $SU(2)$
and let $\bar R$ denotes the $R$-matrix of
the quantum group $U_q (su(2))$ where $su(2)$ denotes the Lie algebra
of $SU(2)$ and $q=e^{\frac{i2\pi}{k+g}}$ where $g=2$. Then there exists a twisting
$F\in U_q(su(2))\otimes U_q(su(2))$ such that
\begin{equation}
\bar R=F^{-1}RF^{-1}
\label{R}
\end{equation}
From this relation we have that the braid group representations obtained from the quantum group $U_q(su(2)$
and obtained from the one-side KZ equation are equivalent.
\end{theorem}
We shall use only the relation (\ref{R}) of this theorem  to derive
the skein relation of the HOMFLY polynomial.

From the property of the quantum group $U_q(su(2))$ we
also have the following formula 
\begin{equation}
\bar R^{ 2} -(q^{\frac12}-q^{-\frac12}) \bar R
 - I
=0
\label{h6}
\end{equation}
Thus we have
\begin{equation}
\bar R -(q^{\frac12}-q^{-\frac12})I - \bar R^{ -1}
=0
\label{h6a}
\end{equation}
By using this formula we have
\begin{equation}
[\bar R -(q^{\frac12}-q^{-\frac12})I - \bar R^{-1 }]
 FW(z_1,w)W(w,z_2)W(z_3,w)W(w,z_4)F^{-1}
=0
\label{h6b}
\end{equation}
Thus by using the relation (\ref{R}) we have
\begin{equation}
\begin{array}{rl}
& Tr F^{-1}R W(z_1,w)W(w,z_2)W(z_3,w)W(w,z_4)F^{-1}
  \\
- &
(q^{\frac12}-q^{-\frac12})
Tr FW(z_1,z_2)W(z_3,z_4)F^{-1} \\
-
& Tr FW(z_1,w)W(w,z_2)W(z_3,w)W(w,z_4)R^{-1}F
=0
\end{array}
\label{h7a}
\end{equation}
Then by using  the formulas
(\ref{m7a}) and (\ref{m8a})
for upcrossing and undercrossing from (\ref{h6a})
we have
\begin{equation}
\begin{array}{rl}
& Tr
F^{-1}W(z_3,w)W(w,z_2)W(z_1,w)W(w,z_4)F^{-1} \\
- &
(q^{\frac12}-q^{-\frac12})
Tr  W(z_1,z_2)W(z_3,z_4) \\
- &
 Tr\langle FW(z_1,w)W(w,z_4)W(z_3,w)W(w,z_2)F\rangle
=0
\end{array}
\label{h7}
\end{equation}
Let us make a further twist that replace $F^2$ by
$F^2x$  where $x$ denotes a nonzero variable. Then from
(\ref{h7}) we have the following skein relation
for the HOMFLY polynomial:
\begin{equation}
xL_{+}+yL_{0}-x^{-1}L_{-}=0 \label{hh7}
\end{equation}
where we define $y=q^{-\frac12}-q^{\frac12}$ and
that $L_{+}$, $L_{0}$
and $L_{-}$ are defined by
\begin{equation}
\begin{array}{rl}
L_{+}=&
Tr F^{-2}x^{-1}
W(z_3,w)W(w,z_2)W(z_1,w)W(w,z_4)
\\
L_{0}=&
Tr W(z_1,z_2)W(z_3,z_4)
\\
L_{-}=&
Tr xF^{2}W(z_1,w)W(w,z_4)W(z_3,w)W(w,z_2)
\end{array}
\label{h8}
\end{equation}
which are as the HOMFLY polynomials for upcrossing, zero crossing
and undercrossing respectively.

 \section{Defining quantum knots and knot invariant}\label{sec10}
 
Now
we have that the quantum Wilson loop $W(z_1, z_1)$ corresponds to a closed
curve in the complex plane with starting and ending
point $z_1$.
Let this quantum Wilson loop $W(z_1, z_1)$ represents the unknot. We shall call $W(z_1, z_1)$ as the quantum unknot. Then from
(\ref{closed2}) we have the following invariant
for the unknot:
\begin{equation}
Tr W(z_1, z_1)= Tr R^{n}A \quad\quad n=0, \pm 1, \pm 2, ...
\label{m6}
\end{equation}
where $A=A_{14}$ is a 2-tensor constant matrix operator. 

In the following let us extend the definition (\ref{m6})
to a knot invariant for nontrivial knots.

Let $W(z_i,z_j)$ represent a piece of curve
with starting point $z_i$ and ending point $z_j$.
Then we let
\begin{equation}
W(z_1,z_2)W(z_3,z_4)
\label{m11}
\end{equation}
represent two pieces of uncrossing curve.
Then by interchanging $z_1$ and $z_3$ we have
\begin{equation}
W(z_3,w)W(w,z_2)W(z_1,w)W(w,z_4)
\label{m12}
\end{equation}
represent the curve specified by $W(z_1,z_2)$ upcrossing the
curve specified by $W(z_3,z_4)$.

Now for a given knot diagram we may cut it into a sum of
parts which are formed by two pieces of curves crossing  each other.
Each of these parts is represented
by  (\ref{m12})( For a knot diagram of the unknot
with zero crossings we simply do not need to cut the
knot diagram).
Then we define the trace of a knot with a
given knot diagram by the following form:
\begin{equation}
 Tr \cdot\cdot\cdot
 W(z_3,w)W(w,z_2)W(z_1,w)W(w,z_4)
\cdot\cdot\cdot
 \label{m14}
\end{equation}
where we use (\ref{m12})  to represent the state of the
two pieces of curves specified by
 $W(z_1,z_2)$ and
$W(z_3,z_4)$. The
 $\cdot\cdot\cdot$ means the product
of a sequence of parts represented by
(\ref{m12}) according to the state of
each part. The ordering of the sequence in (\ref{m14})
 follows the ordering of the parts given by the orientation of the
knot diagram. We shall call the sequence of crossings in
the trace (\ref{m14}) as the generalized Wilson
loop of the knot diagram. For the knot diagram of the unknot with zero crossings we simply
let it be $W(z,z)$ and call it the quantum Wilson loop.

We shall
 show that the generalized Wilson loop of a knot diagram has all the properties of the knot diagram  and that
(\ref{m14}) is  a knot invariant. From this we shall call a generalized Wilson loop as a quantum knot. 

\section{Examples of quantum knots} \label{exampleofknot}

Before the proof that a generalized Wilson loop of a knot diagram has all the properties of the knot diagram 
in the following let us first consider
some examples to illustrate the way to define (\ref{m14}) and
the way of applying the  braiding formulas (\ref{m7a}),
 (\ref{m8a}) and (\ref{m9}) to
equivalently transform (\ref{m14}) to a simple
expression of the form  $Tr R^{-m}W(z,z)$ where $m$
is an integer. 

Let us first consider the knot in Fig.1.
For this knot we have that (\ref{m14}) is given by
\begin{equation}
Tr W(z_2,w)W(w,z_2)W(z_1,w)W(w,z_1)
\label{m15a}
\end{equation}
where the product of quantum Wilson lines  is from the definition (\ref{m12})
represented a crossing at $w$.
In applying (\ref{m12}) we let $z_1$ be the
starting and the ending point.

\begin{figure}[hbt]
\centering
\includegraphics[scale=0.6]{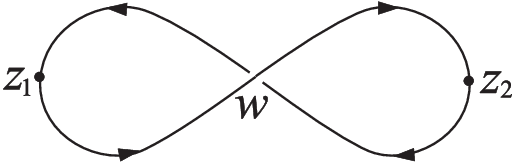}

                Fig.1
\end{figure}
Then  we have that (\ref{m15a}) is equal to
\begin{equation}
\begin{array}{rl}
&Tr W(w,z_2)W(z_1,w)W(w,z_1)W(z_2,w) \\
=&Tr RW(z_1,w)W(w,z_2)R^{-1}
RW(z_2,w)W(w,z_1)R^{-1} \\
=&Tr W(z_1,z_2)W(z_2,z_1) \\
=&Tr W(z_1,z_1)
\end{array}
\label{m16}
\end{equation}
where we have used (\ref{m9}).
We see that (\ref{m16}) is just the knot invariant (\ref{m6}) of
the unknot.
Thus the knot in Fig.1 is with the same knot invariant of the unknot and this
agrees with the fact that this knot is topologically equivalent
to the unknot.

Then let us derive the Reidemeister move 1. Consider
the diagram in Fig.2. We have that by (\ref{m12}) the definition (\ref{m14})
 for this
diagram is given by:
\begin{equation}
\begin{array}{rl}
& Tr\langle W(z_2,w)W(w,z_2)W(z_1,w)W(w,z_3)\rangle\\
=& Tr\langle W(z_2,w)RW(z_1,w)W(w,z_2)
R^{-1}W(w,z_3)\rangle\\
=& Tr\langle W(z_2,w)RW(z_1,z_2)R^{-1}W(w,z_3)\rangle\\
=& Tr\langle R^{-1}W(w,z_3)W(z_2,w)RW(z_1,z_2)\rangle\\
=& Tr\langle W(z_2,w)W(w,z_3)W(z_1,z_2)\rangle\\
=& Tr\langle W(z_2,z_3)W(z_1,z_2)\rangle\\
=& Tr\langle W(z_1,z_3)\rangle\\
\end{array}
\label{m18}
\end{equation}
where $W(z_1,z_3)$ represent a piece of curve with initial
end point $z_1$ and final end point $z_3$ which has no
crossing. When Fig.2 is a part of a knot we can also
derive a result similar to (\ref{m18}) which is for the
Reidemeister move 1. This shows that the Reidemeister move 1 holds.
\begin{figure}[hbt]
\centering
\includegraphics[scale=0.4]{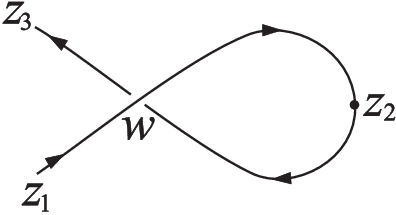}

                    Fig.2
\end{figure}
Then let us derive Reidemeister move 2.
By (\ref{m12}) we have that the definition (\ref{m14})
for the two pieces of curve in 
Fig.3a is given by
\begin{equation}
\begin{array}{rl}
&Tr\langle W(z_5,w_1)W(w_1,z_2)W(z_1,w_1)W(w_1,z_6)\cdot\\
&W(z_4,w_2)W(w_2,z_3)W(z_2,w_2)W(w_2,z_5)\rangle
\end{array}
\label{m19}
\end{equation}
where the two products of $W$ separated by the 
 $\cdot$
are for the two crossings in Fig.3a.
We have that (\ref{m19}) is
equal to
\begin{equation}
\begin{array}{rl}
&Tr\langle W(z_4,w_2)W(w_2,z_3)W(z_2,w_2)W(w_2,z_5)\cdot\\
         & W(z_5,w_1)W(w_1,z_2)W(z_1,w_1)W(w_1,z_6)\rangle \\
=&Tr\langle W(z_1,z_3)W(z_4,z_6)\rangle
\end{array}
\label{m20}
\end{equation}
where we have repeatly used (\ref{m9}). This shows that the diagram in Fig.3a is equivalent to
two uncrossing curves. When Fig.3a is a part of a knot
we can also derive a result similar to (\ref{m20}) for
the Reidemeister
move 2.
This shows that the Reidemeister
move 2 holds.
\begin{figure}[hbt]
\centering
\includegraphics[scale=0.4]{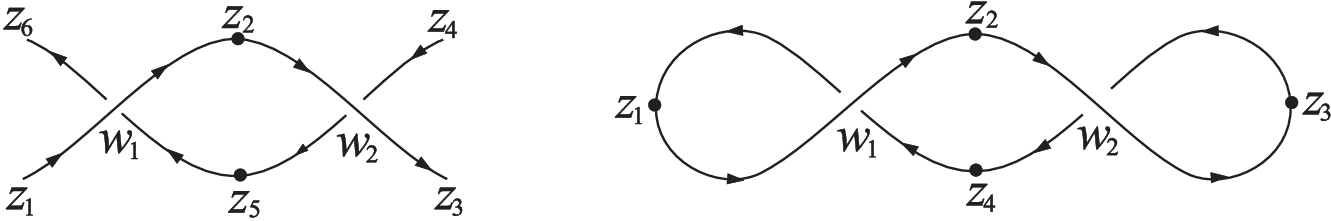}

               Fig.3a   \hspace*{3.6 cm}          Fig.3b
\end{figure}
As an illustration let us consider the knot in Fig.3b which is related to
the Reidemeister move 2. By (\ref{m12}) we have that
the definition (\ref{m14}) for this knot is given by
\begin{equation}
\begin{array}{rl}
&Tr\langle W(z_3,w_2)W(w_2,z_3)W(z_2,w_2)W(w_2,z_4)\cdot\\
&W(z_4,w_1)W(w_1,z_2)W(z_1,w_1)W(w_1,z_1)\rangle \\
=&Tr\langle 
RW(z_2,w_2)W(w_2,z_3)W(z_3,w_2)W(w_2,z_4)\cdot\\
&W(z_4,w_1)W(w_1,z_1)W(z_1,w_1)W(w_1,z_2)R^{-1}\rangle \\
=&Tr\langle W(z_2,w_2)W(w_2,z_3)W(z_3,w_2)W(w_2,z_4)\cdot\\
&W(z_4,w_1)W(w_1,z_1)W(z_1,w_1)W(w_1,z_2)\rangle \\
=&Tr\langle W(z_2,z_2)\rangle
\end{array} 
\label{m20a}
\end{equation}
where we let the curve be with $z_2$ as the initial
and final end point and we have used (\ref{m7a}) and
(\ref{m8a}). This shows that the knot in Fig.3b is with the 
same knot
invariant of a trivial knot. This agrees with the fact
that this knot is equivalent to the trivial knot.
Similar to the above derivations we can derive the
Reidemeister move 3.
                
Let us then consider a trefoil knot in Fig.4a.
By (\ref{m12}) and similar to the above examples
we have that the definition (\ref{m14})
for this knot is given by:
\begin{equation}
\begin{array}{rl}
&Tr\langle W(z_4,w_1)W(w_1,z_2)W(z_1,w_1)W(w_1,z_5)\cdot
W(z_2,w_2)W(w_2,z_6)\\
&W(z_5,w_2)W(w_2,z_3)\cdot 
W(z_6,w_3)W(w_3,z_4)W(z_3,w_3)W(w_3,z_1)\rangle \\
=&Tr\langle 
W(z_6,z_1)W(z_3,z_6)W(z_1,z_3)\rangle
\end{array}
\label{m21}
\end{equation}
where we have repeatly used (\ref{m9}). Then
 we have that (\ref{m21}) is equal to:
\begin{equation}
\begin{array}{rl}
&Tr\langle
W(z_6,w_3)W(w_3,z_1)W(z_3,w_3)W(w_3,z_6)W(z_1,z_3)\rangle
\\
=&Tr\langle
RW(z_3,w_3)W(w_3,z_1)W(z_6,w_3)W(w_3,z_6)W(z_1,z_3)\rangle
\\
=&Tr\langle
RW(z_3,w_3)RW(z_6,w_3)W(w_3,z_1)
R^{-1}W(w_3,z_6)W(z_1,z_3)\rangle\\
=&Tr\langle
W(z_3,w_3)RW(z_6,z_1)
R^{-1}W(w_3,z_6)W(z_1,z_3)R\rangle\\
=&Tr\langle
W(z_3,w_3)RW(z_6,z_3)W(w_3,z_6)\rangle\\
=&Tr\langle W(w_3,z_6)W(z_3,w_3)RW(z_6,z_3)\rangle\\
=&Tr\langle RW(z_3,w_3)W(w_3,z_6)W(z_6,z_3)\rangle\\
=&Tr\langle RW(z_3,z_3)\rangle
\end{array}
\label{m22}
\end{equation}
where we have used (\ref{m7a}) and (\ref{m9}).
We see that (\ref{m22}) is a knot invariant for the trefoil knot in Fig.4a.
\begin{figure}[hbt]
\centering
\includegraphics[scale=0.4]{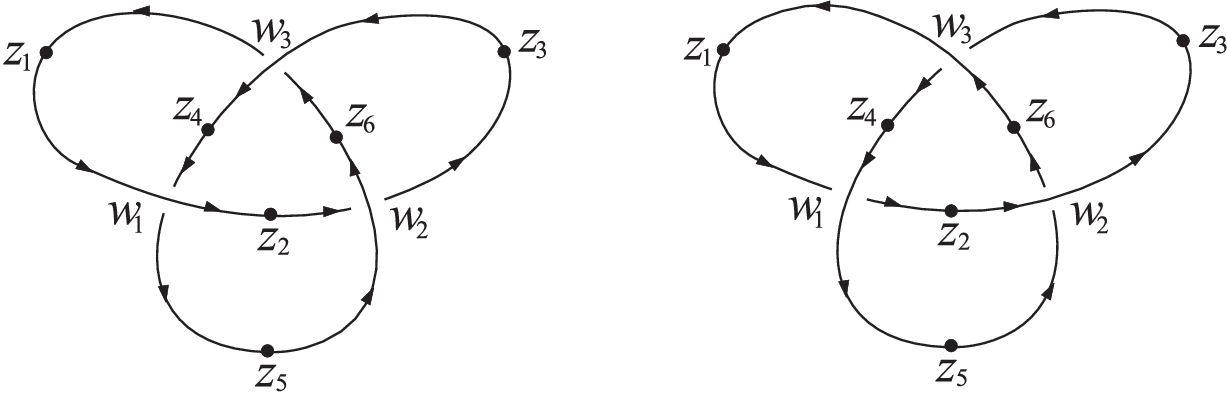}

             Fig.4a  \hspace*{3.5cm}  Fig.4b
\end{figure}
Similarly for the trefoil knot in Fig. 4b which
is the mirror image of the trefoil knot in Fig.4a
we have that the definition (\ref{m14}) for this knot is
equal to 
$Tr\langle W(z_5,z_5)R^{-1}\rangle $ which is a knot invariant for the trefoil knot in
Fig.4b. We notice that this knot invariant is different from
(\ref{m22}). This shows that these two
trefoil knots are not topologically equivalent.

With the above new knot invariants we can now
give a classification of knots.
Let $C_1$ and $C_2$ be two knots. Then since we can derive
the Reidemeister moves by using (\ref{m7a}), (\ref{m8a}) and (\ref{m9}), we have that
$C_1$ and $C_2$ are topologically equivalent if and only
if the $W$-product and the correlation (\ref{m14}) for $C_1$  can be
transformed to the $W$-product and the correlation (\ref{m14}) for $C_2$ by using (\ref{m7a}), (\ref{m8a}) and (\ref{m9}). 
Thus correlations (\ref{m14}) are  knot invariants which
can completely classify
knots.
 In  \cite{Ng},  we use other approach to show that  knots 
can be classfied by the number of product of $R$ and $R^{-1}$
matrices.
 More calculations and examples of the
above knot invariants will be given elsewhere.

\section{ Quantum knots as gluons and quantum Wilson line as gloun propagator} \label{glounpropagator}

Since  the quantum knot based on the $U(1)$ group has  properties of photon with winding numbers as frequencies and with the loop nature giving the spin, it  is as the model of photon. Then we let quantum knots based on $SU(2)$ and $SU(3)$ be as models of glouns. 

When the quantum knots act on complex vectors representing spinors (for electrons or quarks),  the quantum knots together with the complex vectors  are as models of electron, mesons and baryons. For example, the meson $\pi$ is modeled by quantum prime knot  ${\bf4_1}$, and the proton is modeled by the quantum prime knot ${\bf 6_2}$.
More details of modeling elementary particles by quantum knots are given in  \cite{Ng}.

Let us then investigate the quantum Wilson line $W(z_0,z)$ with
$U(1)$ group where $z_0$ is fixed for the photon field. We want to
show that this quantum Wilson line $W(z_0,z)$ may be regarded as
the quantum photon propagator for a photon propagating from $z_0$
to $z$.
As we have shown in the above section on computation of quantum
Wilson line; to compute $W(z_0,z)$ we need to write $W(z_0,z)$ in
the form of two (connected) Wilson lines:
$W(z_0,z)=W(z_0,z_1)W(z_1,z)$ for some $z_1$ point.  Then we have:
\begin{equation}
W(z_0,z_1)W(z_1,z)=e^{-\hat{t}\log [\pm (z_1-z_0)]}Ae^{\hat{t}\log [\pm
(z-z_1)]} \label{graviton6b}
\end{equation}
where $\hat{t}=\frac{e_0^2}{k_0}$ for the $U(1)$ group ($k_0>0$ is a
constant and we may for simplicity let $k_0=1$) where the term
$e^{-\hat{t}\log [\pm (z-z_0)]}$ is obtained by solving the first form
of the dual form of the KZ equation and the term $e^{\hat{t}\log [\pm
(z_0-z)]}$ is obtained by solving the second form of the dual form
of the KZ equation.

Then we may write $W(z_0,z)$ in the following form:
\begin{equation}
W(z_0,z)=W(z_0,z_1)W(z_1,z)=e^{\hat{t}\log \frac{(z-z_1)}{(z_1-z_0)}}A
\label{graviton6}
\end{equation}

Let us fix $z_1$ with $z_0$ such that:
\vspace*{-2pt}
\begin{equation}
\frac{|z-z_1|}{|z_1-z_0|}=\frac{r_1}{n_e^2} \label{egraviton6}
\end{equation}
for some positive integer $n_e $ such that $r_1\geq n_e^2$;
and we let $z$ be a point on a path of connecting $z_0$ and $z_1$
and then a closed loop is formed with $z$ as the starting and
ending point. (This loop can just be the photon-loop of the
electron in this electromagnetic interaction by this photon
propagator (\ref{graviton6}).) Then (\ref{graviton6}) has a factor
$e_0^2\log \frac{r_1}{n_e^2}$ which is the fundamental solution of
the two dimensional Laplace equation and is analogous to the
fundamental solution $\frac{e^2}{r}$ (where $e\,{:=}\,e_0 n_e$
denotes the observed (renormalized) electric charge and $r$
denotes the three dimensional distance) of the three dimensional
Laplace equation for the Coulomb's law. Thus the operator
$W(z_0,z)\,{=}\,W(z_0,z_1)W(z_1,z)$ in (\ref{graviton6}) can be
regarded as the quantum photon propagator propagating from $z_0$
to $z$.

We remark that when there are many photons we may introduce the
space variable $x$ via the Lorentz metric $ds^2=dt^2-dx^2$ as a
statistical variable to obtain the Coulomb's law $\frac{e^2}{r}$
from the fundamental solution $e_0^2\log \frac{r_1}{n_e^2}$ as a
statistical law for electricity (We shall give such a space-time
statistics later).

The quantum photon propagator (\ref{graviton6}) gives a repulsive
effect since it is analogous to the Coulomb's law $\frac{e^2}{r}$.
On the other hand we can reverse the sign of $\hat{t}$ such that
this photon operator can also give an attractive effect:
\begin{equation}
W(z_0,z)=W(z_0,z_1)W(z_1,z)=e^{-\hat{t}\log
\frac{(z_1-z_0)}{(z-z_1)}}A\,, \label{graviton6a}
\end{equation}
where we fix $z_1$ with $z$ such that:
\begin{equation}
\frac{|z_1-z_0|}{|z-z_1|}=\frac{r_1}{n_e^2} \label{egraviton6a}
\end{equation}
for some positive integer $n_e $ such that $r_1\geq n_e^2$; and we
again let $z$ be a point on a path of connecting $z_0$ and $z_1$
and then a closed loop is formed with $z$ as the starting and
ending point. (This loop again can just be the photon-loop of the
electron in this electromagnetic interaction by this photon
propagator (\ref{graviton6}).) Then (\ref{graviton6a}) has a
factor $-e_0^2\log \frac{r_1}{n_e^2}$ which is the fundamental
solution of the two dimensional Laplace equation and is analogous
to the attractive fundamental solution $-\frac{e^2}{r}$ of the
three dimensional Laplace equation for the Coulomb's law.
Thus the quantum photon propagator in (\ref{graviton6}) (and in
(\ref{graviton6a})) can give repulsive or attractive effect
between two points $z_0$ and $z$ for all $z$ in the complex plane.
These repulsive or attractive effects of the quantum photon
propagator correspond to two charges of the same sign and of
different sign respectively.
On the other hand when $z=z_0$ the quantum Wilson line
$W(z_0,z_0)$  in (\ref{graviton6}) which is the quantum photon
propagator becomes a quantum Wilson loop $W(z_0,z_0)$ which is
identified as a photon, as shown in the above sections.

Let us then derive a usual form of photon propagator from the
quantum photon propagator $W(z_0,z)$. Let us choose a path
connecting $z_0$ and $z$. Let us consider the following path:
\begin{equation}
z=z(s)=z_1+ a_0[\theta(s_1-s)e^{-i\beta_1(s_1-s)}+
\theta(s-s_1)e^{i\beta_1(s_1-s)}] 
\label{g1g}
\end{equation}
where $\beta_1>0 $ is a parameter and $z(s_0)=z_0$ for some proper
time $s_0$; and $a_0$ is some complex constant; and $\theta$ is a
step function given by $\theta(s)=0 $ for $s<0$, $\theta(s)=1 $
for $s\geq 0$. Then on this path we have:
\begin{equation}
\begin{array}{rl}
  & W(z_0,z)
= W(z_0,z_1)W(z_1,z) =e^{\hat{t}\log \frac{(z-z_1)}{(z_1-z_0)}}A
=  e^{\hat{t}\log \frac{a_0[\theta(s-s_1)e^{-i\beta_1(s_1-s)}+
\theta(s_1-s)e^{i\beta_1(s_1-s)}]}{(z_1-z_0)}}A \\
&\\
=&  e^{\hat{t}\log b[\theta(s-s_1)e^{-i\beta_1(s_1-s)}+
\theta(s_1-s)e^{i\beta_1(s_1-s)}]}A
=b_0 [\theta(s-s_1)e^{-i\hat{t}\beta_1(s_1-s)}+
\theta(s_1-s)e^{i\hat{t}\beta_1(s_1-s)}]A\\
\end{array}
\label{g2}
\end{equation}
for some complex constants $b$ and $b_0$. From this chosen of the
path (\ref{g1g}) we have that the quantum photon propagator is
proportional to the following expression:
\begin{equation}
\frac{1}{2\lambda_1}[\theta(s-s_1)e^{-i\lambda_1(s-s_1)}+
\theta(s_1-s)e^{i\lambda_1(s-s_1)}] \label{g3}
\end{equation}
where we define $\lambda_1=-\hat{t}\beta_1=e_0^2\beta_1>0$.
 We see that this is the usual propagator of a
particle $x(s)$ of harmonic oscillator with mass-energy parameter
$\lambda_1 >0$ where $x(s)$ satisfies the following harmonic
oscillator equation:
\begin{equation}
\frac{d^2x}{ds^2} =-\lambda_1^2 x(s) \label{g4}
\end{equation}
We regard (\ref{g3}) as the propagator of a photon with
mass-energy parameter $\lambda_1$. Fourier transforming (\ref{g3})
we have the following form of photon propagator:
\begin{equation}
\frac{i}{k_E^2-\lambda_1^2} \label{g5}
\end{equation}

Then when the gauge group is the $SU(2)\otimes U(1)$ (or  $SU(3)\otimes U(1)$) group the analysis is similar to the case of  $U(1)$ group. In the following section we show that the corresponding  quantum gloun propagator gives the nuclear force. 

\section{Quantum gluon propagator as nuclear  potential} 
 \label{nuclearpotential}

Let  the quantum Wilson line $ W(z_0,z_1)W(z_1,z)$ be as a quantum golun propagator. Let
\begin{equation}
{\bf p}_i:=W(K_i)Z_i
\label{baryon}
\end{equation}
$i=1,2$ denotes two protons. Then the quantum golun propagator $W(z_0,z_1)W(z_1,z)$ represents a nuclear force on the interaction of ${\bf p}_i, i=1,2$ by the following formula:
\begin{equation}
{\bf p}_1^{*}  W(z_0,z_1)W(z_1,z){\bf p}_2
\label{baryon2}
\end{equation}
Similarly  let
\begin{equation}
{\bf n}_i:=W(K_i)Z_i
\label{baryons}
\end{equation}
$i=1,2$ denotes two neutrons. Then the quantum golun propagator $ W(z_0,z_1)W(z_1,z)$ represents a nuclear force on the interaction of ${\bf n}_i, i=1,2$ by the following formula:
\begin{equation}
{\bf n}_1^{*}   W(z_0,z_1)W(z_1,z){\bf n}_2
\label{baryon3}
\end{equation}

Then the quantum golun propagator $ W(z_0,z_1)W(z_1,z)$ represents a nuclear force on the interaction of ${\bf p}$ and ${\bf n}$ by the following formula:
\begin{equation}
{\bf n}^{*}   W(z_0,z_1)W(z_1,z){\bf p}
\label{baryon4}
\end{equation}


For the  quantum gloun propagator, for simplicity we consider the $SU(2)\otimes U(1)$ group which is as a subgroup of  $SU(3)\otimes U(1)$.
 We have:
\begin{equation}
\hat{t}:=\frac{e_0^2}{k_0+g_0}\sum_a t^a \otimes t^a
\label{t1}
\end{equation}
where $ t^a=iT^a$ and:
\begin{equation}
T^1= \left ( \begin{array}{cc}
             0 & 1 \\
             1 & 0 \end{array}\right), \qquad
T^2= \left ( \begin{array}{cc}
             0 & i \\
            -i & 0 \end{array}\right), \qquad
T^3= \left ( \begin{array}{cc}
             1 & 0 \\
             0 & -1 \end{array}\right)
\end{equation}
 Let us call this Casimir opeartor $T$ as the mass operator for the $\pi$ mesons.
We have
\begin{equation}
 T=\sum_{a=1}^{3} T^a \otimes T^a =
\left ( \begin{array}{cccc}
             1 &  0 &  0 & 0 \\
             0 & -1 &  2 & 0 \\
             0 &  2 & -1 & 0 \\
             0 &  0 &  0 & 1  \end{array}\right)
\label{casimor} 
\end{equation}
This Casimir operator has eigenvalues $1$ and $-3$ where $1$ is with
multiplicity $3$. Since $1$ is positive it is considered as an eigenvalue for mass and energy. The negative eigenvalue $-3$ is considered to be nonphysical. It can be checked  that
\begin{equation}
v_1=\left ( \begin{array}{c}
             1  \\
             0  \\
             0  \\
             0   \end{array}\right),
\quad
v_3= \left ( \begin{array}{c}
              0  \\
              \frac12 \\
             \frac12  \\
              0  \end{array}\right),
\quad
v_4=\left ( \begin{array}{c}
             0 \\
             0 \\
              0 \\
             1  \end{array}\right),
\label{casimir1a}
\end{equation} 
are  three independent eigenvectors for the eigenvalue $1$.  We show in \cite{Ng} that the igenvalue $1$  gives the mass $135 Mev$ for the $\pi^0$ meson.

Then an eigenvector of $-3$ is given by:
\begin{equation}
v_2= \left ( \begin{array}{c}
              0  \\
              \frac12 \\
             \frac12  \\
              0  \end{array}\right)
\label{casimir1c}
\end{equation} 

We can write $T$ in the diagonalized form:
\begin{equation}
 T=
\left ( \begin{array}{cccc}
             v_1 &  v_2  &   v_3 &  v_4 \end{array}\right)
\left ( \begin{array}{cccc}
             1 &  0 &  0 & 0 \\
             0 &  -3 &  0 & 0 \\
             0 &   0 & 1 & 0 \\
             0 &  0 &  0 & 1  \end{array}\right)
 \left ( \begin{array}{cccc}
             v_1 &  v_2  &   v_3 &  v_4 \end{array}\right)^{-1}
=:VDV^{-1}
\label{casimir1d}
\end{equation}

As the photon case, the  eigenvalue $1$  for $ r_1 > n_e^2$ gives repulsive effect for the interacting elementary particles such as ${\bf p}$ and ${\bf n}$.  This gives the effect of the $\pi^0$ meson for the nuclear force induced by the  $\pi^0$ meson. Thus the quantum gloun propagator includes the repulsive effect of  the scattering  of nucleons such as ${\bf p}$, ${\bf n}$ by the  exchange of the  $\pi^0$ meson.

Then  eigenvalue $-3$  for $ r_1 > n_e^2$ gives attractive effect for the interacting elementary particles such as ${\bf p}$ and ${\bf n}$. Thus the quantum gloun propagator includes the attractive and repulsive effect of  the scattering  of nucleons such as ${\bf p}$, ${\bf n}$.  
Statistically, when  $ r_1 > n_e^2$ and $ r_1- n_e^2$ is small,  the attractive and repulsive effects are respectively given by:
\begin{equation}
 \frac{3e^2}{r}  \quad  \mbox{and} \quad  \frac{-e^2}{r}   
\label{attractive3}
\end{equation}
where $r$ is the usual three dimensional space distance.
Then since $|-3|>1$ we have that the attractive effect by $-3$ is greater than the repulsive effect by $1$.  Thus the nuclear interaction of nucleon such as ${\bf p}$, ${\bf n}$ is  attracive in the forming of nuclei since this is the case of  $ r_1 > n_e^2$ and $ r_1- n_e^2$ is small.

Then we notice that when $ r_1 <n_e^2$ we have that the interacting potential of the  quantum gloun propagator of the eigenvalue $-3$ gives repulsive effect which  is greater than the attractive effect of the eigenvalue $1$. Thus when the distance between the neucleons is very small  (or $ r_1 <n_e^2$) we have that the force between the neuclons is repulsive. This agrees with the experiments that the  neuclear force between the neuclons is repulsive when the distance between the neuclons is very small.


On the other hand when $r_1$ is large such that $ r_1 >>n_e^2$ 
 we have that the effect of quantum gloun propagator of the eigenvalue $-3$ is much smaller than the effect of quantum gloun propagator of the eigenvalue $1$ since we have:
\begin{equation}
e^{-3\log \frac{r_1}{n_e^2} } \to 0 \quad \mbox{as} \quad r_1\to \infty.
\label{attractive}
\end{equation}
Thus  the net effect of  the quantum gloun propagator 
 is repulsive when  $r_1$ is large.  This agrees with the experiments of the neuclear force.

Then because of the mass term of the quantum gloun propagator is larger than that of the  quantum photon propagator we have that the repulsive effect of the quantum gloun propagator is of short range  while the repulsive effect of the electromagetic interaction is of long range. Thus the quantum gloun propagator gves the nuclear force for the scattering interaction of nucleons and for the forming of nuclei.


The  quantum gloun propagator gives the nuclear force  for the scattering interaction of nucleons and for the forming of nuclei. Thus this  quantum gloun propagator can be used to replace the Paris potential for the scattering interection of nucleons.

Let us consider the knot model of an elementary particle such as the following knot model of proton:
\begin{equation}
{\bf p}=W(K)Z(s)
\label{asymptotic}
\end{equation}
where the quark components of $Z(s)$ are at proper time $s$ forming the vector $Z(s)$. In this case the quantum knot $W(K)$ is  a gloun propagator in closed form with $z_0=z$.  Because of the closed form we have that  the net interaction effect from the gloun propagator is zero.  Thus  the net interaction effect of the quark components of $Z(s)$ is zero. Thus this knot model of proton  gives the asymptotic freedom property of strong interaction of the quarks of  the proton ${\bf p}$.


Let us then derive (and define)  a classical gloun  propagator from the
quantum gloun propagator.  As the photon case let us choose a path
connecting $z_0$ and $z$. Let us consider the  path  (\ref{g1g}).
From this chosen of the
path (\ref{g1g}) we have that the quantum gloun propagator is
proportional to the following expression:
\begin{equation}
V\frac{1}{2\lambda_1}[\theta(s-s_1)e^{-i D_1(s-s_1)}+
\theta(s_1-s)e^{i D_1(s-s_1)}]V^{-1}
 \label{ga3}
\end{equation}
where we define $D_1=-D e_0^2\beta_1$.

\section{Renormalization group method for critical phenomena of QED} \label{reQED}

In this section and the following sections we shall show that the QED model (\ref{1.1}) with the parameter $\frac{1}{ih+\kappa}$  is a
model of superconductivity and magnetism. We shall continue to show that the Planck constant parameter $h\neq 0$ basically gives dynamical wave effects and conductivity while the parameter $\kappa>0$ basically gives magnetic effects. For this QED model (\ref{1.1}) to be a
model of superconductivity and magnetism we introduce many variables $Z_a$ for electrons indexed by $a$ and many variables $A_{1b},
A_{2b}$ for photons indexed by $b$  which are coupled together with the following seagull vertex:
\begin{equation}
-\frac{1}{ih+\kappa}e^2 Z_a^{*}Z_a A_{b}^2
\label{u1}
\end{equation}
and the usual vertex for QED (giving the Coulomb potential) in (\ref{1.1}) with a single
$\frac{dZ_a}{ds}$.
We shall show that the coupling term
(\ref{u1}) gives critical phenomena which gives the phenomena of superconductivity and magnetism.

From the interaction term
(\ref{u1}) the one-loop mass-energy term of electron and photon
are respectively given by:
\begin{equation}
\begin{array}{rl}
&\frac{e^2}{(ih+\kappa)}\frac{\frac{1}{ih+\kappa}}{2\pi}\int\frac{dk}{\frac{1}{(ih+\kappa)^2}k^2
+\lambda_1^2}
=\frac{e^2}{(ih+\kappa)}
\frac{(ih+\kappa)}{2\pi}\int\frac{dk}{k^2+(ih+\kappa)^2\lambda_1^2}
=
\frac{e^2}{2(ih+\kappa)\lambda_1}
=:\frac{\omega^2}{(ih+\kappa)}
=:\omega^2(-i\bar{h}+\bar{\kappa})=:r \label{u2a}
\end{array}
\label{u2b1}
\end{equation}
and
\begin{equation}
\begin{array}{rl}
&\frac{e^2}{(ih+\kappa)}\frac{\frac{1}{ih+\kappa}}{2\pi}\int\frac{dq}{\frac{1}{(ih+\kappa)^2}q^2+\mu^2}
= \frac{e^2}{(ih+\kappa)}
\frac{(ih+\kappa)}{2\pi}
\int\frac{dq}{q^2+(ih+\kappa)^2\mu^2}
=
\frac{e^2}{2(ih+\kappa)\mu}
=:\frac{\lambda_2^2}{(ih+\kappa)}=:-i\bar{\lambda}_2^2\bar{h}+
\bar{\lambda}_2^2\bar{\kappa}
\end{array}
\label{u2b}
\end{equation}
In terms of Feymann diagram the above formula (\ref{u2a}) shows that a photon line in the forming of a closed loop gives a photon propagator of the form $\frac{ih+\kappa}{k^2+(ih+\kappa)^2\lambda_1^2}$. Also the above formula (\ref{u2b}) shows that an electron line in the forming of a closed loop gives an electron propagator of the form $\frac{ih+\kappa}{q^2+(ih+\kappa)^2\mu^2}$. In these two formulae each $e^2$ gives a factor $\frac{1}{ih+\kappa}$.

Similarly for the two-loop case (with one photon loop and one electron
loop where the photon loop consists of two photon lines and the electron
loop consists of one electron line) constructed by (\ref{u1})  the mass-energy term of electron
is given by:
\begin{equation}
\begin{array}{rl}
&-\frac{e^4}{(ih+\kappa)^2}\frac{(ih+\kappa)^3}{(2\pi)^2}
\int\frac{dkdq}{(k^2+(ih+\kappa)^2\lambda_1^2)^2
(q^2+(ih+\kappa)^2\mu^2)}
=
-\frac{e^4}{8(ih+\kappa)^3\lambda_1^3\mu}
=-\frac{e^4((\bar{\kappa}^2-\bar{h}^2)-2i\bar{\kappa}\bar{h})(-i\bar{h}+\bar{\kappa})}
{8\lambda_1^3\mu}
\end{array}
\label{u3}
\end{equation}
where the factor $(ih+\kappa)^3$ is from the two photon lines and one electron line.

Then for the three-loop case (with two photon loops and one electron
loop where one photon loop consists of two photon lines and the electron
loop consists of two electron lines while the other photon loop consists of one photon line, as the case of (\ref{u2a})) constructed by (\ref{u1}) we have that the mass-energy term
is given by:
\begin{equation}
\begin{array}{rl}
&\frac{e^6}{(ih+\kappa)^3}\frac{(ih+\kappa)^5}{(2\pi)^3}
\int\frac{dkdqdl}
{(k^2+(ih+\kappa)^2\lambda_1^2)^2(q^2+(ih+\kappa)^2\mu^2)^2(l^2+(ih+\kappa)^2\lambda_1^2)}
\\
&\\
= &\frac{e^6}{32(ih+\kappa)^5\lambda_1^4\mu^3}
=\frac{e^6[(\bar{\kappa}^3-3\bar{\kappa}\bar{h}^2)-i(3\bar{\kappa}^2\bar{h}-\bar{h}^3)](-ih+\kappa)^2}
{32\lambda_1^4\mu^3}
\end{array}
\label{u4}
\end{equation}
where the factor $(ih+\kappa)^5$ is from the the three photon lines and the two electron lines.

Then we consider the four point vertex with four external electron
lines. From (\ref{u1}) the one-photon-loop interaction term
is given by:
\begin{equation}
\begin{array}{rl}
&I_1(p_1+p_2) = \frac{e^4(ih+\kappa)^2}{(ih+\kappa)^2(2\pi)} \int\frac{dk}{(k^2+(ih+\kappa)^2\lambda_2^2)((p_1+p_2-k)^2+(ih+\kappa)^2\lambda_2^2)}\\
\end{array}
\label{u5}
\end{equation}
where $p_1, p_2$ are the momentums of two external electron lines
respectively, and the two photon propagators are obtained by the Dyson geometry series.

Setting $p_1, p_2 =0$ we get an interaction constant
$u$ given by
\begin{equation}
\begin{array}{rl}
u & := 2I_1(0)= \frac{e^4 }{2(ih+\kappa)^3\lambda_2^3}
=\frac{e^4[(\bar{\kappa}^3-3\bar{}\bar{h}^2)-i(3\bar{\kappa}^2\bar{h}-\bar{h}^3)]}{2\lambda_2^3}
\end{array}
\label{u6}
\end{equation}
where the coefficient 2 comes from permutation of $I_1(p_1+p_2)$.

Similarly from (\ref{u1}) the three-loop interaction term (two
photon-loops and one electron-loop) of the four point vertex is
given by
\begin{equation}
\begin{array}{rl}
&I_2(p_1+p_2)
=\frac{e^8(ih+\kappa)^6}{(ih+\kappa)^4(2\pi)^3}\int\frac{dkdqdl}
{(k^2+(ih+\kappa)^2\lambda_2^2)}
\frac{1}
{((p_1+p_2-k)^2+(ih+\kappa)^2\lambda_2^2)(q^2+(ih+\kappa)^2\omega^2)}\times\\
&\\
&\times\frac{1}{((p_1+p_2-q)^2+(ih+\kappa)^2\omega^2)((p_1+p_2-q)^2+(ih+\kappa)^2\omega^2)}
\frac{1}{
(l^2+(ih+\kappa)^2\lambda_2^2)((p_1+p_2-l)^2+(ih+\kappa)^2\lambda_2^2)}
\end{array}
\label{u7}
\end{equation}
where the photon propagators and electron propagators are obtained by the Dyson geometry series.
Setting $p_1$,$p_2 =0$ we have the following interaction constant
which will be added to (\ref{u6}):
\begin{equation}
2I_2(0)=\frac{e^8}{32(ih+\kappa)^7\lambda_2^6\omega^3} \label{u8}
\end{equation}
By re-scaling $\lambda_1, \mu$ to $\lambda_1^{\prime},
\mu^{\prime}$ from
(\ref{u2a})-(\ref{u8}) we re-scale $r$ to $r'$ and $u$ to $u'$
given by
\begin{equation}
\begin{array}{rl}
r'&=
\frac{e^2}{2(ih+\kappa)\lambda_1^{\prime}}-\frac{e^4}{8(ih+\kappa)^3\lambda_1^{\prime
3} \mu^{\prime}}
    +\frac{e^6}{32(ih+\kappa)^5\lambda_1^{\prime 4} \mu^{\prime 3}}
\\
u'&= \frac{e^4}{2(ih+\kappa)^3\lambda_2^{\prime 3}}
   +\frac{e^8}{32(ih+\kappa)^7\lambda_2^{\prime 6}\omega^{\prime 3}}
\end{array}
\label{u10}
\end{equation}
From (\ref{u10}) we get the following renormalization-group
recursion equation:
\begin{equation}
\begin{array}{rl}
r'&= b^2 r-Ab^{6-\frac23 \epsilon}r^3 u^{-\frac23}
    +Cb^{8-2\epsilon}r^4 u^{-2} \\
u'&=b^{\epsilon}u + Bb^{2\epsilon-3}u^2 r^{-\frac32}
\end{array}
\label{u11}
\end{equation}
where we define
\begin{equation}
\frac{\lambda_1}{\lambda_1^{\prime}}=b^2 \quad \mbox{and}\quad
\frac{\mu^{\prime}}{\mu}=b^{\frac23 \epsilon} \label{u12}
\end{equation}
where $b$ is regarded as a scale factor and $\epsilon >0$ is a
parameter
which is analogous to the
$\epsilon$ in the $\epsilon$-expansion of the
Ginzburg-Landau-Wilson (GLW) model in statistical physics
\cite{Wilson}-\cite{Onsager};
and
\begin{equation}
A=\frac{2^{\frac13}}{e^{\frac43}(ih+\kappa)^2}, \quad B=
\frac{1}{8(ih+\kappa)^{\frac{5}{2}}}, \quad
C=\frac{1}{(ih+\kappa)^7} \label{u13}
\end{equation}

As analogous to the GLW model the parameter $\epsilon$ is related
to the dimension of the space \cite{Wilson}-\cite{Onsager}.
Thus as a space-time statistics we consider that $\epsilon$ is
related to the space-time statistics that $d=3-\epsilon$ is the
(fractional) dimension of the space, which is related to the
ultraviolet limit, as we have considered in the above sections;
while the parameter $b$ is related to the infrared limit.

 We remark that
(\ref{u11}) is analogous to the
renormalization-group equation of the Ginzburg-Landau-Wilson model
in statistical physics. We may regard the GLW model as an
approximation of this QED statistical model.

Following the renormalization group method in statistical physics we set up
the following equation to find fixed points of (\ref{u11}):
\begin{equation}
\begin{array}{rl}
 r&= b^2 r-Ab^{6-\frac23 \epsilon}r^3 u^{-\frac23}
+Cb^{8-2\epsilon}r^4 u^{-2} \\
u&= b^{\epsilon}u +Bb^{2\epsilon-3}u^2 r^{-\frac32}
\end{array}
\label{u14}
\end{equation}
The nontrivial fixed point is given by:
\begin{equation}
\begin{array}{rl}
r^{*}& = B^{-\frac23}\frac {(1-b^{\epsilon})^{\frac23}
[b^2-1+\frac{CB^2b^{2+2\epsilon}}{(1-b^{\epsilon})^2}]}
{Ab^{4+\frac23 \epsilon}}
\\
u^{*}& = B^{-2}\frac {(1-b^{\epsilon})^2 [b^2-1
+\frac{CB^2b^{2+2\epsilon}}{(1-b^{\epsilon})^2}]^{\frac32}}
{A^{\frac32}b^{3+3\epsilon}}
\end{array}
\label{u15}
\end{equation}
Linearizing  the renormalization-group transformation  near
the fixed point (\ref{u15}) we get the linearized transformation
matrix. We consider two cases:

{\bf Case 1)}.
$\epsilon \ll 1$ and $b$ is large such that $\epsilon\log b$
is very small;

{\bf Case 2)}. $\epsilon=1$, and $b<1$, $b\to 1$ such that $\epsilon|\log
b|=|\log b|$ is very small.

For these two cases the eigenvalues of this matrix are
approximately given by:
\begin{equation}
\begin{array}{rl}
\lambda_{R1} & \approx b^2[-2+
CB^2b^{2\epsilon}(b^{\epsilon}-1)^{-2}]\\
&\\
            & \approx b^2[-2+CB^2(3-2\epsilon \mbox{ln}b]\\
            &\\
\lambda_{R2} & \approx 0
\end{array}
\label{u16}
\end{equation}

Let us first consider the Case 1).  From (\ref{u16}) a rough
estimate for $\lambda_{R1}
>1$ is:
\begin{equation}
-2+ \mbox{Re}CB^2 >0
\label{u17}
\end{equation}
and
\begin{equation}
\mbox{Im}CB^2 =0 \label{u17c}
\end{equation}

From (\ref{u16}) we have that when (\ref{u17}) and (\ref{u17c})
hold the critical exponent $\nu$ of correlation length is given
by:
\begin{equation}
\begin{array}{rl}
\nu = & \frac{\ln b}{\ln\lambda_{R1}}\approx\frac12
+\frac{2\mbox{Re}CB^2\epsilon}{2(-2+3\mbox{Re}CB^2)}\\
&\\
= &\frac12 +\frac{\epsilon}{2}
\end{array}
\label{u19}
\end{equation}
Thus for this Case 1) (which means that the dimension of the
space is 3) we have that $\frac12 <\nu <1$.

Let us then consider the Case 2). For this case the dimension of
space is $3-\epsilon=2$. From (\ref{u16}) a rough estimate for
$\lambda_{R1}
>1$ is:
\begin{equation}
\mbox{Re}CB^2 >0 \label{u17a}
\end{equation}
and
\begin{equation}
\mbox{Im}CB^2 =0 \label{u17b}
\end{equation}

From (\ref{u16}) we have that when (\ref{u17a}) and (\ref{u17b})
hold the critical exponent $\nu$ of correlation length is given
by:
\begin{equation}
\begin{array}{rl}
\nu = & \frac{\ln b}{\ln\lambda_{R1}} \approx \frac12
+\frac{2\mbox{Re}CB^2\epsilon}{4\mbox{Re}CB^2}\\
&\\
= &\frac12 +\frac{\epsilon}{2}=1
\end{array}
\label{u20}
\end{equation}
Thus for this  Case 2) (which means that the dimension of the
space is 2) we have that $\nu =1$.
It is interesting that this is just the result given by
the exact solution of the two dimensional Ising model
\cite{Huang,Onsager}. 

\section{Renormalization group method for critical phenomena of QCD} \label{reQCD}

In this section we consider the QCD model (\ref{n1}) with the parameter $\frac{1}{ih+\kappa}$. We shall continue to show that the Planck constant parameter $h\neq 0$ basically gives dynamical wave effects and conductivity while the parameter $\kappa>0$ basically gives magnetic effects. 

For this QCD model (\ref{n1}), as the  QED model (\ref{1.1}),  we introduce many variables $Z_a$ for quarks  indexed by $a$ and many variables $A_{1b},
A_{2b}$ for gloun indexed by $b$  which are coupled together with the following seagull vertex:
\begin{equation}
-\frac{1}{ih+\kappa}e^2 Z_a^{*}Z_a A_{b}^2
\label{uc1}
\end{equation}

For the component of  gloun propagator with respect to the eigenvalue $1$, the result of the renormalization group method is exactly the same as the QED case.
Let us then consider  the component of  gloun propagator with respect to the eigenvalue $-3$. In this case the classical  component of  gloun propagator is given by:
\begin{equation}
\frac{i}{k_E^2-3^2\lambda_1^2} \label{2g5}
\end{equation}
when $h=1, \kappa \to 0$.

Then, for general $h, \kappa $,  as the QED case, we have the  following  
 renormalization group equation:
\begin{equation}
\begin{array}{rl}
r'&= b^2 r- \frac{1}{3^2}Ab^{6-\frac23 \epsilon}r^3 u^{-\frac23}
    + \frac{1}{3^3} Cb^{8-2\epsilon}r^4 u^{-2} \\
u'&=b^{\epsilon}u + Bb^{2\epsilon-3}u^2 r^{-\frac32}
\end{array}
\label{uc11}
\end{equation}
where $r$ is defined in (\ref{u2b1}) and the equation of $r'$ is obtained by multiplying $3$ to get the right hand side of the  equation of $r'$.

Thus for the QCD case we have two renormalization group equations ( \ref{uc11}) and (\ref{u11}).

Then we set up
the following equation to find fixed points of (\ref{uc11}):
\begin{equation}
\begin{array}{rl}
 r&= b^2 r- \frac{1}{3^2} Ab^{6-\frac23 \epsilon}r^3 u^{-\frac23}
+ \frac{1}{3^3} Cb^{8-2\epsilon}r^4 u^{-2} \\
u&= b^{\epsilon}u +Bb^{2\epsilon-3}u^2 r^{-\frac32}
\end{array} 
\label{uc14}
\end{equation}
The nontrivial fixed point is given by:
\begin{equation}
\begin{array}{rl}
r^{*}& = B^{-\frac23}\frac {(1-b^{\epsilon})^{\frac23}
[b^2-1+\frac{ 3^{-3} CB^2b^{2+2\epsilon}}{(1-b^{\epsilon})^2}]}
{3^{-2} Ab^{4+\frac23 \epsilon}}
\\
u^{*}& = B^{-2}\frac {(1-b^{\epsilon})^2 [b^2-1
+\frac{ 3^{-3} CB^2b^{2+2\epsilon}}{(1-b^{\epsilon})^2}]^{\frac32}}
{ 3^{-3} A^{\frac32}b^{3+3\epsilon}}
\end{array}
\label{uc15}
\end{equation}
Linearizing  the renormalization-group transformation  near
the fixed point (\ref{uc15}) we get the linearized transformation
matrix. We consider two cases:

{\bf Case 1)}.
$\epsilon \ll 1$ and $b$ is large such that $\epsilon\log b$
is very small;

{\bf Case 2)}. $\epsilon=1$, and $b<1$, $b\to 1$ such that $\epsilon|\log
b|=|\log b|$ is very small.

For these two cases the eigenvalues of this matrix are
approximately given by:
\begin{equation}
\begin{array}{rl}
\lambda_{R1} & \approx b^2[-2+
 \frac{1}{3^3} CB^2b^{2\epsilon}(b^{\epsilon}-1)^{-2}]\\
&\\
            & \approx b^2[-2+ \frac{1}{3^3} CB^2(3-2\epsilon \mbox{ln}b]\\
            &\\
\lambda_{R2} & \approx 0
\end{array}
\label{uc16}
\end{equation}

Let us first consider the Case 1).  From (\ref{u16}) a rough
estimate for $\lambda_{R1}
>1$ is:
\begin{equation}
-2+ \mbox{Re} \frac{1}{3^3} CB^2 >0
\label{uc17}
\end{equation}
and
\begin{equation}
\mbox{Im} \frac{1}{3^3} CB^2 =0 \label{uc17c}
\end{equation}


Let us then consider the Case 2). For this case the dimension of
space is $3-\epsilon=2$. From (\ref{u16}) a rough estimate for
$\lambda_{R1}
>1$ is:
\begin{equation}
\mbox{Re} \frac{1}{3^3} CB^2 >0 \label{uc17a}
\end{equation}
and
\begin{equation}
\mbox{Im} \frac{1}{3^3} CB^2 =0 \label{uc17b}
\end{equation}

Thus for QCD case, we have two nontrivial fixed points. This gives some special properties on phase diagram of QCD.

\section{Ferromagnetism,  Bose-Einstein condensation and  supeconductivity of QED }

Let us then investigate the phenomena of superconductivity and
magnetism by the above renormalization group method. Let us
first consider the Case 1) with $0<\epsilon <<1$ for the three dimensional space statistics. For satisfying the inequality (\ref{u17}) of this
case from the formula for $B$ and $C$ we have $h=0$ or
$\kappa=0$. Let us first consider the case $h=0$. For this subcase
we have that the equation (\ref{u17c}) also trivially holds. Then
the inequality  (\ref{u17}) is simplified to the following
equality:
\begin{equation}
-2+\frac{1}{64\kappa^{12}} >0 \label{u39}
\end{equation}
This gives an estimate of the critical value $\kappa_c$ of
$\kappa$:
\begin{equation}
2^{-\frac{7}{12}} >\kappa_c \label{u39a}
\end{equation}
Let us identify the constant $\kappa$  as a constant proportional to the absolute temperature $T$ such that $\kappa=k_BT$ where $k_B$ denotes the Baltzmann constant.
Since $h=0$ gives the ground state of electrons we have that $\kappa>0$ gives the statc magnetic property of the ground state of electrons. Moreover since $h\neq 0$ gives dynamical effect we have that $h=0$ and $\kappa>0$ gives the static magnetic property. Thus this subcase can be
identified as a static magnetic phase transition where $\kappa_c=k_BT_c$ and $T_c$
is as the critical temperature for the magnetic phase
transition.

Moreover from the winding numbers of the photon loops we have that the sign $\pm$ of $\pm h\neq 0$ gives the two opposite spins of photons and of electrons. Thus the dynamical interaction of spins is from $\pm h\neq 0$. Then when $h=0$ there is no dynamical interaction among the spins. In this case the antiferromagnetic phase does not appear since antiferromagnetism is formed by the dynamical interaction of two spins with opposite signs. Thus this static magnetic phase
transition given by this case $h=0$ and $\kappa>0$ is the ferromagnetic phase transition where the paramagnetic phase is transited to the ferromagnetic phase when the temperature $T$ is decreased from $T>T_c$ to $T<T_c$.

Let us then consider the subcase that $h\neq 0$ and $\kappa=0$. We shall show that this
subcase can be identified as the Bose-Einstein condensation and is related to the Type I superconductivity. 
For this subcase $h \neq 0$ and $\kappa=0$,
 the equation (\ref{u17c}) always
holds. Then the inequality  (\ref{u17}) is simplified to the
following equality:
\begin{equation}
-2+\frac{1}{64 h^{12}} >0 \label{u40}
\end{equation}
This gives a condition for the critical value $h_c$ of $h$:
\begin{equation}
2^{-\frac{7}{12}} >h_c \label{u41}
\end{equation}
Then
since the critical temperature $\kappa_c=\kappa= 0$ the  following subcase: 
\begin{equation}
2^{-\frac{7}{12}} >h_c>0, \quad  \kappa_c= 0
\label{uu41}
\end{equation}
can be identified as the Bose-Einstein condensation.

 Let us extend this case   (\ref{uu41})
 to that $\kappa \approx 0$ but is not equal to zero.
Let us find out such a range for $h\neq 0$ and $\kappa \approx 0$. This range is then as the 
 phase of (Type I) superconductivity. First we must have that $h>\kappa$. Then we have:
\begin{equation}
\begin{array}{rl}
& CB^2=2^{-6}\frac{1}{(ih+\kappa)^{12}}= 2^{-6} 
 [(\bar{\kappa}^3-
3\bar{\kappa}\bar{h}^2 + 3\bar{\kappa}^2\bar{h}-\bar{h}^3)^2 \times \\
&\\
& \times (\bar{\kappa}^3-3\bar{\kappa}\bar{h}^2-3\bar{\kappa}^2\bar{h}+\bar{h}^3)^2
 -4\bar{\kappa}^2\bar{h}^2(\bar{h}^2-3\bar{\kappa}^2)^2(\bar{\kappa}^2-3\bar{h}^2)^2 ]\\
&\\
&
+ i2^{-6}\cdot 4\bar{\kappa}\bar{h}(\bar{\kappa}^3-
3\bar{\kappa}\bar{h}^2+3\bar{\kappa}^2\bar{h}-\bar{h}^3)
 (\bar{\kappa}^3-3\bar{\kappa}\bar{h}^2-3\bar{\kappa}^2\bar{h}+\bar{h}^3)
(\bar{h}^2-3\bar{\kappa}^2)(\bar{\kappa}^2-3\bar{h}^2)
\end{array}
\label{ssu42s}
\end{equation}
Thus for (\ref{u17b}) holds we must have the condition $\kappa=0$ or the following condition:
\begin{equation}
h^2=3\kappa^2
\label{supercond}
\end{equation}
The condition (\ref{supercond}) will be as the basic state of superconductivity.

These two conditions are of the same kind for (\ref{u17b}) holds.
Then when $h_c$ is very small we have that these two states $\kappa_c=0$ and $h_c^2=3\kappa_c^2$ are closed to each other. Thus when $h_c$ is very small, the range for $\kappa_c\approx 0$ is $0\leq \kappa_c\leq 3^{-\frac12}h_c$ where
$h_c<2^{-\frac{7}{12}}$.
Thus the range of the phase of  superconductivity (which will be of Type I) is included in the following range:
\begin{equation}
0\leq \kappa_c\leq \frac{1}{\sqrt{3}}h_c \quad \mbox{and} \quad  0<h_c<2^{-\frac{7}{12}}
\label{type1}
\end{equation}
where $\kappa_c=k_B T_c$ and $T_c$ is as the critical temperature for (Type I) superconductivity. We shall later find out  
 the exact range of Type I superconductivity.


For a fixed $\kappa_c>0$, the condition $h^2=3\kappa_c^2$ will be as the basic state of superconductivity and the condition $h^2>3\kappa_c^2$ will be as the normal metallic state. Then we shall show that the Cooper pair of electrons with different spins is as a mechanism of reducing $h$ for the transition from the normal metallic state to the state of superconductivity.
Briefly, in the normal metallic state the two electrons of a Cooper pair are described by two phases $e^{\pm ih}$. Then the different sign of $\pm ih$ gives the effect of reducing $h$ for the transition from the normal metallic state to the state of superconductivity. We shall give more details of the mechanism of Cooper pair for reaching the state of superconductivity.

We remark that it may be the case that the state $h^2=3\kappa_c^2$ can be reached without the forming of  Cooper pairs. In this case the state $h^2=3\kappa_c^2$ gives phenomenon which is similar but is different from the usual phenomenon of superconductivity with the Cooper pair mechanism. In this paper we shall only investigate the case of the usual superconductivity with the Cooper pair mechanism.


Thus when $0<h_c<2^{-\frac{7}{12}}$ and $h_c^2=3\kappa_c^2$ (such that $h_c$ and $\kappa_c$ are very small) we have the following formula for the critical temperature $T_c$:
\begin{equation}
k_BT_c=\kappa_c=\frac{1}{\sqrt{3}}h_c=:\frac{1}{\sqrt{3}}\Delta(0)=:\frac{1}{\sqrt{3}}\bar{h}_c\omega_D
\label{ssu42s1}
\end{equation}
where $\Delta(0):=h_c$ is as the energy gap of a Cooper pair of electrons and $\omega_D$ denotes the Debye frequency of the phonons and $\bar{h}_c$ denotes a value of the Planck constant parameter $h$ such that $\bar{h}_c\omega_D=h_c$. Thus this formula (\ref{ssu42s1}) represents the critical temperature $T_c$ in terms of the energy gap $\Delta(0)$.

We notice that when $h_c$ is further restricted such that $\bar{h}$ is the usual Planck constant (devided by $2\pi$) this formula (\ref{ssu42s1}) of the critical temperature $T_c$ is of the same form as the formula of the critical temperature $T_c$ in the BCS theory. Here the coefficient $\frac{1}{\sqrt{3}}=\frac{2}{3.4641}$ is close to the corresponding coefficient $\frac{2}{3.5}$ of the formula of the critical temperature $T_c$ in the BCS theory \cite{Bar,Abr,Bog,Hub,Gor,Tin}.

As the BCS theory because $\omega_D$ is proportional to $M^{-\frac12}$ where $M$ denotes the mass of the isotope of the material in the phase of superconductivity, the formula (\ref{ssu42s1}) gives the isotope effect \cite{Bar,Abr,Bog,Hub,Gor,Tin}.

\section{Conventional Type I and Type II superconductivity of QED}

Let us then consider the subcase that $h\neq 0$ and $\kappa>0$ of the Case 1) with the three dimensional space statistics ($0<\epsilon <<1$).
We shall show that this subcase gives the conventional Type II superconductivity (When $\kappa>0$ is very small,  this subcase also gives the Type I superconductivity).
For this subcase $h\neq 0$ and $\kappa>0$, from (\ref{ssu42s}), the condition (\ref{u17}) can be written in the following form:
 \begin{equation}
 \begin{array}{rl}
& 
[(\bar{\kappa}^3-
3\bar{\kappa}\bar{h}^2+3\bar{\kappa}^2\bar{h}-\bar{h}^3)^2
(\bar{\kappa}^3-3\bar{\kappa}\bar{h}^2-3\bar{\kappa}^2\bar{h}+\bar{h}^3)^2 \\
&\\
&
-4\bar{\kappa}^2\bar{h}^2(\bar{h}^2-3\bar{\kappa}^2)^2(\bar{\kappa}^2-3\bar{h}^2)^2 ]>2\cdot  2^{6}
\end{array}
 \label{u17hi}
\end{equation}
Then from (\ref{ssu42s}),  the condition (\ref{u17b}) can be written in the following form:
 \begin{equation}
(\bar{h}^2-3\bar{\kappa}^2)(\bar{\kappa}^2-3\bar{h}^2)=0
 \label{su43}
\end{equation}
Then from (\ref{su43}) we have
that the inequality (\ref{u17}) (or (\ref{u17hi}))  is of the following form:
\begin{equation}
(\bar{\kappa}^3-3\bar{\kappa}\bar{h}^2+3\bar{\kappa}^2\bar{h}-\bar{h}^3)^2
(\bar{\kappa}^3-3\bar{\kappa}\bar{h}^2-3\bar{\kappa}^2\bar{h}+\bar{h}^3)^2>2^{7}
 \label{su42}
\end{equation}
where $0<\bar{\kappa}<1$ and $0<\bar{h}<1$.
Then from (\ref{su43}) we have the following two subcases:
\begin{equation}
h^2-3\kappa^2=0
 \label{su42c}
\end{equation}
\begin{equation}
\kappa^2-3h^2=0
 \label{su42d}
\end{equation}

Since the subcase (\ref{su42c}) gives $h^2>\kappa^2$, this subcase (\ref{su42c}) is identified as the state of the conventional Type II (and Type I) superconductivity.
On the other hand the subcase (\ref{su42d}) gives $h^2<\kappa^2$.
Thus this subcase is for magnetic phase transition, as the above case of ferromagnetism. Now since $h^2\neq 0$ gives the dynamical interaction of up and down spins specified by the $\pm$ of $\pm ih$, we have that this subcase can be identified as the state of the (conventional) antiferromagnetism. 

For the subcase (\ref{su42c}), the inequality (\ref{su42}) is simplified to the following  inequality:
\begin{equation}
\begin{array}{rl}
2^{-6}2^{12}\bar{\kappa}^{12}=\frac{2^{6}}{2^{24}\kappa^{12}}= \frac{1}{2^{18}\kappa^{12}}>2
\end{array}
 \label{su42e}
\end{equation}
This gives the condition on the parameter $\kappa$ (and the parameter $h^2=3\kappa^2$)
 for the phase transition of the conventional Type II  (and Type I) superconductivity. We can write this condition in the following form:
\begin{equation}
\begin{array}{rl}
\frac{\sqrt{3}}{2}\cdot 2^{-\frac{7}{12}}>h_c
\end{array}
 \label{su42e2}
\end{equation}
where we let $h=h_c$ and $h_c^2=3\kappa^2$.

For the subcase (\ref{su42d}), the inequality (\ref{su42}) is simplified to the following inequality:
\begin{equation}
\begin{array}{rl}
2^{-6}2^{12}\bar{h}^{12}=2^{6}\frac{2^{6}}{2^{24}h^{12}}= \frac{1}{2^{18}h^{12}}>2
\end{array}
 \label{su42f2}
\end{equation}
This gives the condition on the parameter $h$
 for the (conventional)  antiferromagnetism.
We can write this condition in the following form (we let $h=h_c$):
\begin{equation}
\begin{array}{rl}
 h_{cAmax}:=\frac{1}{2}\cdot 2^{-\frac{7}{12}}>h_c
\end{array}
 \label{su42f}
\end{equation}

Let us then consider the following condition for satisfying (\ref{u17b}):
\begin{equation}
\begin{array}{rl}
(\bar{\kappa}^3-3\bar{\kappa}\bar{h}^2+3\bar{\kappa}^2\bar{h}-\bar{h}^3)
(\bar{\kappa}^3-3\bar{\kappa}\bar{h}^2-3\bar{\kappa}^2\bar{h}+\bar{h}^3)=0
\end{array}
 \label{su42g}
\end{equation}
This condition can be simplified to the following conditions:
\begin{equation}
\begin{array}{rl}
\kappa^2=h^2; \quad \mbox{or} \quad \kappa^2=(2\pm \sqrt{3})^2h^2
\end{array}
 \label{su42h}
\end{equation}
With each of these conditions we have that the eigenvalue $\lambda_{R1}<-1$.
The conditions (\ref{su42h}) can also be written in the following form:
\begin{equation}
\begin{array}{rl}
\kappa^2=h^2; \quad \mbox{or} \quad h^2=(2\pm \sqrt{3})^2\kappa^2
\end{array}
 \label{su42k}
\end{equation}
These conditions are as further states of phase transitions.
We remark that because these further states of phase transitions are different from the above states of phase transitions of antiferromagnetism and conventional 
 superconductivity we shall also call these further states of phase transitions as the states of cross-over.


 The two states $h>0, \kappa=0$ and
$h^2=3\kappa^2>0$ are as states of superconductivity when $h$ is restricted by the condition (\ref{su42e2}).
Thus for a fixed $h_c$ such that $h_{c}<\frac{\sqrt{3}}{2}\cdot 2^{-\frac{7}{12}}$,
the range $0<\kappa^2\leq \frac13 h_c^2$ is between two states of superconductivity.

The state $\frac{1}{(2+\sqrt{3})^2}h^2=\kappa^2$ is a state of phase transition. Thus, when the temperature is decreasing (or when the field is cooled) the material is from the state of superconductivity $\frac13 h^2=\kappa^2$ tending to the state $\frac{1}{(2+\sqrt{3})^2}h^2=\kappa^2$. 

 Let $h_{cI}$ and $h_{cII}$ be two $h$ values related by $ \frac13 h_{cI}^2=\frac{1}{(2+\sqrt{3})^2}h_{cII}^2$. From this relation
let us define a critical value $h_{cImax}$ by:
\begin{equation}
h_{cImax} := \frac{\sqrt{3}}{2+\sqrt{3}}\frac{\sqrt{3}}{2}\cdot 2^{-\frac{7}{12}}
\label{I5}
\end{equation}
Then 
the region of Type I superconductivity is the following region:
\begin{equation}
 \frac{1}{(2+\sqrt{3})^2}h_{cI}^2<\kappa^2\leq\frac13 h_{cI}^2;
 \,\, 0<h_{cI}< h_{cImax}.
\label{meI2}
\end{equation}


Then let us define a critical value $h_{cIImax}$ by:
\begin{equation}
h_{cIImax}:= \frac{\sqrt{3}}{2}\cdot 2^{-\frac{7}{12}}
\label{I4}
\end{equation}
Then from the definitions of $h_{cAmax}$, $h_{cImax}$ and $h_{cIImax}$ we have  $h_{cImax}<h_{cAmax}<h_{cIImax}$.
Then
the region of Type II superconductivity is the following region:
\begin{equation}
\frac{1}{(2+\sqrt{3})^2}h_{cII}^2<\kappa^2\leq\frac13 h_{cII}^2;
 \,\, h_{cImax}\leq h_{cII}< h_{cIImax}.
\label{meI3}
\end{equation}
where the range of $\kappa$ is connected with the corresponding range of $\kappa$ of the Type I superconductivity by the relation of $h_{cI}$ and $h_{cII}$.
We remark that the critical values $ h_{cImax}$ and  $ h_{cIImax}$ correspond to the critical magnetic field strength $H_{c1}$ and  $H_{c2}$ of the Ginzberg-Landau model
of superconductivity.

\subsection{\bf Meissner effect}

In the above sections we show that the magnetic properties are from
 the photon loops.
When a material is in the phase of Type I superconductivity
with the Cooper pair mechanism,
the spins (which are determined by the photon loops with the phases $e^{\pm ih}$ where $h\neq 0$) of the two electrons of a Cooper pair are opposite to each other.  Thus the net magnetic fluxes of the two photon loops of the two electrons of a Cooper pair is equal to $0$.
Further, when a material is in the state of Type I superconductivity,
the photons (or loops of photon) in this material
had to be attached to (or absorbed by) the Cooper pairs of electrons to form pairs with spins opposite to each other. Thus, in the state of Type I superconductivity, the net numbers of
 loops of photons
(and electrons) is equal to $0$.
Then since an external magnetic field acting on a material is as 
 a set of loops of photons in this material, when this material is in the state of Type I superconductivity we have that the net  winding numbers of
 loops of photons
 (and electrons) is equal to $0$.  Thus an external magnetic field acting on a material gives the zero effect of the net  winding numbers of
 loops of photons
in this material when this material is remained in the state of Type I superconductivity under the action of the external field.
Thus the external magnetic field is expelled from this material in the state of Type I superconductivity. Thus we have the Meissner effect.

We remark that in the case that the state $h^2=3\kappa_c^2$ is reached without the Cooper pair mechanism then we may not have the Meissner effect. In this case we have the phenomenon of persistent current (which is a phenomenon of superconductivity) appearing in a normal metal which is without the Meissner effect. Such phenomenon of persistent current 
was experimentally observed \cite{Lev7,Cha6,Mai6,Bar6,Blu6}.

The region of Type I and Type II superconductivity is the following region:
\begin{equation}
\frac{1}{(2+\sqrt{3})^2}h_{c}^2<\kappa^2\leq\frac13 h_{c}^2\quad \mbox{and}
 \quad 0<h_{c}< h_{cIImax}.
\label{meI5}
\end{equation}
We notice that in this region the parameter $\kappa>0$ gives static spins for magnetic effect as that of the state for antiferromagnetism. Then this region
is with larger $h_c^2$ (for a fixed $\kappa$) than the state for antiferromagnetism. Since $h_c^2$ gives dynamical effect for spins we have that this region is with
more dynamical effect for spins than the state for antiferromagnetism. This dynamical effect for spins
weakens the effect of antiferromagnetism which supresses the 
net  winding number of loops of photons.
This gives the Meissner effect in the region of Type I superconductivity where $h_c$ is restricted by the condition $h_c\leq h_{cImax}<h_{cAmax}$ that the effect of weakening the supression of 
the net  winding number of loops of photons.
 is not strong to give the Meissner effect.

On the other hand in the region of Type II superconductivity where $h_c$ is given by the condition $h_{cImax}\leq h_c< h_{cIImax}$ that the effect of weakening the supression of
the net  winding number of loops of photons
 is stronger. In particular $h_c$ can be large such that $h_{cAmax}<h_c<h_{cIImax}$ since $h_{cAmax}< h_{cIImax}$.

Thus by the effect of weakening the supression of 
the net  winding number of loops of photons,
the 
the effect of  winding number of loops of photons  can appear in the region of Type II superconductivity when an external magnetic field is applied.
This gives the phenomenon of the mixed state of Meissner effect with the appearing of 
 loops of photons in this region of Type II superconductivity (We shall give more analysis on this mixed state in the following section on high-$T_c$ superconductivity where the Type II superconductivity and the high-$T_c$ superconductivity are unified).
Since photons are in loop form,
in the mixed state when external magnetic field penetrating the material in the state of superconductivity,
the loops of photon  with the same orientation (or spin) of the magnetic field when attaching together form the macroscopic observable vortices (or magnetic fluxes).
These quantum  loops of photons  give the following formula
of the winding numbers which are as the quantization of energy of photons:
 \begin{equation}
n eq_{min}:= n\cdot n_e e_0 \cdot \frac{n_m e_0 \pi }{k_0},
\label{diracs}
\end{equation}
for $n=0,1,2, 3,\dots$.
When these  loops of photons  (with the same winding number) attaching together giving the macroscopic vortices, these vortices are with the same winding numbers 
and thus are with the following quantum of magnetic fluxes:
 \begin{equation}
n\Phi_0:=n\frac{h_0}{2e}:=n eq_{min}:= n\cdot n_e e_0 \cdot \frac{n_m e_0 \pi }{k_0} \label{diracss}
\end{equation}
for $n=0,1,2, 3,\dots$,
where $\Phi_0:=\frac{h_0}{2e}:=\frac{q_{min}}{2}$ is the unit of magnetic flux.

Thus from the quantization
of  loops of photons, 
 we have 
the  phenomenon of quantization of the vortices.
We shall give more description of the possibility of appearance of  vortices in the following section on high-$T_c$ superconductivity.  

\subsection{\bf Paramagnetic Meissner effect}

From the above region of Type I superconductivity,
 as the temperature decreasing (or the field is cooled) to
the state $\frac{1}{(2+\sqrt{3})^2}h_{cI}^2=\kappa^2$ we have a phase transition (or cross-over). Crossing  this state we have that $h_{cI}^2>(2+\sqrt{3})^2\kappa^2$. Thus crossing this state the spins of electrons are with more dynamical effect than the region of Meissner effect. 
In this case the antiferromagnetic state of the spins of the electrons becomes
a state which is more dynamical. Since the state of paramagnetic Meissner effect is with more dynamical effect than that of the antiferromagnetic state and the Type I superconductivity we have that the range of $\kappa$ of  paramagnetic Meissner effect is  the following range:
\begin{equation}
0 <\kappa^2\leq\frac{1}{(2+\sqrt{3})^2}h_{cI}^2
\label{me2}
\end{equation}
In addition, we have the following condition:
\begin{equation}
0<h_{cI}<h_{cImax}
\label{me2I}
\end{equation}
These two conditions give a region of paramagnetic Meissner effect which is connected to the region of Type I superconductivity. This region  of paramagnetic Meissner effect is only for the  Type I superconductivity.
We shall in the next section give the phase figures of the regions of the Type I and Type II superconductivity and this region of the paramagnetic Meissner effect when we investigate the phase diagram of high-$T_c$ superconducivity.
Experimentally the phenomenon of this Type I paramagnetic Meissner effect was  observed \cite{Gei}.
We shall also consider another (conventional Type II and unconventional) paramagnetic Meissner effect which was observed before this conventional Type I paramagnetic Meissner effect \cite{Tho,Kos,Sve,Bra5,Kho,Sig,Rie,Li}. 

\subsection{\bf Paramagnetism and antiferromagnetism}

Symmetric to  the range $\frac{1}{(2+\sqrt{3})^2}h_c^2<\kappa^2\leq \frac13h_c^2$
(the range of Meissner effect with vortex phenomena) which is a diamagnetic effect, 
the  following range of $\kappa$:
\begin{equation}
(2+\sqrt{3})^2h_c^2>\kappa^2\geq 3h_c^2
\label{me4a}
\end{equation}
is the range of antiferromagnetism, which is a property similar to that of
diamagnetism,
for a fixed
$h_c$ such that $\frac{1}{2}\cdot 2^{-\frac{7}{12}}>h_c$.

Then  symmetric to  the range $0<\kappa^2\leq\frac{1}{(2+\sqrt{3})^2}h_{cI}^2$
 (the range of paramagnetic Meissner effect) which is of paramagnetic effect,
the following range of $\kappa$:
\begin{equation}
\kappa^2\geq (2+\sqrt{3})^2h_c^2
\label{me4b}
\end{equation}
is the range of paramagnetism, for a fixed $h_c$ such that
$\frac{1}{2}\cdot 2^{-\frac{7}{12}}>h_c$.
Thus the state $\kappa^2=(2+\sqrt{3})^2h_c^2$ is the state for the phase transition of paramagnetism transited to antiferromagnetism
as $\kappa^2$ decreasing, for a fixed $h_c$ such that
$\frac{1}{2}\cdot 2^{-\frac{7}{12}}>h_c$.

\subsection{\bf Spin density waves and charge density waves in spin glass form}

Then as $\kappa$ decreasing the next state is $\kappa^2=3h_c^2$ which is the state for the phase transition of antiferromagnetism transited to another phase.
Following this state is the state $\kappa^2=h_c^2$ which is a symmetric state of $\kappa$ and $h$ where $\kappa$ is for static spin effect and $h$ is for dynamical wave effect.
Thus the following range of $\kappa$:
\begin{equation}
3h_c^2>\kappa^2> h_c^2
\label{me5}
\end{equation}
is identified as the phase of spin density waves where $\kappa$ gives static spin effect, $h$ gives wave and (charge) effect and that $\kappa^2> h_c^2$ indicates that the ratio of spin effect is greater than charge effect, and thus is the property of spin density waves (in contrast to charge density waves),
for a fixed $h_c$ such that $\frac{1}{2}\cdot 2^{-\frac{7}{12}}>h_c$.
In this 3D Case 1), this phase of  spin density waves  is identified as the phase of  spin density waves in spin glass form.
Then as a symmetry the following range of $\kappa$:
\begin{equation}
h_c^2>\kappa^2 >\frac{1}{3}h_c^2
\label{me7c}
\end{equation}
is identified as the phase of charge density waves, where $h_c^2>\kappa^2$ indicates that the ratio of spin effect is less than charge effect, and thus is the property of charge density waves (in contrast to spin density waves),
for a fixed $h_c$ such that $\frac{\sqrt{3}}{2}\cdot 2^{-\frac{7}{12}}>h_c$.
(Since this phase of charge density waves is a phase of more charge effect, the condition $\frac{\sqrt{3}}{2}\cdot 2^{-\frac{7}{12}}>h_c$ which is with larger $h_c$ for charge effect and is the condition for superconductivity is used instead of the condition $\frac{1}{2}\cdot 2^{-\frac{7}{12}}>h_c$ for antiferromagnetism).
In this 3D Case 1),  this phase of  charge density waves  is identified as the phase of charge density waves in spin glass form.
Then the state $\kappa^2=h_c^2$ is the state of phase transition (or cross-over) for the phase of spin density waves transited to the phase of charge density waves, as $\kappa$ decreasing.

On the other hand since the range (\ref{me5}) of spin density waves can be identified as a range of insulator and the range (\ref{me7c}) of charge density waves can be identified as a range of semiconductor, the state $\kappa^2=h_c^2$ is also the state of phase transition for the phase of insulator transited to the phase of conductor, as $\kappa$ decreasing.

\subsection{\bf Vortex core filled with normal metallic electrons}

Since the state of conventional Type II superconductivity is with the conditions that $\frac13h_c^2>\kappa^2>(2-\sqrt{3})^2h_c^2$ and $h_{cImax}\leq h_c< h_{cIImax}$ we have  that  the state of a vortex core in the region of conventional Type II superconductivity is with the condition  $\frac13h_c^2>\kappa^2>(2-\sqrt{3})^2h_c^2$ while the condition $h_{cImax}\leq h_c< h_{cIImax}$ can not be satisfied. Then since the state of Type I superconductivity is with the conditions that $\frac13h_c^2>\kappa^2>(2-\sqrt{3})^2h_c^2$ and $0<h_c< h_{cImax}$
and that the state of a vertex core is not the state of Type I superconductivity
we have  that  the state of a vortex core can not satisfy the condition $0<h_c< h_{cImax}$.

 Thus the state of a vortex core can only satisfy the condition $h_c \geq h_{cIImax}$. This condition means that the dynamical effect is relatively stronger (since the parameter $h_c$ is relatively larger) and is as a condition for the normal metallic state since for a given three dimensional meterial the normal metallic state is the only ordering state which can satisfy this condition $h_c \geq h_{cIImax}$.

Thus the state of a vertex core is the state of normal metallic state. This agrees with the experiments on the conventional Type II superconductivity that the vortex core is filled with normal metallic electrons.

On the other hand the state  of antiferromagnetism is with the condition $0<h_c< h_{cAmax}$ where $h_{cAmax}<h_{cIImax}$. Thus a vertex core in the region of the conventional Type II superconductivity can not be in the state of antiferromagnetism. We shall show that for high-$T_c$ superconductivity a vertex core can be in the state of antiferromagnetism.

\subsection{\bf Ferromagnetic superconductivity}

As explained in the above section on photon loops the spins of photons and of electrons are determined by the sign $\pm$ of the factor $e^{\pm ih\nu}$. Thus the spins of the two electrons of a Cooper pair are determined by the signs of $\pm h$ for $\nu >0$. Then we notice that the above state equations are of the form
$h^2=c\kappa^2$ for some constant $c>0$ (In particular for the state of superconductivity we have $c=3$). Then since $h^2=(-h)^2$ we have that it is possible for the two electrons of a Cooper pair have the same spin or have different spins when the Cooper pair is in the state of superconductivity. In this case when the two electrons are with different spins we have that the Cooper pair is a singlet pair and is for the conventional (antiferromagnetic) superconductivity. 

On the other hand when the two electrons are with the same spin we have that the Cooper pair is a triplet pair and we have the  ferromagnetic superconductivity. Thus it is possible for the existence of the
ferromagnetic superconductivity. Experimentally such ferromagnetic superconductivity has been observed for materials such as $UGe_2$ \cite{Sax,Aok,Huy,Nav,Mac}.

Then we notice that for the singlet Cooper pair because of the zero sum of
$\pm ih\nu$ of the factor $e^{ih\nu}e^{-ih\nu}$ of the singlet Cooper pair we have that the singlet Cooper pair is in a state of reduced energy (We shall consider in more detail of the singlet Cooper pair). On the other hand because of the sum  $ i2h\nu$ of $e^{ih\nu}e^{ih\nu}$ of the triplet Cooper pair we have that the triplet Cooper pair is in a state of increased energy.

For the ferromagnetic state we have $\kappa_c^2>h^2$ where $h^2\approx 0$ (which means that $h^2$ is very small while $h^2>0$) and superconductivity is with the condition
$\kappa_c^2=\frac13 h_c^2$.  Then from pairing we have $h_c=2h> 0$. Then we have
$\frac13 h_c^2= \frac43 h^2>h^2$. Thus the increasing of energy of the state of the triplet Cooper pair gives a possible way from the ferromagnetic state to the state of superconductivity. This shows that it is possible to have the state of  ferromagnetic superconductivity.


Further, since the ferromagnetic state is with $h^2\approx 0$
we have that the ferromagnetic superconductivity corresponds to the case that $h_c=2h> 0$ is small.
Thus for the range $\frac13 h_c^2\geq \kappa_c^2> \frac{1}{(2+\sqrt{3})^2} h_c^2$ of ferromagnetic superconductivity the critical temperature $T_c=\frac{1}{k_B}\kappa_c$ should be low. This agrees with the experiments such as the experiment of superconductivity of $UGe_2$ which shows that the ferromagnetic superconductivity is with a low critical temperature $T_c$ of about $1 K$ \cite{Sax,Aok,Huy,Nav,Mac}.

\section{High-$\bf{T_c}$ superconductivity of QED}

Let us then consider the 2D Case 2) with $\epsilon=1$. For this case we have the
subcase that $h_c\neq 0$ and $\kappa_c > 0$. We have that (\ref{u17a})
is of the following form:
\begin{equation}
\mbox{Re}(i\bar{h}+\bar{\kappa})^{12}>0
 \label{u42}
\end{equation}
where $0\leq \bar{h}, \bar{\kappa}\leq 1$. Also we have that (\ref{u17b})
is of the following form:
\begin{equation}
\mbox{Im}(i\bar{h}+\bar{\kappa})^{12}=0
 \label{u43}
\end{equation}

As similar to the above Case 1) from (\ref{u43}) we have the following condition:
\begin{equation}
(h^2-3\kappa^2)(\kappa^2-3h^2)=0
 \label{u44}
\end{equation}
This gives the following conditions:
\begin{equation}
h^2-3\kappa^2=0
 \label{u45}
\end{equation}
\begin{equation}
\kappa^2-3h^2=0
 \label{u46}
\end{equation}

{\bf Remark}. In this Case 2) for the two subcases $h=0$ and $\kappa=0$ we also get conditions similar to (\ref{u45}) and (\ref{u46}). However these conditions do not give restirctions on $h$ or $\kappa$ to give conditions of phases. Thus the two subcases $h=0$ and $\kappa=0$ do not give more phases and thus can be omitted. $\P$


From (\ref{u45}) we have that (\ref{u42}) is simplified to the following condition:
\begin{equation}
\frac{1}{2^{18}\kappa^{12}}>0
 \label{u47}
\end{equation}
Unlike the corresponding condition in the Case 1) this condition always hold for $\kappa>0$. Thus for $\epsilon=1$ the parameter $\kappa>0$ can be large such that the critical 
 temperature $T_c$ can be very high. This corresponds to the high-$T_c$ superconductivity phase transition.
On the other hand from (\ref{u46}) we have that (\ref{u42}) is simplified to the following condition:
\begin{equation}
\frac{1}{2^{18}h^{12}}>0
 \label{u48}
\end{equation}
Unlike the corresponding condition in the Case 1) this condition always holds for $h\neq 0$. Thus for $\epsilon=1$ the parameter $h\neq 0$ (and the parameter $\kappa^2=3h^2$) can be large such that the critical 
 temperature $T_N$ for antiferromagnetism can be high.

As  Case 1) let us then consider the following condition from which  (\ref{u42}) holds:
\begin{equation}
\begin{array}{rl}
(\bar{\kappa}^3-3\bar{\kappa}\bar{h}^2+3\bar{\kappa}^2\bar{h}-\bar{h}^3)
(\bar{\kappa}^3-3\bar{\kappa}\bar{h}^2-3\bar{\kappa}^2\bar{h}+\bar{h}^3)=0
\end{array}
 \label{u49}
\end{equation}
This condition can be simplified to the following conditions:
\begin{equation}
\begin{array}{rl}
h^2=\kappa^2; \quad \mbox{or} \quad h^2=(2\pm \sqrt{3})^2\kappa^2
\end{array}
 \label{u50}
\end{equation}
With each of these conditions we have that the eigenvalue $\lambda_{R1}<-1$.
Each of these conditions gives a state of phase transition.

We notice that the two states (\ref{u45}) and (\ref{u46}) of superconductivity and antiferromagnetism are as phases of superconductivity and antiferromagnetism which are one dimensional lines (and not as two dimensional regions) in the two dimensional phase plane $(\kappa, h)$. Thus these two phases will be unstable. 
However when the two dimensional phase plane $(\kappa, h)$
is transformed to
the two dimensional phase plane $(T, x)$ of the cuprates where $x$ denotes the doping parameter,  a part of the phase line (\ref{u45}) of superconductivity is bifurcated into two lines which encompass a region (We shall show the existence of this bifurcation and the forming of the region). 
Similarly the region of antiferromagnetism is also formed by such bifurcation of phase line.  These two regions are then the two well known regions of superconductivity and antiferromagnetism in the well known phase diagram in the phase plane $(T, x)$ of the cuprates \cite{Lia,Bra,Tim}.
In forming these two regions
the width of each of the two regions is the length of the coordinate $x$. In forming the region of superconductivity when $\kappa$ (or $T$) is fixed we have that the parameter $h$ is also fixed by the equation $h^2=3\kappa^2$. Then as the doping parameter $x$ ranging from one end of the region to another end of the region,  the corresponding parameter $h$ (which corresponds to the variation of $x$) is kept to be fixed by the equation $h^2=3\kappa^2$. This then has the effect that the phase line $h^2=3\kappa_c^2$ is bifurcated into two lines to form the region of superconductivity in the phase plane $(T, x)$.

 Similarly the region of antiferromagnetism in the phase plane $(T, x)$ is formed in this way. This then give the existence of the two phases of superconductivity and antiferromagnetism in the phase plane $(T, x)$ of the cuprates \cite{Lia,Bra,Tim}.

Let us then show the existence of the bifurcation of phase line and the forming of 2D phase region. Let us consider the case of superconductivity. 
For the conventional Type II superconductivity we have the following range of
superconductivity:
 \begin{equation}
\frac13h^2\geq\kappa^2>(2-\sqrt{3})^2h^2
\label{supercon5}
\end{equation}
for a fixed $h$ such that $h_{cImax}\leq h< h_{cIImax}$.
Then since the conventional Type II superconductivity is as a part of the complete region of superconductivity we have that this range give a range in the region of superconductivity. This range is as a part of the maximal range in the axis of $x$ of the complete region of superconductivity. This thus give a nonempty region of superconductivity and the existence of bifurcation for forming this region. This shows the existence of bifurcation and the forming of the 2D region of high-$T_c$ (and Type II) superconductivity.

Similarly we can show the existence of the bifurcation and the formation of the 2D region of antiferromagnetism.

\subsection{\bf Meissner effect with quantized vortex phenomenon}

We have shown that the existence of the conventional Type II superconductivity gives the bifurcation of the phase line $h^2=3\kappa_c^2$ of superconductivity and this bifurcation gives the stable existence of the region of high-$T_c$ superconductivity. When bifurcation forming the region of high-$T_c$ superconductivity,  this region of high-$T_c$ superconductivity is from the left bifurcation phase line $h^2=3\kappa_c^2$ of superconductivity  to the right bifurcation phase line $h^2=3\kappa_c^2$ of superconductivity. This region of high-$T_c$ superconductivity is just above (and contains) the region of the conventional Type II superconductivity. Thus the region of high-$T_c$ superconductivity and the conventional Type II superconductivity are both formed by the phase line $h^2=3\kappa_c^2$. 
Further for the high-$T_c$ superconductivity since the region of high-$T_c$ superconductivity is above the region of conventional Type II superconductivity there exists $h$ of the phase line $h^2=3\kappa_c^2$ such that $h_{cAmax}<h$. Thus as the conventional Type II superconductivity the Meissner effect with vortex phenomenon may appear in the region of high-$T_c$ superconductivity.
Moreover, since the unified region of high-$T_c$ superconductivity and the conventional Type II superconductivity is formed by the bifurcation of the phase line $h^2=3\kappa_c^2$, this unified region is easier to be penetrated by the magnetic fluxes of external magnetic field. Thus the mixed state of the Meissner effect with vortex phenomenon is easier to appear in this unified region.
Then as the case of the conventional Type II superconductivity,  the vortices are from  photon loops and are with the winding number as quantization and thus are with the following quantum of magnetic fluxes:
 \begin{equation}
n\Phi_0:=n\frac{h_0}{2e}:=n eq_{min}:= n\cdot n_e e_0 \cdot \frac{n_m e_0 \pi }{k_0}, \label{diracss2}
\end{equation}
for $n=0,1,2,3, \dots$,
where $\Phi_0:=\frac{h_0}{2e}:=\frac{q_{min}}{2}$ is the unit of magnetic flux.
Thus we have the quantization of the vortex phenomenon. 

\subsection{\bf  Bifurcation and quantum critical point of superconductivity}


When  the phase line $h^2=3\kappa_c^2$ of  superconductivity intersects with the bottom line  $\kappa_c=0$ (or $T=0$) of the phase plane $(\kappa,h )$, the intersection point is as the quantum critical point of superconductivity since it is the instability point on the base line  $T=0$ of the bifurcation of the phase line $h^2=3\kappa_c^2$.

\subsection{\bf Type II Paramagnetic Meissner effect}

For a fixed $\kappa_c$  the range $3\kappa_c^2< h^2$ is
separated into two ranges
$3\kappa_c^2<h^2<(2+\sqrt{3})^2\kappa_c^2$ and
$(2+\sqrt{3})^2\kappa_c^2\leq h^2$. Crossing the state $h^2=3\kappa_c^2$ 
with the Meissner effect,  the spins of electrons are with more dynamical effect than the range of Meissner effect. In this case the antiferromagnetic state of the spins of the electrons become the paramagnetic state which is more chaotic. In this case,  the Meissner effect becomes a paramagnetic Meissner effect. 
Thus the following range:
\begin{equation}
3\kappa_c^2<h^2<(2+\sqrt{3})^2\kappa_c^2
\label{mei}
\end{equation}
is the range of paramagnetic Meissner effect.

Then by the bifurcation of the phase line $h^2=3\kappa_c^2$  the region of  (conventional Type II and  high-$T_c$) superconductivity at the left side of  the (original non-bifurcated) phase line $h^2=3\kappa_c^2$ (or the quantum critical point)
 is  formed  as the overlapping of the region  (\ref{mei}) with the region $3\kappa_c^2> h^2$. Thus this  region of  (conventional Type II and  high-$T_c$) superconductivity at the left side of  the  (original non-bifurcated) phase line $h^2=3\kappa_c^2$ is also as a region of  paramagnetic Meissner effect. This paramagnetic Meissner effect in the region of  (conventional Type II and  high-$T_c$) superconductivity is observed  
 if we let $h^2$ be fixed and let $\kappa_c$ decreasing to $0$ (or the field is cooled) 
\cite{Tho,Kos,Sve,Bra5,Kho,Sig,Rie,Li}.

\subsection{\bf Range of normal metallic state}

Then
the next state $h^2=(2+\sqrt{3})^2\kappa_c^2$ is a state of phase transition.
Crossing  this state we have  $(2+\sqrt{3})^2\kappa_c^2<h^2$. Thus crossing this state the spins of electrons are with more dynamical effect than the range of paramagnetic Meissner effect. In this case the antiferromagnetics state of the spins of the electrons become the normal metallic state which is more chaotic. In this case,  the paramagnetic Meissner effect becomes the normal metallic effect. Thus the following range:
\begin{equation}
(2+\sqrt{3})\kappa_c^2<h^2
\label{mei2}
\end{equation}
is the range of normal metallic state.

\subsection{\bf Paramagnetism}

As a symmetry to that the range $3\kappa_c^2<h^2< (2+\sqrt{3})^2\kappa_c^2$ is the range of paramagnetic Meissner effect  we have that the following range of $h^2$:
 \begin{equation}
 \frac{1}{(2+\sqrt{3})^2}\kappa_c^2<h^2<\frac{1}{ 3}\kappa_c^2
\label{mei4}
\end{equation}
 is the range of paramagnetism.

\subsection{\bf Diamagnetism}

Then as a symmetry to that the range $h^2>(2+\sqrt{3})^2\kappa_c^2$ is the range of normal metallic state, the following range of $h^2$:
\begin{equation}
0<h^2<\frac{1}{(2+\sqrt{3})^2}\kappa_c^2
\label{mei4a}
\end{equation}
is the range of diamagnetism, for an (arbitrary) fixed $\kappa_c$. Thus the state $h^2=\frac{1}{(2+\sqrt{3})^2}\kappa_c^2$ is the state of phase transition (or cross-over) from diaamagnetism transited to  pamamagnetism as $h^2$ increasing.

\subsection{\bf Spin density waves and charge density waves}

Then as $h^2$ increasing the next state is $\kappa_c^2=3h^2$ which is the state  of phase transition for the phase (or state) of antiferromagnetism transited to the phase of spin density waves. Following this state is the state $\kappa_c^2=h^2$. Thus the following range of $h^2$:
\begin{equation}
\frac13 \kappa_c^2<h^2<\kappa_c^2
\label{mei5}
\end{equation}
is the range of spin density waves, for 
a fixed $\kappa_c$.
Then
the following range of $h^2$:
\begin{equation}
\kappa_c^2<h^2 <3\kappa_c^2
\label{mei7c}
\end{equation}
is as the (first range) of charge density waves, for 
a fixed $\kappa_c$ (We remark that the above range of paramagnetic Meissner effect can be regarded as the second range of charge density waves).
Thus the state $h^2=\kappa_c^2$ is as the state of phase transition for the phase of spin density waves transited to the phase of charge density waves,
where $\kappa_c>0$ is for static spins while $h\neq 0$ is for dynamical charge.
We remark that the range of spin density waves and charge density waves is also called the range of pseudogap. This identification of the range of spin density waves and charge density waves agrees with
experimental results \cite{Wis3,Deg,Ber3,Seo}.

 As the Case 1)
the state $h^2=\kappa_c^2$ is also as the state of phase transition for the phase of Mott insulator transited to the phase of semiconductor, where the phase (\ref{mei5}) of spin density waves is a phase of Mott insulator and the  phase (\ref{mei7c}) of charge density waves is a phase of semiconductor.

\subsection{\bf The condition (\ref{u17hi}) as condition of  3D phase
 in the complete phase diagram of high-${\bf T_c}$ superconductivity}

 We notice that the condition (\ref{u17hi}) in the above section on conventional Type II superconductivity divides the phase diagram of material with high-$T_c$ superconductivity into two phases. The phase satisfies the condition (\ref{u17hi}) is for the phase diagram in the above section including the conventional Type II superconductivity which is identified as the low-$T_c$ phase appearing in the bottom of the phase diagram of cuprates such as the cuprate $La_{2-x}Sr_{x}CuO_{4}$ \cite{Lia,Bra,Tim}. 

Then the other phase of the phase diagram of high-$T_c$ superconductivity satisfies the following condition which is complementary to the condition (\ref{u17hi}):
 \begin{equation}
 \begin{array}{rl}
&
[(\bar{\kappa}^3-
3\bar{\kappa}\bar{h}^2+3\bar{\kappa}^2\bar{h}-\bar{h}^3)^2
(\bar{\kappa}^3-3\bar{\kappa}\bar{h}^2-3\bar{\kappa}^2\bar{h}+\bar{h}^3)^2 \\
&\\
&
-4\bar{\kappa}^2\bar{h}^2(\bar{h}^2-3\bar{\kappa}^2)^2(\bar{\kappa}^2-3\bar{h}^2)^2 ]\leq  2^{7}
\end{array}
 \label{u17high2}
\end{equation}
This phase is the main phase in the upper high-$T_c$ part of the phase diagram of material with high-$T_c$ superconductivity  such as the cuprate $La_{2-x}Sr_{x}CuO_{4}$.

We notice that the conventional Type II superconductivity is of three dimensional space nature (with $0<\epsilon <<1$) while high-$T_c$ superconductivity is of two dimensional space nature (with $\epsilon=1$). Thus the two phases specified by the condition (\ref{u17high2}) and the condition (\ref{u17hi}) give a dimensional crossover. Such a phase diagram agrees with the fact that the cuprates such as $La_{2-x}Sr_{x}CuO_{4}$ are quisi 2D materails which are 3D materails but are with 2D properties in the experiments of superconductivity.
 Thus the conventional Type II superconductivity is a phase in the 3D phase and is as the 3D phase part 
of the complete phase of high-$T_c$ superconductivity. 

\subsection{\bf The upper limit of ${\bf T_c}$ of high-${\bf T_c}$ superconductivity}

 Since any point $(\kappa,h)$ of the phase line $h^2=3\kappa^2$ of  superconductivity can be the bifurcation point and that there is no upper bound on $h$ for the high-$T_c$ superconductivity we have that there is no upper limit on the critical temperature $T_c$ of high-$T_c$ superconductivity.
From the above phase diagram we see that the increasing of $h$ corresponds to the increasing of the parameter $x$ of cuprates \cite{Lia,Bra,Tim}. Thus the changing of the Planck constant parameter $h$ can be achieved by methods such as doping of cuprates.

\subsection{\bf Phase diagram  of $\bf{La_{2-x}Sr_{x}CuO_{4}}$}

As an example let us consider the experiments of the cuprate $La_{2-x}Sr_{x}CuO_{4}$ \cite{Lia,Bra,Tim}.
Let us consider the $(T,x)$ phase diagram of $La_{2-x}Sr_{x}CuO_{4}$ in Fig.1 (This phase diagram is derived from the well known experimental $(T,x)$ phase diagram of $La_{2-x}Sr_{x}CuO_{4}$ \cite{Lia,Bra,Tim}).

\begin{figure}[hbt]
\centering
\includegraphics[scale=0.2]{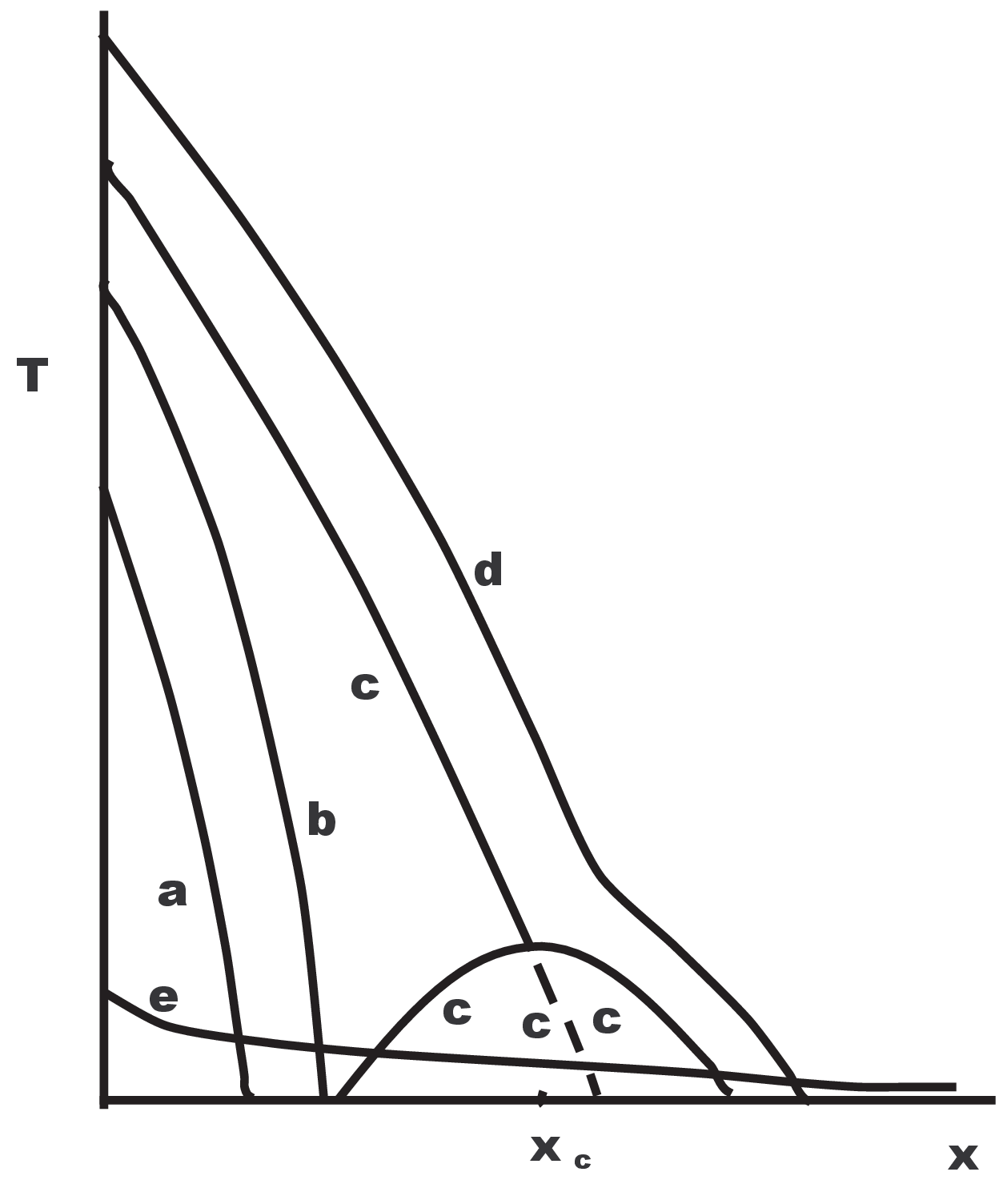}
\begin{minipage}[t]{8cm}
\begin{center}
Fig.5. $(T,x)$ phase diagram of $La_{2-x}Sr_{x}CuO_{4}$. \end{center}
In this phase diagram, the line {\bf a} is the (bifurcation) phase line $\kappa^2=3h^2$ of antiferromagnetism, the line {\bf b} is the cross-over line $h^2=\kappa^2$ of insulator to semiconductor, the line {\bf c} is the (bifurcation) phase line $h^2=3\kappa^2$ of superconductivity, the line {\bf d} is the cross-over line $h^2=(2+\sqrt{3})^2\kappa^2$ of state of pseudogap to normal metallic state; the line {\bf e} is the cross-over line from the phase of 2D to the phase of 3D (which coexists with the phase of 2D).
\end{minipage}
\end{figure}

 In this phase diagram the parameter $\kappa$ corresponds to $T$ by the relation $\kappa=k_B T$. The line {\bf a} denotes the phase line $\kappa^2= 3h^2$ (of the phase plane
$(\kappa, h)$) in the phase plane $(T, x)$. The line {\bf b} denotes the cross-over line $\kappa^2= h^2$ (of the phase plane
$(\kappa, h)$) in the phase plane $(T, x)$. The line {\bf c} denotes the phase line $h^2=3\kappa^2$ (of the phase plane
$(\kappa, h)$) of superconductivity in the phase plane $(T, x)$. Then the line {\bf d} denotes the cross-over line $h^2=(2+\sqrt{3})^2\kappa^2$ (of the phase plane
$(\kappa, h)$) in the phase plane $(T, x)$. In this phase diagram the original phase line $h^2=3\kappa^2$ with the part denoted by the dash line is bifurcated into the phase line $h^2=3\kappa^2$ denoted by the line {\bf c} in the phase plane $(T, x)$.

The two lines  {\bf a} and {\bf c} are the two basic phase lines for antiferromagnetism and superconductivity respectively in the well known $(T, x)$ phase diagram of $La_{2-x}Sr_{x}CuO_{4}$ \cite{Lia,Bra,Tim}). 
In the phase diagram in Fig.1 we have two more phase lines {\bf b} and {\bf d}. The region between {\bf a} and {\bf b} is the region of spin density waves and region of insulator. The region between {\bf b} and {\bf c} is the region of charge density waves and region of semiconductivity. The region between {\bf a} and {\bf c} is usually called the region of pseudogap.
 The region between {\bf c} and {\bf d} is the region of paramagnetic Messiner effect and is usually called the region of non-Fermi liquid (We shall also call this region as the extended region of pseudogap. We shall explain more on this extended region of pseudogap). 
The right side of {\bf d} is the region of normal metallic state. Then the bifurcation of the phase line {\bf c} gives the region of high-$T_c$ superconductivity.  The bifurcation  region of superconductivity at the left side of the original non-bifurcated (the dash-line)  {\bf c}  {\bf c} is also a  region of paramagnetic Messiner effect. The intersection point of  the original non-bifurcated (the dash-line)  {\bf c} and the base line $T=0$ is the quantum critical point of superconductivity.

Then the left side of {\bf a} is the region of  antiferromagnetism (The region of  antiferromagnetism is also obtained by the bifurcation of the phase line {\bf a} where a part of this region of  antiferromagnetism obtained by bifurcation is at the outside of this phase diagram in Fig.1).

Then let us consider the special phase line {\bf e}.
This phase line {\bf e} is the cross-over line from the phase of two dimensional phenomenon (i.e. the Case 2) with $\epsilon=1$) to the phase of three dimensional phenomenon (i.e. the Case 1) with $0<\epsilon<<1$). Above this line is the phase of two dimensional phenomenon described in this section on high-$T_c$ superconductivity and below this line is the
phase of three dimensional phenomenon described in the above section on conventional Type II  superconductivity (This three dimensional phenomenon coexists with the basic two dimensional phenomenon). In this 3D phase the region in the bifurcation region of superconductivity is the region of conventional Type II superconductivity while the region between the line {\bf c} and the line {\bf d} is the region of Type I superconductivity (In the 2D phase this region between the line {\bf c} and the line {\bf d} is the extended region of pseudogap).
We remark that in Fig.1 the phase line {\bf e} crosses the phase line {\bf d} with the equation $h^2=(2+\sqrt{3})^2\kappa^2$ and decreases to $0$ at the right side of the phase line {\bf d}.
This is because that the right side of the phase line {\bf d} is the normal metallic state which is at the outside of the phase diagram of conventional Type II superconductivity described in the above section.
Furthermore the crossing point is at $\kappa=\frac{1}{(2+\sqrt{3})^2}h^2$, as shown in Fig.1.
This is because that in the phase diagram of the conventional Type II superconductivity in the above section, for a given $h$ such that $h_{cImax}\geq h$, the range of the paramagnetic Meissner effect is the range $\frac{1}{(2+\sqrt{3})^2}h^2>\kappa^2>0$. This gives an exact fitting of the phase diagram of the 3D conventional Type II superconductivity in the phase diagram of the quasi-2D high-$T_c$ superconductivity.
Furthermore, since the region of normal metallic state in the vortex core of the conventional superconductivity is the region: $\frac13h^2>\kappa^2>\frac{1}{(2+\sqrt{3})^2}h^2$ and $h\geq h_{cIImax}$, and the phase diagram of the conventional superconductivity is in the below of the line {\bf e}, we have that in Fig.1 the region of the normal metallic state of the conventional superconductivity is the region enclosed by the (left) line {\bf c} and the line {\bf e} (and is at the right side of the (left) line {\bf c} and the above of the line {\bf e}).
Thus this normal metallic state of the conventional superconductivity in the phase diagram in Fig.1 contains the region of the normal metallic state of the high-$T_c$ superconductivity enclosed by the line {\bf d} and the line {\bf e} (and is at the right side of the line {\bf d} and the above of the line {\bf e}).
Thus this phase diagram in Fig.1 is a complete phase diagram on magnetism and superconductivity that it includes the Type I superconductivity, the conventional Type II superconductivity and the high-$T_c$ superconductivity.  

\subsection{\bf Computation of $\bf{T_c}$ of $\bf{La_{2-x}Sr_x CuO_4}$}
 Let us consider the mechanism of  forming  Cooper pair
 for the superconductivity of $La_{2-x}Sr_x CuO_4$.
Let us first consider the possibility of superconductivity of the transition element $Cu$. This element is with the outer shells $3d^{10}4s^1$ of valence electrons. Since $Cu$ is a metal we have that $Cu$ is in the state $3\kappa^2<h^2$. 
For $Cu$ in the state $3\kappa^2=h^2$ of superconductivity we have that $h^2$ is needed to be reduced such that $3\kappa^2=h^2$.
Let us then consider the Cooper pair of electrons as such a mechanism. Since an atom $Cu$ has only one $s$ electron in the outer shells we have that this $s$ electron can not be used to form a Cooper pair.
Thus we can only consider the $d$ electrons in the outer shells of $Cu$. There are five $d$ electrons with up spins and five $d$ electrons with down spins in the outer shells of $Cu$. Thus the five $d$ electrons with up spins of one $Cu$ may be with the five $d$ electrons with down spins of another $Cu$ to form five pairs of electrons.
Then since the $d$ electrons are in the inner of the outer shells that the binding energy on the $d$ electrons is stronger that such pairings of the $d$ electrons can not be formed stabely. Thus for the $Cu$ metal there is no mechanism for the forming of stable pairings of electrons and thus it is difficult for the metal $Cu$ in the state of superconductivity. This agrees with the experiments that $Cu$ is not a superconductor (We remark that $Cu$ being not in the state of superconductivity gives the antiferromagnetic state of $La_2 CuO_4$).

Then let us consider the case of $La_{2-x}Sr_x CuO_4$. 
We shall show that the interaction of $Cu$ with $La(Sr)$ opens a channel for the freedom of $d$ valence electrons of $Cu$ by doping $La(Sr)$. With this channel opening the Cooper pair of $d$ valence electrons between two $Cu$ atoms can be formed by the attractive effect of the seagull vertex interaction.
Let the two $d$ electrons for the forming of Cooper pair be with the wave factors $e^{\pm ih_{Cu}}$ respectively in their wave functions where $|h_{Cu}|$ is the energy of the $d$ electrons for the interaction of forming the Cooper pair (The $\pm$ sign is for the up and down spins of the two $d$ electrons.
This forms a wave factor $e^{i(h_{Cu}-h_{Cu})}$ of the Cooper pair  where $h_{Cu}-h_{Cu}=0$. This cancelation of energy for the forming of Cooper pair gives the reduction of the energies of the $Cu$ atoms. Thus the forming of Cooper pair is a mechanism for the reduction of energy.

Let us specify the parameter $h$ of $h^2=3\kappa^2$. We notice that when there is no mechanism of pairing there is no reduction of $h$ of $Cu$ to let the metal $Cu$ from the normal metallic state to the state of superconductivity. Thus more reduction of energy gives higher critical temperature $T_c$ of superconductivity (when the reduced energy is not exceeded the remained energy). Thus the reduction of energy is proportional to the critical temperature $T_c$ of superconductivity.
Let us find out the exact proportion.
Each $d$ electron of the metal $Cu$ gives an energy $|h_{Cu}|$.
The total energy of $d$ electrons of the metal $Cu$ is $10|h_{Cu}|$. Let $h_1$ be the energy reduced from this total energy by the mechanism of Cooper pairing. Then since $h_1$ is proportional to $T_c$ from the state equation of superconductivity we have $h_1^2=3\kappa_{1}^2=3(k_BT_{c1})^2$.
On the other hand
let $h_2=10|h_{Cu}|-h_1$ be the energy after the reduction of energy. Then
this remained energy must also be proportional to the critical temperature $T_c$ of superconductivity because of the general state equation $h_2^2=3\kappa_{2}^2$ of superconductivity. Thus from the state equation of superconductivity we have $h_2^2=3\kappa_{2}^2=3(k_BT_{c2})^2$.
Thus from the formula $h^2=3\kappa^2$ we have that the maximum of $T_c$ is at the case that $h_1$ is increased to equal to $h_2$. Thus for the maximum of $T_c$ we have $h_1=h=h_2$.
Thus we have the following formula of $T_c$ of $La_{2-x}Sr_x CuO_4$:
\begin{equation}
k_BT_c=\kappa=\frac{1}{\sqrt{3}}\Delta_{LSCO}
\label{CuO4}
\end{equation}
where $\Delta_{LSCO}=h=5|h_{Cu}|$ is the energy gap of the high-$T_c$ superconductivity of $La_{2-x}Sr_x CuO_4$ and $T_c$ is the highest critical temperature of superconductivity.

We remark that in the above derivation the energy gap $\Delta_{LSCO}$ is from the $Cu$ atom. We shall show that there is a case that a small energy gap $|h_{La(Sr)}|$ from two $La(Sr)$ atoms may appear such that $\Delta_{LSCO}=h=5|h_{Cu}|+2|h_{La(Sr)}|$.

Let us then find a way to determine the energy $|h_{Cu}|$. Since the $d$ valence electrons are approximately as free electrons we have that the energy $|h_{Cu}|$ is proportional to the ionization energy of the $d$ electrons.
Then since an atom $Cu$ in the $CuO_2$ plane is vertically interacted with two $O$ atoms in two $La(Sr)O$ planes in opposite directions we have that  the $d$ valence electrons of $Cu$ can be in three states of ionization energies. Thus the $d$ valence electrons of $Cu$ can be in the state of third ionization energy of $Cu$
(We shall explain this point in more detail when we consider the effect of $La(Sr)$ on $Cu$).
 Thus the maximum value of the  energy $|h_{Cu}|$ is proportional to the third ionization energy of $Cu$
(It will be more precise that  the maximum value of the  energy $|h_{Cu}|$ is proportional to
the third ionization energy of the electrons of $Cu$ in the $CuO_2$ plane. This ionization energy includes the effect of the oxygen atoms on the $d$ electrons).

From the existing table of ionization energy of $Cu$ we have that the
third ionization energy of $Cu$ is approximately equal to
$3555$ kJ/mol.
Let $|h_{Cu}|\approx \xi 3555$ kJ/mol where $\xi$ is a proportional constant. We shall show that $\xi\approx 2.83133971\cdot 10^{-5}$ when we compute the $T_c$ of the transition element $Nb$.
Then from (\ref{CuO4}) we can compute the highest critical temperature $T_c$ of $La_{2-x}Sr_x CuO_4$:
\begin{equation}
T_c \approx 34.95 K \,\, \mbox{(Computed value of $T_c$ of $La_{2-x}Sr_x CuO_4$)}
\label{CuO5}
\end{equation}
This agrees with  the experimental value of $T_c\approx 35 K$ of $La_{2-x}Sr_x CuO_4$ at the optimal doping $x_c\approx 0.15$.

In a $La(Sr)O$ plane
we have that one $O$ atom can interact with two $La(Sr)$ atoms. This gives an interaction between two $La(Sr)$ atoms. With this additional interaction we have that the seagull vertex interaction gives the attractive effect for the forming of the Cooper pair of $s$ valence electrons where the $La$ atom is with outer shells $5d^16s^2$ and the $Sr$ atom is with outer shell $6s^2$.
The two $s$ valence electrons for the forming of Cooper pair are with wave factors $e^{\pm ih_{La(Sr)}}$ respectively in their wave functions where $|h_{La(Sr)}|$ is the energy of the $s$ electrons for the interaction of forming the Cooper pair.
Then the $Cu$ atom and the two $La(Sr)$ atoms on the centered $c$-axis of $La_{2-x}Sr_x CuO_4$ are  regarded as a clustering of atoms and are regarded as a larger atom.
Thus in the case of the appearing of Cooper $ss$-pairs of $La(Sr)$ atoms we have that the total energy gap of this clustering of atoms is $\Delta_{LSCO}=h=5|h_{Cu}|+2|h_{La(Sr)}|$.
Then  since a $La(Sr)$ atom in a $La(Sr)O$ plane can interact with only one $O$ atom at the outside of this $La(Sr)O$ plane and next to this $La(Sr)$ atom we have that the $s$ valence electrons of $La(Sr)$ can have at most two states of ionization energies.
Thus
the maximum value of the energy parameter $|h_{La(Sr)}|$ is proportional to the second ionization energy of $La$ (or $Sr$). From the existing table of ionization energies we have that the
second ionization energy of $La$ is approximately equal to
$1067$ kJ/mol (The second ionization energy of $Sr$ is approximately equal to
$1064.2$ kJ/mol and can be regarded as equal to that of the $La$ atom).
Let $|h_{La(Sr)}|\approx \xi 1067$ kJ/mol. 
Then from (\ref{CuO4}), with $\Delta_{LSCO}=5|h_{Cu}|+2|h_{La(Sr)}|$, we can compute the highest critical temperature $T_c$ of $La_{2-x}Sr_x CuO_4$:
\begin{equation}
T_c \approx 39.14 K \quad \mbox{(Computed  $T_c$ of $La_{2-x}Sr_x CuO_4$)}
\label{CuO8}
\end{equation}
This agrees with  another experimental value of $T_c\approx 39 K$ of $La_{2-x}Sr_x CuO_4$ at the optimal doping $x\approx 0.15$. Thus from this high-$T_c$ theory of superconductivity we can explain the two experimental values of $T_c$ of $La_{2-x}Sr_x CuO_4$.

\subsection{\bf Doping mechanism and Jahn-Teller effect of degenerate states}

In the literature there are two main approaches 
 for the mechanism of high-$T_c$ superconductivity. One approach
 is that the pairing mechanism is a purely  electron-electron interaction \cite{And,And5}.
Another 
approach is that  the pairing mechanism is an electron-phonon  interaction 
\cite{Mul,Bus}.

In this new gauge model  the pure electron-electron interaction and the  electron-phonon  interaction 
 are unified consistently to give the mechanism of  high-$T_c$ superconductivity. As a pure electron-electron interaction an attractive potential between electrons is derived from the seagull vertex term  for the interaction of forming Cooper pairs of electrons.  Then as an  electron-phonon  interaction 
 let us derive the doping mechanism of degenenrate states of electrons to give the Jahn-Teller effect and the Jahn-Teller polarons of  high-$T_c$ superconductivity.

To this end we invesitigate the mechanism of the increasing of the doping parameter $x$ giving the increasing of $h$ for $La_{2-x}Sr_x CuO_4$.
 Let us consider the following function:
 \begin{equation}
 f(x)=1850.3(1-\frac{x}2)+4138\frac{x}2
 \label{LSCO63}
\end{equation}
 where $1850.3$ kJ/mol and $4138$ kJ/mol are approximately the
third ionization energies of $La$ and $Sr$ respectively. Since $1850.3<4138$ and $h$ is proportional to the ionization energy we have that this function gives the relation that the increasing of $x$ giving the increasing of $h$ for $La_{2-x}Sr_x CuO_4$.

We have $1850.3<1957.9<4138$ where $1957.9$ kJ/mol is  the state of second ionization energy of $Cu$.
Let us establish the following relation (A main point of this relation is that $1957.9-1850.3$ is small and $4138>1957.9$):
\begin{equation}
f(x_0)=1850.3(1-\frac{x_0}2)+4138\frac{x_0}2=1957.9
\label{LSCO6}
\end{equation}
for some $x_0$ ($0<x_0<1$).
This then relates the state of third ionization energy of $La$ and the  state of second ionization energy of $Cu$. 

When  (\ref{LSCO6})  holds (or approximately holds)  the channel connecting the  state of second ionization energy of $Cu$ and   the  state of third ionization energy of $La$
is opened (We notice that the relation   (\ref{LSCO6})  gives the degeneration of electron states. The resulting degenerate state gives a Jahn-Teller electron-phonon interaction effect which is described as a channel opening \cite{Mul,Bus,Ber}.  Let us call   (\ref{LSCO6}) and the  resulting  Jahn-Teller  effect of degeneration as  the degenerate state of channel opening).. 

This channel opening gives
a freedom of electric current with
a direction  orthogonal to the $CuO_2$ plane.
From this freedom of electric current, the Cooper pairs of $d$ valence electrons of $Cu$ can be formed
(In other words, through this channel opening the $d$ valence electrons of $Cu$ can transit from the valence band to the conduction band from which Cooper pairs of these $d$ valence electrons can be formed).
Thus this channel opening  together with the $CuO_2$ plane gives the region of 3D superconductivity. From this region of 3D superconductivity we have the existence of quasi-2D bifurcation region of the unconventional high-$T_c$ superconductivity given by the $CuO_2$ plane.


We notice that $x_0\approx 0.094$. Thus the doping range $x_0<x<x_1$ for some $x_1<1$ is a range for the state of high-$T_c$ superconductivity. This agrees with the experiment that the state of high-$T_c$ superconductivity of $La_{2-x}Sr_x CuO_4$ is in the doping range $0.06 \leq x\leq 0.3$.
When (\ref{LSCO6}) holds we have that the $d$ valence electrons of $Cu$ (forming Cooper pairs) are in the basic state of second ionization energy, and other states are to be reached from this basic state. Also when (\ref{LSCO6}) holds the $d$ valence electrons of $Cu$ and $La(Sr)$ are unified to occupy the sequence of state of ionization energies.
Since the $Cu$ atoms in the $CuO_2$ plane are with one more ionization direction than the $La(Sr)$ atoms in the $La(Sr)O$ plane we have that the $d$ valence electrons of $Cu$ occupy the higher states while the $d$ valence electrons of $La(Sr)$ occupy the lower states.

 Then since there are no $d$ valence  electrons of $La(Sr)$ forming Cooper pairs and the $d$ valence electrons of $Cu$ have three states of ionization energies
we have that the $d$ valence electrons of $Cu$ (forming Cooper pairs) have the states of first, second and third ionization energies.
Thus the maximum value of the energy parameter $|h_{Cu}|$ is proportional to the third ionization energy of $Cu$. This gives a more detail explanation for that the energy parameter $|h_{Cu}|$ of the $d$ valence electrons of $Cu$ in $La_{2-x}Sr_x CuO_4$ is proportional to the third ionization energy of $Cu$.

\subsection{\bf Computation of $\bf{T_c}$ of $\bf{YBa_2Cu_3O_{7}}$}

Let us then investigate the mechanism of the  high-$T_c$ superconductivity of the cuprate $YBa_2Cu_3O_{6+x}$ ($0\leq x\leq 1$).
Similar to  $La_{2-x}Sr_x CuO_4$
let us consider again the Cooper pairs of electrons of $Cu$ as a mechanism for reducing $h^2$ to reach the state $3\kappa^2=h^2$ of superconductivity.
As  $La_{2-x}Sr_x CuO_4$, we consider the $d$ valence electrons in the outer shells of valence electrons of $Cu$.
It is known that $YBa_2Cu_3O_{6+x}$ has the multi-layer form $BaO/CuO/BaO/CuO_2/Y/CuO_2/BaO$ \cite{Chu}. The two $CuO_2$ planes, one $Y$ plane and one $CuO$ plane are as a quantum system.
The $BaO$ planes are as external to this quantum system.  Let us suppose that the multi-layer $CuO/BaO/CuO_2/Y/CuO_2$ can still be regarded as a quasi 2D system such that the above theory of high-$T_c$ superconductivity can be applied.
In this case the three $Cu$ atoms
at the (centered) $c$-axis of (a unit cell) joining the two $CuO_2$ planes
and the $CuO$ plane in the multi-layer $CuO/BaO/CuO_2/Y/CuO_2$ form a cluster. Thus the three $Cu$ atoms can be regarded as a larger atom.

Let us first consider the $CuO_2$ planes.
In a $CuO_2$ plane one $Cu$ atom can interact with two $O$ atoms. This gives an interaction between two $Cu$ atoms. With this additional interaction we have that the seagull vertex interaction gives the attractive effect for the forming of the Cooper pair of $d$ electrons.
As  $La_{2-x}Sr_x CuO_4$,
the two $d$ electrons for the forming of Cooper pair are with the wave factors $e^{\pm ih_{Cu}}$ respectively in their wave functions where $|h_{Cu}|$ is the energy of the $d$ electrons for the interaction of forming the Cooper pair.
Further since the two $CuO_2$ planes are separated by a $Y$ plane there is only one direction of ionization orthogonal to the two $CuO_2$ planes we have that the $d$ valence electrons of $Cu$ of these two $CuO_2$ planes have only two states of ionization energies.
We shall show that the maximum value of the energy parameter $|h_{Cu}|$ of the $d$ valence electrons of $Cu$ of these two $CuO_2$ planes is proportional to the third ionization energy of $Cu$.

Let us then consider the $CuO$ plane.
In the $CuO$ plane one $O$ atom can interact with two $Cu$ atoms in chain form.
This gives an interaction between two $Cu$ atoms. With this additional interaction we have that the seagull vertex interaction gives the attractive effect for the forming of the Cooper pair of $d$ electrons.
The two $d$ electrons for the forming of Cooper pair are with wave factors $e^{\pm ih_{CuO}}$ respectively in their wave functions where $|h_{CuO}|$ is the energy of the $d$ electrons for the interaction of forming the Cooper pair in the $CuO$ plane.
This energy parameter $|h_{CuO}|$ gives a wave factor $e^{i(h_{CuO}-h_{CuO})}$ of the Cooper pair of $d$ electrons of $Cu$ in the $CuO$ plane. This cancelation of energy for the forming of Cooper pair gives the reduction of the energies of the $Cu$ atoms. Thus the forming of Cooper pair of $d$ electrons of $Cu$ in the $CuO$ plane is a mechanism for the reduction of energy.
Since there is only one direction of ionization orthogonal to the $CuO$ plane we have that the $d$ valence electrons of $Cu$ of the $CuO$ plane have only two states of ionization energies.
We shall show that
the maximum value of the energy parameter $|h_{CuO}|$ of the $d$ valence electrons of $Cu$ in the $CuO$ plane is proportional to the second ionization energy parameter of $Cu$.
For  the doping mechanism of superconductivity,
let us consider the following function:
 \begin{equation}
 f(x)=1980-3600(1-x)
 \label{YBCO63}
\end{equation}
where $1980$ kJ/mol and $3600$ kJ/mol are the third ionization energies of $Y$ and $Ba$ respectively. This function gives the relation that the increasing of $x$ giving the increasing of $h$ ($0\leq x\leq 1$).
Then let us set  the following relation of  the degenerate state of channel opening (A main point for setting up this relation is that $1980-1957.9$ is small and $3600>1980$):
\begin{equation}
f(x_0)=1980-3600(1-x_0) =1957.9x_0
\label{YBCO6}
\end{equation}
for some $x_0$ ($0<x_0<1$).  
Equivalently we have the following relation:
\begin{equation}
1980 =1957.9x_0+3600(1-x_0)
\label{YBCO62}
\end{equation}
When this relation (\ref{YBCO6}) holds
 (or approximately holds), the channel connecting the two states of  ionization energy of $Cu$ and  $Y$ 
 is opened. 
This channel opening gives a freedom of electric current with
a direction orthogonal to the  two $CuO_2$ planes. 
From this freedom of electric current, the Cooper pairs of the $3d$-level electrons of $Cu$  can be formed.
Thus this channel opening gives the 3D region of conventional superconductivity. From this 3D conventional superconductivity we have the existence of quasi-2D bifurcation region of unconventional high-$T_c$ superconductivity given by the  two $CuO_2$ planes.


We notice that $x_0\approx 0.9865$.
Thus the cuprate $YBa_2Cu_3O_{6+x}$ comes into the range  of high-$T_c$ superconductivity when $x$ is just greater than $x_0$. Since $1-x_0$ is very small we have that $x_c\approx 1$
is the optimal doping with highest $T_c$ of high-$T_c$ superconductivity. This agrees with the experiment that the cuprate $YBa_2Cu_3O_{7}$ is in optimal doping with highest $T_c$ of high-$T_c$ superconductivity.

When (\ref{YBCO6}) holds, the $d$ valence electrons of $Cu$ in the $CuO$ plane and the two $CuO_2$ planes are in the basic state of second ionization energy, and other states are to be reached form this basic state.
Also when (\ref{YBCO6}) holds the $d$ valence electrons of $Cu$ in the $CuO$ plane and the two $CuO_2$ planes are unified to occupy the sequence of state of ionization energies.
Since the $Cu$ and $O$ are interacted in chain form in the $CuO$ plane we have that the $Cu$ atoms in the $CuO_2$ plane are with one more ionization direction than the $Cu$ atoms in the $CuO$ plane. Thus the $d$ valence electrons of $Cu$ in the $CuO_2$ planes occupy the higher states of ionization energies while the $d$ valence electrons of the $Cu$ atoms in the $CuO$ plane occupy the lower states.
Thus the $d$ valence electrons of $Cu$ in the $CuO_2$ planes have the states of second and third ionization energies and the $d$ valence electrons of $Cu$ in the $CuO$ plane have the state of first and second ionization energies.
Thus the maximum value of the energy parameter $|h_{CuO}|$ of the $d$ valence electrons of $Cu$ in the $CuO$ plane is proportional to the second ionization energy of $Cu$ and the maximum value of the energy parameter $|h_{Cu}|$ of the $d$ valence electrons of $Cu$ in the two $CuO_2$ planes is proportional to the third ionization energy of $Cu$.
Then we have the following formula of $T_c$ of $YBa_2Cu_3O_{6+x}$:
\begin{equation}
k_BT_c=\kappa=\frac{1}{\sqrt{3}}\Delta_{YBCO}
\label{YBCO}
\end{equation}
where $\Delta_{YBCO}=h=10|h_{Cu}|+5|h_{CuO}|$ is the energy gap of 
$YBa_2Cu_3O_{6+x}$, where $10=\frac{20}2$ is from the 20  $d$-electrons of two $Cu$ atoms of two $CuO_2$ plane and $5=\frac{10}2$ is from the 10  $d$-electrons of one $Cu$ atom of the $CuO$ plane of a unit cell of $YBa_2Cu_3O_{6+x}$. 
Then from the  table of ionization energy, the
second and third ionization energies of $Cu$ are approximately equal to $1957.9$ kJ/mol and
$3555$ kJ/mol.
Let $|h_{Cu}|\approx \xi 3555$ kJ/mol and $|h_{CuO}|\approx \xi 1957.9$ kJ/mol.
Then from (\ref{YBCO}) we can compute the highest critical temperature $T_c$ of $YBa_2Cu_3O_{6+x}$:
\begin{equation}
T_c \approx 89.14 K \quad \mbox{(Computed value of $T_c$ of $YBa_2Cu_3O_{6+x}$)}
\label{YBCO5}
\end{equation}
This agrees with  the experimental value of $T_c\approx 90 K$ of $YBa_2Cu_3O_{6+x}$ \cite{Chu}.

Similar to the $La$ atom for $La_{2-x}Sr_x CuO_4$, when the channel through $Y$ is open,  and the $s$ valence electrons of $Y$ also form a Cooper $ss$-pair, then this gives an energy parameter $|h_{Y}|$ which is proportional to the first ionization energy of $Y$.
Then $\Delta_{YBCO}=h=10|h_{Cu}|+5|h_{CuO}|+|h_{Y}|$ where $|h_{Y}|\approx \xi 600$ kJ/mol ($ 600$ kJ/mol is approximately the
first ionization energy of $Y$). Then from (\ref{YBCO}) we have:
\begin{equation}
T_c \approx 90.32 K \quad \mbox{(Computed value of $T_c$ of $YBa_2Cu_3O_{6+x}$)}
\label{YBCO9}
\end{equation}
This also agrees with  the experimental value of $T_c\approx 90 K$ of $YBa_2Cu_3O_{6+x}$ \cite{Chu}.

\subsection{\bf Computation of $\bf{T_c}$ of $\bf{MgB_2}$}

The intermetallic material $MgB_2$ is found to have the phenomenon of two-gap superconductivity \cite{Nag2,Cha,Liu,Cho,Iav,Bou,Sol}.
 Let us show that this phenomenon can also be derived from the above theory of high-$T_c$ superconductivity.
This material $MgB_2$ has two layers consisting of $Mg$ and $B$ respectively. It is known that $MgB_2$ has a structure that the $B$ atoms in a layer of $B$ form a honeycomb. Then each $Mg$ atom in the layer of $Mg$ is located in the center of each hexagon of the honeycomb. Thus a hexagon of six $B$ atoms is enclosed by a hexagon of six $Mg$ atoms and is separated from other hexagons of $B$ atoms by the hexagon of six $Mg$ atoms. 
Let us call such a hexagon of six $B$ atoms as a unit element in a layer of $B$ atoms.

The two
 $s$-valence electrons and the two $s$-electrons of the inner shell of a  $B$ atom give an energy $2|h_{B1}|+2|h_{B2}|$. Then since the $B$ atoms of a unit element are clustered together that these $B$ atoms can be regarded as a larger atom. Thus the total $s$-electrons of a unit element of $B$ atoms give the total energy $12|h_{B1}|+12|h_{B2}|$.
Let us consider the intermetallic $Mg_{1-x}Be_{x}B_2$ ($0\leq x\leq 1$). This intermetallic is of the $AlB_2$-type. 
For the doping mechanism of of superconductivity,
let us consider  consider the following function:
 \begin{equation}
 f(x)=737.7(1-x)+899.5x
 \label{MgB263}
\end{equation}
 where $737.7$ kJ/mol and $899.5$ kJ/mol are approximately the
first ionization energies of $Mg$ and $Be$ respectively. Since $737.7<899.5$ and $h$ is proportional to the ionization energy we have that this function gives the relation that the increasing of $x$ giving the increasing of $h$ for $Mg_{1-x}Be_{x}B_2$.
We have $737.7<800.6<899.5$ where $800.6$ kJ/mol is  the first ionization energy of $B$.
Let us set the following relation of  the degenerate state of channel opening: 
\begin{equation}
f(x_0)=737.7(1-x_0)+899.5x_0=800.6    
\label{MgB269}
\end{equation}
for some $x_0$ ($0\leq x_0\leq 1$).

 When this relation holds (or approximately holds), 
 the channel connecting the two states of first ionization energy of $B$ and  $Mg$
is opened. This channel opening gives
a freedom of electric current with
a direction orthogonal to the $B$ plane.
By this freedom of electric current, the electrons of a cluster of atoms in a unit cell can be coupled to the electrons of other clusters to form Cooper pairs. In this way
the Cooper pairs of $s$-valence (and nonvalence) electrons of $B$, the Cooper pairs of $s$-valence electrons of $Mg$ 
 can be formed.
Thus this channel opening gives the 3D region ($\pi$-band) of conventional superconductivity. From  this 3D conventional superconductivity we have the existence of quasi-2D  bifurcation region ($\sigma$-band) of  unconventional high-$T_c$ superconductivity given by the $B$ plane.


We notice that $x_0\approx 0$. 
This agrees with the experiment that the state of superconductivity of $Mg_{1-x}Be_{x}B_2$ is in the doping range $0 \leq x<1$ and that $MgB_2$ is in the state of superconductivity \cite{Fel}. 

When (\ref{MgB269}) holds we have that the $s$ valence electrons of $B$ are in the state of first ionization energy of $B$. Then the $s$ nonvalence electrons of $B$ are in the state of second ionization energy of $B$. Also the $s$ valence electrons of $Mg$ are in the state of first ionization energy of $Mg$. The $s$ valence electrons of $B$ and $Mg$ in the state of first ionization energy are in the opened channel of 3D superconductivity
while the $s$ nonvalence electrons of $B$ in the state of second ionization energy of $B$ are for the quasi-2D high-$T_c$ superconductivity of the $B$ plane.
Thus the maximum value of the energy parameter $|h_{B1}|$ and $|h_{B2}|$ are proportional to the first and second  ionization energies of $B$ respectively, and the maximum value of the energy parameter $|h_{Mg}|$ of the $s$ valence electrons of $Mg$ is proportional to the first  ionization energy of $Mg$.
Then we have the following formula of $T_c$ of $MgB_2$:
\begin{equation}
k_BT_c=\kappa=\frac{1}{\sqrt{3}}\Delta_{MgB2}
\label{MgB2}
\end{equation}
where 
$\Delta_{MgB2}=h=6|h_{B1}|+6|h_{B2}| +|h_{Mg}|$ is the energy gap of 
   superconductivity of $MgB_2$, where $6=\frac{12}2$ is from the $12$  valence or nonvalence  $s$  electrons of the $B$ atoms in a hexagonal unit cell of $MgB_2$. From the existing table of ionization energy of $B$, the first and second ionization energies of $B$
are approximately equal to $800.6$ kJ/mol and
$2427.1$ kJ/mol respectively.
Let $|h_{B1}|\approx \xi 800.6$ kJ/mol
and $|h_{B2}|\approx \xi 2427.1$ kJ/mol, 
 $|h_{Mg}|\approx \xi 737.7$ kJ/mol. Thus from (\ref{MgB2}) we have:
\begin{equation}
T_c \approx 39.53K \quad \mbox{(Computed value of $T_c$ of $MgB_2$)}
\label{MgB25}
\end{equation}
This agrees with  the experimental value of $T_c\approx 39 K$ of $MgB_2$ \cite{Nag2}.

We notice that the energy gap $\Delta_{MgB2}=h=6|h_{B1}|+6|h_{B2}|+|h_{Mg}|$ contains the sum of two energy gaps $\Delta_{1}=h=6|h_{B1}|$ and $\Delta_{2}=h=6|h_{B2}|$.
We have that
$\Delta_{1}$  correspond to the conventional 3D superconductivity 
while 
$\Delta_{2}$ correspond to the unconventional quasi-2D superconductivity.
In this case we have the phenomeon of two-gap superconductivity \cite{Nag2,Cha,Liu,Cho,Iav,Bou,Sol}.

Then, from experiments we have two energy gaps: $\Delta_{\sigma}(4.2K)\approx 7.1 meV$ and $\Delta_{\pi}(4.2K)\approx 2.3 meV$ where $\Delta_{\sigma}$ is a 2D energy gap and  $\Delta_{\pi}$ is a 3D energy gap \cite{Nag2,Cha,Liu,Cho,Iav,Bou,Sol}.
Thus, $\Delta_{\sigma}$ corresponds to $\Delta_{2}$ and  $\Delta_{\pi}$ corresponds to
$\Delta_{1}$. Thus we have:
\begin{equation}
\frac{\Delta_{2}}{\Delta_{1}}=\frac{2427.1}{800.6}\approx 3.03
\, (\mbox{Computed value of}\, \frac{\Delta_{\sigma}(4.2K)}{\Delta_{\pi}(4.2K)})
\label{MgB27}
\end{equation}
This agrees with the experimental value of  $\frac{\Delta_{\sigma}(4.2K)}{\Delta_{\pi}(4.2K)}=\frac{7.1}{2.3}\approx 3.08$ \cite{Nag2,Cha,Liu,Cho,Iav,Bou,Sol}.

\subsection{\bf Energy gap and pseudogaps}

For the high-$T_c$ superconductivity we
 have the phase line
$h^2=3\kappa^2$ of  superconductivity. In this case
the energy  $\Delta(0)=|h|>0$ is also a superconducting energy gap when this energy is the remained energy after the forming of Cooper pairs.

On the other hand we also have other states of phase transitions of the form $h^2=c\kappa^2$ for some constants $c>0$ such as $c=\frac13$ for antiferromagnetic state. As the state of superconductivity the parameters $h$ of these states may also be as energy gaps. However since these states are not the state of superconductivity we have that these energy gaps are identified as pseudogaps in the sense that these energy gaps are not the energy gap of the state of superconductivity.
These pseudogaps are in the region of pseudogap of phase diagram such as the phase diagram in Fig.1.

 This agrees with experiments that there are phenomena of pseudogaps which are not directly related to superconductivity \cite{Tra,Fuj,Ich,Eme3,Zaa5,Mac3,Bia,Kiv},\cite{Tan,Val, Ma, Man3}. 


\subsubsection{\bf Pseudogap of superconductivity: Coexistence of pseudogap and superconducting gap}

Let us then investigate the pseudogap phenomena which are directly related to superconductivity and giving the temperature ${\bf T^{*}}>T_c$ 
\cite{Tra,Fuj,Ich,Eme3,Zaa5,Mac3,Bia,Kiv},\cite{Tan,Val, Ma, Man3}. 

We first consider the case of the coexistence of pseudogap and superconducting gap.
Let us find out the mechanism of this coexistence.

As an example to illustrate this mechanism
let us consider again the cuprate $La_{2-x}Sr_xCuO_4$ ($0\leq x\leq 1$) (or $LSCO$). We have that, together with the forming of Cooper pair of the $d$-electrons of $Cu$ in the quisi-2D $CuO_2$ plane, the third state of ionization energy of $Cu$ in the quisi-2D $CuO_2$ plane
gives the energy gap of the high-$T_c$ superconductivity of $LSCO$ by the channel opening connecting to the 3D region of conventional superconductivity.
This is as the mechanism of forming the energy gap of high-$T_c$ superconductivity of $LSCO$.

Thus the states of ionization energies of $Cu$ in the $CuO_2$ plane are related to the gap phenomenon of
high-$T_c$ superconductivity of $LSCO$.

 Let us then consider the fourth state of ionization energy of $Cu$ in the $CuO_2$ plane. If there were no the third state of ionization energy of $Cu$, then since this fourth state of ionization energy of $Cu$ is the higher state most closed to the state of ionization energy of $Cu$ of channel opening, this fourth state of ionization energy of $Cu$ would play the role of the third state of ionization energy of $Cu$. In this case the fourth state of ionization energy of $Cu$ in the quisi-2D $CuO_2$ plane would give the quisi-2D high-$T_c$ superconductivity of $LSCO$ by the channel opening connecting to the 3D region of conventional superconductivity.

Now since the third state of ionization energy of $Cu$ exists and is next to the state of channel opening we have that the superconducting $d$-electrons of $Cu$ in the quisi-2D $CuO_2$ plane are gathered in the third state
to give the quisi-2D high-$T_c$ superconductivity of $LSCO$ and that there are only a few $d$-electrons of $Cu$ in higher states such as the fourth state of ionization energy of $Cu$.

This means that
the fourth state (and higher states) of ionization energy of $Cu$ are potential states of giving  superconducting gap to $LSCO$ in the sense that these states exist for giving superconductivity (and giving the superconducting energy gap) but only a few superconducting $d$-electrons have enough energy to reach to these states.
In this case the stable critical temperature $T_c$ is mainly given by the energy gap of the third state of ionization energy of $Cu$, as computed in the example of $LSCO$.

Thus, together with the Cooper pairing mechanism and the doping mechanism of giving high-$T_c$ superconductivity, the fourth state (and higher states) of ionization energy of $Cu$  constitute the mechanism of giving pseudogap phenomenon to $LSCO$.

Further, we notice that since the fourth ionization energy of $Cu$ is approximately equal to $5536$ kJ/mol which is larger than the third ionization energy $\approx 3555$ kJ/mol of $Cu$ we have that the pseudogap is larger than the superconducting gap and the coresponding temperature $T^{*}$ of pseudogap is larger than $T_c$ (We shall give more details of this point). This agrees with experiments that the temperature $T^{*}$ of pseudogap is larger than $T_c$.

Furthermore,
at the underdoped range,
because of the shape of the phase line $h^2=3\kappa^2$ in the $(T,x)$-phase plane and that $5536>3555$,
 the maximal $T^{*}$ appears in the underdoped range while the $T_c>0$ begins to appear by the bifurcation of the phase line $h^2=3\kappa^2$ in the underdoped range.

 Then as the doping parameter $x$ further increases, because of the shape of the phase line $h^2=3\kappa^2$ in the $(T,x)$-phase plane,
 the $T^{*}$ begins to decrease while $T_c$ increases in the underdoped range.

 Then as the doping parameter $x$ increases to the optimal doping
 we have that $T_c$ increases to the maximal value by the bifurcation of the phase line $h^2=3\kappa^2$ while $T^{*}$ decreases and approaches to $T_c$. This then gives the nodal-antinodal phenomenon of gap and pseudogap.
Thus the energy gap-pseudogap phenomenon of $LSCO$ is a two-gap phenomenon where the two gaps are with similar structure.

We remark that this two-gap phenomenon is similar to the two-gap phenomenon of $MgB_2$ of which the two gaps both contribute to the superconductivity of $MgB_2$ while in the energy gap-pseudogap phenomenon of $LSCO$ one gap is for superconductivity and the other gap is a (virtually) potential state for superconductivity in a sense as described in the above.

We can then extend this pseudogap mechanism of $LSCO$ as the general mechanism of pseudogaps of meterials with high-$T_c$ superconductivity. We have that this general mechanism of pseudogaps  generally consists of the states of ionization energies (of an element $M$) following the state of ionization energy (of the element $M$) giving the high-$T_c$ superconductivity; the doping mechanism; and the Cooper pairing mechanism of high-$T_c$ superconductivity.

\subsubsection{\bf The temperature ${\bf T^{*}}$ and ${\bf T_P}$}

In the above we have shown that the range of paramagnetic Meissner effect is given by the following range:
\begin{equation}
3\kappa_c^2\leq h^2 <(2+\sqrt{3})^2\kappa_c^2
\label{pseudo}
\end{equation}
where $3\kappa_c^2= h^2$ is the state and the phase of high-$T_c$ superconductivity and the state $h^2 =(2+\sqrt{3})^2\kappa_c^2$ is the state for the cross-over to the phase of normal metallic state. In the phase diagram of the cuprates such as $La_{2-x}Sr_{x}CuO_{4}$ this range (\ref{pseudo}) is narrow \cite{Lia,Bra,Tim}.

Thus there is the phenomenon that the phase of normal metallic state is in competition with the phase $3\kappa_c^2= h^2$ of high-$T_c$ superconductivity  \cite{Lia,Bra,Tim}.
Thus the state $h^2 =(2+\sqrt{3})^2\kappa_c^2$ is for the above pseudogap related to the energy gap of the phase $3\kappa_c^2= h^2$ of high-$T_c$ superconductivity.

Furthermore below the line {\bf e} the region defined by the range (\ref{pseudo}) is the region of Type I superconductivity.  Thus the state $h^2 =(2+\sqrt{3})^2\kappa_c^2$ can be for the above pseudogap related to the energy gap of superconductivity.

Thus from the point of view of Type I superconductivity this state
$h^2 =(2+\sqrt{3})^2\kappa_c^2$ is also a state for an energy gap while the state $h^2 =3\kappa_c^2$ is already a state for an energy gap.  Thus the state $h^2 =(2+\sqrt{3})^2\kappa_c^2$ can be as the state of the above pseudogap related to the energy gap of high-$T_c$ superconductivity.

Let us also call the region specified by (\ref{pseudo}) as the extended region of pseudogap. This region is usually called the region of non-Fermi liquid.

From the range (\ref {pseudo}) of extended pseudogap
we have the following formula of temperature $T^{*}$ given by:
\begin{equation}
k_BT^{*}=\frac{1}{c_0} h^{*}
\label{pseudo2}
\end{equation}
where $\sqrt{3}\leq c_0\leq 2+\sqrt{3}$.

Let the pseudogap $\Delta^{*}(0)=h^{*}$. Then the temperature $T^{*}$ is given by:
\begin{equation}
\frac{2\Delta^{*}(0)}{k_BT^{*}}=2c_0
\label{pseudo3}
\end{equation}
where $\sqrt{3}\leq c_0\leq 2+\sqrt{3}$.

This agrees with the experimental value of $T^{*}$ which is given by \cite{Mei2,Ido}:
\begin{equation}
\frac{2\Delta^{*}(0)}{k_BT^{*}}\approx 4\sim 8 \quad (\mbox{Experimental value})
\label{pseudo4}
\end{equation}

Then, for $LSCO$, we have that $\frac{\Delta^{*}(0)}{\Delta(0)}\approx \frac{5536}{3555}$ gives a range of $T^{*}>T_c$. This agrees with the experimental results \cite{Mei2,Ido}.

On the other hand let a temperature $T_P$ related to pseudogap be given by the following equation:
\begin{equation}
\frac{2\Delta_P(0)}{k_BT_P}=2(2+\sqrt{3})
\label{pseud5}
\end{equation}
where the pseudogap $\Delta_P(0)=\Delta(0)$. Then we have $T_P<T_c$. This agrees with experimental value
\cite{Tra,Fuj,Ich,Eme3,Zaa5,Mac3,Bia,Kiv},
 \cite{Tan,Val}:
\begin{equation}
\frac{2\Delta(0)}{k_BT_P}\approx 4\sim 8 \quad (\mbox{Experimental value})
\label{pseud6}
\end{equation}

\subsubsection{\bf Coexistence of pseudogap and superconducting gap}

The pseudogap $\Delta^{*}(0)$ of the state (\ref{pseudo3}) and the superconducting energy gap $\Delta(0)$ of the state $h^2=3\kappa^2$ of superconductivity are both from the forming of Cooper pairs and the doping mechanism of superconductivity. Thus the two gaps coexist. This agrees with experimental results \cite{Ma,Man3}.

\subsubsection{\bf Pseudogap of superconductivity: Existence of pseudogap when superconductivity suppressed}

Let us then consider the case of existence of pseudogap  when the superconductivity is suppressed in the process of doping. Let us consider the example $LBCO$ which is similar to $LSCO$ with $Ba$ replacing $Sr$. For this example it is known that the state of superconductivity is firstly suppressed in the process of doping at $x=\frac18$ \cite{Tra}. Since the state of high-$T_c$ superconductivity is 
suppressed at $x=\frac18$ we have that the opening channel is just closed at $x=\frac18$.
 This channel closed gives no supporting 3D superconducting region for the stable existence of the quasi-2D high-$T_c$ superconductivity.

Then, 
though the opening channel is  closed, the pairing mechanism can still hold (In this case we shall not call this pairing mechanism as the Cooper pairing mechanism where we reserve the name Cooper pairing for this pairing mechanism when this pairing mechanism gives high-$T_c$ superconductivity).

This
pairing mechanism thus gives the existence of pseudogap given by the fourth state (and higher states) of ionization energy of $Cu$. Also since the state of superconductivity is just suppressed we have that the superconducting gap given by the third state of ionization energy of $Cu$ becomes a pseudogap as that of the fourth state (and higher states) of ionization energy of $Cu$. This gives the existence of pseudogap when the state of superconductivity is suppressed where the mechanism of forming pseudogap is just the mechanism of forming superconducting gap. This agrees with experimental results \cite{Val}.

This resolves the problem of the relation between pseudogaps and the onset of high-$T_c$ superconductivity.

\subsection{\bf Stripe phenomena}

The phase of the charge density waves is with the range $\kappa_c^2<h^2<3\kappa_c^2$ (which is a part of the phase of pseudogap). At the state $\kappa_c^2=h^2$ because of the condition (\ref{u17b}), the wave effect of $h$ is suppressed. Then for the range $\kappa_c^2<h^2<3\kappa_c^2$ the wave effect of $h$ appears while the static magnetic effect is remained with $\kappa_c>0$. This wave effect of $h$ and the static magnetic effect of $\kappa_c$ gives the phenomenon of wave interference.

Then in the two dimensional case this phenomenon of wave interference gives the stripe phenomenon \cite{Tra,Fuj,Ich,Eme3,Zaa5,Mac3,Bia,Kiv}; where the charge stripe is from the case that the wave effect of $h$ is not suppressed by the interference while the static magnetic stripe is from the case that the wave effect of $h$ is suppressed by the interference.

\subsection{\bf Checkerboard phenomena}

As a symmetry to the phase of charge density waves we have that the phase of spin density waves (SDW) is with the range $\frac13\kappa_c^2<h^2<\kappa_c^2$.
Then for this range $\frac13\kappa_c^2<h^2<\kappa_c^2$ the wave effect of $h$ appears weakly while the static magnetic effect is remained with $\kappa_c>0$. This wave effect of $h$ and the static magnetic effect of $\kappa_c$ gives  interference wave phenomenon.

Since $h$ gives attractive interaction between electrons (by the seagull vertex term) when $h^2<\kappa_c^2$ we have that the attractive interaction between the two pairing electrons
is weaker than the case $h^2>\kappa_c^2$ for charge density waves. Thus in this case of $h^2<\kappa_c^2$ the two pairing electrons
separately give stripe phenomena. Then the superposition  of these two stripe phenomena  gives the checkerboard phenomenon.

We notice that this phase of SDW is a phase of pseudogaps given by
$\kappa_c^2=h^2$ and $\kappa_c^2=3h^2$. Experimentally this phase of pseudogaps is with the checkerboard phenomenon \cite{Han2,Mce,Che,Pol,Zhu2,Che2,Pol2}. Thus the range $\frac13\kappa_c^2<h^2<\kappa_c^2$ identified as the phase of SDW with the checkerboard phenomenon agrees with the experiments that the checkerboard phenomenon appears in this phase of pseudogaps.

\subsection{\bf Vortex with checkerboard phenomena}

When vortex appears in the range of superconductivity, the vortex core must be with the condition $h^2<c^2\kappa_c^2$ for some $c^2<3$ or with the condition
$h^2>c^2\kappa_c^2$ for some $c^2>3$;
 and there are no further conditions on the states in the vortex core. Thus it is possible for the vortex core to be filled with electrons in the state of antiferromagnetism which is with the condition
$h^2=\frac13\kappa_c^2$ (This phase line is bifurcated to form the region of antiferromagnetism). Then since the normal metallic state is with the condition $h^2>(2+\sqrt{3})^2\kappa_c^2$ which is more dynamical than the antiferromagnetic state we have that for stability the antiferromagnetic state appears more frequently in the vortex core than the normal metallic state. This agrees with the experiments of high-$T_c$ superconductivity \cite{And3,And6,Boe,Hil}.

In case that the vortex core is in the region  of the antiferromagnetic phase, the region surrounding the vortex is then a region of spin density waves since the region of spin density waves is next to the region of antiferromagnetism. In this case the checkerboard phenomenon appears in the region surrounding the vortex, as observed in the experiments of high-$T_c$ superconductivity \cite{Han2,Mce,Che,Pol,Zhu2,Che2,Pol2}.

\subsection{\bf For fixed $\bf{\kappa_c}$ the pressure $\bf{P}$ corresponding to Planck parameter $\bf{h}$}

 It is shown that at a fixed doping in the phase plane $(T,P)$ ($P$ denotes pressure) a cuprate has a phase diagram
which is similar to the phase diagram of the cuprate in the phase plane $(T,x)$ \cite{Cuk}.

In this phase diagram in the phase plane $(T,P)$ as $P$ increasing we have the effect of the increasing of $h$, as similar to that the increasing of $x$ has the effect of the increasing of $h$. In particular the phase line $h^2=3\kappa^2$ of superconductivity appears in this phase diagram, where a bifurcation of forming the region of superconductivity appears at the pressure $P=21 GPa$ \cite{Cuk}.

Each phase in this phase diagram in the phase plane $(T,x)$ also corresponds to a crystal structure of the cuprate.
Thus the relation between pressure $P$ and the parameter $h$ is that the varying of $P$ giving the various phases of crystal structure of the cuprate (or a general material) and each phase of crystal structure of the cuprate (or the general material) corresponds to a range of the parameter $h$. Thus by varying pressure the cuprate (or the general material) can go through all the phases of crystal structure of the cuprate (or the general material).

\subsection{\bf Ferromagnetic superconductivity}

 Similar to the Case 1) of Type II conventional superconductivity, because of the term $h^2$ in the above equations of states and that the spins of the two electrons of a Cooper pair are given by the signs of $\pm h$ of $h^2$, besides antiferromagnetic superconductivity we also have ferromagnetic superconductivity for a material of unconventional (quai-2D) superconductivity.
As the Case 1),
a triplet Cooper pair has the effect of increased energy. Thus, as the Case 1), the triplet Cooper pair is a mechanism for reaching the state of ferromagnetic superconductivity with the state condition $\kappa_c^2=\frac13 h_c^2$
from the ferromagnetic state with the condition $\kappa^2> h^2$ where $h^2\approx 0$ (which means that $h^2$ is very small while $h^2>0$).
Then, as the Case 1), since the state $\kappa_c^2=\frac13 h_c^2$ of ferromagnetic superconductivity is reached from the condition
$h^2\approx 0$ by the mechanism of triplet Cooper pair, we have that $h_c^2$ is also small. Thus
the critical temperature $T_c=\frac{1}{k_B}\kappa_c$ of ferromagnetic superconductivity should be low.

\section{\bf Examples of  electron-doped high-$\bf{T_c}$ superconductivity}

In the above section the superconductors are of hole-doped. In this section let us consider some electron-doped examples.

\subsection{\bf Computation of $\bf{T_c}$ of $\bf{Nd_{2-x}Ce_xCuO_4}$}

Let us consider the electron-doped cuprate $Nd_{2-x}Ce_xCuO_4$ ( $0 \leq x \leq 1$).
For the doping mechanism of  superconductivity,
let us consider   the following function:
 \begin{equation}
f(x)=2130(1-2\cdot\frac{x}2)+1949\cdot \frac{x}2
 \label{e1}
\end{equation}
 where $2130$ kJ/mol and $1949$ kJ/mol are approximately the
third ionization energies of $Nd$ and $Ce$ respectively. Since $1949< 2130$ and $h$ is proportional to the ionization energy we have that this function gives the relation that the increasing of $x$ gives the decreasing of $h$ for $Nd_{2-x}Ce_xCuO_4$. 
We have $1949<1957.9< 2130$ where $1957.9$ kJ/mol is  the second ionization energy of $Cu$.

We notice that there is a factor $2$ of $2\cdot\frac{x}2$. This factor $2$ comes from that there are two kinds of doping on $Nd$: the  $d$-electron doping and the  $f$-electron doping.
We shall show that this factor $2$  gives a shift of the dome shape region of superconductivity such that this region of superconductivity overlaps with the region of antiferromagnetism in the phase diagram of   $Nd_{2-x}Ce_xCuO_4$.


Then we let: 
\begin{equation}
f(x_0)=2130(1-2\cdot\frac{x_0}2)+1949\cdot \frac{x_0}2=1957.9  
\label{e2}
\end{equation}
for some $x_0$ ($0< x_0< 1$).

 When this relation holds (or approximately holds), 
 the channel connecting the state of second ionization energy of  $d$-valence electrons of $Cu$ and   the state of third ionization energy of  $d$-valence electrons of $Ce$
is open. This channel opening gives
a freedom of electric current with
a direction orthogonal to the $CuO_2$ plane.
By this freedom of electric current, the electrons of a cluster of atoms in a unit cell can be coupled to the electrons of other clusters to form Cooper pairs. In this way
the Cooper pairs of $d$-valence electrons of $Cu$ are formed.
Thus this channel opening gives the 3D region of conventional superconductivity. From  this 3D conventional superconductivity we have the existence of quasi-2D  bifurcation region  of  unconventional high-$T_c$ superconductivity given by the $CuO_2$ plane.

We notice that $x_0\approx 0.148  <1$. 
Since  the function $f$ is a decreasing function of $x$,
 we have that $x_0$ is as a doping point  such that the range $0<x_1<x\leq x_0$ for some $x_1$  is a range of superconductivity and the range  $ x_0<x<1$ is as the range that $T_c$ as a function of $x$ is decreasing for $ x_0<x<1$. This means that   $x_0$ is the optimal doping point of superconductivity of $Nd_{2-x}Ce_xCuO_4$ . This agrees with the experiments that the optimal doping of superconductivity of  $Nd_{2-x}Ce_xCuO_4$ is at $x_0\approx  0.15$.

Also we notice that the  factor $2$ gives a shift of the dome shape region of superconductivity such that the optimal doping point  $x_0$ is shifted to the value $0.148$.

In this channel opening of electron-doped case, the $d$-valence electron of $Ce$ are not free to be in the   state of third  ionization energy  due to the  $f$-valence electron of $Nd$. Thus
 the $d$-valence electrons of $Cu$ are not pushed to higher state  of  ionization energy (i.e. the state of third ionization energy)  by the  the $d$-valence electron of $Ce$, as the case of $La_{2-x}Sr_xCuO_4$. Thus the $d$-valence electrons of $Cu$ are remained in the basic  state of second ionization energy.

Thus we have the following formula of $T_c$ of  $Nd_{2-x}Ce_xCuO_4$:
\begin{equation}
k_BT_c=\kappa=\frac{1}{\sqrt{3}}\Delta_{NCCO}
\label{e3}
\end{equation}
where 
$\Delta_{NCCO}=h=5|h_{Cu1}| +2|h_{Nd}|$ is the energy gap of 
   superconductivity of $Nd_{2-x}Ce_xCuO_4$, where $5=\frac{10}2$ is from the $10$  valence  $3d$-electrons of the $Cu$ atom in a  unit cell of $Nd_{2-x}Ce_xCuO_4$, and 
$|h_{Cu1}|\approx \xi 1957.9$ kJ/mol,
 $|h_{Nd}|\approx \xi 1040$ kJ/mol  where 1040  kJ/mol is the second  ionization energy of  $Nd$. The gap $2|h_{Nd}|$ is from the $s$-valence electrons of the two $Nd$ atoms in a unit cell of  $Nd_{2-x}Ce_xCuO_4$.
Thus from (\ref{e3}) we have:
\begin{equation}
T_c \approx 23.6 K \quad \mbox{(Computed value of $T_c$ of $Nd_{2-x}Ce_xCuO_4$)}
\label{e4}
\end{equation}
This agrees with  the experimental value of $T_c\approx 24 K$ of $Nd_{2-x}Ce_xCuO_4$ \cite{Nai,Nai2,Saw}.

\subsection{\bf  Computation of $\bf{T_c}$ of $\bf{Ln_{2-x}Ce_xCuO_4}$}

Let us consider the electron-doped cuprates $Ln_{2-x}Ce_xCuO_4$ ( $0 \leq x \leq 1$) where $Ln= Pr, Sm, Eu$.
For  $Ln= Pr, Sm, Eu$, the doping mechanism of  superconductivity is the same as that of $Nd_{2-x}Ce_xCuO_4$,  with $2130$ replaced by $2086,2260, 2404$ where $2086,2260, 2404$  kJ/mol are the third ionization energies of   $ Pr, Sm, Eu$ respectively.
Then the results are similar to the case of $Nd_{2-x}Ce_xCuO_4$. 

For $Ln=Pr$ we have that
the optimal doping $x_0\approx 0.152$. This agrees with  the experimental value  $x_0\approx 0.15$ of  $Pr_{2-x}Ce_xCuO_4$.
Then we have the following formula of $T_c$ of  $Nd_{2-x}Ce_xCuO_4$:
\begin{equation}
k_BT_c=\kappa=\frac{1}{\sqrt{3}}\Delta_{PCCO}
\label{ep3}
\end{equation}
where 
$\Delta_{NCCO}=h=5|h_{Cu1}| +2|h_{Pr}|$ is the energy gap of 
   superconductivity of $Pr_{2-x}Ce_xCuO_4$, 
$|h_{Cu1}|\approx \xi 1957.9$ kJ/mol,
 $|h_{Pr}|\approx \xi 1020$ kJ/mol  where 1020  kJ/mol is the second  ionization energy of  $Pr$. The gap $2|h_{Pr}|$ is from the $s$-valence electrons of the two $Pr$ atoms in a unit cell of  $Pr_{2-x}Ce_xCuO_4$.
Thus from (\ref{ep3}) we have:
\begin{equation}
T_c \approx 23.6 K \quad \mbox{(Computed value of $T_c$ of $Pr_{2-x}Ce_xCuO_4$)}
\label{ep6}
\end{equation} 
This agrees with  the experimental value of $T_c\approx 24 K$ of $Pr_{2-x}Ce_xCuO_4$ \cite{Nai,Nai2,Saw}.

 For $Ln=Sm$ we have that
the optimal doping $x_0\approx 0.14$. This agrees with  the experimental value  $x_0\approx 0.15$ of  $Sm_{2-x}Ce_xCuO_4$.
Then we have the following formula of $T_c$ of  $Sm_{2-x}Ce_xCuO_4$:
\begin{equation}
k_BT_c=\kappa=\frac{1}{\sqrt{3}}\Delta_{SCCO}
\label{em3}
\end{equation}
where 
$\Delta_{SCCO}=h=5|h_{Cu1}| $ is the energy gap of 
   superconductivity of $Sm_{2-x}Ce_xCuO_4$. This is due to 
the Cooper pair of  the $s$-valence electrons of the two $Sm$ atoms  have not been formed, because of the stronger antiferromagnetic effect of $Sm_{2-x}Ce_xCuO_4$ on the shifted region of superconductivity of $Sm_{2-x}Ce_xCuO_4$.
Thus from (\ref{em3}) we have:
\begin{equation}
T_c \approx 19.5 K \quad \mbox{(Computed value of $T_c$ of $Sm_{2-x}Ce_xCuO_4$)}
\label{em6}
\end{equation} 
This agrees with  the experimental value of $T_c\approx 19 K$ of $Sm_{2-x}Ce_xCuO_4$ \cite{Nai,Nai2,Saw}.

There is another experimental value of $T_c\approx 16 K$ of $Sm_{2-x}Ce_xCuO_4$. This can be explained by that,  due to  the stronger antiferromagnetic effect of $Sm_{2-x}Ce_xCuO_4$ on the shifted region of superconductivity of $Sm_{2-x}Ce_xCuO_4$,  the Cooper  pairing numbers of $d$-electrons of $Cu$ is changed from 5 to 4. In this case we have $\Delta_{SCCO}=h=4|h_{Cu1}| $. Then we have:
\begin{equation}
T_c \approx 15.6 K \quad \mbox{(Computed value of $T_c$ of $Sm_{2-x}Ce_xCuO_4$)}
\label{em7}
\end{equation} 
This agrees with  the experimental value of $T_c\approx 16 K$ of $Sm_{2-x}Ce_xCuO_4$.

 For $Ln=Eu$ we have that
the optimal doping $x_0\approx 0.12$. This agrees with  the experimental value  $x_0\approx 0.15$ of  $Eu_{2-x}Ce_xCuO_4$.
Then we have the following formula of $T_c$ of  $Eu_{2-x}Ce_xCuO_4$:
\begin{equation}
k_BT_c=\kappa=\frac{1}{\sqrt{3}}\Delta_{ECCO}
\label{en3}
\end{equation}
where 
$\Delta_{ECCO}=h=3|h_{Cu1}| $ is the energy gap of 
   superconductivity of $Eu_{2-x}Ce_xCuO_4$. This is due to that the Cooper  pairing numbers of $d$-electrons of $Cu$ is changed from 5 to 3, and that the Cooper pair of  the $s$-valence electrons of the two $Eu$ atoms  have not been  formed, due to the stronger antiferromagnetic effect of $Eu_{2-x}Ce_xCuO_4$ on the shifted region of superconductivity of $Eu_{2-x}Ce_xCuO_4$.
Thus from (\ref{en3}) we have:
\begin{equation}
T_c \approx 12 K \quad \mbox{(Computed value of $T_c$ of $Eu_{2-x}Ce_xCuO_4$)}
\label{en6}
\end{equation} 
This agrees with  the experimental value of $T_c\approx 13 K$ of $Eu_{2-x}Ce_xCuO_4$ \cite{Nai,Nai2,Saw}.

\subsection{\bf  Computation of $\bf{T_c}$ of $\bf{La_{2-x}Ce_xCuO_4}$}

Let us consider the electron-doped cuprate $La_{2-x}Ce_xCuO_4$ ( $0 \leq x \leq 1$). 
Consider   the following function:
 \begin{equation}
f(x)=1850(1-\frac{x}2)+3547\cdot \frac{x}2
 \label{ea1}
\end{equation}
 where $1850$ kJ/mol and $3547$ kJ/mol are approximately the
third  and fourth ionization energies of $La$ and $Ce$ respectively. 
Let us set: 
\begin{equation}
f(x_0)=1850(1-\frac{x_0}2)+3547\cdot\frac{x_0}2=1957.9  
\label{ea2}
\end{equation}
for some $x_0$ ($0\leq x_0\leq 1$).

 When this relation holds (or approximately holds), 
 the channel connecting the state of second ionization energy of  $d$-valence electrons of $Cu$ and   the state of third ionization energy of  $d$-valence electrons of $La$
is open. This channel opening gives
a freedom of  $d$-valence electrons of $Cu$ with
a direction orthogonal to the $CuO_2$ plane.
By this freedom of electric current,
the Cooper pairs of $d$-valence electrons of $Cu$ are formed.
Thus this channel opening gives the 3D region of conventional superconductivity. From  this 3D conventional superconductivity we have the existence of quasi-2D  bifurcation region  of  unconventional high-$T_c$ superconductivity given by the $CuO_2$ plane.

We notice that $x_0\approx 0.12$. 
This agrees with the experiment that the state of superconductivity of $La_{2-x}Ce_xCuO_4$  is in the doping range $0.075 \leq x<0.15$ \cite{Zue}.

Then  in the channel opening of this electron-doped case,
due to the doping of the $f$-valence electron of $Ce$,  the 
$T$-structure of $La_{2}CuO_4$ is changed  to the $T'$-structure of $La_{2-x}Ce_xCuO_4$ \cite{Nai,Nai2}. 
 In this case 
the $d$-valence electron of $La$ is
 free in the
state of third ionization energy and  in the direction of the  $CuO_2$ plane.
Thus  the $d$-valence electrons of $Cu$  in the
state of second ionization energy and  in the direction of the  $CuO_2$ plane can
 be pushed to the higher state of third ionization energy by the $d$-valence electron of $La$ in the state of third ionization energy, as the case of $La_{2-x}Sr_xCuO_4$,  for the quasi-2D superconductivity. Then the $d$-valence electrons of $Cu$  in the direction orthogonal to the $CuO_2$ plane are remained in the state of second ionization energy for the conventional 3D superconductivity.

Then we have the following formula of $T_c$ of $La_{2-x}Ce_xCuO_4$:
\begin{equation}
k_BT_c=\kappa=\frac{1}{\sqrt{3}}\Delta_{LCCO}
\label{ea3}
\end{equation}
where 
$\Delta_{LCCO}=h=3|h_{Cu1}|+2|h_{Cu}| +2|h_{La}|$ is the energy gap of 
   superconductivity of $La_{2-x}Ce_xCuO_4$,  and  we determine that  two pairs of $d$-valence electrons of $Cu$ are in the state of third ionization energy  and three  pairs  of $d$-valence electrons of $Cu$ are in the state of second ionization energy.

Thus from (\ref{ea3}) we have:
\begin{equation}
T_c \approx 29.5 K \quad \mbox{(Computed value of $T_c$ of $La_{2-x}Ce_xCuO_4$)}
\label{ea4}
\end{equation}
This agrees with  the experimental value of $T_c\approx 30 K$ of $La_{2-x}Ce_xCuO_4$ \cite{Nai,Nai2,Saw,Zue}.

\subsection{\bf Computation of $\bf{T_c}$ of $\bf{Sr_{1-x}La_xCuO_2}$}

Let us consider the electron-doped cuprate $Sr_{1-x}La_xCuO_2$ ( $0 \leq x \leq 1$).
For the doping mechanism of  superconductivity,
let us consider   the following function:
 \begin{equation}
f(x)=5500(1-x)+5940x
 \label{er1}
\end{equation}
 where $5500$ kJ/mol and $5940$ kJ/mol are approximately the
fourth ionization energy of $Sr$ and the fifth ionization energy of $La$ respectively. 

Then we let: 
\begin{equation}
f(x_0)=5500(1-x_0)+5940x_0=5535  
\label{er2}
\end{equation}
for some $x_0$ ($0< x_0< 1$) where $5535$ kJ/mol  is  approximately the
fourth ionization energy of $Cu$.

We notice that $x_0\approx 0.0795 $. 
This agrees with the experiments that the range of superconductivity of  $Sr_{1-x}La_xCuO_2$ is at $0\leq x_0 \leq  0.125$ \cite{Tan9}.

 When this relation holds (or approximately holds), 
 the channel connecting the state of fourth ionization energy of  $3s, 3p$-electrons of $Cu$ and   the state of  fourth ionization energy of  $4s, 4p$-electrons of $Sr$
is open. This channel opening gives
a freedom of  $3s, 3p$-electrons of $Cu$ with
a direction orthogonal to the $CuO_2$ plane.
By this freedom of electric current, the $3s, 3p$-electrons of $Cu$ in a unit cell can be coupled to the  $3s, 3p$-electrons of $Cu$ of other unit cell to form Cooper pairs. In this way
the Cooper pairs of $3s,3p$-electrons of $Cu$ are formed.
Thus this channel opening gives the 3D region of conventional superconductivity. From  this 3D conventional superconductivity we have the existence of quasi-2D  bifurcation region  of  unconventional high-$T_c$ superconductivity given by the $CuO_2$ plane.

Then we notice that, in the channel opening of this electron doped case, due to the $d$-valence electron of $La$, the  $4s, 4p$-electrons of $Sr$ are still 
not yet free to form Cooper pairs.
 In this case, the $3s,3p$-electrons of $Cu$ are not pushed to the higher state of fifth  ionization energy of $Cu$, by the $4s, 4p$-electrons of $Sr$. Then, in this case, the $3s$-electrons of $Cu$ are lowered to the state of third  ionization energy of $Cu$ for the 3D superconductivity while the
$3p$-electrons of $Cu$ are remained in the state of fourth  ionization energy of $Cu$ for the quasi 2D superconductivity.

Thus we have the following formula of $T_c$ of  $Sr_{1-x}La_xCuO_2$:
\begin{equation}
k_BT_c=\kappa=\frac{1}{\sqrt{3}}\Delta_{SLCO}
\label{er3}
\end{equation}
where 
$\Delta_{SLCO}=h=3|h_{Cu4}| + |h_{Cu}| +|h_{Sr}|$ is the energy gap of superconductivity of $Sr_{1-x}La_xCuO_2$, where $3=\frac{6}2$ is from the $6$  $3p$-electrons of the $Cu$ atom in a  unit cell of $Sr_{1-x}La_xCuO_2$, and 
$|h_{Cu4}|\approx \xi 5535$ kJ/mol,
 $|h_{Sr}|\approx \xi 1064.2$ kJ/mol  where 1064.2  kJ/mol is the second  ionization energy of  $Sr$ for $s$-valence electrons of $Sr$. The gap $|h_{Sr}|$ is from the $s$-valence electrons of the  $Sr$ atom in a unit cell of  $Sr_{1-x}La_xCuO_2$.
Thus from (\ref{er3}) we have:
\begin{equation}
T_c \approx 42.2 K \quad \mbox{(Computed value of $T_c$ of $Sr_{1-x}La_xCuO_2$)}
\label{er4}
\end{equation}
This agrees with  the experimental value of $T_c\approx 43 K$ of $Sr_{1-x}La_xCuO_2$  \cite{Smi,Tan9,Koj}.

\section{More examples of high-${\bf T_c}$ superconductivity}

Let us in this section  consider more hole-doped examples of high-$T_c$ superconductivity.

\subsection{\bf Computation of ${\bf T_c}$ of $\bf{HgBa_2CuO_{4+x}}$}

Let us first consider  the cuprate $HgBa_2CuO_{4+x}$ ($0\leq x\leq 1$).
Similar to $La_{2-x}Sr_x CuO_4$
the Cooper pairs of valence electrons of $Cu$ (and $Hg$) is as a mechanism for reducing $h^2$ to reach the state $\kappa^2=3h^2$ of superconductivity.
It is known that $HgBa_2CuO_{4+x}$ has the multi-layer form $BaO/HgO_{x}/BaO/CuO_2/BaO$.
In a unit cell of $HgBa_2CuO_{4+x}$ the $Cu$ atom and $Hg$ atom at the same (centered) $c$-axis
joining the  $CuO_2$ plane and the $HgO_{x}$ plane in the multi-layer $HgO_{x}/BaO/CuO_2$ form a cluster. This cluster of atoms may be regarded as a larger atom.
For the multi-layer $HgO_{x}/BaO/CuO_2$,
let us first consider the $HgO_{x}$ plane. The interaction between the $Hg$ atoms and $O$ atoms gives an interaction between two $Hg$ atoms.
With this additional interaction the seagull vertex interaction gives an attractive effect for the forming of the Cooper pair of $d$ valence electrons (and $s$ valence electrons) of $Hg$ where $Hg$ is with outer shells $5d^{10}6s^2$.
The two $d$ electrons for the forming of Cooper pair are with wave factors $e^{\pm ih_{Hg2}}$ respectively in their wave functions where $|h_{Hg2}|$ is the energy of the $d$ electrons for the interaction of forming the Cooper pair in the $HgO_{x}$ plane (Similarly two $s$ electrons for the forming of Cooper pair are with wave factors $e^{\pm ih_{Hg1}}$).
Then,
the $d$ valence electrons of $Hg$ have only two basic states of first and second ionization energies. Thus the maximum value of the energy parameter $|h_{Hg2}|$ of the $d$ valence electrons of $Hg$ is proportional to the second ionization energy of $Hg$ (Similarly the energy parameter $|h_{Hg1}|$ is proportional to the first ionization energy parameter of $Hg$).
This energy parameter $|h_{Hg2}|$ gives a wave factor $e^{i(h_{Hg2}-h_{Hg2})}$ of the Cooper pair of $d$ electrons of $Hg$. This cancelation of energy for the forming of Cooper pair gives the reduction of energies of the $Hg$ atoms. Thus the forming of Cooper pair of $d$ electrons of $Hg$ is a mechanism for the reduction of energy.

Let us then consider the $CuO_2$ plane. As the cuprate $La_{2-x}Sr_x CuO_4$, the interaction between $Cu$ atoms and $O$ atoms
gives an interaction between two $Cu$ atoms. With this additional interaction the seagull vertex interaction gives an attractive effect for the forming of the Cooper pair of $d$ electrons.
As $La_{2-x}Sr_x CuO_4$,
the two $d$ electrons for the forming of Cooper pair are with the wave factors $e^{\pm ih_{Cu}}$ respectively in their wave functions where $|h_{Cu}|$ is the energy of the $d$ electrons for the interaction of forming the Cooper pair.

As $YBa_2Cu_2O_{6+x}$ let us find out the mechanism of doping giving superconductivity.
 Let us consider the following function:
 \begin{equation}
 f(x)=3555-3600(1-x)
\label{HgBaCO83}
\end{equation}
where $3555$ kJ/mole and $3600$ kJ/mole are the third ionization energies of $Cu$ and $Ba$ respectively.

 This function gives the relation that the increasing of $x$ giving the increasing of $h$. Then we set up the following relation
 (A main point for setting up this relation is that $3555-3300$ is small and $3600>3555$):
\begin{equation}
f(x_0)=3555-3600(1-x_0) =3300x_0
\label{HgBaCO8}
\end{equation}
for some $x_0$ ($0<x_0< 1$) where $3300$ kJ/mol is the third ionization energy of $Hg$.  

When this relation holds  (or approximately holds) we have that the channel for the high-$T_c$ superconductivity given by the $CuO_2$ plane is open. This channel opening gives a freedom of electric current with
a direction orthogonal to the $CuO_2$ plane.
From this freedom of electric current, the Cooper pairs of $d$ valence electrons of $Cu$ and the Cooper pairs of $d$ valence electrons of $Hg$ can be formed.
Thus this channel opening  together with the $CuO_2$ plane gives the region of 3D superconductivity. From this region of 3D superconductivity we have the existence of quasi-2D bifurcation region of the unconventional high-$T_c$ superconductivity given by the $CuO_2$ plane.

We notice that $x_0\approx 0.15$.
Thus
$HgBa_2CuO_{4+x}$ comes into the range  of high-$T_c$ superconductivity when $x_0\leq x\leq x_1$ for some $x_1$ such that $x_0< x_1\leq 1$.  This agrees with the experiment that
 $HgBa_2CuO_{4+x}$ is in the state of high-$T_c$ superconductivity for $x_2\leq x\leq x_1$ for some $x_2\geq 0.15$ \cite{Put}.

When (\ref{HgBaCO8}) holds we have that the $d$ valence electrons of $Cu$ are in the basic state of third ionization energy, and other states are to be reached from this state.
Also when (\ref{HgBaCO8}) holds the $d$ valence electrons of $Cu$ and $Hg$ are unified to occupy the sequence of state of ionization energies. Since, comparing the $CuO_2$ plane and the $HgO_x$ plane, the $Cu$ atoms in the $CuO_2$ plane are with one more ionization direction than the $Hg$ atoms in the $HgO_x$ plane we have that the $d$ valence electrons of $Cu$ occupy the higher states while the $d$ valence electrons of $Hg$ occupy the lower states.
Then since the $d$ valence  electrons of $Hg$ has two states of ionization energies (and the states of first and second ionization energies can be reached from the third state) and that the $d$ valence electrons of $Cu$ have three states of ionization energies
we have that the $d$ valence electrons of $Hg$ have the states of first and second ionization energies while the $d$ valence electrons of $Cu$ have the states of third, fourth and fifth ionization energies.
Thus the maximum value of the energy parameter $|h_{Hg2}|$ of the $d$ valence electrons of $Hg$ is proportional to the second ionization energy of $Hg$ and the maximum value of the energy parameter $|h_{Cu}|$ is proportional to the fifth ionization energy of $Cu$.

On the other hand the $s$ valence electrons of $Hg$ have the state of first ionization energy. Thus the maximum value of the energy parameter $|h_{Hg1}|$ of the $s$ valence electrons of $Hg$ is proportional to the first ionization energy of $Hg$.

We can now specify the parameter $h$ of $h^2=3\kappa^2$.
As $La_{2}Sr CuO_4$,
each $d$ valence electron of a (centered) $Cu$ atom (of a unit cell) in the $CuO_2$ planes gives an energy $|h_{Cu}|$.
Then each $d$ valence electron of a (centered) $Hg$ atom gives an energy $|h_{Hg2}|$.
Then each $s$ valence electron of a (centered) $Hg$ atom gives an energy $|h_{Hg1}|$.
Thus the total energy of forming Cooper pairs of $d$ (and $s$) electrons of a cluster of $Cu$ atom and $Hg$ atom is $10|h_{Cu}|+10|h_{Hg2}|+2|h_{Hg1}|$.

Let $h_1$ be the energy reduced from this total energy by the mechanism of Cooper pairing.
On the other hand
let $h_2=10|h_{Cu}|+10|h_{Hg2}|+2|h_{Hg1}|-h_1$ be the energy after the reduction of energy.
Then as $La_{2-x}Sr_x CuO_4$ we have that $h_1=h=h_2$. Thus we have the following formula of $T_c$ of $HgBa_2CuO_{4+x}$:
\begin{equation}
k_BT_c=\kappa=\frac{1}{\sqrt{3}}\Delta_{HgBCO}
\label{HgBCO}
\end{equation}
where $\Delta_{HgBCO}=h=5|h_{Cu}|+5|h_{Hg2}|+|h_{Hg1}|$ is the energy gap of
 superconductivity of $HgBa_2CuO_{4+x}$ and $T_c$ is the highest critical temperature of superconductivity.

From the existing table of ionization energies we have that the
 fifth ionization energy of $Cu$, and the
 first and second ionization energies of $Hg$ are approximately equal to $7700$ kJ/mol, $1007.1$ kJ/mol and $1810$ kJ/mol respectively.
Thus we let $|h_{Cu}|\approx \xi 7700$ kJ/mol, $|h_{Hg1}|\approx \xi 1007.1$ kJ/mol and $|h_{Hg2}|\approx \xi 1810$ kJ/mol.
Then from (\ref{HgBCO}) we can compute the highest critical temperature $T_c$ of $HgBa_2CuO_{4+x}$:
\begin{equation}
T_c \approx 95.41 K \,\, \mbox{(Computed $T_c$ of $HgBa_2CuO_{4+x}$)}
\label{HgBCO5}
\end{equation}
This agrees with  the experimental value of $T_c\approx 94 K$ of $HgBa_2CuO_{4+x}$ \cite{Put}.

In case that the Cooper $ss$-pair of $Hg$ does not occur, we have $\Delta_{HgBCO}=h=5|h_{Cu}|+5|h_{Hg2}|$.
Then from (\ref{HgBCO}) the highest
$T_c$ of $HgBa_2CuO_{4+x}$ is given by:
\begin{equation}
T_c \approx 93.49 K \,\, \mbox{(Computed value of $T_c$ of $HgBa_2CuO_{4+x}$)}
\label{HgBCO6}
\end{equation}
This also agrees with  the experimental value of $T_c\approx 94 K$ of $HgBa_2CuO_{4+x}$.

\subsection{\bf Computation of $\bf{T_c}$ of $\bf{HgBa_2CaCu_2O_{6+x}}$}

Let us then
 consider
$HgBa_2CaCu_2O_{6+x}$ ($0\leq x\leq 1$) with a multi-layer
 $HgO_{x}/BaO/CuO_2/Ca/CuO_2$. The $BaO$ planes are as external to this multi-layer.
 We suppose that the layer $CuO_2/Ca/CuO_2$
of this cuprate can still be regarded as a quasi 2D material such that the above theory of high-$T_c$ superconductivity can be applied.
In a unit cell the $Cu$ atoms, the $Ca$ atom  and $Hg$ atom at the centered $c$-axis joining the two $CuO_2$ planes, the $Ca$ plane and the $HgO_{x}$ plane in this multi-layer form a cluster. This cluster of atoms may be regarded as a larger atom.

Let us consider the two $CuO_2$ planes. Because of the existence of $Ca$ between the two $CuO_2$ planes the $d$ valence electrons of $Cu$ in a $CuO_2$ plane can only have one direction of ionization  orthogonal to this $CuO_2$ plane. Thus the $d$ valence electrons of $Cu$ can only have two states of ionization energies.
Then, as analogous to the cuprate $HgBa_2CuO_{4+x}$,
the doping mechanism for $HgBa_2CaCu_2O_{6+x}$ is as that of $HgBa_2CuO_{4+x}$ with the same function $f$ and the relation (\ref{HgBaCO8}).

When this relation (\ref{HgBaCO8}) holds we have that the channel for the high-$T_c$ superconductivity given by the two $CuO_2$ planes is open. In this case of channel opening
the Cooper pairs of $d$ valence electrons of $Cu$, the Cooper pairs of $d$ valence electrons of $Hg$
and the bifurcation region of high-$T_c$ superconductivity can be formed.
When (\ref{HgBaCO8}) holds we have that the $d$ valence electrons of $Hg$ and the $d$ valence electrons of $Cu$ in the $CuO_2$ plane next to
the $BaO$ plane (which is next to the $HgO_{x}$ plane) are in the basic state of third ionization energy and other states are to be reached from this state.

Further the $d$ valence electrons of $Hg$ and the $d$ valence electrons of $Cu$ in the $CuO_2$ plane next to
the $BaO$ plane are unified to occupy the sequence of state of ionization energies.
Since the $Cu$ atoms in the $CuO_2$ plane are with one more ionization direction than the $Hg$ atoms in the $HgO_x$ plane we have that the $d$ valence electrons of $Cu$ occupy the higher states while the $d$ valence electrons of $Hg$ occupy the lower states.
Thus the $d$ valence electrons of $Cu$ in this $CuO_2$ plane
have the state of fourth ionization energy since these $d$ valence electrons can have two states of ionization energies.
Thus the maximum value of the energy parameter $|h_{Cu}|$ is proportional to the fourth ionization energy of $Cu$. In this case let us denote this energy parameter $|h_{Cu}|$ by $|h_{Cu4}|$.

On the other hand the $CuO_2$ plane which is not next to the $BaO$ plane is separated from the $HgO_{x}$ and $BaO$ planes by a $Ca$ plane.
Because of this separation the $d$ valence electrons of $Cu$ in this $CuO_2$ plane is without the ionization direction connecting to the $HgO_{x}$ and $BaO$ planes. Thus the $d$ valence electrons of $Cu$ in this $CuO_2$ plane can have the basic state of third ionization energy and the lower state of second ionization energy. Thus the maximum value of the energy parameter $|h_{Cu}|$ is proportional to the third ionization energy of $Cu$.

Then since the $d$ valence electrons of $Hg$ in the $HgO_{x}$ plane have two
directions of ionization with $O$ (with one direction of ionization orthogonal to the $HgO_{x}$ plane) we have that the $d$ valence electrons of $Hg$ have two states of
ionization energies.
Then the $d$ valence electrons of $Hg$ can not jump from the state of
second ionization energy of $1810$ kJ/mol to the state of first ionization energy of $1007.1$ kJ/mol (or vise versa) when a $s$ valence electron of $Ca$ is in the state of second ionization energy of $1145.4$ kJ/mol. Thus the $d$ valence electrons of $Hg$ have to occupy the states of second and third ionization energies. Thus the maximum value of the energy parameter $|h_{Hg2}|$ of the $d$ valence electrons of $Hg$ is proportional to the third ionization energy of $Hg$.

We can now specify the parameter $h$ of $h^2=3\kappa^2$.
As $La_{2-x}Sr_x CuO_4$,
each $d$ electron of a (centered) $Cu$ atom  gives an energy parameter $|h_{Cu}|$ (For the $d$ electron of a $Cu$ atom of the $CuO_2$ plane which is next to the $BaO$ plane, we denote the energy parameter by $|h_{Cu4}|$).
Then each $d$ electron of a (centered) $Hg$ atom gives an energy parameter $|h_{Hg2}|$.
Thus the total energy of $d$
electrons of a cluster of the two (centered) $Cu$ atoms and the $Hg$ atom is $10|h_{Cu}|+10|h_{Cu4}|+10|h_{Hg2}|$ (We consider the case that the Cooper $ss$-pairs of $Hg$ of $Hg$ have not been formed).

Let $h_1$ be the energy reduced from this total energy by the mechanism of Cooper pairing.
On the other hand
let $h_2=10|h_{Cu}|+10|h_{Cu4}|+10|h_{Hg2}|-h_1$ be the energy after the reduction of energy.
Then as $La_{2-x}Sr_x CuO_4$ we have that $h_1=h=h_2$. Thus we have the following formula of $T_c$ of $HgBa_2CaCu_2O_{6+x}$:
\begin{equation}
k_BT_c=\kappa=\frac{1}{\sqrt{3}}\Delta_{HgBCCO}
\label{HgBCCO}
\end{equation}
where $\Delta_{HgBCCO}=h=5|h_{Cu}|+5|h_{Cu4}|+5|h_{Hg2}|$ is the energy gap of
 $HgBa_2CaCu_2O_{6+x}$.
Let $|h_{Cu}|\approx \xi 3555$ kJ/mol, $|h_{Cu4}|\approx \xi 5536$ kJ/mol and $|h_{Hg2}|\approx \xi 3300$ kJ/mol where $5536$ kJ/mol is approximately the fourth ionization energies of $Cu$.
Then from (\ref{HgBCCO}) we can compute the highest critical temperature $T_c$ of $HgBa_2CaCu_2O_{6+x}$:
\begin{equation}
T_c \approx 121.63 K \,\, \mbox{(Computed $T_c$ of $HgBa_2CaCu_2O_{6+x}$)}
\label{HgBCCO5}
\end{equation}
This agrees with  the experimental value of $T_c\approx 120 K$ of $HgBa_2CaCu_2O_{6+x}$ \cite{Sch5,Gao}.

\subsection{\bf Computation of $\bf{T_c}$ of $\bf{Tl_2Ba_2CuO_{6+x}}$}

We then consider the cuprate $Tl_2Ba_2CuO_{6+x}$ ($-1\leq x\leq 1$).
Let us first find out the mechanism of doping giving superconductivity.
 Similar to
 $HgBa_2CuO_{4+x}$,
 let us consider the following function:
 \begin{equation}
f(x)=1971(1+\frac{x}{2}) +965.2(-\frac{x}{2})
 \label{Tl2BaCO83}
\end{equation}
 where $1971$ kJ/mol and $965.2$ kJ/mole are approximately the second ionization energies of $Tl$ and $Ba$ respectively.

 This function gives the relation that the increasing of $x$ giving the increasing of $h$. Then we set  the following relation:
\begin{equation}
f(x_0)=1971(1+\frac{x_0}{2}) +965.2(-\frac{x_0}{2}) =1957.9
\label{Tl2BaCO8}
\end{equation}
for some $x_0$ ($-1<x_0< 1$) where $1957.9$ kJ/mol is the second
ionization energy of $Cu$.

 When this relation holds  (or approximately holds) we have that the channel for the high-$T_c$ superconductivity given by the $CuO_2$ plane is open. In this case the Cooper pairs of $d$ valence electrons of $Cu$,  the Cooper pairs of $d$ valence electrons of $Tl$,
and the bifurcation region of high-$T_c$ superconductivity can be formed.

We notice that $x_0\approx -0.026$.
Thus $Tl_2Ba_2CuO_{6+x}$ comes into the range  of high-$T_c$ superconductivity when $x_0\leq x\leq x_1$ for some $x_1$ such that $x_0< x_1< 1$.  This agrees with the experiment that
 $Tl_2Ba_2CuO_{6+x}$ is in the state of high-$T_c$ superconductivity for $x_2\leq x\leq x_1$ for some $x_2\geq -0.026$ \cite{She5,Tor,Kan,Par,Par5}.

When (\ref{Tl2BaCO8}) holds we have that the $d$ valence electrons of $Cu$ and the $d$ valence electrons of $Tl$ are in the basic state of second ionization energy, and that other states are to be reached from this state. 

Further the $d$ valence electrons of $Cu$ and the $d$ valence electrons of $Tl$ are unified to occupy the sequence of state of ionization energies. Since the $Cu$ atoms in the $CuO_2$ plane are with one more ionization direction than the $Tl$ atoms in the $TlO_{1+\frac{x}{2}}$ plane we have that the $d$ valence electrons of $Cu$ occupy the higher states while the $d$ valence electrons of $Tl$ occupy the lower states.
Thus, as $La_{2-x}Sr_x CuO_4$, the $d$ valence electrons of $Cu$ have the states of second, third, and fourth ionization energies and the $d$ valence electrons of $Tl$ have the states of first and second ionization energies. 
Thus the maximum value of the energy parameter $|h_{Cu}|$ is proportional to the fourth ionization energy of $Cu$ and the maximum value of the energy parameters $|h_{Tl2}|$ for the $d$ valence electrons of $Tl$ is proportional to the second ionization energy of $Tl$.

Then it is clear that the energy parameters $|h_{Tl1}|$ for the $s$ valence electrons of $Tl$ is proportional to the first ionization energy of $Tl$.
As $HgBa_2CuO_{4+x}$, we have the following formula of $T_c$ of $Tl_2Ba_2CuO_{6+x}$:
\begin{equation}
k_BT_c=\kappa=\frac{1}{\sqrt{3}}\Delta_{Tl2BCO}
\label{Tl2BCO}
\end{equation}
where $\Delta_{Tl2BCO}=h=5|h_{Cu}|+10|h_{Tl2}|+2|h_{Tl1}|$ is the energy gap of superconductivity of $Tl_2Ba_2CuO_{6+x}$.
We have $|h_{Cu}|\approx \xi 5536$ kJ/mol where $5536$ kJ/mol is approximately the
fourth ionization energy of $Cu$,
$|h_{Tl2}|\approx \xi 1971$ kJ/mol and $|h_{Tl1}|\approx \xi 589.4$ kJ/mol.
Then from (\ref{Tl2BCO}) we can compute the highest critical temperature $T_c$ of $Tl_2Ba_2CuO_{6+x}$:
\begin{equation}
T_c \approx 95.2 K \,\, \mbox{(Computed $T_c$ of $TlBa_2CuO_{6+x}$)}
\label{Tl2BCO5}
\end{equation}
This roughly agrees with the experimental value of $T_c\approx 90K$ of $TlBa_2CuO_{6+x}$ \cite{She5,Tor,Kan,Par,Par5}.

In case of the nonappearing of the Cooper $ss$-pair of $Tl$ we have $\Delta_{Tl2BCO}=h=5|h_{Cu5}|+10|h_{Tl2}|$. Then we have:
\begin{equation}
T_c \approx 92.96 K \,\, \mbox{(Computed $T_c$ of $Tl_2Ba_2CuO_{6+x}$)}
\label{Tl2BCO9}
\end{equation}
This agrees with the experimental value of $T_c\approx 90 K$ of $Tl_2Ba_2CuO_{6+x}$.

\subsection{\bf Computation of $\bf{T_c}$ of $\bf{TlBa_2CaCu_2O_{7+x}}$}.
Let us then consider the cuprate $TlBa_2CaCu_2O_{7+x}$ ($-1\leq x\leq 1$).
As the cuprate $Tl_2Ba_2CuO_{6+x}$,
let us consider the follwing function:
 \begin{equation}
f(x)=1971(1+x) +965.2(-x)
 \label{TlBaCO83}
\end{equation}
 where $1971$ kJ/mol and $965.2$ kJ/mole are approximately the second ionization energies of $Tl$ and $Ba$ respectively.

 This function gives the relation that the increasing of $x$ giving the increasing of $h$. Then we set  the following relation:
\begin{equation}
f(x_0)=1971(1+x_0) +965.2(-x_0) =1957.9
\label{TlBaCO8}
\end{equation}
for some $x_0$ ($-1<x_0< 1$) where $1957.9$ kJ/mol is the second
ionization energy of $Cu$.  When this relation holds  (or approximately holds) we have that the channel for the high-$T_c$ superconductivity given by the $CuO_2$ plane is open.
In this case the Cooper pairs of $d$ valence electrons of $Cu$,  the Cooper pairs of $d$ valence electrons of $Tl$
and the bifurcation region of high-$T_c$ superconductivity can be formed.

We notice that $x_0\approx -0.013$.
Thus
 $TlBa_2CaCu_2O_{7+x}$ comes into the range  of high-$T_c$ superconductivity when $x_0\leq x\leq x_1$ for some $x_1$ such that $x_0< x_1< 1$.  This agrees with the experiment that $TlBa_2CaCu_2O_{7+x}$ is in the state of high-$T_c$ superconductivity for $x_2\leq x\leq x_1$ for some $x_2\geq -0.013$ \cite{She5,Tor,Kan,Par,Par5}.

When (\ref{TlBaCO8}) holds we have that the $d$ valence electrons of $Tl$ and the $d$ valence electrons of $Cu$ in the $CuO_2$ plane next to
the $BaO$ plane (which is next to the $TlO_{1+x}$ plane) are in the basic state of second ionization energy and other states are to be reached from this state.

Further the $d$ valence electrons of $Tl$ and the $d$ valence electrons of $Cu$ in the $CuO_2$ plane next to
the $BaO$ plane are unified to occupy the sequence of state of ionization energies. Since the $Cu$ atoms in the $CuO_2$ plane are with one more ionization direction than the $Tl$ atoms in the $TlO_{1+x}$ plane we have that the $d$ valence electrons of $Cu$ occupy the higher states while the $d$ valence electrons of $Tl$ occupy the lower states.
Thus, the $d$ valence electrons of $Tl$ have the states of first and second ionization energies, and
because there is a $Ca$ plane between the two $CuO_2$ planes, the $d$ valence electrons of $Cu$ in this $CuO_2$ plane have the states of second and third ionization energies.
Thus the maximum value of the energy parameter $|h_{Cu}|$ for the $d$ valence electrons of $Cu$ is proportional to the third ionization energy of $Cu$ and the maximum value of the energy parameters $|h_{Tl2}|$ for the $d$ valence electrons of $Tl$ is proportional to the second ionization energy of $Tl$.

On the other hand the $CuO_2$ plane which is not next to the $BaO$ plane is separated from the $TlO_{1+x}$ and $BaO$ planes by a $Ca$ plane.
The $d$ valence electrons of $Cu$ in this $CuO_2$ plane have the basic state of second ionization energy and the state of first or third ionization energies are to be reached from this state.

Then the $d$ valence electrons of $Cu$ can not jump from the state of
second ionization energy of $1957.9$ kJ/mol to the state of first ionization energy of $745.5$ kJ/mol when a $s$ valence electron of $Ca$ is in the state of second ionization energy of $1145.4$ kJ/mol. Thus the $d$ valence electrons of $Cu$ in this $CuO_2$ plane also have to occupy the states of second and third ionization energies.
Thus the maximum value of the energy parameter $|h_{Cu}|$ is proportional to the third ionization energy of $Cu$.

Then it is clear that the energy parameters $|h_{Tl1}|$ for the $s$ valence electrons of $Tl$ is proportional to the first ionization energy of $Tl$.
 Then, as $HgBa_2CaCu_2O_{6+x}$,
we have the following formula of $T_c$ of $TlBa_2CaCu_2O_{7+x}$:
\begin{equation}
k_BT_c=\kappa=\frac{1}{\sqrt{3}}\Delta_{TlBCCO}
\label{TlBCCO}
\end{equation}
where $\Delta_{TlBCCO}=h=10|h_{Cu}|+5|h_{Tl2}|+|h_{Tl1}|$ is the energy gap of $TlBa_2CaCu_2O_{7+x}$.
We have $|h_{Cu}|\approx \xi 3555$ kJ/mol,
$|h_{Tl2}|\approx \xi 1971$ kJ/mol and $|h_{Tl1}|\approx \xi 589.4$ kJ/mol where $589.4$ kJ/mol is approximately the
first ionization energy of $Tl$.
Then from (\ref{TlBCCO}) we can compute the highest critical temperature $T_c$ of $TlBa_2CaCu_2O_{7+x}$:
\begin{equation}
T_c \approx 90.40K \,\, \mbox{(Computed $T_c$ of $TlBa_2CaCu_2O_{7+x}$)}
\label{TlBCCO5}
\end{equation}
This agrees with the experimental value of $T_c\approx 90K$ of $TlBa_2CaCu_2O_{7+x}$ \cite{She5,Tor,Kan,Par,Par5}.

In case of the nonappearing of the Cooper $ss$-pair of $Tl$ we have $\Delta_{TlBCCO}=h=10|h_{Cu5}|+5|h_{Tl2}|$. Then we have:
\begin{equation}
T_c \approx 89.28 K \,\, \mbox{(Computed $T_c$ of $TlBa_2CaCu_2O_{7+x}$)}
\label{TlBCCO9}
\end{equation}
This also agrees with the experimental value of $T_c\approx 90 K$ of $TlBa_2CaCu_2O_{7+x}$.

\subsection{\bf Computation of $\bf{T_c}$ of $\bf{Tl_2Ba_2CaCu_2O_{8+x}}$}

Let us then consider the cuprate $Tl_2Ba_2CaCu_2O_{8+x}$ ($-1\leq x\leq 1$).
The doping mechanism of this cuprate
 is the same as that of $Tl_2Ba_2CuO_{6+x}$ with the same relation (\ref{Tl2BaCO8}).
 When this relation holds  (or approximately holds) we have that the channel for the high-$T_c$ superconductivity given by the $CuO_2$ planes is open. In this case the Cooper pairs of $d$ valence electrons of $Cu$, the Cooper pairs of $d$ valence electrons of $Tl$
and the bifurcation region of high-$T_c$ superconductivity can be formed.
Thus  $Tl_2Ba_2CaCu_2O_{8+x}$ comes into the range  of high-$T_c$ superconductivity when $x_0\leq x\leq x_1$ for some $x_1$ such that $x_0< x_1< 1$.  This agrees with the experiment that  $Tl_2Ba_2CaCu_2O_{8+x}$ is in the state of high-$T_c$ superconductivity for $x_2\leq x\leq x_1$ for some $x_2\geq -0.026$ \cite{She5,Tor,Kan,Par,Par5}.
When (\ref{Tl2BaCO8}) holds we have that the $d$ valence electrons of $Cu$ and the $d$ valence electrons of $Tl$ are in the basic state of second ionization energy, and that other states are to be reached from this state.

 Further the $d$ valence electrons of $Cu$ and the $d$ valence electrons of $Tl$ are unified to occupy the sequence of states of ionization energies. Since the $Cu$ atoms in the two $CuO_2$ planes are with one more ionization direction than the $Tl$ atoms in the $TlO_{1+\frac{x}{2}}$ plane we have that the $d$ valence electrons of $Cu$ occupy the higher states while the $d$ valence electrons of $Tl$ occupy the lower states.
Thus, the $d$ valence electrons of $Tl$ have the states of first and second ionization energies, and
because there is a $Ca$ plane between the two $CuO_2$ planes, the $d$ valence electrons of $Cu$ in the two $CuO_2$ planes have the states of second and third ionization energies.
Thus the maximum value of the energy parameter $|h_{Cu}|$ for the $d$ valence electrons of $Cu$ is proportional to the third ionization energy of $Cu$ and the maximum value of the energy parameters $|h_{Tl2}|$ for the $d$ valence electrons of $Tl$ is proportional to the second ionization energy of $Tl$.

Then it is clear that the energy parameters $|h_{Tl1}|$ for the $s$ valence electrons of $Tl$ is proportional to the first ionization energy of $Tl$.
Thus, as the above cuprates,
we have the following formula of $T_c$ of $Tl_2Ba_2CaCu_2O_{8+x}$:
\begin{equation}
k_BT_c=\kappa=\frac{1}{\sqrt{3}}\Delta_{Tl2BCCO}
\label{Tl2BCCO}
\end{equation}
where $\Delta_{Tl2BCCO}=h=10|h_{Cu}|+10|h_{Tl2}|+2|h_{Tl1}|$ is the energy gap 
of 
 $Tl_2Ba_2CaCu_2O_{8+x}$.

We have $|h_{Cu}|\approx \xi 3555$ kJ/mol,
$|h_{Tl2}|\approx \xi 1971$ kJ/mol and $|h_{Tl1}|\approx \xi 589.4$ kJ/mol. Then from (\ref{Tl2BCCO}) we can compute the highest critical temperature $T_c$ of $Tl_2Ba_2CaCu_2O_{8+x}$:
\begin{equation}
T_c \approx 110.92K \,\, \mbox{(Computed $T_c$ of $Tl_2Ba_2CaCu_2O_{8+x}$)}
\label{Tl2BCCO5}
\end{equation}
This agrees with the experimental value of $T_c\approx 110K$ of $Tl_2Ba_2CaCu_2O_{8+x}$ \cite{She5,Tor,Kan,Par,Par5}.

In case of the nonappearing of the Cooper $ss$-pair of $Tl$ we have $\Delta_{Tl2BCCO}=h=10|h_{Cu5}|+10|h_{Tl2}|$. Then we have:
\begin{equation}
T_c \approx 108.69K \,\, \mbox{(Computed $T_c$ of $Tl_2Ba_2CaCu_2O_{8+x}$)}
\label{Tl2BCCO9}
\end{equation}
This also agrees with the experimental value of $T_c\approx 110 K$ of $Tl_2Ba_2CaCu_2O_{8+x}$
\cite{She5,Tor,Kan,Par,Par5}.

\subsection{\bf Computation of $\bf{T_c}$ of $\bf{Tl_2Ba_2Ca_2Cu_3O_{10+x}}$}

We then consider the cuprate $Tl_2Ba_2Ca_2Cu_3O_{10+x}$ ($-1\leq x\leq 1$).
While the layers of this material is increased let us suppose that the conducting layer of this cuprate can still be regarded as a quasi 2D material such that the above theory of high-$T_c$ superconductivity can be applied.
In a unit cell the $Cu$ atoms and $Tl$ atoms at the centered $c$-axis joining the three $CuO_2$ planes and the two $TlO_{1+\frac{x}{2}}$ planes
 form a cluster. Then the three $Cu$ atoms and the two $Tl$ atoms  may be regarded as a larger atom.
As $Tl_2Ba_2CaCu_2O_{8+x}$,
the doping mechanism giving superconductivity of $Tl_2Ba_2Ca_2Cu_3O_{10+x}$ is the same as that of $Tl_2Ba_2CuO_{6+x}$ with the same relation (\ref{Tl2BaCO8}).
 When this relation holds  (or approximately holds) we have that the channel for the high-$T_c$ superconductivity given by the $CuO_2$ plane is open. In this case the Cooper pairs of $d$ valence electrons of $Cu$, the Cooper pairs of $d$ valence electrons of $Tl$
and the bifurcation region of high-$T_c$ superconductivity can be formed.
Thus  $Tl_2Ba_2Ca_2Cu_3O_{10+x}$ comes into the range  of high-$T_c$ superconductivity when $x_0\leq x\leq x_1$ for some $x_1$ such that $x_0< x_1< 1$.  This agrees with the experiment that $Tl_2Ba_2Ca_2Cu_3O_{10+x}$ is in the state of high-$T_c$ superconductivity for $x_2\leq x\leq x_1$ for some $x_2\geq -0.026$ \cite{She5,Tor,Kan,Par,Par5}.
As $Tl_2Ba_2CaCu_2O_{8+x}$, when (\ref{Tl2BaCO8}) holds we have that the $d$ valence electrons of $Cu$ in the two outer $CuO_2$ planes have the states of second and third ionization energies and the $d$ valence electrons of $Tl$ have the states of first and second ionization energies. Thus the maximum value of the energy parameter $|h_{Cu}|$ of the $d$ valence electrons of $Cu$ in the two outer $CuO_2$ planes is proportional to the third ionization energy of $Cu$ and the maximum value of the energy parameter $|h_{Tl2}|$ of the $d$ valence electrons of $Tl$ is proportional to the second ionization energy of $Tl$.

 Then for the inner $CuO_2$ plane between two $Ca$ planes, because there is no orthogonal ionization direction for $d$ valence electrons of $Cu$ in this $CuO_2$ plane, we have that the $d$ valence electrons of $Cu$ in this $CuO_2$ plane have only one state of ionization energy which is the basic state of second ionization energy. Thus the maximum value of the energy parameter $|h_{Cu}|$ for these $d$ valence electrons of $Cu$ is proportional to the second ionization energy of $Cu$. Let us denote this energy parameter $|h_{Cu}|$ by $|h_{Cu2}|$.

 Then the energy parameters $|h_{Tl1}|$ for the $s$ valence electrons of $Tl$ is proportional to the first ionization energy of $Tl$.
Then, as the above cuprates,
we have the following formula of $T_c$ of $Tl_2Ba_2Ca_2Cu_3O_{10+x}$:
\begin{equation}
k_BT_c=\kappa=\frac{1}{\sqrt{3}}\Delta_{Tl2BC2CO}
\label{Tl2BC2CO}
\end{equation}
where $\Delta_{Tl2BC2CO}=h=10|h_{Cu}|+5|h_{Cu2}|+10|h_{Tl2}|+2|h_{Tl1}|$ is the energy gap of high-$T_c$
 superconductivity of $Tl_2Ba_2Ca_2Cu_3O_{10+x}$.

We have $|h_{Cu}|\approx \xi 3555$ kJ/mol, $|h_{Cu2}|\approx \xi 1957.9$ kJ/mol,
 $|h_{Tl2}|\approx \xi 1971$ kJ/mol and $|h_{Tl1}|\approx \xi 589.4$ kJ/mol.
 Then from (\ref{Tl2BC2CO}) we can compute the highest critical temperature $T_c$ of $Tl_2Ba_2Ca_2Cu_3O_{10+x}$:
\begin{equation}
T_c \approx 130.15K \,\, \mbox{(Computed $T_c$ of $Tl_2Ba_2Ca_2Cu_3O_{10+x}$)}
\label{Tl2BC2CO5}
\end{equation}
This roughly agrees with the experimental value of $T_c\approx 127K$
of $Tl_2Ba_2Ca_2Cu_3O_{10+x}$ \cite{She5,Tor,Kan,Par,Par5}.

In case of the nonappearing of the Cooper $ss$-pair of $Tl$ we have $\Delta_{Tl2BC2CO}=h=10|h_{Cu}|+5|h_{Cu2}|+10|h_{Tl2}|$. Then we have:
\begin{equation}
T_c \approx 127.91K \,\, \mbox{(Computed $T_c$ of $Tl_2Ba_2Ca_2Cu_3O_{10+x}$)}
\label{Tl2BC2CO9}
\end{equation}
This agrees with the experimental value of $T_c\approx 127 K$ of $Tl_2Ba_2Ca_2Cu_3O_{10+x}$ \cite{She5,Tor,Kan,Par,Par5}.

\subsection{\bf Computation of $\bf{T_c}$ of $\bf{HgBa_2Ca_2Cu_3O_{8+x}}$}

Let us then consider the cuprate $HgBa_2Ca_2Cu_3O_{8+x}$ ($0\leq x\leq 1$).
While the layers of this material is increased let us suppose that the conducting layer of this cuprate can still be regarded as a quasi 2D material such that the above theory of high-$T_c$ superconductivity can be applied.
In a unit cell the $Cu$ atoms and $Hg$ atom at the centered $c$-axis joining the three $CuO_2$ planes and the $HgO_{x}$ plane
 form a cluster. Then the three $Cu$ atoms and the $Hg$ atom may be regarded as a larger atom.
Let us consider the three $CuO_2$ planes. Because of the existence of $Ca$ between the three $CuO_2$ planes the $d$ valence electrons of $Cu$ in the two outer $CuO_2$ planes can only have one direction of ionization  orthogonal to these two $CuO_2$ planes. Thus the $d$ valence electrons of $Cu$ can only have two states of ionization energies. Then since the inner $CuO_2$ plane has no orthogonal direction of ionization with $O$ we have that the $d$ valence electrons of $Cu$ of this inner $CuO_2$ plane has only one state of ionization energy.
Then, as $HgBa_2CuO_{4+x}$,
the doping mechanism of superconductivity for $HgBa_2Ca_2Cu_3O_{8+x}$ is as that of $HgBa_2CuO_{4+x}$ with the same function $f$ and the relation (\ref{HgBaCO8}).
When this relation (\ref{HgBaCO8}) holds we have that the channel for the high-$T_c$ superconductivity given by the three $CuO_2$ planes is open. In this case the Cooper pairs of $d$ valence electrons of $Cu$, the Cooper pairs of $d$ valence electrons of $Hg$
and the bifurcation region of high-$T_c$ superconductivity can be formed.
When (\ref{HgBaCO8}) holds we have that the $d$ valence electrons of $Cu$ in the outer $CuO_2$ plane next to
a $BaO$ plane (which is next to the $HgO_{x}$ plane) are in the basic state of third ionization energy, and other states are to be reached from this state.
Then, as similar to $HgBa_2CaCu_2O_{6+x}$, the $d$ valence electrons of $Cu$ in this $CuO_2$ plane had to occupy the states of third and  fourth ionization energies. Thus the maximum value of the energy parameter $|h_{Cu}|$ is proportional to the fourth ionization energy of $Cu$. In this case let us denote this parameter $|h_{Cu}|$ by $|h_{Cu4}|$.

On the other hand the other outer $CuO_2$ plane which is not next to the $BaO$ plane (which is next to the $HgO_{x}$ plane) is separated from the $HgO_{x}$ and the $BaO$ plane by a $Ca$ plane. Because of this separation the $d$ valence electrons of $Cu$ in this $CuO_2$ plane can jump from the state of third ionization energy to the state of second ionization energy. Thus the $d$ valence electrons of $Cu$ in this $CuO_2$ plane can only have the two states of second and third ionization energies. Thus the maximum value of the energy parameter $|h_{Cu}|$ is proportional to the third ionization energy of $Cu$.

Then for the inner $CuO_2$ plane since it is independent of the $HgO_{x}$ plane and has only one state of ionization energy we have that the $d$ valence electrons of $Cu$ in this $CuO_2$ plane is in the basic state of second ionization energy. Thus the maximum value of the energy parameter $|h_{Cu}|$ is proportional to the second ionization energy of $Cu$. In this case let us denote this parameter $|h_{Cu}|$ by $|h_{Cu2}|$.

As $HgBa_2CaCu_2O_{6+x}$, the highest state of $d$ valence electrons of $Hg$ in
 $HgBa_2Ca_2Cu_3O_{8+x}$ is the state of third ionization energy of $Hg$.
Thus the maximum value of the energy parameter $|h_{Hg2}|$ of the $d$ valence electrons of $Hg$ is proportional to the third ionization energy of $Hg$.
Thus, as the above cuprates, we have the following formula of $T_c$ of $HgBa_2Ca_2Cu_3O_{8+x}$ :
\begin{equation}
k_BT_c=\kappa=\frac{1}{\sqrt{3}}\Delta_{HgBCCO}
\label{HgBC2CO}
\end{equation}
where $\Delta_{HgBCCO}=h=5|h_{Cu}|+5|h_{Cu4}|+5|h_{Cu2}|+5|h_{Hg2}|$ is the energy gap of
of $HgBa_2Ca_2Cu_3O_{8+x}$ (We omit the effect of the Cooper $ss$-pairs of $Hg$).

We have $|h_{Cu}|\approx \xi 3555$ kJ/mol, $|h_{Cu2}|\approx \xi 1957.9$ kJ/mol, $|h_{Cu4}|\approx \xi 5536$ kJ/mol; and
  $|h_{Hg2}|\approx \xi 3300$ kJ/mol.
Then from (\ref{HgBCCO}) we can compute the highest critical temperature $T_c$ of $HgBa_2Ca_2Cu_3O_{8+x}$:
\begin{equation}
T_c \approx 141.06 K \,\, \mbox{(Computed $T_c$ of $HgBa_2Ca_2Cu_3O_{8+x}$ )}
\label{HgBC2CO5}
\end{equation}
This agrees with the experimental 
 highest $T_c\approx 140 K$ of $HgBa_2Ca_2Cu_3O_{8+x}$ \cite{Sch5,Gao}.

We remark that the descrepancy
$141.06-140=1.06 K$
may  due to that $HgBa_2Ca_2Cu_3O_{8+x}$ is just deviated from as a meterial having a quasi-2D conducting layer.

\section{Composition of high-${\bf T_c}$ superconductors}

Based on the above theory of  high-$ T_c$ superconductivity let us compose some  materials which are predicted to be high-${\bf T_c}$ superconductors.

\subsection{$\bf{CaCu_5}$-type  intermetallic ${\bf La_{1-x}Ca_{x}Cu_5}$}

As analogous to  the intermetallic $MgB_2$ of $AlB_2$-type,
let us  compose more  intermetallic high-$T_c$ superconductors.
Let us first
consider  an intermetallic $La_{1-x}Ca_{x}Cu_5$ ($0\leq x\leq 1$). 
We have that the intermetallics $LaCu_{5}$ and $CaCu_{5}$ are of the $CaCu_5$-type \cite{Chak, Sub}.
Thus 
 the crystal structure of $La_{1-x}Ca_{x}Cu_5$ can be formed in the $CaCu_5$-type with $Ca$  corresponding to $La_{1-x}Ca_{x}$. 
 
This $CaCu_5$-type is similar to the $AlB_2$-type with the hexagon 
of six $B$ atoms replaced by an enlarged hexagon of twelve $Cu$ atoms and a hexagon of six $Cu$ atoms is intercalated in one of the two $Ca$ planes which replaces one of the two $Al$ planes of $AlB_2$. There are eighteen $Cu$ atoms near the faces of a unit cell of the hexagonal structure of $CaCu_5$ (where twelve $Cu$ atoms are from the two hexagons of six $Cu$ atoms intercalated in the two $Ca$ planes and six $Cu$ atoms are from the enlarged hexagon of $Cu$). FIG.6 
 shows the crystal structure of two of the three layers of a unit cell of $La_{1-x}Ca_xCu_5$ of $CaCu_5$-type. FIG.6a is the middle layer and FIG.6b is the lower or upper layer. The $Cu$ atoms are depicted by the black balls and the $La(Ca)$ atoms are depicted by the white balls.

 \begin{figure}[hbt]
\centering
\includegraphics[scale=0.35]{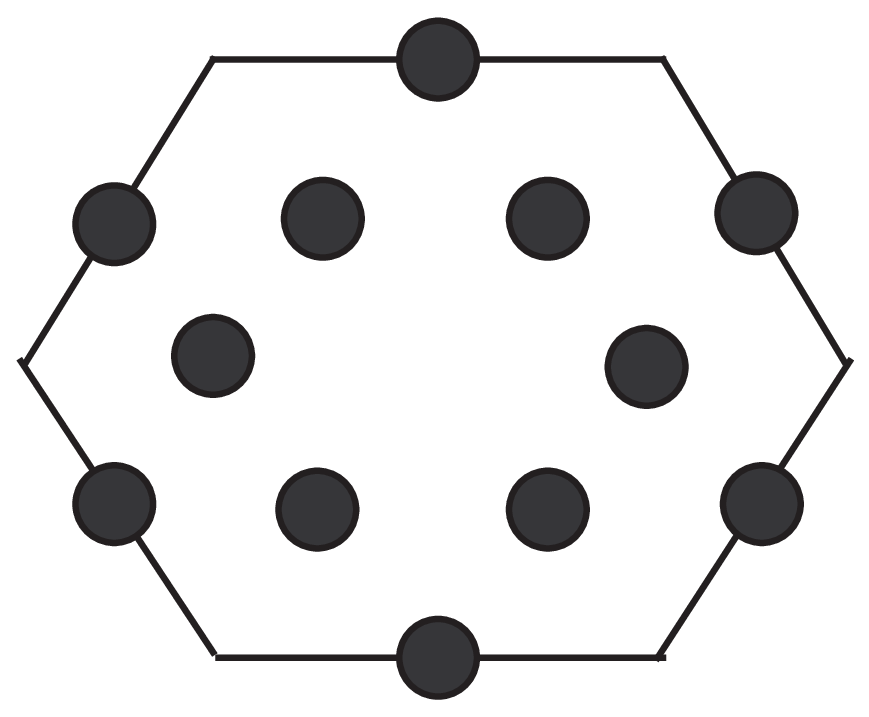}


FIG. 6a
 \end{figure}

 \begin{figure}[hbt]
\centering
\includegraphics[scale=0.35]{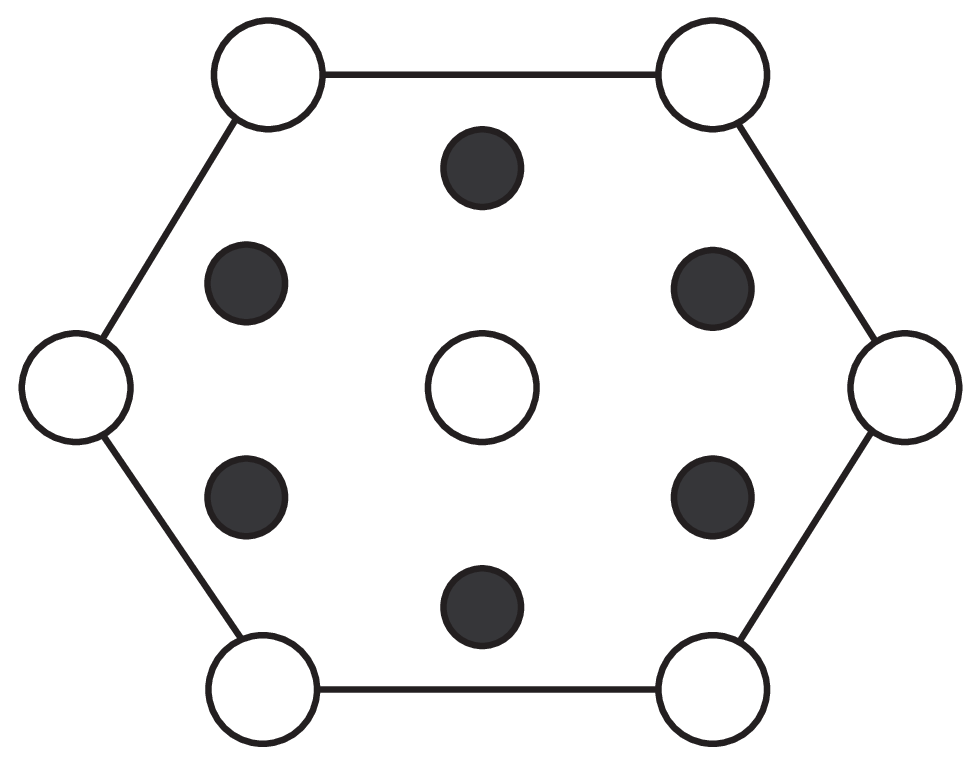}

            FIG. 6b
 \end{figure}
 



Thus in a unit cell of $CaCu_5$
nine $Cu$ atoms of the twelve $Cu$ atoms of the enlarged hexagon of $Cu$ and six $Cu$ atoms of the two hexagons of $Cu$ and one face-centered $Ca$ atom are as a cluster of atoms.
 
For the doping mechanism of superconductivity let us consider the following function:
\begin{equation}
f(x)=1850(1-x)+4912.4x
 \label{LaCaCu83}
\end{equation}
 where $1850$ kJ/mol and $4912.4$ kJ/mole are approximately the third ionization energies of $La$ and $Ca$ respectively. 
 This function gives the relation that the increasing of $x$ gives the increasing of $h$. Then we set the following relation of  the degenerate state of channel opening:
\begin{equation}
f(x_0)=1850(1-x_0)+4912.4x_0=1957.9
\label{LaCaCu8}
\end{equation}
for some $x_0$ ($0<x_0< 1$) where $1957.9$ kJ/mol is approximately the second
ionization energy of $Cu$ (A main point is that  $1850$ is close to $1957.9$).  When this relation holds (or approximately holds), the channel connecting the state of third ionization energy of $La$ and the state of second ionization energy of $Cu$ 
can be opened. This channel opening gives a freedom of electric current with
a direction orthogonal to the $Cu$ plane. 
By this freedom of electric current, the valence electrons of a cluster of atoms in a unit cell can be coupled to the valence electrons of  clusters of other unit cells to form Cooper pairs by the attractive electron-electron interaction
derived from the seagull vertex term. In this way 
the Cooper pairs of $d$-valence electrons of $Cu$ and the Cooper pair of the $s$-valence electrons of $La$ can be formed.
Thus this channel opening  together with the $Cu$ plane gives the region of 3D conventional superconductivity (as similar to the $\pi$-band of $MgB_2$). From this region of 3D superconductivity we have the existence of quasi-2D bifurcation region of the high-$T_c$ superconductivity given by the $Cu$ plane (as similar to the $\sigma$-band of $MgB_2$).

We notice that $x_0\approx 0.035$.
Thus $La_{1-x}Ca_{x}Cu_5$ comes into the range  of high-$T_c$ superconductivity when $ x_0<x<x_1$ for some $x_1$ such that $x_0< x_1\leq 1$.  

When (\ref{LaCaCu8}) holds  (or approximately holds) giving channel opening, the $d$-valence electrons of $Cu$ in the $Cu$ plane are in the basic state of second ionization energy and the  $d$ valence electron of $La$ is in the basic state of third ionization energy, and that other states are to be reached from these two states. Further the $d$-valence electrons of $Cu$ in the $Cu$ plane and the $d$-valence electrons of $La$ are unified to occupy a sequence of states such that the $d$-valence electrons of $Cu$ in the $Cu$ plane occupy the higher states while the $d$-valence electrons of $La$ occupy the lower state. 

Thus, as $MgB_2$, 
the $d$-valence electrons of $Cu$ are in the  state of third ionization energy, the $d$-valence electron of $La$ is in the state of third ionization energy;
and  
the $s$-valence electrons of $La$ are in the state of first ionization energy.

Thus,
the maximum value of the energy parameter $|h_{Cu}|$  of the $d$-valence electrons of $Cu$  are proportional to the third ionization energy of $Cu$, 
and the energy parameter $|h_{La}|$ of the $s$-valence electrons of $Y$ are proportional to the first ionization energy of $La$.
Then, as $MgB_2$, if we omit the effect of the hexagon of six $Cu$ atoms intercalated in the $La$ plane, then we have the following formula of $T_c$ of $La_{1-x}Ca_{x}Cu_5$:
\begin{equation}
k_BT_c=\kappa=\frac{1}{\sqrt{3}}\Delta_{LaCaCu} 
\label{LaCaCu}
\end{equation}
where $\Delta_{LaCaCu}=h=45|h_{Cu}|+|h_{La}|$ is the energy gap of the  superconductivity of $La_{1-x}Ca_{x}Cu_5$
where the coefficient $45=5\cdot9$ are from the $5=\frac{10}2$ for the ten $d$-valence electrons of the $Cu$, and the $9$ for the nine $Cu$ atoms in the enlarged hexagon of $Cu$.
(The effect of $s$-valence electrons of $Ca$ is omitted).

We have $|h_{Cu}| \approx \xi 3555$ kJ/mol, and
 $|h_{La}|\approx \xi 538.1$ kJ/mol where $3555$ kJ/mol is approximately the third ionization energy of $Cu$ and $538.1$ kJ/mol is approximately the first ionization energy of $La$.
Then from (\ref{LaCaCu}) we can compute the highest critical temperature $T_c$ of $La_{1-x}Ca_{x}Cu_5$ which can be upped to:
\begin{equation}
T_c \approx 315.6 K \,\, \mbox{(Computed $T_c$ of $La_{1-x}Ca_{x}Cu_5$)}
\label{LaCaCu5}
\end{equation}

If we include the effect of the hexagon of six $Cu$ atoms intercalated in the $La$ plane and suppose that this effect is the same as the enlarged hexagon of twelve $Cu$ atoms. Then the term $45|h_{Cu}|$ in $\Delta_{LaCaCu}$ becomes $75|h_{Cu}|$. Then the highest critical temperature $T_c$ of $La_{1-x}Ca_{x}Cu_5$ can be upped to:
\begin{equation}
T_c \approx 525.3 K \,\, \mbox{(Computed $T_c$ of $La_{1-x}Ca_{x}Cu_5$)}
\label{LaCaCu55}
\end{equation}

We remark that we may use other alkaline elements such as $Sr$ (with the third ionization energy $\approx 4138 $ kJ/mol)
to replace the element $Ca$ to form the intermetallic $La_{1-x}Sr_{x}Cu_5$ (We have that $SrCu_5$ is in the $CaCu_5$-type \cite{Chak, Sub}). 
The doping mechanism is similar to $La_{1-x}Ca_{x}Cu_5$ with $x_0\approx 0.047$. The highest critical temperature can be upped to $T_c \approx 525.3 K$. We have that $La_{1-x}Sr_{x}Cu_5$ is easier to be synthesized than $La_{1-x}Ca_{x}Cu_5$ since $4138<4912.4 $.

Further, we may use  a mixture of $Ca$, $Sr$ and $Ba$, denoted by $A$, to replace the element $Ca$ to form the intermetallic $La_{1-x}A_{x}Cu_5$ ($0\leq x\leq 1$).

Also we may use  a mixture of $Cu$ and $Zn$ to replace the element $Cu$ to form the intermetallic $La_{1-x}A_{x}Cu_5(1-y)Zn_{5y}$ ($0\leq y\leq 1$).

More generally  the rare earth element  $La$ 
may be replaced by  rare earth elements  or $Y$, or the mixture of rare earth elements such as $Mm$ (Misch metal). Also the $Cu$ and $Zn$ may be replaced by mixture of transition elements  to form  $CaCu_5$-type superconductors $R_{1-x}A_x TM_5$ where $R$ denotes a mixture of rare earth elements and $Y$ and  $TM$ denotes a mixture of  transition elements.


\subsection{ $\bf{CaCu_5}$-type intermetallic ${\bf Gd_{x}Ca_{1-x}Cu_5}$}
We may use other rare earth elements such as $Gd$ to replace the element $La$ to form the intermetallic $Gd_{x}Ca_{1-x}Cu_5$ ($0\leq x\leq 1$) (We have that $GdCu_5$ can be formed in the $CaCu_5$-type \cite{Chak, Sub}).
For the doping mechanism we consider the following function:
\begin{equation}
f(x)=1957.9- 1145.4(1-x)
 \label{GdCaCu83}
\end{equation}
 where $1957.9$ kJ/mol and $1145.4$ kJ/mole are approximately the second ionization energies of $Cu$ and $Ca$ respectively. 
 This function gives the relation that the increasing of $x$ gives the increasing of $h$. Then we set the following relation  of  the degenerate state of channel opening:
\begin{equation}
f(x_0)=1957.9- 1145.4(1-x_0)=1990x_0
\label{GdCaCu8}
\end{equation}
for some $x_0$ ($0<x_0< 1$) where $1990$ kJ/mol is approximately the third
ionization energy of $Gd$. 
Then we have $x_0\approx 0.96$.
Then $Gd_{x}Ca_{1-x}Cu_5$ comes into the range  of high-$T_c$ superconductivity when $ x_2<x<x_0$ for some $x_2$ such that $0\leq x_2<x_0< 1$.  
Then, as $La_{1-x}Ca_xCu_5$, we have the following formula of $T_c$ of $Gd_{x}Ca_{1-x}Cu_5$:
\begin{equation}
k_BT_c=\kappa=\frac{1}{\sqrt{3}}\Delta_{GdCaCu} 
\label{GdCaCu}
\end{equation}
where $\Delta_{GdCaCu}=h=75|h_{Cu}|+|h_{Gd}|$ is the energy gap of the  superconductivity of $Gd_{x}Ca_{1-x}Cu_5$,
 $|h_{Cu}| \approx \xi 3555$ kJ/mol, and
 $|h_{Gd}|\approx \xi 592.5$ kJ/mol where  $592.5$ kJ/mol is approximately the first ionization energy of $Gd$.
Then from (\ref{GdCaCu}) we can compute the highest critical temperature $T_c$ of $Gd_{x}Ca_{1-x}Cu_5$ which can be upped to:
\begin{equation}
T_c \approx 525.4 K \,\, \mbox{(Computed $T_c$ of $Gd_{x}Ca_{1-x}Cu_5$)}
\label{GdCaCu5}
\end{equation}

We remark that we may use other alkali elements such as $Sr$ (with the second ionization energy $\approx 1064.2 $ kJ/mol)
to replace the element $Ca$ to form the intermetallic $Gd_{x}Sr_{1-x}Cu_5$ ($0\leq x\leq 1$). The doping mechanism is similar to $Gd_{x}Ca_{1-x}Cu_5$ with $x_0\approx 0.96$. The critical temperature can be upped to $T_c \approx 525.4 K$.

\subsection{ $\bf{CaCu_5}$-type intermetallic ${\bf Y_{x}Ca_{1-x}Cu_5}$}
We may also use other rare earth elements such as $Y$ to replace the element $La$ to form the intermetallic $Y_{x}Ca_{1-x}Cu_5$ ($0\leq x\leq 1$) (We have that $YCu_5$ is in the $CaCu_5$-type \cite{Chak, Sub}).
For the doping mechanism we consider the following function:
\begin{equation}
f(x)=1957.9- 1145.4(1-x)
 \label{YCaCu83}
\end{equation}
 where $1957.9$ kJ/mol and $1145.4$ kJ/mole are approximately the second ionization energies of $Cu$ and $Ca$ respectively. 
 This function gives the relation that the increasing of $x$ gives the increasing of $h$. Then we set  the following relation  of  the degenerate state of channel opening:
\begin{equation}
f(x_0)=1957.9- 1145.4(1-x_0)=1980x_0
\label{YCaCu8}
\end{equation}
for some $x_0$ ($0<x_0< 1$) where $1980$ kJ/mol is approximately the third
ionization energy of $Y$. 
Then we have $x_0\approx 0.96$.
Then $Y_{x}Ca_{1-x}Cu_5$ comes into the range  of high-$T_c$ superconductivity when $ x_2<x<x_0$ for some $x_2$ such that $0\leq x_2<x_0< 1$.  
Then, as $La_{1-x}Ca_xCu_5$, we have the following formula of $T_c$ of $Y_{x}Ca_{1-x}Cu_5$:
\begin{equation}
k_BT_c=\kappa=\frac{1}{\sqrt{3}}\Delta_{YCaCu} 
\label{YCaCu}
\end{equation}
where $\Delta_{YCaCu}=h=75|h_{Cu}|+|h_{Y}|$ is the energy gap of the  superconductivity of $Y_{x}Ca_{1-x}Cu_5$,
 $|h_{Cu}| \approx \xi 3555$ kJ/mol, and
 $|h_{Y}|\approx \xi 600$ kJ/mol where  $600$ kJ/mol is approximately the first ionization energy of $Y$.
Then from (\ref{YCaCu}) we can compute the highest critical temperature $T_c$ of $Y_{x}Ca_{1-x}Cu_5$ which can be upped to:
\begin{equation}
T_c \approx 525.43 K \,\, \mbox{(Computed $T_c$ of $Y_{x}Ca_{1-x}Cu_5$)}
\label{YCaCu5}
\end{equation}

We remark that we may use other alkaline elements such as $Sr$ 
to replace the element $Ca$ to form the intermetallic $Y_{x}Sr_{1-x}Cu_5$ ($0\leq x\leq 1$). The doping mechanism is similar to $Y_{x}Ca_{1-x}Cu_5$ with $x_0\approx 0.96$. The critical temperature can be upped to $T_c \approx 525.43 K$.

\subsection{ $\bf{CaCu_5}$-type intermetallic ${\bf LaNi_{5-5x}Cu_{5x}}$}
Let us 
consider an intermetallic $LaNi_{5(1-x)}Cu_{5x}$ ($0\leq x\leq 1$). 
We have that  $LaCu_{5}$ and $LaNi_{5}$ are of the $CaCu_5$-type \cite{Chak, Sub}.
Thus 
 the crystal structure of $LaNi_{5(1-x)}Cu_{5x}$ can be formed in the $CaCu_5$-type with $Ca$  corresponding to $La$ and $Cu$  corresponding to $Ni_{(1-x)}Cu_{x}$. 
 For the doping mechanism of superconductivity let us 
consider the following function:
\begin{equation}
f(x)=1753(1-x)+1957.9x
 \label{LaNiCu83}
\end{equation}
 where $1753$ kJ/mol and $1957.9$ kJ/mole are approximately the second ionization energies of $Ni$ and $Cu$ respectively. 
 This function gives the relation that the increasing of $x$ gives the increasing of $h$. Then we set  the following relation of  the degenerate state of channel opening:
\begin{equation}
f(x_0)=1753(1-x_0)+1957.9x_0=1850
\label{LaNiCu8}
\end{equation}
for some $x_0$ ($0<x_0< 1$) where $1850$ kJ/mol is approximately the third
ionization energy of $La$.  When this relation holds (or approximately holds) we have that the channel connecting the state of third ionization energy of $La$ and the state of second ionization energy of $Ni(Cu)$ 
can be opened. 
In this case of channel opening 
the Cooper pairs of $d$-valence electrons of $Ni(Cu)$, the Cooper pairs of $s$-valence electrons of $La$,
and the bifurcation region of high-$T_c$ superconductivity can be formed.

We notice that $x_0\approx 0.4734$.
Thus $LaNi_{5(1-x)}Cu_{5x}$ comes into the range  of high-$T_c$ superconductivity when 
$x_0\leq x\leq x_1$ for some $x_1$ such that $ x_0<x_1\leq 1$. Further we notice that $1753$ is close to $1850$. Thus it is possible that $LaNi_{5(1-x)}Cu_{5x}$ comes into the range  of high-$T_c$ superconductivity for  $0\leq x\leq 0.4734$. 

We remark that for this intermetallic $LaNi_{5(1-x)}Cu_{5x}$ we have another doping relation giving superconductivity. Let us consider the following function:
\begin{equation}
g(x)=1753(1-x)+3555x
 \label{YNiCu83a}
\end{equation}
 where $3555$ kJ/mole is approximately the third ionization energy of $Cu$. 
 This function also gives the relation that the increasing of $x$ gives the increasing of $h$. Then we set  the following relation of  the degenerate state of channel opening:
\begin{equation}
g(x_0^{\prime})=1753(1-x_0^{\prime})+3555x_0^{\prime}=1850
\label{YNiCu8a}
\end{equation}
for some $x_0^{\prime}$ ($0<x_0^{\prime}< 1$) where $1850$ kJ/mol is again approximately the third
ionization energy of $La$ (The ionization energy of $Cu$ can be the second or the third ionization energies).  When this relation holds (or approximately holds),  the channel connecting the state of third ionization energy of $La$ and the state of second ionization energy of $Ni(Cu)$ can be opened. Then we have $x_0^{\prime}\approx 0.0538$.
Thus $LaNi_{5(1-x)}Cu_{5x}$ comes into the range  of high-$T_c$ superconductivity when 
$x_0^{\prime}\leq x\leq x_1$ for some $x_1$ such that $ x_0^{\prime}<x_1\leq 1$. We notice that $x_0^{\prime}$ is more close to $0$ than $x_0$.

Since $x_0^{\prime}\approx 0$, we have that $LaNi_5$ can be formed in the degenerate state of channel opening. Thus  $LaNi_5$  can be formed as a superconductor.

When (\ref{LaNiCu8}) holds (or when (\ref{YNiCu8a}) holds) giving channel opening, the $d$-valence electrons of $Ni$ in the $Ni(Cu)$ plane are in the basic state of second ionization energy and the  $d$-valence electrons of $La$ are in the basic state of third ionization energy, and that other states are to be reached from these two states. Further the $d$-valence electrons of $Ni$ in the $Ni(Cu)$ plane and the $d$-valence electrons of $La$ are unified to occupy a sequence of states such that the $d$-valence electrons of $Ni$ in the $Ni(Cu)$ plane occupy the higher states while the $d$-valence electrons of $La$ occupy the lower state. 
Thus, as $MgB_2$, the $d$-valence electrons of $Ni$ are in the state of third ionization energy, the $d$-valence electron of $La$ is in the state of third ionization energy;
and the $s$-valence electrons of $Ni$ are in the state of first ionization energy, 
the $s$-valence electrons of $La$ are in the state of first ionization energy.

Thus,
the maximum value of the energy parameter $|h_{Ni}|$  of the $d$-valence electrons of $Ni$ is proportional to the third ionization energy $\approx 3395$ kJ/mol of $Ni$, the energy parameter $|h_{Ni1}|$ of the $s$-valence electrons of $Ni$ is proportional to the first ionization energy $\approx 737.1$ kJ/mol of $Ni$,
and the energy parameter $|h_{La}|$ of the $s$-valence electrons of $La$ is proportional to the first ionization energy of $La$.
Then, as $La_{1-x}Ca_{x}Cu_{5}$, 
we have the following formula of $T_c$ of $LaNi_{5(1-x)}Cu_{5x}$:
\begin{equation}
k_BT_c=\kappa=\frac{1}{\sqrt{3}}\Delta_{LaNiCu} 
\label{LaNiCu}
\end{equation}
where $\Delta_{LaNiCu}=h=60|h_{Ni}|+15|h_{Ni1}|+|h_{La}|$ is the energy gap of 
 $LaNi_{5(1-x)}Cu_{5x}$ and $60=4\cdot 15$ where $4=\frac82$ is from the $8$ $d$-valence electrons of $Ni$, and the $15$ for the fifteen $Ni(Cu)$ atoms of the cluster of $Ni(Cu)$ atoms in a unit cell
(For simplicity we have simplified  the effect of $d$-valence electrons of $Cu$).
Then from (\ref{LaNiCu}) we can compute the critical temperature $T_c$ of $LaNi_{5(1-x)}Cu_{5x}$ which can be upped to:
\begin{equation}
T_c \approx 423.321 K \,\, \mbox{(Computed $T_c$ of $LaNi_{5(1-x)}Cu_{5x}$)}
\label{LaNiCu5}
\end{equation}

Since $x_0^{\prime}\approx 0$, we have that $LaNi_5$ can be formed in the degenerate state of channel opening and thus $LaNi_5$  can be formed as a superconductor with  $T_c$ upped to $ 423.321 K$.  
This  intermetallic $LaNi_5$ had been used for hydrogen storage since  $LaNi_5$ is easy to be activated and can store a large amount of hydrogen under ambient pressure and in the room temperature \cite{Vuc,Vin}. 
Since the activation and the  hydrogen storage of a material is from the activity of electrons of this material,  this property of  $LaNi_5$ shows that  the $d$-valence electrons of $N_i$ and $La$ in $LaNi_5$ are  in the degenerate state of channel opening.

We remark that the rare earth elements  $La$ of the above $CaCu_5$-type superconductors may be replaced by  rare earth elements  or $Y$, or the mixture of rare earth elements such as $Mm$ (Misch metal). Also the $Cu$ and $Ni$ of  $LaNi_{5(1-x)}Cu_{5x}$ may be replaced by mixture of transition elements and metals such as $Al$ to form  $CaCu_5$-type superconductors.

\subsubsection{More evidences of superconductivity of ${\bf LaNi_{5}}$}
 Let us find more evidences for that $LaNi_{5}$ is a superconductor.  In the theory of  $NiMH$  rechargeable batteries ($M$ denotes a metal with the function of storing hydrogen $H$ and we consider $M=LaNi_{5}$ or $M=RTM_5$ where $R$ is a mixture of rare earth elements with $La$ as a main part such as  $Mm$ and $TM$ denotes a mixture of metals with $Ni$ as the main part), we have the following reaction of charging at the positive electrode
\cite{Vuc,Vin}:
\begin{equation}
Ni(OH)_2+OH^- \to Ni(OOH) + H_2O+e^-
 \label{charge}
\end{equation}
where $Ni(OH)_2$ is as the positive electrode. 
This is as a part of charging. The inverse reaction is as a discharging.

Then we have the following reaction of charging at the negative electrode:
\begin{equation}
 M+H_2O+e^- \to MH + OH^- 
 \label{charge2}
\end{equation}
where $M$ is as the negative electrode. This is as another part of charging. The inverse reaction is as a discharging. Combining the above two reactions we have the following overall reaction of charging:
\begin{equation}
 M + Ni(OH)_2 \to MH +Ni(OOH)
 \label{charge3}
\end{equation}
This reaction is as an exchange of a hydrogen atom $H$. The inverse reaction is as a discharging.
 To express the charging effect, this reaction is usually regarded as an exchange of  a proton $H^+$.

Here to express the charging effect, let us write the reaction (\ref{charge3}) in the following form:
\begin{equation}
 M^+ + Ni(OH)_2 +e^-\to M H +Ni(OOH)
\label{charge4}
\end{equation}
where $M$ is written as $M^+ +e^-$ to express the charging effect of charging $e^-$ to $M^+$ to give $M=M^+ +e^-$ (The inverse reaction gives the discharging effect of discharging $M=M^+ +e^-$ to $M^+$ and $e^-$).

We notice that the freedom of the discharging of $M$ to $M^+$ for the hydrogen storage $MH$ is an evidence 
 for that the $3d$-electrons of $Ni$ of $LaNi_5$ (or  $MmNi_5$) are in the degenerate state of channel opening. Thus this is also an evidence for that $LaNi_5$  (or  $MmNi_5$) is in the state of superconductivity.

Then it is observed that for the $NiMH$  batteries
 a large electric current can be discharged  and that $NiMH$  batteries are with 
a battery capacity which is more than  two times of that of the $ NiCd$ batteries which are of the same structure as  $NiMH$  batteries except that the negative electrode is composed with the element $Cd$ instead of $M$ \cite{Vuc,Vin,Tan2}. This is an evidence 
 for that the $3d$-electrons of $Ni$ of $LaNi_5$ (or  $MmNi_5$) are in the degenerate state of channel opening. Thus this is also an evidence for that $LaNi_5$  (or  $MmNi_5$) is in the state of superconductivity.

 \subsubsection{ $\bf{CaCu_5}$-type  negative electrode of NiMH battery}
Besides $LaNi_5$ and $MmNi_5$,
a general $CaCu_5$-type  negative electrode of the commercial   $NiMH$ batteries is the intermetallic of the form $Ml(Ni_{5-x-y}Co_{x}Z_{y)})_{1+z}$,  $(0\leq x\leq 1, 0\leq y\leq 1,  0\leq x+y\leq 1+a)$ and $a, z$ are small numbers,  where
$Ml$ is a mixture of rare earth elements with $La$ as the main part, $Z$ is a mixture of metals and $Z$ is usually a mixture of $Mn$ and $ Al$ \cite{Lir}.
The intermetallics $LaNi_5$, $MmNi_5$ and $Ml(Ni_{5-x-y}Co_{x}Z_{y)})_{1+z}$ are used as alloys for hydrogen storage, and are extended as the  negative electrode of $NiMH$.

Similar to $LaNi_5$, we expect that the  negative electrode $Ml(Ni_{5-x-y}Co_{x}Z_{y)})_{1+z}$ of the  $NiMH$  is a superconductor with $T_c>300 K$. Let us show, from experiments and the gauge model of superconductivity, that $Ml(Ni_{5-x-y}Co_{x}Z_{y)})_{1+z}$ is a superconductor with $T_c>300 K$.

The maximum of discharge capacity of  $NiMH$  is 
at a temperature about  $30^oC$ \cite{Lir}. 
In the range  $-30 \sim 30^oC$, the increasing of temperature gives the  increasing of  the discharge capacity (and  the increasing of internal resistence of $NiMH$).  Then because of the corrosion of oxidation resistance during discharging, 
in the range  $30 \sim 80^oC$,  the increasing of temperature gives the  decreasing of  the discharge capacity (and  the continuous increasing of internal resistence of $NiMH$).  

In the experiment  on the alloy $Ml(Ni_{5-x-y}Co_{x}Z_{y)})_{1.08}$ with $Co$ as a substituent for $Ni$, in which the content of $Co$ is $0, 1.31\%, 2.63\%, 3.94\%, 5.25\%,$ and $6.56\%$,  it is shown that the discharge capacity of  
$NiMH$  increases at higher temperature ($40  \sim 80^o C$) and decreases at lower  temperature ($-30  \sim 0^o C$) with an increasing $Co$ content \cite{Lir}. 

This experiment can be interpreted as follows. The substution of $Co$ has the effect of anti-corrosion of oxidation resistance during discharge.  When the temperature increases in the range  $40  \sim 80^o C$, more $Co$ contents gives more anti-corrosion of oxidation resistance. This gives lesser internal resistence of the $NiMH$  and thus give more discharge capacity than the case with fewer $Co$ contents. On the other hand when  the temperature decreases in the range  $-30  \sim 0^o C$, since the temperature is low the effect of corrosion of oxidation resistance is small and thus the effect of increasing of discharge capacity with more $Co$ content is negative. This means that  in the range  $-30  \sim 0^o C$ the discharge capacity of $NiMH$  with more  $Co$ substituent  is less than  the discharge capacity of $NiMH$  with fewer  $Co$ substituent.
 Thus in the range  $-30  \sim 0^o C$ (or when the temperature is very low) the internal resistence of $NiMH$  with fewer  $Co$ substituent is less than the  internal resistence of  $NiMH$  with more  $Co$ substituent. This means that the resistance of the alloy $Ml(Ni_{5-x-y}Co_{x}Z_{y)})_{1.08}$ with fewer $Co$ substituent for $Ni$ is less than  the resistance of the alloy $Ml(Ni_{5-x-y}Co_{x}Z_{y)})_{1.08}$ with more $Co$ substituent for $Ni$.

Then we notice that the $Co$ atom is with seven $3d$-valence electrons while the $Ni$ atom is with eight  $3d$-valence electrons. Thus, from the gauge model of superconductivity,  the $Co$ atom can only give three Cooper pairs of $3d$-valence electrons (and with one $3d$-valence electron does not form  Cooper pairs) while  the $Ni$ atom can  give four Cooper pairs of $3d$-valence electrons.  Thus, since the experiment result shows that   in the range  $-30  \sim 0^o C$ the resistance of the alloy $Ml(Ni_{5-x-y}Co_{x}Z_{y)})_{1.08}$ with fewer $Co$ substituent for $Ni$ is less than  the resistance of the alloy $Ml(Ni_{5-x-y}Co_{x}Z_{y)})_{1.08}$ with more $Co$ substituent  for $Ni$, we have that  in the range  $-30  \sim 0^o C$ the  eight  $3d$-valence electrons of $Ni$ must form Cooper pairs to give lesser resistance. Thus  in the range  $-30  \sim 0^o C$ the $CaCu_5$-type alloy $Ml(Ni_{5-x-y}Co_{x}Z_{y)})_{1.08}$ is in the state of superconductivity. 

Then from the graph of the discharge capacity this state of superconductivity can actually be extended to the range $0  \sim 30^0C$ \cite{Lir}.  Thus  the $CaCu_5$-type alloy $Ml(Ni_{5-x-y}Co_{x}Z_{y)})_{1.08}$ is a superconductor with $T_c>300 K$.

 \subsection{ $\bf{CaCu_5}$-type  ${\bf La_{1-x}Ce_{x}Ni_{5-5c}Cu_{5c}}$}
Let us
consider an intermetallic $La_{1-x}Ce_{x}Ni_{5(1-c)}Cu_{5c}$ ($0\leq x, c\leq 1$).
This intermetallic may also be of the $CaCu_5$-type.
 For the doping mechanism let us 
consider the following function:
\begin{equation}
f(x)=1850(1-x)+1949x
 \label{MmNiCu83}
\end{equation}
 where $1850$ kJ/mol and $1949$ kJ/mole are approximately the third ionization energies of $La$ and $Ce$ respectively. 
 This function gives the relation that the increasing of $x$ gives the increasing of $h$. Then we set  the following relation of  the degenerate state of channel opening:
\begin{equation}
f(x_0)=1850(1-x_0)+1949x_0=1753(1-c)+1957.9c
\label{MmNiCu8}
\end{equation}
for some $x_0$ ($0<x_0< 1$) and some $c>0$. When this relation holds (or approximately holds) we have that the channel connecting the state of third ionization energy of $La(Ce)$ and the state of second ionization energy of $Cu(Ni)$ can be opened. 
In this case of channel opening 
the Cooper pairs of $d$ valence electrons of $Cu(Ni)$ 
and the bifurcation region of high-$T_c$ superconductivity can be formed.

When $c<\frac12$, this 
 $La_{1-x}Ce_{x}Ni_{5(1-c)}Cu_{5c}$
is similar to $LaNi_{5(1-x)}Cu_{5x}$ with $Ni$ as the main part and we omit the details. When $c>\frac12$, then $Cu$ is as the main part. Thus, as $La_{1-x}Ca_{x}Cu_{5}$, when $c$ is close to $1$
we have the following  
formula of approximate value of $T_c$ of
 $La_{1-x}Ce_{x}Ni_{5(1-c)}Cu_{5c}$:
\begin{equation}
k_BT_c=\kappa=\frac{1}{\sqrt{3}}\Delta_{LaCeNiCu} 
\label{LaCeNiCu}
\end{equation}
where $\Delta_{LaCeNiCu}=h=75|h_{Cu}|+|h_{La}|$ is the energy gap of 
$La_{1-x}Ce_{x}Ni_{5(1-c)}Cu_{5c}$.
(For simplicity we have simplified the effect of $s,d, f$-valence electrons of $Ni$ and $Ce$).
Then from (\ref{LaCeNiCu}) we can compute the critical temperature $T_c$ of $La_{1-x}Ce_{x}Ni_{5(1-c)}Cu_{5c}$ which can be upped to:
\begin{equation}
T_c \approx 525.3 K \,\, \mbox{(Computed $T_c$ of $La_{1-x}Ce_{x}Ni_{5(1-c)}Cu_{5c}$)}
\label{LaCeNiCu5}
\end{equation}

We remark that we may replace $Ce$ with a mixture of rare earth elements. When $x_0\geq \frac25$ this family of $CaCu_5$-type superconductors includes the intermetallic $MmNi_{5(1-c)}Cu_{5c}$ where $Mm$ denotes a Mischmetal which is a mixture of rare earth elements with the fractional part of $Ce$ more than $\frac25$. As $LaNi_5$ this intermetallic $MmNi_{5(1-c)}Cu_{5c}$ is good for hydrogen storage and rechargeble batteries \cite{Vin}.

\subsection{Self-Discharge of ${\bf AB_5}$-NiMH as Experimental Evidence of Superconductivity of ${\bf LaNi_5}$}\label{sec001}

Let us give more experimental evidences for the superconductivity of $LaNi_5$.
It is well known that the $AB_5$-type  $NiMH$ batteries (with $LaNi_5$ and its variations as the negative electrodes) are with very high self-discharge rate comparing to other rechargeable batteries such as NiCd batteries and the $Li$ batteries, and that the self-discharge rate of $AB_5$-$NiMH$ is the highest among the non-$NiMH$ batteries such as NiCd batteries and  $Li$ batteries.

Then a new kind of $NiMH$ batteries called as the $ENELOOP$ is developed \cite{Ter}. The negative electrode of this $ENELOOP$ is an intermetallic of the $Ce_2Ni_7$-type. This crystal structure $Ce_2Ni_7$ of $ENELOOP$ is a mixed-type of $AB_5$ and $AB_2$ with $AB_5$ of $\frac23$ part and  $AB_2$ of $\frac13$ part  
where $A$ is a mixture of $R$ (where $R$ is a mixture of $La$ with rare earth elements such as $Ce$) and $Mg$, and $B$ is $Ni$ (which is mixed with a little amount of $Al$)\cite{Ter}.
In this section we show that the property of very high self-discharge rate of $AB_5$-$NiMH$ 
and the very low self-discharge rate of $ENELOOP-NiMH$ is an experimental evidence for that $LaNi_5$ and its variations are in the state of superconductivity with critical temperature $T_c> 300 K$.

Let the number of $Ni$ elements be proportional to  1 in three unit cells of $AB_5$ (i.e. $LaNi_5$ and its variations with $CaCu_5$-type) corresponding to a unit cell of $Ce_2Ni_7$. Then the
number of $Ni$ elements is proportional to $\frac{16}{18}$ in a unit cell of $Ce_2Ni_7$\cite{Ter}. 
Then we let the doping be $x=0,\frac{1}{18}, \frac{2}{18}$. Then we have that the doping gives $Ce_2Ni_7$-$ENELOOP$ transited to $CaCu_5$-$AB_5$-$NiMH$ by the doping of the factor
$Ni_xMg_{a(\frac{1}{9}-x)}$ where $Ni$ may be mixed with a little amount of $Mn, Co$, $0<a<1$ is a constant which gives the mixture 
$A=R_{1-a(\frac{1}{9}-x)}Mg_{a(\frac{1}{9}-x)}$. When $x=0$, this gives the 
 $Ce_2Ni_7$-type crystal structure of $ENELOOP$. When $x=\frac19$, this gives the  $CaCu_5$-type crystal structure of $AB_5$-$NiMH$.
We have that the doping parameter $x$ is small. In this case, the doping should give a little change of property if there is no phase transition 
appearing in this doping.

We remark that this doping is similar to the doping of the cuprate $ La_2CuO_4 $ transited to the  form $ La_{2-x}Sr_xCuO_4$. 
Further, in this phase transition of $ La_2CuO_4 $ there is a sudden change of the crystal structure from the orthorhombic phase to the tetragonal phase during the phase transition. Thus this is another form of the superconductivity phase transition of $ La_2CuO_4 $ to $ La_{2-x_0}Sr_{x_0}CuO_4$ for some $x_0>0$.

Now let us consider the doping of $Ce_2Ni_7$-$ENELOOP$ to $CaCu_5$-$AB_5$-$NiMH$.  
We want to show that there is a sudden change as that in the Fig.5 for the doping from  $La_2CuO_4 $ to $La_{2-x_0}Sr_{x_0}CuO_4$. Then, as the doping in Fig.5, this sudden change is a phase transition.
Then we show that this phase transition is the appearance of the phase of superconductivity, as that in the Fig.5 for $La_2CuO_4 $.

For this doping of $Ce_2Ni_7$-$ENELOOP$ to $CaCu_5$-$AB_5$-$NiMH$,
we have that the small doping of $x=0$ to $x=\frac19$ gives the change of crystal structure changing from 
$Ce_2Ni_7$-type to $CaCu_5$-type. Thus this is a sudden change of the crystal structure. Thus this is a phase transition.

On the other hand, let us consider another phenomenon of this phase transition.
It is known that $ENELOOP$ is a $NiMH$ battery with a very low self-discharge (about 2 to 3 $\%$ per month  from the beginning of self-discharge) comparing to that of the $AB_5$-$NiMH$ (about 30 $\%$ per month  from the beginning of self-discharge, and about 20 $\%$ on the first day. c.f. Nickel-Metal Hydride, Wikimedia).

Thus we have that the small doping of $x=0$ to $x=\frac19$ gives that the very low self-discharge rate of $ENELOOP$ is changed to a very high self-discharge rate of $AB_5$-$NiMH$. 

Then let us consider the relation of self-discharge with the resistance.
We have that the resistance of the negative electrode is an important factor of the self discharge of batteries, when other factors of discharge are not changed, we have that if the  self-discharge is higher, then the resistance of  the negative electrode is lower.

Then, since other factors of self-discharge of $ENELOOP$ and $AB_5$-$NiMH$ are basically of the same in the beginning period of self-discharge, we have that this small doping of $x=\frac19$ gives the sudden change of resistance of the negative electrode of $ENELOOP$ to a state with very low resistance of the negative electrode of $AB_5$-$NiMH$ (or the $LaNi_5$ and its variations). Thus this is a phase transition.  

Thus we have two phenomena of this phase transition in this small doping of $x=\frac19$. One is the change of the crystal structure and the other is the change of a state of resistance of a metal to a state of metal with very low resistance. 

Thus, since this phase transition is about the resistance, we have that this phase transition is from the phase of conductivity to the phase of superconductivity. 



Thus, as a phase transition, the change of 
self-discharge rates and the change of crystal structures of 
$Ce_2Ni_7$-$ENELOOP$ to $CaCu_5$-$AB_5$-$NiMH$, 
is an experimental evidence that the material $LaNi_5$ and its variations are in the state of superconductivity with critical temperature $T_c> 300 K$.

\subsection{Discharge Capacity of ${\bf AB_5}$-NiMH as Experimental Evidence of Superconductivity of ${\bf LaNi_5}$}

In this subsection we show that the property of  high discharge capacity of $NiMH$ is  an experimental evidence  for that $LaNi_5$ and its variations are in
the state of superconductivity with critical temperature $T_c> 300 K$.

It is well known that the $AB_5$-$NiMH$ batteries (with $LaNi_5$ and its variations as the negative electrodes) are with high discharge capacity comparing to other rechargeable batteries such as $AB_2$-$NiMH$ batteries. 
Let us first
 use this fact to show that this is an experimental evidence that $LaNi_5$ and its variations are in
the state of superconductivity with critical temperature $T_c> 300 K$.

It is known that the ability of hydrogen storage of $AB_2$ material (such as  $VNi_2, TiNi_2, ZrNi_2, VCr_2$) is higher than that of $LaNi_5$ and its variations (or called as $AB_5$ material). This ability of hydrogen storage is for the discharge ability of $NiMH$. Then the $AB_5$ material is with higher discharge capacity than the $AB_2$ material means that  the $AB_5$ materials are with better conductivity such that they are with 
higher discharge capacity than that of $AB_2$ materials even though they are with lower ability of hydrogen storage. Then since the $AB_2$ materials are a kind of conductors we have that the $AB_5$ materiasl are in a higher level of conductivity. This is then an experimental evidence that $AB_5$ materials are 
in
the state of superconductivity with critical temperature $T_c> 300 K$.

In the following of this subsection  we then show that this property of high discharge capacity of $AB_5$-$NiMH$ (with $LaNi_5$ and its variations as the negative electrodes) and the lower discharge capacity $NiMH$ batteries ($ENELOOP$) is an experimental evidence for that $LaNi_5$ and its variations are in the state of superconductivity with critical temperature $T_c> 300 K$.

It is known that $ENELOOP-NiMH$ battery is with the discharge capacity 2000 mAh comparing to the discharge capacity 2500 mAh of the $AB_5$-$NiMH$ battery\cite{Ter}. 
In the following of this section let us show that this is an experimental evidence that $LaNi_5$ and its variations are in
the state of superconductivity with critical temperature $T_c> 300 K$.

Now let us consider the doping of $Ce_2Ni_7$-$ENELOOP$ to $AB_5$-$NiMH$. We want to show that there is a sudden change as that in the Fig.5 for the doping from  $La_2CuO_4 $ to $La_{(2-x_0)}Sr_{x_0}CuO_4$. Then, as the doping in Fig.5, this sudden change is a phase transition. 
Then we show that this phase transition is the appearance of the phase of superconductivity, as that in the Fig.5 for $La_2CuO_4 $.

Then, we have that the small doping of $x=0$ to $x=\frac19$ gives the change of crystal structure changing from 
$Ce_2Ni_7$-type to $CaCu_5$-type. Thus this is a sudden change of the crystal structure. Thus this is a phase transition.

On the other hand, as stated in the above, this $ENELOOP-NiMH$ battery is with the discharge capacity 2000 mAh comparing to the discharge capacity 2500 mAh of the $AB_5$-$NiMH$ battery \cite{Ter}. 

Thus we have that the small doping of $x=0$ to $x=\frac19$ gives that the discharge capacity 2000 mAh of $ENELOOP$ is changed to a  higher discharge capacity 2500 mAh of 
$AB_5$-$NiMH$. 

Then let us consider the relation of discharge with the resistance.
We have that the resistance of the negative electrode is an important factor of the discharge capacity of batteries, when other factors of discharge capacity are not changed, we have that if the discharge capacity is higher, then the resistance of  the negative electrode is lower.

Then, since other factors of discharge capacity of $ENELOOP$ and $AB_5$-$NiMH$ are basically of the same in the beginning period of discharge, we have that this small doping of $x=\frac19$ gives the change of resistance of $ENELOOP$ to a state with very low resistance of $AB_5$-$NiMH$.

Let us show that this change of resistance is a sudden change. 
Then this is a phenomenon of phase transition.

We have that for the $ENELOOP$, each factor of the number of $N_i$ gives the discharge capacity
$\frac{2000}{16}=125$ mAh.  Thus the doping from $\frac{16}{18}$ to 1 should give the increasing of the discharge capacity equal to $125\times 2=250$ mAh. Then, from the discharge rate 2500 mAh of $AB_5$-$NiMH$,  we have that the actual increasing of the discharge capacity is $2500-2000=500$ mAh. This is two times of the increasing of the discharge capacity 250 mAh of $ENELOOP$. Thus this increasing of discharge capacity  is a sudden increasing of the discharge capacity. Thus this increasing of discharge capacity is a phenomenon of phase transition.

Thus we have two phenomena of this phase transition in this small doping of $x=\frac19$. One is the change of the crystal structure and the other is the change of a state of resistance of a metal to a state of metal with much lower resistance. 

Then, since this phase transition is about the resistance, we have that this phase transition is from the phase of conductivity to the phase of superconductivity.


Thus, as a phase transition, the sudden change of 
discharge capacity and the sudden change of crystal structures of $Ce_2Ni_7$-$ENELOOP$ to $CaCu_5$-$AB_5$-$NiMH$ is an experimental evidence that the material $LaNi5$ and its variation are in the state of superconductivity with critical temperature $T_c> 300 K$.

\subsection{Activation of ${\bf AB_5}$-NiMH as Experimental Evidence of Superconductivity of ${\bf LaNi_5}$}

In this subsection we show that the property of  easy activation of $NiMH$ is  an experimental evidence  for that $LaNi_5$ and its variations are in
the state of superconductivity with critical temperature $T_c> 300 K$.

It is shown experimentally that $AB_5$-$NiMH$ is of easy activation in the sense that 
about 2 to 3 decrepitation cycles can give a full activation such that the experimental discharge capacity approximates the theoretical discharge capacity. 

Also it is  shown experimentally that $AB_2$-$NiMH$ can be activated but is not as easy as the $AB_5$-type and is of about 40 decrepitation cycles for activation\cite{Er}. 
This experimentally  shows 
that the $AB_2$-type material are in 
the degenerate state of channel opening.
Thus the negative electrode of $AB_2$-$NiMH$ is in the state of conductivityty and
$AB_2$-$NiMH$ is in the state of activation but is not in the state of easy activation. Let us call this kind of activation with more than 4 decrepitation cycles for a full activation as an ordinary activation.

 Since conductivity is a main factor for activation, in this section we shall show that the fact that $AB_5$-$NiMH$ is easy to be activated is an experimental evidence that $LaNi_5$ and its variations are in the state of superconductivity.

Now let us consider the doping of $Ce_2Ni_7$-$ENELOOP$ to $AB_5$-$NiMH$. We want to show that there is a sudden change as that in the Fig.5 for the doping from  $La_2CuO_4 $ to $La_{2-x_0}Sr_{x_0}CuO_4$. Then, as the doping in Fig.5, this sudden change is a phase transition.  
Then we show that this phase transition is the appearance of the phase of superconductivity, as that in the Fig.5 for $La_2CuO_4 $.

 Then, we have that the small doping of $x=0$ to $x=\frac19$ gives the change of crystal structure changing  from 
$Ce_2Ni_7$-type to $CaCu_5$-type. Thus this is a sudden change of the crystal structure. Thus this is a phase transition.

On the other hand, 
this $ENELOOP-NiMH$ battery is with the discharge capacity 2000 mAh comparing to the discharge capacity 2500 mAh of the $AB_5$-$NiMH$ battery\cite{Ter}. 

Let us consider the activation of $NiMH$ batteries. Let us consider $Ce_2Ni_7$-$ENELOOP$. 
It is known that  $Ce_2Ni_7$-$ENELOOP$ is with higher hydrogen storage than the $AB_5$-$NiMH$ of $CaCu_5$-type. 
Thus, if $Ce_2Ni_7$-$ENELOOP$ and  $AB_5$-$NiMH$ are of the same level of easy activation, then $Ce_2Ni_7$-$ENELOOP$ should be with higher discharge capacity than  $AB_5$-$NiMH$ since the level  of easy activation also means that the experimental discharge capacity approximates the theoretical discharge capacity. 
However it is found that $Ce_2Ni_7$-$ENELOOP$ is with discharge capacity $2000$ mAh while $AB_5$-$NiMH$ is with discharge capacity $2500$ mAh in the beginning of discharge. 
Then it is known that $AB_5$-$NiMH$ is of easy activation in the sense that 
about 2 to 3 decrepitation cycles can give a full activation such that the experimental discharge capacity approximates the theoretical discharge capacity. 
Thus the activation of $Ce_2Ni_7$-$ENELOOP$
is not of the same level as the activation of $AB_5$-$NiMH$ of the $CaCu_5$-type.
Thus  $Ce_2Ni_7$-$ENELOOP$ is in the state of ordinary activation but is not in the state of easy activation.

Thus this transition from $Ce_2Ni_7$-$ENELOOP$ to $CaCu_5$-$AB_5$-$NiMH$
is a transition from the state of ordinary activation 
to the state of easy activation.

Thus we have two phenomena of this phase transition in this small doping of $x=\frac19$. One is the change of the crystal structure and the other is the change of a state of activation 
to the state of 
easy activation. 
Then, since 
easier activation means better conductivity, and that easy activation 
is the best kind of conductivity,
we have that this phase transition is from the phase of conductivity to the phase of superconductivity.


Thus, as a phase transition, the easy activation of $AB_5$-$NiMH$ and the change of crystal structures of $ENELOOP$ and $AB_5$-$NiMH$ is an experimental evidence that the material $LaNi_5$ and its variation are in the state of superconductivity with critical temperature 
$T_c>300 K$.

\subsection{High-Rate Discharge of ${\bf AB_5}$-NiMH as Experimental Evidence of Superconductivity of ${\bf LaNi_5}$}

In this subsection we show that the property of  high-rate discharge of $NiMH$ is also an experimental evidence  for that $LaNi_5$ and its variations are in
the state of superconductivity with critical temperature $T_c> 300 K$.

It is shown from the generations of products of $ENELOOP$ 
that $Ce_2Ni_7$-$ENELOOP$ battery with the discharge capacity 2000 mAh can be in  the state of high-rate discharge of  up to 1C where C means the capacity and it is not good for $Ce_2Ni_7$-$ENELOOP$ in  the state of high-rate discharge of 2C.

Then it is known that $AB_5$-$NiMH$ batteries can discharge a large current in a short time.
It is shown that $AB_5$-$NiMH$ with discharge capacity 2500 mAh can be in  the state of high-rate discharge of 
up to 10C (c.f. GP Betteries, Nickel Metal Technical Handbook, and SUPP Betteries, Nickel Metal Technical Handbook).

Thus this transition from $Ce_2Ni_7$-$ENELOOP$ to $CaCu_5$-$AB_5$-$NiMH$
is a transition from the state of  high-rate discharge of roughly up to 2C
to the state of  high-rate discharge of up to 10C.

Thus we have two phenomena of this phase transition in this small doping of $x=\frac19$. One is the change of the crystal structure and the other is the change of a state of  high-rate discharge of roughly up to 2C
to the state of  high-rate discharge of up to 10C.

Then since high-rate discharge is of a short time discharge (such as the discharge of digital camera) that the electric current during discharge is stored in the 
negative electrode. Thus the appearing of high-rate discharge of up to 10C in the negative electrode means that the resistance of the $LaNi_5$ and its variations (as the negative electrode) is very low. Thus this phase transition  from the phase of high-rate discharge of up to 2C to the phase of high-rate discharge of up to 10C is a phase transition from the phase of conductivity to the phase of  superconductivity.


Thus, as a phase transition, the transition from the phase of high-rate discharge of up to 2C of $Ce_2Ni_7$-$ENELOOP$ to the phase of the high-rate discharge of up to 10C of $AB_5$-$NiMH$ and the transition of crystal structures of $Ce_2Ni_7$ to $CaCu_5$ 
is an experimental evidence that the material $LaNi_5$ and its variations are in the state of superconductivity with critical temperature $T_c> 300 K$.

\subsection{ $\bf{CaCu_5}$-type superconductor ${\bf YNi_{5-5x}Cu_{5x}}$}
Let us
consider an intermetallic $YNi_{5(1-x)}Cu_{5x}$ ($0\leq x\leq 1$). 
We have that the intermetallics $YCu_{5}$ and $YNi_{5}$ are of the $CaCu_5$-type \cite{Chak, Sub}.
Thus we may suppose that
 the crystal structure of $YNi_{5(1-x)}Cu_{5x}$ can be formed in the $CaCu_5$-type with $Ca$  corresponding to $Y$ and $Cu$  corresponding to $Ni_{(1-x)}Cu_{x}$. 
 
For the doping mechanism giving superconductivity let us consider the following function:
\begin{equation}
f(x)=1980-3393(1-x)
 \label{YNiCu83}
\end{equation}
 where $1980$ kJ/mol and $3393$ kJ/mole are approximately the third ionization energies of $Y$ and $Ni$ respectively. 
 This function gives the relation that the increasing of $x$ gives the increasing of $h$. Then we set  the following relation of  the degenerate state of channel opening:
\begin{equation}
f(x_0)=1980-3393(1-x_0)=1957.9x_0
\label{YNiCu8}
\end{equation}
for some $x_0$ ($0<x_0< 1$) where $1957.9$ kJ/mol is approximately the second
ionization energy of $Cu$.  When this relation holds we have that the channel connecting the state of third ionization energy of $Y$ and the state of second ionization energy of $Cu$ (for the high-$T_c$ superconductivity given by the $Cu$ plane) can be opened. 
In this case of channel opening 
the Cooper pairs of $d$-valence electrons of $Cu$, the Cooper pairs of $s$-valence electrons of $Y$,
and the bifurcation region of high-$T_c$ superconductivity can be formed.

We notice that $x_0\approx 0.9846$.
Thus $YNi_{5(1-x)}Cu_{5x}$ comes into the range  of high-$T_c$ superconductivity when $ x_2<x<x_0$ for some $x_2$ such that $0\leq x_2<x_0< 1$.  

When (\ref{YNiCu8}) holds giving channel opening, the $d$-valence electrons of $Cu$ in the $Cu$ plane are in the basic state of second ionization energy and the  $d$-valence electrons of $Y$ are in the basic state of third ionization energy, and that other states are to be reached from these two states. Further the $d$ valence electrons of $Cu$ in the $Cu$ plane and the $d$ valence electron of $Y$ are unified to occupy a sequence of states such that the $d$ valence electron of $Cu$ in the $Cu$ plane occupy the higher states while the $d$ valence electron of $Y$ occupy the lower state. 

Thus, similar to $MgB_2$, the $d$-valence electrons of $Cu$ are in the state of third ionization energy, the $d$-valence electron of $Y$ is in the state of third ionization energy;
and the $s$-valence electrons of $Cu$ are in the state of first ionization energy, 
the $s$-valence electrons of $Y$ is in the state of first ionization energy.

Thus,
the maximum value of the energy parameter $|h_{Cu}|$  of the $d$ valence electrons of $Cu$ is proportional to the third ionization energy of $Cu$, 
and the energy parameter $|h_{Y}|$ of the $s$-valence electrons of $Y$ is proportional to the first ionization energy of $Y$.

Then, similar to $La_{1-x}Ca_{x}Cu_{5}$, 
we have the following formula of $T_c$ of $YNi_{5(1-x)}Cu_{5x}$:
\begin{equation}
k_BT_c=\kappa=\frac{1}{\sqrt{3}}\Delta_{YNiCu} 
\label{YNiCu}
\end{equation}
where $\Delta_{YNiCu}=h=75|h_{Cu}|+|h_{Y}|$ is the energy gap of the  superconductivity of $YNi_{5(1-x)}Cu_{5x}$
(We have simplified the effect of $s$, $d$ valence electrons of $Ni$).

Then from (\ref{YNiCu}) we can compute the highest critical temperature $T_c$ of $YNi_{5(1-x)}Cu_{5x}$ which can be upped to:
\begin{equation}
T_c \approx 525.43 K \,\, \mbox{(Computed $T_c$ of $YNi_{5(1-x)}Cu_{5x}$)}
\label{YNiCu5}
\end{equation}

We remark that we may use other transition elements such as $Co$ 
to replace the element $Ni$. In this case we have the intermetallic $YCo_{5(1-x)}Cu_{5x}$ with $x_0\approx 0.9$. The critical temperature can be upped to $T_c \approx 525.43 K$.  

Let us then consider the case of $Ni$ as the main part (instead of the above case of $Cu$ as the main part) for $YNi_{5(1-x)}Cu_{5x}$. To this end let us wirte this intermetallic in the form $YNi_{5y}Cu_{5(1-y)}$ ($0\leq y\leq 1$). Then we consider the following function for doping mechanism:
\begin{equation}
g(y)=1980-3555(1-y)
 \label{YNiCu83b}
\end{equation}
 where $1980$ kJ/mol and $3555$ kJ/mole are approximately the third ionization energies of $Y$ and $Cu$ respectively. 
 This function gives the relation that the increasing of $y$ gives the increasing of $h$. Then we set the following relation of  the degenerate state of channel opening:
\begin{equation}
g(y_0)=1980-3555(1-y_0)=1753y_0
\label{YNiCu8b}
\end{equation}
for some $y_0$ ($0<y_0< 1$) where $1753$ kJ/mol is approximately the second
ionizationy energy of $Ni$.  When this relation holds we have that the channel connecting the state of third ionization energy of $Y$ and the state of second ionization energy of $Ni$ (for the high-$T_c$ superconductivity given by the $Ni(Cu)$ plane) can be opened. 

Then, the rest is similar to 
$ LaNi_{5(1-x)}Cu_{5x}$ ($0\leq x\leq 1$) with $La$ replaced by $Y$ and $x$ replaced by $y$. Thus we have the following formula of $T_c$ of $YNi_{5y}Cu_{5(1-y)}$:
\begin{equation}
k_BT_c=\kappa=\frac{1}{\sqrt{3}}\Delta_{YNiCu} 
\label{YNiCub}
\end{equation}
where $\Delta_{YNiCu}=h=60|h_{Ni}|+15|h_{Ni1}|+|h_{Y}|$ is the energy gap of superconductivity of $YNi_{5y}Cu_{5(1-y)}$ (when $Ni$ instead of $Cu$ is the main part).

Then from (\ref{YNiCub}) we can compute the critical temperature $T_c$ of $YNi_{5y}Cu_{5(1-y)}$ which can be upped to:
\begin{equation}
T_c \approx 423.21 K \,\, \mbox{(Computed $T_c$ of $YNi_{5y}Cu_{5(1-y)}$)}
\label{YNiCu5b}
\end{equation}

\subsection{ $\bf{CaCu_5}$-type intermetallic ${\bf YNi_{5(1-x)}Co_{5x}}$}
Let us consider an intermetallic $YNi_{5(1-x)}Co_{5x}$ ($0\leq x\leq 1$). 
We have that the intermetallics $YCo_{5}$ and $YNi_{5}$ are of the $CaCu_5$-type 
Thus we may suppose that
 the crystal structure of $YNi_{5(1-x)}Co_{5x}$ can be formed in the $CaCu_5$-type with $Ca$  corresponding to $Y$ and $Cu$  corresponding to $Ni_{5(1-x)}Co_{5x}$. 
 Then we consider the following function for doping mechanism:
\begin{equation}
f(x)=1980-3232(1-x)
 \label{YNiCo83}
\end{equation}
 where $1980$ kJ/mol and $3232$ kJ/mole are approximately the third ionization energies of $Y$ and $Co$ respectively. 
 This function gives the relation that the increasing of $x$ gives the increasing of $h$. Then we set  the following relation of  the degenerate state of channel opening:
\begin{equation}
f(x_0)=1980-3232(1-x_0)=1753x_0
\label{YNiCo8}
\end{equation}
for some $x_0$ ($0<x_0< 1$) where $1753$ kJ/mol is approximately the second
ionization energy of $Ni$.  When this relation holds we have that the channel connecting the state of third ionization energy of $Y$ and the state of second ionization energy of $Ni$ (for the high-$T_c$ superconductivity given by the $Ni(Co)$ plane) can be opened. 
In this case of channel opening 
the Cooper pairs of $s$, $d$ valence electrons of $Ni$, the Cooper pairs of $s$ valence electrons of $Y$,
and the bifurcation region of high-$T_c$ superconductivity can be formed.

We have that $x_0\approx 0.846$.
Thus $YNi_{5x}Co_{5(1-x)}$ comes into the range  of high-$T_c$ superconductivity when $ x_2<x<x_0$ for some $x_2$ such that $0\leq x_2<x_0< 1$.

When (\ref{YNiCo8}) holds giving channel opening we have that the $d$ valence electrons of $Ni$ in the $Ni(Co)$ plane are in the basic state of second ionization energy of $Ni$ and the  $d$ valence electrons of $Y$ are in the basic state of third ionization energy of $Y$, and that other states are to be reached from these two states. Further the $d$ valence electrons of $Ni$ in the $Ni(Co)$ plane and the $d$ valence electron of $Y$ are unified to occupy a sequence of states such that the $d$ valence electrons of $Ni$ in the $Ni(Co)$ plane occupy the higher states while the $d$ valence electron of $Y$ occupy the lower state. 
Thus, as $MgB_2$, the $d$ valence electrons of $Ni$ are in the basic state of third ionization energy of $Ni$, the $d$ valence electron of $Y$ is in the state of third ionization energy of $Y$;
and the $s$ valence electrons of $Ni$ are in the basic state of first ionization energy of $Ni$, 
the $s$ valence electrons of $Y$ are in the basic state of first ionization energy.

Thus,
the maximum value of the energy parameter $|h_{Ni}|$  of the $d$ valence electrons of $Ni$ is proportional to the third ionization energy of $Ni$,
 the maximum value of the energy parameter $|h_{Ni1}|$  of the $s$ valence electrons of $Ni$ is proportional to the first ionization energy of $Ni$; 
and the energy parameter $|h_{Y}|$ of the $s$ valence electrons of $Y$ is proportional to the first ionization energy of $Y$.

Then, similar to $YNi_{5(1-x)}Cu_{5x}$, 
we have the following formula of $T_c$ of $YNi_{5x}Co_{5(1-x)}$:
\begin{equation}
k_BT_c=\kappa=\frac{1}{\sqrt{3}}\Delta_{YNiCo} 
\label{YNiCo}
\end{equation}
where $\Delta_{YNiCo}=h=60|h_{Ni}|+15|h_{Ni1}|+|h_{Y}|$ is the energy gap of the  superconductivity of $YNi_{5x}Co_{5(1-x)}$
where the coefficients $15$, $60=4\cdot15$ are from the $4=\frac82$ for the eight $d$ valence electrons of the $Ni$, and the $15$ for the fifteen $Ni(Co)$ atoms of the cluster of $Ni(Co)$ atoms in a unit cell.
(For simplicity we have simplified the effect of $s$, $d$ valence electrons of $Co$ for $\Delta_{YNiCo}$).
We have $|h_{Ni}| \approx \xi 3393$ kJ/mol,  $|h_{Ni1}|\approx \xi 736.7$ kJ/mol and
 $|h_{Y}|\approx \xi 600$ kJ/mol where $736.7$ kJ/mol and $3393$ kJ/mol are approximately the first and third ionization energies of $Ni$.
Then from (\ref{YNiCo}) we can compute the highest critical temperature $T_c$ of $YNi_{5x}Co_{5(1-x)}$ which can be upped to:
\begin{equation}
T_c \approx 423.21 K \,\, \mbox{(Computed $T_c$ of $YNi_{5x}Co_{5(1-x)}$)}
\label{YNiCo5}
\end{equation}

We remark that the rare earth elements or $Y$ of all the above $CaCu_5$-type superconductors can be replaced by the mixture of rare earth elements such as $Mm$ (Misch metal) with the outmost shell of the form $5d^16s^2$.

\subsection{$\bf{CaCu_5}$-type superconductor ${\bf Sr_{1-x}Ca_{x}Cu_5}$}
Let us consider an intermetallic $Sr_{1-x}Ca_{x}Cu_5$ ($0\leq x\leq 1$). This intermetallic is a special case of the  $CaCu_5$-type intermetallic $R_{1-x}A_{x}Cu_5$ with $x=1$ ($A$ denotes a mixture of $Ca$ and $Sr$).
We have that $CaCu_5$ and $SrCu_5$ can be formed in the $CaCu_5$-type phase. 
Thus we may suppose that the intermetallic $Sr_{1-x}Ca_{x}Cu_5$  can also be formed in the $CaCu_5$-type phase.
For the doping mechanism of superconductivity let us consider the following function:
\begin{equation}
f(x)=5500(1-x)+6491x 
 \label{CaSrCu83}
\end{equation}
 where $6491$ kJ/mol and $5500$ kJ/mole are approximately the fourth ionization energies of $Ca$ and $Sr$ respectively. 
 This function gives the relation that the increasing of $x$ gives the increasing of $h$. Then we set  the following relation of  the degenerate state of channel opening:
\begin{equation}
f(x_0)= 5500(1-x_0)+6491x_0  =5536
\label{CaSrCu8}
\end{equation}
for some $x_0$ ($0<x_0< 1$) where $5536$ kJ/mol is approximately the fourth
ionization energy of $Cu$.  When this relation holds we have that the channel connecting the two states of fourth ionization energy of $Sr$ and  $Cu$ (for the high-$T_c$ superconductivity given by the $Cu$ plane) can be opened. 
This channel opening gives a freedom of electric current with
a direction orthogonal to the $Cu$ plane. 
From this freedom of electric current, the Cooper pairs of the $3s,3p$-level electrons of $Cu$ and the $4s,4p$-level electrons of $Sr$ can be formed.
Thus this channel opening gives the 3D region of conventional superconductivity. From this 3D conventional superconductivity we have the existence of quasi-2D bifurcation region of unconventional high-$T_c$ superconductivity given by the $Cu$ plane.

We notice that $x_0\approx 0.0363$.
Thus $Sr_{1-x}Ca_{x}Cu_5$ comes into the range  of high-$T_c$ superconductivity when $x_0\leq x\leq x_1$ for some $x_1$ such that $x_0< x_1< 1$. 
 
When (\ref{CaSrCu8}) holds giving channel opening we have that the $3s,3p$-level electrons of $Cu$ in the $Cu$ plane are in the basic state of fourth ionization energy and the $4s,4p$-level electrons of $Sr$ are in the basic state of fourth ionization energy, and that other states are to be reached from these two states. Further the $3s,3p$-level electrons of $Cu$ and the $4s,4p$-level electrons of $Sr$ are unified to occupy a sequence of states such that the $3s,3p$-level electrons of $Cu$ in the $Cu$ plane occupy the higher states while the $4s,4p$-level electrons of $Sr$ occupy the lower states. 

Thus, the $3s,3p$-level electrons of $Cu$ are in the state of fifth ionization energy of $Cu$;
and the $4s,4p$-level electrons of $Sr$ are in the basic state of fourth ionization energy. 
The $4s,4p$-level electrons of $Sr$ in the basic state of fourth ionization energy are in the opened channel of 3D superconductivity,
while the $3s,3p$-level electrons of $Cu$ in the state of fifth ionization energy are for the quasi-2D high-$T_c$ superconductivity of the $Cu$ plane. 

Then when (\ref{CaSrCu8}) holds giving channel opening the Cooper pairs of the 
$s,p$-valence electrons of $Sr$ can also be formed. These $s$-valence electrons of $Sr$ are in the basic state of first ionization energy (and can be in the states of third and second ionization energies of valence electrons respectively).

Thus, the energy parameter $|h_{Cu5}|$  of the $3s,3p$-level electrons of $Cu$ is proportional to the fifth ionization energy of $Cu$; 
 the  energy parameter $|h_{Sr4}|$ of the $4s,4p$-level electrons of $Sr$ is proportional to the fourth ionization energy of $Sr$; and the energy parameter $|h_{Sr}|$ of the $s$-valence electrons of $Sr$ is proportional to the first ionization energies $Sr$ respectively (when the $s$-valence electrons of $Sr$ is in the basic state of first ionization energy).

Then, 
we have the following formula of $T_c$ of $Sr_{1-x}Ca_{x}Cu_5$:
\begin{equation}
k_BT_c=\kappa=\frac{1}{\sqrt{3}}\Delta_{CaSrCu} 
\label{CaSrCu}
\end{equation}
where $\Delta_{CaSrCu}=h=60|h_{Cu5}|+4|h_{Sr4}|+|h_{Sr}|$ is the energy gap of superconductivity of $Sr_{1-x}Ca_{x}Cu_5$
where the coefficients $4$, $60=4\cdot15$ are from the $4=\frac82$ for the eight $3s,3p$-level electrons of the $Cu$ or the eight $4s,4p$-level electrons of the $Sr$, and the $15$ for the fifteen $Cu$ atoms of the $CaCu_5$-type structure, as the above case $A_xR_{1-x}Cu_5$ with $x<1$.

We have  $|h_{Cu5}|\approx \xi 7700$ kJ/mol, $|h_{Sr}|\approx \xi 549.5$ kJ/mol and
 $|h_{S4r}|\approx \xi 5500$ kJ/mol where $7700$ kJ/mol is approximately the fifth ionization energy of $Cu$, and  $549.5$ kJ/mol is approximately the first ionization energy of $Sr$.
Then from (\ref{CaSrCu}) we can compute the highest critical temperature $T_c$ of $Sr_{1-x}Ca_{x}Cu_5$ which can be upped to:
\begin{equation}
T_c \approx 952.7 K \,\, \mbox{(Computed $T_c$ of $Sr_{1-x}Ca_{x}Cu_5$)}
\label{CaSrCu5}
\end{equation}

Similar to the  intermetallic $R_{1-x}A_{x}Cu_5$ in the previous examples, we may generally replace $Cu$ with $TM=Ni, Co, Zn$ and the mixture of $Cu, Ni, Co, Zn$ to construct the intermetallic $Sr_{1-x}Ca_{x}TM_5$.

\subsection{$\bf{CaCu_5}$-type intermetallic ${\bf R_{1-x}A_{x}Cu_5}$}

The intermetallic $R_{1-x}A_{x}Cu_5$ ($0<x < 1$)  where $R$ denotes a  rare earth element including the element $Y$ and the mixture thereof  and $A$ denotes a mixture of alkaline elements has two cases of high-$T_c$ superconductivity: one case is the $Cu$ $d$-electron case 
(from the effect of both $A$ and $R$) 
and the other case is  the $Cu$ $s,p$-electron case.  (from the effect of $A$ only or from the effect of both $A$ and $R$).
When  both cases of superconductivity appear at some doping,  the effects of superconductivity of these two cases can be combined. From this combination of effects of superconductivity we can have a larger critical current consisting of $d$-electrons and  $s,p$-electrons and a higher critical temperature $T_c$.
As an example let us consider the  intermetallic $La_{1-x}Sr_{x(1-y)}Ca_{xy}Cu_5$ ($0<x, y < 1$) where we let $R=La$ and  $A$ be a mixture of alkaline elements of the form $A=Sr_{1-y}Ca_{y}$. Let the doping parameters $x,y$ be such that 
\begin{equation}
 \begin{array}{rl}
1850.3\cdot \frac{(1-x)}{1-x+xy} +4912.4\cdot\frac{xy}{1-x+xy}   &=1957.9 \\
&\\
5500(1-y) +6491y  &= 5536 \\
\end{array}
\label{LaCaSrCu5}
\end{equation}
 When this relation of  the degenerate state of channel opening holds (or approximately holds),  the $s,p$-channel connecting the two states of fourth ionization energy of $Sr$ and  $Cu$  of the $3s,3p$-electrons of $Cu$ and the $4s,4p$-electrons of $Sr$ can be opened and also  the $d$-channel connecting the two states of  ionization energy of $La$ and  $Cu$  of the $3d$-electrons of $Cu$ and the $4d$-electrons of $La$ can be opened
This two-channel-opening gives a freedom of electric current with
a direction orthogonal to the $Cu$ plane. 
From this freedom of electric current, the Cooper pairs of the $3s,3p, 3d$-electrons of $Cu$ and the $4s,4p$-electrons of $Sr$  can be formed.

Then, this two-channel-opening can give two types of high-$T_c$ superconductivity: the $s,p$-electron superconductivity and the $d$-electron superconductivity.  This combination of $s,p$-electron superconductivity and $d$-electron superconductivity can give  larger critical current and higher  critical temperature $T_c$.

\subsection{Critical  current density of  ${\bf CaCu_5}$-type superconductors}
We notice that the hexagonal  unit cells of the above examples of  $CaCu_5$-type  intermetallics are with a large cluster of superconducting electrons. Thus these $CaCu_5$-type  intermetallics are with high critical  current density of superconducting current.

Also we notice that the doping of the above  $CaCu_5$-type  intermetallics gives the degenerate state of channel opening for superconductivity.  Thus the doping  of the above  $CaCu_5$-type intermetallics has the effect of  introducing  pinning centers of  superconducting current to these $CaCu_5$-type  intermetallics. 

We may apply more dopings to form more $CaCu_5$-type intermetallics.  For superconductivity these dopings must 
give the degenerate state of  channel opening.
 By more dopings we have the following form of intermetallic which includes the above examples of $CaCu_5$-type  intermetallics:
\begin{equation}
R_{1-x}A_x(TM_{5(1-y)}M_{y})_{1+z}
 \label{pin}
\end{equation}
where $0\leq y\leq y_0=0.02$, 
$0\leq z\leq z_0=0.01$ for some small parameters $y_0, z_0$ 
(The doping parameters $y_0, z_0$ are small to keep the $CaCu_5$-type phase); $TM=Cu, Ni, Co, Zn$,
and $M$ is a mixture of elements $Ti, Cr, Mn, Fe, Co, Zr$, $Nb, Mo, Hf, Ta, W$, $Ga, In, Al, Si, Ge, Sn$. These elements are with states of second or third ionization energies quite close to that of $R$ or $TM$. Thus for some chosen $M$, $y$ and $z$, this intermetallic can 
give the degenerate state of channel opening.

\section{Examples of conventional superconductivity}

As examples of conventional superconductivity let us first consider the transition metals.

\subsection{\bf Transition element $\bf{Ti}$}

Let us first consider the element $Ti$. This element is with the outer shells $3d^24s^2$ of valence electrons. Since $Ti$ is a metal we have that $Ti$ is in the state $\kappa^2<\frac13 h^2$.

For $Ti$ in the state of superconductivity we have that $h^2$ is needed to be reduced such that $\kappa^2\geq \frac{1}{(2+\sqrt{3})^2}h^2$ (It is not that $\kappa^2$ is increased such that $\kappa^2\geq (2+\sqrt{3})^2h^2$ because in this case the $|h|$ will be higher than the least upper bound
 for conventional superconductivity and thus will not be in the state of superconductivity).

 For the decreasing of $h$ let the $\kappa$ be decreased (or let the temperature $T$ be decreased). In this case we need a mechanism for the decreasing of $h$ to reach the state of superconductivity as $T$ decreasing.

Let us then consider the Cooper pair as such a mechanism. Let the two $s$ electrons be with the factor $e^{\pm ih_{Ti}}$ in their wave functions where $h_{Ti}$ is the energy of the states of $s$ valence electrons of $Ti$ (The $\pm$ is for the up and down spins of the two $s$ electrons respectively).

When the temperature $T$ is low enough we have that a $s$ electron of an atom $Ti$ with wave factor
$e^{ih_{Ti}}$ is coupled with the wave factor $e^{-ih_{Ti}}$ of a $s$ electron of another $Ti$. This forms a wave fator $e^{i(h_{Ti}-h_{Ti})}$ of the Cooper pair of the two corresponding $s$ electrons.

This cancelation of energy for the forming of Cooper pair gives the reduction of the energies of the $Ti$ atoms. Thus the forming of Cooper pair is a mechanism for the reduction of energy.
In this case the forming of Cooper pair of $s$ electrons of $Ti$ has reduced the state of normal metallic state of $Ti$ to the state $h^2=3\kappa^2$ of superconductivity.

\subsubsection{\bf Computation of ${\bf T_c}$ of ${\bf Ti}$}

Let us specify the parameter $h$ of $h^2=3\kappa^2$. We notice that when there is no mechanism of pairing there is no reduction of $h$ of $Ti$ to let the metal $Ti$ from the state of normal metal to the state of superconductivity. Thus more reduction of energy gives higher critical temperature $T_c$ of superconductivity. Thus the reduction of energy is proportional to the critical temperature $T_c$ of superconductivity.

Let us find out the exact proportion.
Each $s$ electron of the metal $Ti$ gives an energy $|h_{Ti}|$. The total energy of $s$ electrons of the metal $Ti$ is $2|h_{Ti}|$. Let $h_1$ be the energy reduced from this total energy by the mechanism of Cooper pairing. Then since $h_1$ is proportional to $T_c$ from the state equation of superconductivity we have $h_1^2=3\kappa_{1}^2=3(k_BT_{c1})^2$.

On the other hand
let $h_2=2|h_{Ti}|-h_1$ be the energy after the reduction of energy. Then
this remaining energy must also proportional to the critical temperature $T_c$ of superconductivity because of the general state equation $h_2^2=3\kappa_{2}^2$ of superconductivity. Thus from the state equation of superconductivity we have $h_2^2=3\kappa_{2}^2=3(k_BT_{c2})^2$.

Thus from the formula $h^2=3\kappa^2$ we have that $h_1=h=h_2$. Thus we have the following formula of $T_c$ of $Ti$:
\begin{equation}
k_BT_c=\kappa=D_{Ti}\Delta_{Ti}
\label{Ti}
\end{equation}
where $D_{Ti}$ is a coefficient of $Ti$ such that $D_{Ti}$ is in the range of Type I superconductivity: $\frac{1}{2+\sqrt{3}}\leq D_{Ti}\leq \frac{1}{\sqrt{3}}$ (We shall see that $s$ valence electrons are  for Type I superconductivity),
$\Delta_{Ti}=h=|h_{Ti}|$ is the energy gap of superconductivity of $Ti$ and $T_c$ is the critical temperature of superconductivity of $Ti$.

Let us then find a way to determine the energy $|h_{Ti}|$. Since the $s$ valence electrons are approximately as free electrons we have that the energy $|h_{Ti}|$ is proportional to the ionization energy of the $s$ electrons.

Then
since the $s$ electrons are in the outer shell of the atom $Ti$ we have that the $s$ valence electrons are in the state with the first ionization energy of $Ti$. Thus the maximal value of the  energy parameter $|h_{Ti}|$ is proportional to the first ionization energy of $Ti$.

From the existing table of ionization energies we have that the
first ionization energy of $Ti$ is approximately equal to
$661$ kJ/mol.
Let $|h_{Ti}|\approx \xi 661$ kJ/mol where $\xi$ is the proportional constant (We shall show that $\xi\approx 2.83133971\cdot 10^{-5}$ when we consider the element $Nb$).

Then from (\ref{Ti}), by adoping the coefficient $D_{Ti}=\frac{1}{2+\sqrt{3}}$, we can compute the critical temperature $T_c$ of $Ti$:
\begin{equation}
T_c \approx 0.604 K \quad \mbox{(Computed value of $T_c$ of $Ti$)}
\label{Ti5}
\end{equation}
This roughly agrees with  the experimental value of $T_c\approx 0.39K$ of $Ti$ under ambient pressure.

The descrepancy is due to that $Ti$ has some magnetic properties which from the above unified theory of magnetism and superconductvity can reduce the parameter $h$ and thus reduce the proportional constant $\xi$
(We also notice that the experimental value of the highest critical temperature $T_c$ approximately equals to $T_c\approx 3.35K$ under high pressure $56 Gpa$. Thus there is a variance of the experimental value of $T_c\approx 0.39K$  due to pressure \cite{Sch}).

We identify this superconductivity induced by the $s$ valence electrons as the Type I superconductivity.

\subsection{\bf Transition elements $\bf{Zr}$ and $\bf{Hf}$}

Similarly since the transition elements $Zr$ and $Hf$ are also with the outer shells of the form $4d^25s^2$  and $5d^26s^2$ respectively which are similar to that of $Ti$ we have that $Zr$ and $Hf$ can be in superconductivity where the superconductivity is induced by the $s$ valence electrons and is as the Type I superconductivity.

Also the formulas of temperature $T_c$ are given by:
\begin{equation}
k_BT_c=D_{e}\Delta_{e}
\label{ZrandHf}
\end{equation}
where $D_{e}$ is a coefficient of $e=Zr$ and $e=Hf$ respectively such that  $\frac{1}{2+\sqrt{3}}\leq D_{e}\leq \frac{1}{\sqrt{3}}$.

Since $e=Ti, Zr, Hf$ are of the same form of outer shells of Type I superconductivity we have that $\Delta_{e}$ are of the same order for $e=Ti, Zr, Hf$ with a variance due to the effect of pressure on $e=Ti, Zr, Hf$.

From the existing table of ionization energies  we have that the first ionization energies of $Zr$ and $Hf$ are approximately given by $670$ and $531$ respectively. Let
$|h_{Zr}|\approx \xi 670$ kJ/mol and $|h_{Hf}|\approx \xi 531$ kJ/mol.

Then from (\ref{ZrandHf}), by adoping the coefficient $D_{Zr}=D_{Hf}=D_{Ti}=\frac{1}{2+\sqrt{3}}$, we have that the  computed $T_c$ of $Zr$ and $Hf$ are respectively given by $0.611K$ and $0.484K$.
This roughly
agrees with  the experimental values of $T_c\approx 0.546K$ of $Zr$ and
$T_c\approx 0.12K$ of $Hf$ under ambient pressure \cite{Sch}.

As $Ti$, the descrepancy is due to that $Zr$ and $Hf$ have some magnetic properties (We notice that the experimental value of the highest critical temperatures $T_c$ of $Zr$ and $Hf$ approximately equal to $T_c\approx 11K$ under pressure $30 Gpa$ and $T_c\approx 8.6K$ under pressure $62 Gpa$ respectively. Thus there are variances of the experimental values of $T_c$ of $Zr$ and $Hf$ due to pressure \cite{Sch}).

\subsection{\bf Transition element $\bf{V}$}

Let us then consider the element $V$. This element is with the outer shells $3d^34s^2$ of valence electrons. Since $V$ is a metal we have that $V$ is in the state $\kappa^2<\frac13 h^2$.

Similar to the element $T_i$ for $V$ in the state of superconductivity we have that $h^2$ is needed to be reduced such that $\kappa^2\geq \frac{1}{(2+\sqrt{3})^2}h^2$. Similar to the element $T_i$ for the decreasing of $h$ the $\kappa$ is to be decreased (or the temperature $T$ is to be decreased). In this case we need a mechanism for the decreasing of $h$ to reach the state of superconductivity.

Again let the Cooper pair be such a mechanism. Let the three $d$ valence electrons be with the wave factor $e^{ih_{V}}$ where $h_{V}$ is the energy of the $d$ electrons for the forming of Cooper pairing.

When the temperature $T$ is low enough we have the phase transition of Bose-Einstein condensation. This phase transition is with the changing of the spin of one of $d$ electrons such that the state of this $d$ electron is changed from the state $e^{ih_{V}}$ to the state $e^{-ih_{V}}$.

This gives the state $e^{i(-h_{V}+h_{V})}$ where $e^{ih_{V}}$ is the state of another $d$ electron of another element $V$. This is thus as a Cooper pair of $d$ electrons.
(When the $d$ elecrons forming Cooper pairs we have that the two $s$ electrons may or may not form pair.

 Let us consider the case that the two $s$ electrons do not form a Cooper pair. This case agrees with the phenomenon that the high-$T_c$ superconductivity is usually with the $d$ waves while the high-$T_c$ superconductivity is extended from the conventional Type II superconductivity and thus the conventional Type II superconductivity is also with the $d$ waves and Cooper $d$-pair).

Thus the forming of Cooper pair reduces the total energy  $2|h_{V}|$ of two $d$ electrons of $V$.
Then as the case of $T_i$, let $h_1$ be the energy reduced from this total energy by the mechanism of Cooper pairing. Then since $h_1$ is proportional to $T_c$ from the state equation of superconductivity we have $h_1^2=3\kappa_{1}^2=3(k_BT_{c1})^2$.

On the other hand
let $h_2=2|h_{V}|-h_1$ be the energy after the reduction of energy. Then
this remaining energy must also proportional to the critical temperature $T_c$ of superconductivity because of the general state equation $h_2^2=3\kappa_{2}^2$ of superconductivity. Thus from the state equation of superconductivity we have $h_2^2=3\kappa_{2}^2=3(k_BT_{c2})^2$.

Thus from the formula $h^2=3\kappa^2$ we have that $h_1=h=h_2$. Thus we have the following formula of $T_c$ of $V$:
\begin{equation}
k_BT_c=\kappa=D_{V}\Delta_{V}
\label{V}
\end{equation}
where $D_{V}$ is a coefficient for $V$ such that $D_{V}=\frac{1}{\sqrt{3}}$ which is the coefficient for the Type II superconductivity and high-$T_c$ superconductivity,
$\Delta_{V}=h=|h_{V}|$ is the energy gap of superconductivity of $V$ and $T_c$ is the critical temperature of superconductivity of $V$.

Let us then find a way to determine the energy $|h_{V}|$.
As similar to the case of high-$T_c$ superconductivity
we have that the $d$ valence electrons of $V$ can be in the state with the third ionization energy of $V$.

From the existing table of ionization energies we have that the
third ionization energy of $V$ is approximately equal to
$2828$ kJ/mol.

Let $|h_{V}|\approx \xi 2828$ kJ/mol where $\xi$ is the proportional constant.
 Then from (\ref{V}) we can compute the critical temperature $T_c$ of $V$:
\begin{equation}
T_c \approx 5.56K \quad \mbox{(Computed value of $T_c$ of $V$)}
\label{V5}
\end{equation}
This agrees with  the experimental value of $T_c\approx 5.38K$ of $V$ under ambient pressure.

It is clear that we should identify the superconductivity of $V$ induced by the  $d$ electrons as the conventional Type II superconductivity.

\subsection{\bf Transition element $\bf{Nb}$}

Let us then consider the element $Nb$. This element is with the outer shells $4d^45s^1$ of valence electrons. Since $Nb$ is a metal we have that $Nb$ is in the state $\kappa^2<\frac13 h^2$. For $Nb$ in the state of superconductivity we have that $h^2$ is needed to be reduced such that $\kappa^2\geq \frac{1}{(2+\sqrt{3})^2}h^2$. In this case we need a mechanism for the decreasing of $h$ to reach the state of superconductivity.

Again let the Cooper pair be such a mechanism. Let the four $d$ electrons be with the wave factor $e^{ih_{Nb}}$ where $h_{Nb}$ is the energy of the $d$ electrons for the forming of Cooper pairing.

When the temperature $T$ is low enough we have the phase transition of Bose-Einstein condensation. This phase transition is with the changing of the spin of two of $d$ electrons such that the state of these two $d$ electron are changed from the state $e^{ih_{Nb}}$ to the state $e^{-ih_{Nb}}$.
This gives the state $e^{i(-h_{Nb}+h_{Nb})}$ where $e^{ih_{Nb}}$ is the state of another two $d$ electrons of another $Nb$ atom. Thus we have two Cooper pairs of $d$ electrons.

Thus the forming of the Cooper pairs reduces the total energy $4|h_{Nb}|$ of the four $d$ electrons of $Nb$.
Then as the case of $V$, let $h_1$ be the energy reduced from this total energy by the mechanism of Cooper pairing. Then since $h_1$ is proportional to $T_c$ from the state equation of superconductivity we have $h_1^2=3\kappa_{1}^2=3(k_BT_{c1})^2$.

On the other hand
let $h_2=4|h_{Nb}|-h_1$ be the energy after the reduction of energy. Then
this remaining energy must also be proportional to the critical temperature $T_c$ of superconductivity because of the general state equation $h_2^2=3\kappa_{2}^2$ of superconductivity. Thus from the state equation of superconductivity we have $h_2^2=3\kappa_{2}^2=3(k_BT_{c2})^2$.

Thus from the formula $h^2=3\kappa^2$ we have that $h_1=h=h_2$. Thus we have the following formula of $T_c$ of $Nb$:
\begin{equation}
k_BT_c=\kappa=D_{Nb}\Delta_{Nb}
\label{Nb}
\end{equation}
where $D_{Nb}$ is a coefficient for $Nb$ such that $D_{Nb}=\frac{1}{\sqrt{3}}$ which is the constant for the Type II superconductivity and high-$T_c$ superconductivity,
$\Delta_{Nb}=h=2|h_{Nb}|$ is the energy gap of superconductivity of $Nb$ and $T_c$ is the critical temperature of superconductivity of $Nb$.

Let us then determine the energy $|h_{Nb}|$.
As similar to the case of high-$T_c$ superconductivity
we have that the $d$ valence electrons of $Nb$ can be in the state with the third ionization energy of $Nb$.

From the existing table of ionization energies we have that the
third ionization energy of $Nb$ is approximately equal to
$2416$ kJ/mol.

Let $|h_{Nb}|\approx \xi 2416$ kJ/mol where $\xi$ is the proportional constant.
 Then from (\ref{Nb}) let us compute the constant $\xi$ by using the experimental value $T_c\approx 9.5K$:
\begin{equation}
k_B9.5=\kappa=\frac{1}{\sqrt{3}}\xi 2\cdot 2416
\label{Nb2}
\end{equation}
This gives $\xi\approx 2.83133971\cdot 10^{-5}$. From this value of $\xi$ we have that the computed $T_c$ of $Nb$ is equal to the experimental value $T_c=9.5K$ of $Nb$ (under ambient pressure).

Since the superconductivity of $Nb$ is induced by the $d$ electrons this is as the Type II superconductivity.

\subsection{\bf Transition element $\bf{Ta}$}

Let us then consider the element $Ta$. This element is with the outer shells $5d^36s^2$ of valence electrons as similar to that of $V$.

Similar to $V$ we have that $Ta$ can be in superconductivity where the superconductivity is induced by  the change of the spin of the $d$ electrons and is  the Type II superconductivity. Also the formula of critical temperature $T_c$ is given by:
\begin{equation}
k_BT_c=D_{Ta}\Delta_{Ta}
\label{Ta}
\end{equation}
where $D_{Ta}=\frac{1}{\sqrt{3}}$ is the coefficient for Type II superconductivity, and the energy gap $\Delta_{Ta}=|h_{Ta}|$ where $h_{Ta}$ is from the wave factor $e^{ih_{Ta}}$ of the $d$ electrons of $Ta$.

Let us then determine the energy $|h_{Ta}|$.
As similar to the case of high-$T_c$ superconductivity
we have that the $d$ valence electrons of $V$ can be in the state with the third ionization energy of $Ta$.

From the existing table of ionization energies we have that the
third ionization energy of $Ta$ is approximately equal to
$2250$ kJ/mol.

Let $|h_{Ta}|\approx \xi 2250$ kJ/mol where $\xi$ is the proportional constant.
Then from (\ref{Ta}) we can compute the critical temperature $T_c$ of $Ta$:
\begin{equation}
T_c \approx 4.424K \quad \mbox{(Computed value of $T_c$ of $Ta$)}
\label{Ta5}
\end{equation}
This agrees with  the experimental value of $T_c\approx 4.483K$ of $Ta$ under ambient pressure.

\subsection{\bf Transition element $\bf{Tc}$}

Let us then consider the element $Tc$. This element is with the outer shells $4d^65s^1$ of valence electrons. Since $Tc$ is a metal we have that $Tc$ is in the state $\kappa^2<\frac13 h^2$.

Similar to the element $T_i$ for $Tc$ in the state of superconductivity we have that $h^2$ is needed to be reduced such that $\kappa^2\geq \frac{1}{(2+\sqrt{3})^2}h^2$. In this case we need a mechanism for the decreasing of $h$ to reach the state of superconductivity.

Again let the Cooper pair be such a mechanism.
Since there are six $d$ electrons we have that one of the six $d$ electrons is with the spin different from the other five $d$ electrons which are with the same spin.

When the temperature $T$ is low enough
the $d$ electron with different spin and with the wave factor $e^{-ih_{Tc1}}$ is coupled with a $d$ electron with wave factor $e^{ih_{Tc1}}$ of another atom $Tc$. 
This forms a Cooper pair with the wave factor $e^{i(-h_{Tc1}+h_{Tc1})}$. This gives a reduction of energy, as explained in the above analysis.

Because the way of forming $d$-pair is similar to the way of forming pair of the $s$ electrons we have that the energy parameter $h_{Tc1}$ is proportional to the first ionization energy of $Tc$, as the case of the $s$ electrons of the atom $Ti$. In this case the superconductivity will be the Type I superconductivity because of the similarity to the $s$ electrons.

Then let us consider the other $d$ electrons. Because of the one pair of Cooper pair there are four
$d$ electrons to be paired or unpaired. Let us consider the case that when the temperature $T$ is low enough the spins of two of the four $d$ electrons are changed such that the state of these two $d$ electron is changed from the state $e^{ih_{Tc2}}$ to the state $e^{-ih_{Tc2}}$.

In this case these two $d$ electrons are coupled with two $d$ electrons of another atom $Tc$ with the state $e^{ih_{Tc2}}$ to form two Cooper pairs. In this case the energy parameter $h_{Tc2}$ is proportional to the second ionization energy of $Tc$, and is not proportional to the third ionization energy of $Tc$. This is because that if the energy parameter $h_{Tc2}$ is proportional to the third ionization energy of $Tc$, then the total energy parameter $h$ will be greater than the upper bound $h_{IImax}$ for conventional Type II superconductivity and thus the state of conventional Type II superconductivity will be destroyed.

Then the forming of Coopers pair reduces the total energy $2|h_{Tc1}|+4|h_{Tc2}|$ of the six $d$ electrons for the forming of Cooper pairing.

Then as the case of $T_i$, let $h_1$ be the energy reduced from this total energy by the mechanism of Cooper pairing and
let $h_2=2|h_{Tc1}|+4|h_{Tc2}|-h_1$ be the energy after the reduction of energy.

Then by the same reason as that for the element $Ti$,
from the formula $h^2=3\kappa^2$ we have that $h_1=h=h_2$. Thus we have the following formula of critical temperature $T_c$ of $Tc$:
\begin{equation}
k_BT_c=\kappa=D_{Tc}\Delta_{Tc}
\label{Tc}
\end{equation}
where $D_{Tc}=\frac{1}{\sqrt{3}}$ and
$\Delta_{Tc}=h=|h_{Tc1}|+2|h_{Tc2}|$ is the energy gap of superconductivity of $Tc$.

 Then from the existing table of ionization energies we have that the
first and second ionization energies of $Tc$ are approximately equal to $702$ kJ/mol
and $1470$ kJ/mol.

Let $|h_{Tc1}|\approx \xi 702$ kJ/mol and $|h_{Tc2}|\approx \xi 1470$ kJ/mol where $\xi$ is the proportional constant.

Then from (\ref{Tc}) we can compute the critical temperature $T_c$ of $Tc$:
\begin{equation}
T_c \approx 7.16K \quad \mbox{(Computed value of $T_c$ of $Tc$)}
\label{Tc5}
\end{equation}
This agrees with  the experimental value of $T_c\approx 7.7K$ of $Tc$ under ambient pressure.

As explained in the above we identify the superconductivity of $Tc$ as the Type I superconductivity.

\subsection{\bf Transition element $\bf{Re}$}

Let us then consider the element $Re$. This element is with the outer shells $6d^57s^2$ of valence electrons. Since $Re$ is a metal we have that $Re$ is in the state $\kappa^2<\frac13 h^2$.

Similar to the element $T_i$ for $Re$ in the state of superconductivity we have that $h^2$ is needed to be reduced such that $\kappa^2\geq \frac{1}{(2+\sqrt{3})^2}h^2$. In this case we need a mechanism for the decreasing of $h$ to reach the state of superconductivity.

Again let the Cooper pair be such a mechanism.
Since there are five $d$ electrons we have that the five $d$ electrons are with the same spin while the two $s$ electrons are with different spins.

When the temperature $T$ is low enough
one $s$ electron of an atom $Re$ with the wave factor $e^{-ih_{Re}}$ is coupled with a $s$ electron with wave factor $e^{ih_{Re}}$ of another atom $Re$. This forms a Cooper pair with the wave factor $e^{i(-h_{Re}+h_{Re})}$. This gives a reduction of energy, as explained in the above analysis.

For this $ss$-pair the energy parameter $h_{Re}$ is proportional to the first ionization energy of $Re$, as the case of the $ss$-pair of the atom $Ti$. In this case the superconductivity will be the Type I superconductivity because of the similarity to the $s$ electrons.

Then let us consider the other five $d$ electrons. Let us consider the case that when the temperature $T$ is low the spins of the five $d$ electrons are remained to be unchanged. In this case there are no more Cooper pairs except the Cooper $ss$-pair.

Then the forming of this only Cooper pair reduces the total energy $2|h_{Re}|$ of interaction of the two $s$ electrons of $Re$.

Then as the case of $T_i$, let $h_1$ be the energy reduced from this total energy by the mechanism of Cooper pairing and
let $h_2=2|h_{Re}|-h_1$ be the energy after the reduction of energy.

Then by the same reason as that for the element $Ti$,
from the formula $h^2=3\kappa^2$ we have that $h_1=h=h_2$. Thus we have the following formula of critical temperature $T_c$ of $Re$:
\begin{equation}
k_BT_c=\kappa=D_{Re}\Delta_{Re}
\label{Re}
\end{equation}
where $D_{Re}=\frac{1}{\sqrt{3}}$ and
$\Delta_{Re}=h=|h_{Re}|$ is the energy gap of superconductivity of $Re$.

 Then from the existing table of ionization energies we have that the
first ionization energy of $Re$ is approximately equal to $760$ kJ/mol.

Let $|h_{Re}|\approx \xi 760$ kJ/mol where $\xi$ is the proportional constant.
Then from (\ref{Re}) we can compute the critical temperature $T_c$ of $Re$:
\begin{equation}
T_c \approx 1.49K \quad \mbox{(Computed value of $T_c$ of $Re$)}
\label{Re5}
\end{equation}
This agrees with  the experimental value of $T_c\approx 1.4K$ of $Re$ under ambient pressure.

As explained in the above we identify the superconductivity of $Re$ as the Type I superconductivity.

\subsection{\bf Transition element $\bf{Mo}$}

Let us then consider the element $Mo$. This element is with the outer shells $5d^56s^1$ of valence electrons. Since $Mo$ is a metal we have that $Mo$ is in the state $\kappa^2<\frac13 h^2$.

Similar to the element $T_i$ for $Mo$ in the state of superconductivity we have that $h^2$ is needed to be reduced such that $\kappa^2\geq \frac{1}{(2+\sqrt{3})^2}h^2$. In this case we need a mechanism for the decreasing of $h$ to reach the state of superconductivity.

Again let the Cooper pair be such a mechanism.
Since there are five $d$ valence electrons and one $s$ valence electron we have that the $s$ valence electron may or may not with the spin different from the five $d$ electrons which are with the same spins. Let us consider the case that the $s$ valence electron is with the spin different from the five $d$ electrons.

When the temperature $T$ is low enough
the $s$ electron with different spin and with the wave factor $e^{-ih_{Mo1}}$ is coupled with a $d$ electron with wave factor $e^{ih_{Mo2}}$ of another atom $Mo$ where the two energy parameters $h_{Mo1}$ and $h_{Mo2}$ are of the same sign. This forms a Cooper pair with the wave factor 
$e^{i(-h_{Mo1}+h_{Mo2})}$. This gives a reduction of energy, as explained in the above analysis.

Because the way of forming the $sd$-pair is similar to the way of forming the $ss$-pair of the $s$ electrons we have that the energy parameter $h_{Mo1}$ is proportional to the first ionization energy of $Mo$ and the energy parameter $h_{Mo2}$ is proportional to the second ionization energy of $Mo$, as the case of the $s$ electrons of the atoms $Ti$ and $Re$. In this case the superconductivity will be the Type I superconductivity because of the similarity to the $s$ electrons.

Then let us consider the other four $d$ electrons. Let us consider the case that when the temperature $T$ is low the spins of the five $d$ electrons are remained to be unchanged. In this case there are no more Cooper pairs except the above Cooper $sd$-pair.

Then the forming of this only Cooper $sd$-pair reduces the total energy $|h_{Mo1}|+|h_{Mo2}|$ of interaction of the $s$ and $d$ electrons.

Then as the case of $T_i$, let $h_1$ be the energy reduced from this total energy by the mechanism of Cooper pairing and
let $h_2=|h_{Mo1}|+|h_{Mo2}|-h_1$ be the energy after the reduction of energy.

Then by the same reason as that for the element $Ti$,
from the formula $h^2=3\kappa^2$ we have that $h_1=h=h_2$. Thus we have the following formula of critical temperature $T_c$ of $Mo$:
\begin{equation}
k_BT_c=\kappa=D_{Mo}\Delta_{Mo}
\label{Mo}
\end{equation}
where $D_{Mo}=\frac{1}{2+\sqrt{3}}$ for Type I superconductivity and
$\Delta_{Mo}=h=\frac12(|h_{Mo1}|+|h_{Mo2}|)$ is the energy gap of superconductivity of $Mo$.

 Then from the existing table of ionization energies we have that the
first and second ionization energies of $Mo$ are approximately equal to $684.3$ kJ/mol and $1559.2$ kJ/mol espectively.

Let $|h_{Mo1}|\approx \xi 684.3$ kJ/mol and $|h_{Mo2}|\approx \xi 1559.2$ kJ/mol
where $\xi$ is the proportional constant.

 Then from (\ref{Mo}) we can compute the critical temperature $T_c$ of $Mo$:
\begin{equation}
T_c \approx 1.023K \quad \mbox{(Computed value of $T_c$ of $Mo$)}
\label{Mo5}
\end{equation}
This agrees with  the experimental value of $T_c\approx 0.92K$ of $Mo$ under ambient pressure.

As explained in the above we identify the superconductivity of $Mo$ as the Type I superconductivity.

\subsection{\bf Transition element $\bf{W}$}

Let us then consider the element $W$. This element is with the outer shells $6d^47s^2$ of valence electrons. Since $W$ is a metal we have that $W$ is in the state $\kappa^2<\frac13 h^2$.

Similar to the element $T_i$ for $W$ in the state
of superconductivity we have that $h^2$ is needed to be reduced such that
$\kappa^2\geq \frac{1}{(2+\sqrt{3})^2}h^2$. In this case we need a mechanism for the decreasing of $h$ to reach the state of superconductivity.

Again let the Cooper pair be such a mechanism. Let us consider the case that when the temperature $T$ is low the spins of the four $d$ electrons are remained to be unchanged.

Then there are four $d$ valence electrons with the same spins and two $s$ valence electrons with different spins. Let us consider the case that the element $W$ is a metal with
paramagnetic
properties (The two $s$ valence electrons with different spins corresponding to antiferromagnetism and the four $d$ valence electrons with the same spins corresponding to ferromagnetism give paramagnetic properties).

Let the two $s$ valence electrons be with the wave factors $e^{\pm ih_{W1}}$ respectively. Then as the above elements we have that the energy parameter $h_{W1}$ is proportional to the first ionization energy of $W$.

Then
because of the paramagnetic
property we have that the corresponding proportional constant $\xi_1$ for $h_{W1}$ is smaller than the proportional constant $\xi$ for the above elements (We remark that the proportional constant $\xi$ is for the case with the highest critical temperature $T_c$ for superconductivity where the effect of magnetic property does not appear). This is because that due to the interaction of the four $d$ valence electrons with the two $s$ valence electrons it needs more energy to separate the two $s$ electrons of a $W$ atom to form Cooper pairs with the two $s$ electrons of other $W$ atoms. From this the remained energy of the formed Cooper pairs are with smaller energy and thus the proportional constant $\xi_1$ is smaller.

When the temperature $T$ is low enough
a $s$ electron  with the wave factor $e^{-ih_{W1}}$ of a $W$ atom is coupled with a $s$ electron with different spin and with wave factor $e^{ih_{W1}}$ of another atom $W$. This forms a Cooper pair with the wave factor $e^{i(-h_{W1}+h_{W1})}$. This gives a reduction of energy, as explained in the above analysis.

Then let us consider the four $d$ electrons. Let us consider the case that when the temperature $T$ is low the spins of the four $d$ electrons are remained to be unchanged. In this case there are no more Cooper pairs except the only Cooper $ss$-pair.

Then the forming of this only Cooper pair reduces the total energy $2|h_{W1}|$ of interaction of the two $s$ electrons of $W$ forming Cooper pairs.

Then as the case of $T_i$, let $h_1$ be the energy reduced from this total energy by the mechanism of Cooper pairing and
let $h_2=2|h_{W1}|-h_1$ be the energy after the reduction of energy.

Then by the same reason as that for the element $Ti$,
from the formula $h^2=3\kappa^2$ we have that $h_1=h=h_2$. Thus we have the following formula of critical temperature $T_c$ of $W$:
\begin{equation}
k_BT_c=\kappa=D_{W}\Delta_{W}
\label{W}
\end{equation}
where $D_{W}=\frac{1}{2+\sqrt{3}}$ and
$\Delta_{W}=h=|h_{W1}|$ is the energy gap of superconductivity of $W$.

 Then from the existing table of ionization energies we have that the
first ionization energy of $W$ is approximately equal to $770$ kJ/mol.
Let $|h_{W}|\approx \xi_1 770$ kJ/mol where $\xi_1$ is the proportional constant.

Then from (\ref{W}) we can compute the critical temperature $T_c$ of $W$:
\begin{equation}
T_c \approx \frac{\xi_1}{\xi}\cdot 0.702K \quad \mbox{(Computed value of $T_c$ of $W$)}
\label{W5}
\end{equation}

Let us use (\ref{W5}) and the experimental value of $T_c\approx 0.012K$ of $W$ to determine the proportional constant $\xi_1$. By substituting the value of $T_c\approx 0.015K$ (which approximates the experimental value of $T_c\approx 0.012K$) into (\ref{W5}) we have:
 \begin{equation}
\xi_1\approx 0.0214\xi
\label{W6}
\end{equation}
Then we have that the computed critical temperature $T_c\approx 0.015K$ which agrees with the experimental value of $T_c\approx 0.012K$ of $W$.

We notice that this proportional constant $\xi_1<\xi$ effectively gives an antiferromagnetic effect to the two $s$ valence electrons with different spins. This constant $\xi_1$ will be used to compute the critical temperature $T_c$ of some other elements.

As explained in the above we identify the superconductivity of $W$ as the Type I superconductivity.

\subsection{\bf Transition element $\bf{Ru}$}

Let us then consider the element $Ru$. This element is with the outer shells $4d^75s^1$ of valence electrons. Since $Ru$ is a metal we have that $Ru$ is in the state $\kappa^2<\frac13 h^2$.

Similar to the element $T_i$ for $Ru$ in the state of superconductivity we have that $h^2$ is needed to be reduced such that $\kappa^2\geq \frac{1}{(2+\sqrt{3})^2}h^2$. In this case we need a mechanism for the decreasing of $h$ to reach the state of superconductivity.

Again let the Cooper pair be such a mechanism.
Since there are seven $d$ electrons we have that two of the seven $d$ electrons are with the spin different from the other five $d$ electrons which are with the same spin.

When the temperature $T$ is low enough one of the two $d$ electrons with different spin and with the wave factor $e^{-ih_{Ru}}$ is coupled with a $d$ electron with wave factor $e^{ih_{Ru}}$ of another atom $Ru$.
This forms a Cooper pair with the wave factor $e^{i(-h_{Ru}+h_{Ru})}$. This gives a reduction of energy, as explained in the above analysis. Because the way of forming $dd$-pair is similar to the way of forming  $ss$-pair we have that the energy parameter $h_{Ru}$ is proportional to the first ionization energy of $Ru$, as the $ss$-pair of electrons of the atom $Ti$. In this case the superconductivity will be the Type I superconductivity because of the similarity to the $s$ electrons.

 Because of this one $dd$-pair there are five
$d$ electrons and one $s$ electron to be paired or unpaired. Let us consider the case that when the temperature $T$ is low the spins of the other five $d$ electrons and the $s$ electron are remained to be unchanged and are not separated to be in Cooper pair form (Other cases would be more difficult to occur).
 In this case there are no more Cooper pairs except the only Cooper $dd$-pair.

Then the forming of this only Cooper pair reduces the total energy $2|h_{Ru}|$ of interaction of the two $d$ electrons of $Ru$ forming Cooper pairs.

Then as the case of $T_i$, let $h_1$ be the energy reduced from this total energy by the mechanism of Cooper pairing and
let $h_2=2|h_{Ru}|-h_1$ be the energy after the reduction of energy.

Then by the same reason as that for the element $Ti$,
from the formula $h^2=3\kappa^2$ we have that $h_1=h=h_2$.

Thus we have the following formula of critical temperature $T_c$ of $Ru$:
\begin{equation}
k_BT_c=\kappa=D_{Ru}\Delta_{Ru}
\label{Ru}
\end{equation}
where $D_{Ru}=\frac{1}{2+\sqrt{3}}$ and
$\Delta_{Ru}=h=|h_{Ru}|$ is the energy gap of superconductivity of $Ru$.

 Then from the existing table of ionization energies we have that the
first ionization energy of $Ru$ is approximately equal to $710.2$ kJ/mol.
Let $|h_{Ru}|\approx \xi 710.2$ kJ/mol where $\xi$ is the proportional constant.

Then from (\ref{Ru}) we can compute the critical temperature $T_c$ of $Ru$:
\begin{equation}
T_c \approx 0.647K \quad \mbox{(Computed value of $T_c$ of $Ru$)}
\label{Ru5}
\end{equation}
This roughly agrees with  the experimental value of $T_c\approx 0.51K$ of $Ru$ under ambient pressure.

The discrepancy between computation and experiment may due to the magnetic property of the Cooper $d$-pair of $Ru$ where the magnetic property can give the effect of reducing the parameter $\xi$.

As explained in the above we identify the superconductivity of $Ru$ as the Type I superconductivity.

\subsection{\bf Transition element $\bf{Os}$}

Let us then consider the element $Os$. This element is with the outer shells $5d^66s^2$ of valence electrons. Since $Os$ is a metal we have that $Os$ is in the state $\kappa^2<\frac13 h^2$.

Similar to the element $T_i$ for $Os$ in the state of superconductivity we have that $h^2$ is needed to be reduced such that $\kappa^2\geq \frac{1}{(2+\sqrt{3})^2}h^2$. Similar to the element $T_i$ we need to have a mechanism for the decreasing of $h$ to reach the state of superconductivity.

Again let the Cooper pair be such a mechanism.
There are two $s$ electrons with the different spins and
six $d$ electrons such that one of the six $d$ electrons is with the spin different from the other five $d$ electrons which are with the same spin.

When the temperature $T$ is low enough
one of the two $s$ electrons with the wave factor $e^{-ih_{Os}}$ is coupled with a $s$ electron with wave factor $e^{ih_{Os}}$ of another atom $Os$.

 This forms a Cooper $ss$-pair with the wave factor $e^{i(-h_{Os}+h_{Os})}$. This gives a reduction of energy, as explained in the above analysis. Because of the $ss$-pairing we have that the energy parameter $h_{Os}$ is proportional to the first ionization energy of $Os$. In this case the superconductivity is the Type I superconductivity.

 The six
$d$ electrons are to be paired or unpaired. Let us consider the case that when the temperature $T$ is low the spins of the six $d$ electrons are remained to be unchanged and are not separated to be in Cooper pair form (Other cases would be more difficult to occur).
 In this case there are no more Cooper pairs except the only Cooper $ss$-pair.

Then the forming of this only Cooper pair reduces the total energy $2|h_{Os}|$ of interaction of the two $d$ electrons of $Os$ forming Cooper pairs.

Then as the case of $T_i$, let $h_1$ be the energy reduced from this total energy by the mechanism of Cooper pairing and
let $h_2=2|h_{Os}|-h_1$ be the energy after the reduction of energy.

Then by the same reason as that for the element $Ti$,
from the formula $h^2=3\kappa^2$ we have that $h_1=h=h_2$.
Thus we have the following formula of critical temperature $T_c$ of $Os$:
\begin{equation}
k_BT_c=\kappa=D_{Os}\Delta_{Os}
\label{Os}
\end{equation}
where $D_{Os}=\frac{1}{2+\sqrt{3}}$ and
$\Delta_{Os}=h=|h_{Os}|$ is the energy gap of superconductivity of $Os$.

Then from the existing table of ionization energies we have that the
first ionization energy of $Os$ is approximately equal to $840$ kJ/mol.
Let $|h_{Os}|\approx \xi 840$ kJ/mol where $\xi$ is the proportional constant.

Then from (\ref{Os}) we can compute the critical temperature $T_c$ of $Os$:
\begin{equation}
T_c \approx 0.766K \quad \mbox{(Computed value of $T_c$ of $Os$)}
\label{Os5}
\end{equation}
This roughly agrees with  the experimental value of $T_c\approx 0.655K$ of $Os$ under ambient pressure.

As the element $Ru$, the discrepancy between computation and experiment may due to the magnetic property of the Cooper $ss$-pair of $Os$ where the magnetic property can give the effect of reducing the parameter $\xi$.

As explained in the above we identify the superconductivity of $Os$ as the Type I superconductivity.

\subsection{\bf Transition element $\bf{Ir}$}

Let us then consider the element $Ir$. This element is with the outer shells $5d^76s^2$ of valence electrons. Since $Ir$ is a metal we have that $Ir$ is in the state $\kappa^2<\frac13 h^2$.

Similar to the element $T_i$ for $Ir$ in the state of superconductivity we have that $h^2$ is needed to be reduced such that $\kappa^2\geq \frac{1}{(2+\sqrt{3})^2}h^2$. In this case we need a mechanism for the decreasing of $h$ to reach the state of superconductivity.

Again let the Cooper pair be such a mechanism.
Since there are seven $d$ electrons we have that two of the seven $d$ electrons is with the spin different from the other five $d$ electrons which are with the same spin. Also the two $s$ electrons are with different spins.

When the temperature $T$ is low enough
one of the two $s$ electrons with the wave factor $e^{-ih_{Ir1}}$ is coupled with a $s$ electron with wave factor $e^{ih_{Ir1}}$ of another atom $Ir$. This forms a Cooper pair with the wave factor $e^{i(-h_{Ir1}+h_{Ir1})}$. This gives a reduction of energy, as explained in the above analysis. Because of $ss$-pair we have that the energy parameter $h_{Ir1}$ is proportional to the first ionization energy of $Ir$, as the case of the $s$ electrons of the atom $Ti$. In this case the superconductivity is the Type I superconductivity.

Then when the temperature $T$ is low enough
one of the two $d$ electrons with spin different from the other five $d$ electrons and with the wave factor $e^{-ih_{Ir2}}$ is coupled with a $d$ electron with wave factor $e^{ih_{Ir2}}$ of another atom $Ir$.
 This forms a Cooper pair with the wave factor $e^{i(-h_{Ir2}+h_{Ir2})}$. This gives a reduction of energy, as explained in the above analysis. Because of $dd$-pair and there is a $ss$-pair we have that the energy parameter $h_{Ir2}$ is proportional to the second ionization energy of $Ir$ (It may also be proportional to the third ionization energy of $Ir$. As we shall show it should be proportional to the second ionization energy of $Ir$). This gives two $dd$-pairs.

Then let us consider the other five $d$ electrons. Because of the two $dd$-pairs there are three
$d$ electrons to be paired or unpaired. Let us consider the case that when the temperature $T$ is low enough the spin of one of the three $d$ electrons is changed such that the state of this $d$ electron is changed from the state $e^{ih_{Ir3}}$ to the state $e^{-ih_{Ir3}}$. In this case this $d$ electron is coupled with a $d$ electron of another atom $Ir$ with the state $e^{ih_{Ir3}}$ to form a Cooper $dd$-pair. In this case the energy parameter $h_{Ir3}$ is proportional to the third ionization energy of $Ir$, since there already have $ss$-pair and $dd$-pair proportional to the first and second ionization energies of $Ir$ respectively.

Then the forming of Cooper pairs reduces the total energy $2|h_{Ir1}|+4|h_{Ir2}|+2|h_{Ir3}| $ of interaction of six $d$ electrons of $Ir$ and two $s$ electrons.

Then as the case of $T_i$, let $h_1$ be the energy reduced from this total energy by the mechanism of Cooper pairing and
let $h_2=2|h_{Ir1}|+4|h_{Ir2}|+2|h_{Ir3}|-h_1$ be the energy after the reduction of energy.

Then by the same reason as that for the element $Ti$,
from the formula $h^2=3\kappa^2$ we have that $h_1=h=h_2$. Thus we have the following formula of critical temperature $T_c$ of $Ir$:
\begin{equation}
k_BT_c=\kappa=D_{Ir}\Delta_{Ir}
\label{Ir}
\end{equation}
where $D_{Ir}=\frac{1}{2+\sqrt{3}}$ and
$\Delta_{Ir}=h=|h_{Ir1}|+2|h_{Ir2}|+|h_{Ir3}|$ is the energy gap of superconductivity of $Ir$.

Then from the existing table of ionization energies we have that the
first, second and third ionization energies of $Ir$ are approximately equal to $880$ kJ/mol,
 $1600$ kJ/mol and $2600$ kJ/mol respectively.

 Then
let $|h_{Ir1}|\approx \xi_1 880$ kJ/mol, $|h_{Ir2}|\approx \xi_1 1600$ kJ/mol and $|h_{Ir2}|\approx \xi_1 2600$ kJ/mol where $\xi_1=0.0214\xi$ is the proportional constant appeared in the case of $W$. Here we need to choose a smaller proportional constant $\xi_1$ instead of the proportional constant $\xi$ because if we choose $\xi$ then the total energy $h$ will be greater than the upper bound $h_{IImax}$ of conventional superconductivity and thus the condition of conventional superconductivity will be violated.

Furthermore, since $\xi_1$ gives the antiferromagnetic effect to the $s$ valence electrons of $W$ we have that $\xi_1$ is also for the antiferromagnetism  of other elements.

Comparing to the case of $Tc$ which is similar to this $Ir$ case we see that $Tc$ is with one $dd$-pair which is formed without the mechanism of changing spin while $Ir$ is with two such $dd$-pairs. Thus $\xi_1$ is from the antiferromagnetism of the additional $dd$-pair of $Tc$.

Such Cooper $dd$-pair is formed by first separating two $d$ electrons in the same atom with different spins such that these two $d$ electrons are coupled to the $d$ electrons of other atoms to form Cooper $dd$-pairs (The original two $d$ electrons in the same atom also form pair which is different from the Cooper pair of superconductivity). Thus with more pairings of $d$ electrons in the same atom the energy of each Cooper pair of superconductivity is lesser because some energies are used to separate the pairings of $d$ electrons in the same atom.

Thus for the atom $Tc$ with the constant $\xi$ for one pairing of $d$ electrons in the same atom we have that the atom $Ir$ is with the smaller constant $\xi_1$ for two pairings of $d$ electrons in the same atom. Thus the replacement of $\xi$ with $\xi_1=0.0214\xi$ corresponds to an additional pairing of $d$ electrons in the same atom. This can be generalized to that for three pairings of $d$ electrons in the same atom the proportional constant $\xi$ is replaced by the smaller proportional constant $\xi_2=(0.0214)^2\xi$, and so on.

Then from (\ref{Ir}) we can compute the critical temperature $T_c$ of $Ir$:
\begin{equation}
T_c \approx 0.13K \quad \mbox{(Computed value of $T_c$ of $Ir$)}
\label{Ir5}
\end{equation}
This agrees with  the experimental value of $T_c\approx 0.14K$ of $Ir$ under ambient pressure.

As explained in the above we identify the superconductivity of $Ir$ as the Type I superconductivity.

\subsection{\bf Transition element $\bf{Rh}$}

Let us then consider the element $Rh$. This element is with the outer shells $4d^85s^1$ of valence electrons. Since $Rh$ is a metal we have that $Rh$ is in the state $\kappa^2<\frac13 h^2$.

Similar to the element $T_i$ for $Rh$ in the state of superconductivity we have that $h^2$ is needed to be reduced such that $\kappa^2\geq \frac{1}{(2+\sqrt{3})^2}h^2$. In this case we need a mechanism for the decreasing of $h$ to reach the state of superconductivity.

Again let the Cooper pair be such a mechanism.
Since there are eight $d$ electrons we have that three of the eight $d$ electrons is with the spin different from the other five $d$ electrons which are with the same spin. Thus there are three pairs of $d$ electrons (with different spins) in the same $Rh$ atom and there are only two $d$ electrons with the same spin and is unpaired (in the sense of the same spin) in the same $Rh$ atom.

In this case, when the temperature is low, the spin sign of one of the two unpaired $d$ electrons of a $Rh$ atom is  changed to form a Cooper pair with  other  $Rh$ atoms. On the other hand,
when the temperature is low, there is no enough energy to separate the three paired $d$ electrons of a $Rh$ atom to form Cooper pair with   other  $Rh$ atoms. Thus there is only one Cooper $dd$-pair.

Then since   there is no Cooper $ss$-pair we have that the energy parameter
$h_{Rh1}$ can be proportional to the first ionization energy of $Rh$.
In this case the Cooper $dd$-pair is similar to Cooper $ss$-pair and thus the superconductivity is the Type I superconductivity.
Then the forming of Cooper pair reduces the total energy $2|h_{Rh1}|$ of two $d$ electrons for forming the Cooper $dd$-pair.

Then as the case of $T_i$, let $h_1$ be the energy reduced from this total energy by the mechanism of Cooper pairing and
let $h_2=2|h_{Rh1}|-h_1$ be the energy after the reduction of energy.

Then by the same reason as that for the element $Ti$,
from the formula $h^2=3\kappa^2$ we have that $h_1=h=h_2$. Thus we have the following formula of critical temperature $T_c$ of $Rh$:
\begin{equation}
k_BT_c=\kappa=D_{Rh}\Delta_{Rh}
\label{Rh}
\end{equation}
where $D_{Rh}=\frac{1}{2+\sqrt{3}}$ and
$\Delta_{Rh}=h=|h_{Rh1}|$ is the energy gap of superconductivity of $Rh$.

Then from the existing table of ionization energies we have that the
first ionization energy of $Rh$ is approximately equal to $719.7$ kJ/mol.

 Then
let $|h_{Rh1}|\approx \xi_2 719.7$ kJ/mol where the proportional constant is the constant $\xi_2=(0.0214)^2\xi$. This is because that $Rh$ is with two more pairings of $d$ electrons in the same $Rh$ atom (which are different from the Cooper $dd$-pairs), as explained in the case of $Ir$.

Then from (\ref{Rh}) we can compute the critical temperature $T_c$ of $Rh$:
\begin{equation}
T_c \approx 0.00030K \quad \mbox{(Computed value of $T_c$ of $Rh$)}
\label{Rh5}
\end{equation}
This agrees with  the experimental $T_c\approx 0.00033K$ of $Rh$ under ambient pressure.

As explained in the above we identify the superconductivity of $Rh$ as the Type I superconductivity.

\subsection{\bf Transition elements $\bf{Pd, Pt, Cu, Au, Ag}$}

Let us then consider the elements $Pd$, $Pt$, $Cu$, $Au$, and $Ag$. These elements are with the outer shells $4d^{10}5s^0$, $5d^{9}6s^1$, $3d^{10}4s^1$, $4d^{10}5s^1$, and $5d^{10}6s^1$ of valence electrons respectively.

Since $Pt$ is a metal with nine $d$ valence electrons we have that there are four pairing of $d$ valence electrons in the same $Pt$ atom. Thus the proportional constant is $\xi_3=(0.0214)^3\xi$. Comparing to $Rh$ we have that the critical temperature $T_c$ is of the order $10^{-6} K$. Thus $T_c\approx 0$. This agrees with the experimental value of $T_c\approx 0$.

Also since $Pd, Cu, Au, Ag$ are metals with ten $d$ valence electrons we have that there are five pairing of $d$ valence electrons in the same atom. Thus the proportional constant is $\xi_4=(0.0214)^4\xi$. Comparing to $Rh$ we have that the critical temperature $T_c$ is of the order $10^{-8} K$. Thus $T_c\approx 0$ for $Pd, Cu, Au, Ag$. This agrees with the experimental value of $T_c\approx 0$ for $Pd, Cu, Au, Ag$.

\subsection{\bf Transition element $\bf{Zn}$}

Let us then consider the element $Zn$. This element is with the outer shells $3d^{10}4s^2$ of valence electrons. Since $Zn$ is a metal we have that $Zn$ is in the state $\kappa^2<\frac13 h^2$.

Similar to the element $T_i$ for $Zn$ in the state of superconductivity we have that $h^2$ is needed to be reduced such that $\kappa^2\geq \frac{1}{(2+\sqrt{3})^2}h^2$. In this case we need a mechanism for the decreasing of $h$ to reach the state of superconductivity.

Again let the Cooper pair be such a mechanism.
Since there are ten $d$ electrons we have that five of the ten $d$ electrons are with the spin different from the other five $d$ electrons which are with the same spin. Also the two $s$ electrons are with the different spins.

When the temperature $T$ is low enough one of the two $s$ electrons with the wave factor $e^{-ih_{Zn}}$ is coupled with a $s$ electron with wave factor $e^{ih_{Zn}}$ of another atom $Zn$.
This forms a Cooper pair with the wave factor $e^{i(-h_{Zn}+h_{Zn})}$. This gives a reduction of energy, as explained in the above analysis. Because of the $ss$-pairing we have that the energy parameter $h_{Zn}$ is proportional to the first ionization energy of $Zn$, as the case of the $s$ electrons of the atom $Ti$. In this case the superconductivity will be the Type I superconductivity because of the similarity to the $s$ electrons.

Because there are ten $d$ electrons which completely fill the $d$ shell we have that
 when the temperature $T$ is low the spins of the $d$ electrons are remained to be unchanged.
In this case there are no more Cooper pairs except the only Cooper $ss$-pair.
Then the forming of this only Cooper pair reduces the total energy $2|h_{Zn}|$ of the two $s$ electrons for forming the Cooper pair.

Then as the case of $T_i$, let $h_1$ be the energy reduced from this total energy by the mechanism of Cooper pairing and
let $h_2=2|h_{Zn}|-h_1$ be the energy after the reduction of energy. Then by the same reason as that for the element $Ti$,
from the formula $h^2=3\kappa^2$ we have that $h_1=h=h_2$.

Thus we have the following formula of critical temperature $T_c$ of $Zn$:
\begin{equation}
k_BT_c=\kappa=D_{Zn}\Delta_{Zn}
\label{Zn}
\end{equation}
where $D_{Zn}=\frac{1}{2+\sqrt{3}}$ and
$\Delta_{Zn}=h=|h_{Zn}|$ is the energy gap of superconductivity of $Zn$.

Then from the existing table of ionization energies we have that the
first ionization energy of $Zn$ is approximately equal to $906.4$ kJ/mol.
Let $|h_{Zn}|\approx \xi 906.4$ kJ/mol where $\xi$ is the proportional constant.

 Then from (\ref{Zn}) we can compute the critical temperature $T_c$ of $Zn$:
\begin{equation}
T_c \approx 0.827K \quad \mbox{(Computed value of $T_c$ of $Zn$)}
\label{Zn5}
\end{equation}
This agrees with  the experimental value of $T_c\approx 0.875K$ of $Zn$ under ambient pressure.

As explained in the above we identify the superconductivity of $Zn$ as the Type I superconductivity.

\subsection{\bf Transition element $\bf{Cd}$}

Let us then consider the element $Cd$. This element is with the outer shells $4d^{10}5s^2$ of valence electrons. Since $Cd$ is a metal we have that $Cd$ is in the state $\kappa^2<\frac13 h^2$.

Similar to the element $T_i$ for $Cd$ in the state of superconductivity we have that $h^2$ is needed to be reduced such that $\kappa^2\geq \frac{1}{(2+\sqrt{3})^2}h^2$. In this case we need a mechanism for the decreasing of $h$ to reach the state of superconductivity.

Again let the Cooper pair be such a mechanism. Since the valence electrons of $Cd$ is completely similar to that of $Zn$, we have that $Cd$ has only one Cooper pair which is a $ss$-pair.

Then as similar to the case of $Zn$, we have the following formula of critical temperature $T_c$ of $Cd$:
\begin{equation}
k_BT_c=\kappa=D_{Cd}\Delta_{Cd}
\label{Cd}
\end{equation}
where $D_{Cd}=\frac{1}{2+\sqrt{3}}$ and
$\Delta_{Cd}=h=|h_{Cd}|$ is the energy gap of superconductivity of $Cd$.

Then from the existing table of ionization energies we have that the
first ionization energy of $Cd$ is approximately equal to $867.8$ kJ/mol.
Let $|h_{Cd}|\approx \xi 867.8$ kJ/mol where $\xi$ is the proportional constant.

Then from (\ref{Cd}) we can compute the critical temperature $T_c$ of $Cd$:
\begin{equation}
T_c \approx 0.791K \quad \mbox{(Computed value of $T_c$ of $Cd$)}
\label{Cd5}
\end{equation}
This roughly agrees with  the experimental value of $T_c\approx 0.56K$ of $Cd$ under ambient pressure.

As explained in the above we identify the superconductivity of $Cd$ as the Type I superconductivity.

\subsection{\bf Transition element $\bf{Hg}$}

Let us then consider the element $Hg$. This element is with the outer shells $5d^{10}6s^2$ of valence electrons. Since $Hg$ is a metal we have that $Hg$ is in the state $\kappa^2<\frac13 h^2$.

Similar to the element $T_i$ for $Hg$ in the state of superconductivity we have that $h^2$ is needed to be reduced such that $\kappa^2\geq \frac{1}{(2+\sqrt{3})^2}h^2$. In this case we need a mechanism for the decreasing of $h$ to reach the state of superconductivity.

Again let the Cooper pair be such a mechanism. Since the valence electrons of $Hg$ is completely similar to that of $Zn$, we have that $Hg$ has only one Cooper pair which is a $ss$-pair.

Then since the crystal structure of $Hg$ is of the rlombohedral type that in each unit cell of $Hg$ there is a cluster of two $Hg$ atoms. Let a $s$ electron of $Hg$ be with the wave factor $e^{ih_{Hg}}$.

 Then for the two clustering $Hg$ atoms, the two clustering $s$ electrons give the wave factor $e^{i2h_{Hg}}$. These two clustering $s$ electrons then are coupled with another two clustering $s$ electrons with the wave factor $e^{-i2h_{Hg}}$ to form a $ss$-pair. In this case the total energy of the $ss$-pair is from $2|h_{Hg}|$ increased to $4|h_{Hg}|$ (We shall generalize such a clustering property to synthesize room temperature superconductor). Then the forming of Cooper pair reduces the total energy $4|h_{Hg}|$ to $h=h_2=2|h_{Hg}|$.

Then as similar to the case of $Zn$ (or the case of $Ti$), we have the following formula of critical temperature $T_c$ of $Hg$:
\begin{equation}
k_BT_c=\kappa=D_{Hg}\Delta_{Hg}
\label{Hg}
\end{equation}
where $D_{Hg}=\frac{1}{\sqrt{3}}$ and
$\Delta_{Hg}=h=h_2=2|h_{Hg}|$ is the energy gap of superconductivity of $Hg$. Here we have $D_{Hg}=\frac{1}{\sqrt{3}}$ and $D_{Hg}$ is not equal to the usual value $\frac{1}{2+\sqrt{3}}$ for Type I superconductivity. This is because that $Hg$ under ambient pressure is with the highest critical temperature $T_c$ such that the critical temperature $T_c$ is decreasing when the pressure is increasing.

Then from the existing table of ionization energies we have that the
first ionization energy of $Hg$ is approximately equal to $1007.1$ kJ/mol.
Let $|h_{Hg}|\approx \xi 1007.1$ kJ/mol where $\xi$ is the proportional constant.

Then from (\ref{Hg}) we can compute the critical temperature $T_c$ of $Hg$ (or $\alpha-Hg$):
\begin{equation}
T_c \approx 3.96K \quad \mbox{(Computed value of $T_c$ of $Hg$)}
\label{Hg5}
\end{equation}
This agrees with  the experimental value of $T_c\approx 4.2K$ of $Hg$ (or $\alpha-Hg$) under ambient pressure.

As explained in the above we identify the superconductivity of $Hg$ as the Type I superconductivity.

On the other hand there is the singly ionized $Hg II$ state of $Hg$. For this state because of the tetragonal modification of the rhombohedral type crystal structure of $Hg$ the clustering effect is not as significant as the $\alpha-Hg$ case \cite{Schi,Sch2}. For this state let the $s$ electron be with the wave  factor $e^{ih_{Hg2}}$.

Then as the above case we have the following formula of critical temperature $T_c$ of $Hg II$:
\begin{equation}
k_BT_c=\kappa=D_{Hg}\Delta_{Hg2}
\label{Hg2}
\end{equation}
where $D_{Hg}=\frac{1}{\sqrt{3}}$ and
$\Delta_{Hg2}=h=h_2=|h_{Hg2}|$ is the energy gap of superconductivity of $HgII$.

Then since $Hg II$ is in the singly ionized state we have that the energy parameter $h_{Hg2}$ is proportional to the second ionization energy of $Hg$. From the existing table of ionization energy of $Hg$ we have that the
second ionization energy of $Hg$ is approximately equal to $1810$ kJ/mol.
Let $|h_{Hg2}|\approx \xi 1810$ kJ/mol where $\xi$ is the proportional constant.

 Then from (\ref{Hg2}) we can compute the critical temperature $T_c$ of $Hg II$ (or $\beta-Hg$) :
\begin{equation}
T_c \approx 3.56K \quad \mbox{(Computed value of $T_c$ of $\beta-Hg$)}
\label{Hg6}
\end{equation}
This agrees with  the experimental value of $T_c\approx 3.99K$ of $\beta-Hg$ under ambient pressure where $Hg II$ is to be with the crystal structure of $\beta-Hg$.

\subsection{\bf Nontransition element $\bf{Al}$}

Let us then consider the element $Al$. This element is with the outer shell $3s^{2}3p^1$ of valence electrons. Since $Al$ is a metal we have that $Al$ is in the state $\kappa^2<\frac13 h^2$.

Similar to the element $T_i$ for $Al$ in the state of superconductivity we have that $h^2$ is needed to be reduced such that $\kappa^2\geq \frac{1}{(2+\sqrt{3})^2}h^2$. In this case we need a mechanism for the decreasing of $h$ to reach the state of superconductivity.

Again let the Cooper pair be such a mechanism.
Since there are two $s$ electrons we have that the two $s$ electrons are with different spins.

When the temperature $T$ is low enough one of the two $s$ electrons with the wave factor $e^{-ih_{Al}}$ is coupled with a $s$ electron with wave factor $e^{ih_{Al}}$ of another atom $Al$.

This forms a Cooper $ss$-pair with the wave factor $e^{i(-h_{Al}+h_{Al})}$. This gives a reduction of energy, as explained in the above analysis.

Because of the $ss$-pairing it is possible that the energy parameter $h_{Al}$ is proportional to the first ionization energy of $Al$, as the $ss$-pair of the atom $Ti$. In this case the superconductivity will be the Type I superconductivity.
Then because there is only one $p$ electron left after the Cooper $ss$-pairing there are no more Cooper pairing.

Then the forming of this only Cooper $ss$-pair reduces the total energy $2|h_{Al}|$ of the two $s$ electrons for the forming of the Cooper pair.

Then as the case of $T_i$, let $h_1$ be the energy reduced from this total energy by the mechanism of Cooper pairing and
let $h_2=2|h_{Al}|-h_1$ be the energy after the reduction of energy.

Then by the same reason as that for the element $Ti$,
from the formula $h^2=3\kappa^2$ we have that $h_1=h=h_2$.
Thus we have the following formula of critical temperature $T_c$ of $Al$:
\begin{equation}
k_BT_c=\kappa=D_{Al}\Delta_{Al}
\label{Al}
\end{equation}
where $D_{Al}=\frac{1}{\sqrt{3}}$ (as the case of $Hg$) and
$\Delta_{Al}=h=|h_{Al}|$ is the energy gap of superconductivity of $Al$.

 Then from the existing table of ionization energies we have that the
first ionization energy of $Al$ is approximately equal to $577.4$ kJ/mol.

Let $|h_{Al}|\approx \xi 577.4$ kJ/mol where $\xi$ is the proportional constant.
Then from (\ref{Al}) we can compute the critical temperature $T_c$ of $Al$:
\begin{equation}
T_c \approx 1.135K \quad \mbox{(Computed value of $T_c$ of $Al$)}
\label{Al5}
\end{equation}
This agrees with  the experimental value of $T_c\approx 1.14K$ of $Al$ under ambient pressure.

As explained in the above we identify the superconductivity of $Al$ as the Type I superconductivity.

\subsection{\bf Nontransition element $\bf{Ga}$}

Let us then consider the element $Ga$. This element is with the outer shell $4s^{2}4p^1$ of valence electrons. Since $Ga$ is a metal we have that $Ga$ is in the state $\kappa^2<\frac13 h^2$.

Similar to the element $T_i$ for $Ga$ in the state of superconductivity we have that $h^2$ is needed to be reduced such that $\kappa^2\geq \frac{1}{(2+\sqrt{3})^2}h^2$. In this case we need a mechanism for the decreasing of $h$ to reach the state of superconductivity.

Again let the Cooper pair be such a mechanism. Since the configuration of valence electrons of $Ga$ is completely similar to that of $Al$, we have that $Ga$ has only one Cooper pair.  Let us consider the case that the Cooper pair is a $ss$-pair where the $s$ electrons are with the wave factors $e^{\pm ih_{In}}$ respectively.

Then the forming of this only Cooper pair reduces the total energy $2|h_{Ga}|$ of the two $s$ electrons for the forming of the Cooper pair.

Then as the case of $T_i$ (or as the case of $Al$), we have the following formula of critical temperature $T_c$ of $Ga$:
\begin{equation}
k_BT_c=\kappa=D_{Ga}\Delta_{Ga}
\label{Ga}
\end{equation}
where $D_{Ga}=\frac{1}{\sqrt{3}}$ (as the case of $Hg$) and
$\Delta_{Ga}=h=|h_{Ga}|$ is the energy gap of superconductivity of $Ga$.

Then from the existing table of ionization energies we have that the
first ionization energy of $Ga$ is approximately equal to $578.8$ kJ/mol.
Let $|h_{Ga}|\approx \xi 578.8$ kJ/mol where $\xi$ is the proportional constant.

Then from (\ref{Ga}) we can compute the critical temperature $T_c$ of $Ga$:
\begin{equation}
T_c \approx 1.134K \quad \mbox{(Computed value of $T_c$ of $Ga$)}
\label{Ga5}
\end{equation}
This roughly agrees with  the experimental value of $T_c\approx 1.09K$ of $Ga$ under ambient pressure.

As explained in the above we identify the superconductivity of $Ga$ as the Type I superconductivity.

\subsection{\bf Nontransition element $\bf{In}$}

Let us then consider the element $In$. This element is with the outer shell $5s^{2}5p^1$ of valence electrons. Since $In$ is a metal we have that $In$ is in the state $\kappa^2<\frac13 h^2$.

Similar to the element $T_i$ for $In$ in the state of superconductivity we have that $h^2$ is needed to be reduced such that $\kappa^2\geq \frac{1}{(2+\sqrt{3})^2}h^2$. In this case we need a mechanism for the decreasing of $h$ to reach the state of superconductivity.

Again let the Cooper pair be such a mechanism. Since the configuration of valence electrons of $In$ is completely similar to that of $Al$, we have that $In$ has only one Cooper pair. Let us consider the case that the Cooper pair is a $ss$-pair where the $s$ electrons are with the wave factors $e^{\pm ih_{In}}$ respectively.

Then the forming of this only Cooper pair reduces the total energy $2|h_{In}|$ of the two $s$ electrons for the forming of the Cooper pair.

Then as the case of $T_i$ (or as the case of $Al$), we have the following formula of critical temperature $T_c$ of $In$:
\begin{equation}
k_BT_c=\kappa=D_{In}\Delta_{In}
\label{In}
\end{equation}
where $D_{In}=\frac{1}{\sqrt{3}}$ (as the case of $Hg$) and
$\Delta_{In}=h=|h_{In}|$ is the energy gap of superconductivity of $In$.

 Let us then determine the energy parameter $|h_{In}|$.
 Because of  the $p$ electron in the outer shell, it is possible for the $s$ valence electrons to have the first ionization energy of $In$ or the second ionization energy of $In$. For $In$ let us consider the case that the $s$ valence electrons are with the second ionization energy of $In$.

 Then from the existing table of ionization energies we have that the
second ionization energy of $In$ is approximately equal to $1820.2$ kJ/mol.

Let $|h_{In}|\approx \xi 1820.2$ kJ/mol where $\xi$ is the proportional constant.
Then from (\ref{In}) we can compute the critical temperature $T_c$ of $In$:
\begin{equation}
T_c \approx 3.578K \quad \mbox{(Computed value of $T_c$ of $In$)}
\label{In5}
\end{equation}
This agrees with  the experimental value of $T_c\approx 3.404K$ of $In$ under ambient pressure.

As explained in the above we identify the superconductivity of $In$ as the Type I superconductivity.

\subsection{\bf Nontransition element $\bf{Tl}$}

Let us then consider the element $Tl$. This element is with the outer shell $6s^{2}6p^1$ of valence electrons. Since $Tl$ is a metal we have that $Tl$ is in the state $\kappa^2<\frac13 h^2$.

Similar to the element $T_i$ for $Tl$ in the state of superconductivity we have that $h^2$ is needed to be reduced such that $\kappa^2\geq \frac{1}{(2+\sqrt{3})^2}h^2$. In this case we need a mechanism for the decreasing of $h$ to reach the state of superconductivity.

Again let the Cooper pair be such a mechanism. Since the configuration of valence electrons of $Tl$ is completely similar to that of $Al$, we have that $Tl$ has only one Cooper pair.

Then because of the $p$ valence electron in the outer shell, it is possible that the Cooper pair is a $sp$-pair. Let us consider this $sp$-pair case where the $s$ electron is with the wave factor
$e^{ih_{Tl1}}$ and the $p$ electron is with the wave factor
$e^{-ih_{Tl2}}$.
Then the forming of this only Cooper pair reduces the total energy $|h_{Tl1}|+|h_{Tl2}|$ of the $s$ electron and the $p$ electron for forming the Cooper pair.

Then as the case of $Mo$, we have the following formula of critical temperature $T_c$ of $Tl$:
\begin{equation}
k_BT_c=\kappa=D_{Tl}\Delta_{Tl}
\label{Tl}
\end{equation}
where $D_{Tl}=\frac{1}{\sqrt{3}}$ (as the case of $Hg$) and
$\Delta_{Tl}=h=\frac12(|h_{Tl1}|+|h_{Tl2}|)$ is the energy gap of superconductivity of $Tl$.

Then from the existing table of ionization energies we have that the
first and second ionization energies of $Tl$ are approximately equal to $589.1$ kJ/mol and $1970.5$ kJ/mol respectively.

Let $|h_{Tl1}|\approx \xi 589.1$ kJ/mol and $|h_{Tl2}|\approx \xi 1970.5$ kJ/mol where $\xi$ is the proportional constant.

Then from (\ref{Tl}) we can compute the critical temperature $T_c$ of $Tl$:
\begin{equation}
T_c \approx 2.516K \quad \mbox{(Computed value of $T_c$ of $Tl$)}
\label{Tl5}
\end{equation}
This agrees with  the experimental value of $T_c\approx 2.39K$ of $Tl$ under ambient pressure.

As explained in the above we identify the superconductivity of $Tl$ as the Type I superconductivity.

\subsection{\bf Nontransition element $\bf{Sn}$}

Let us then consider the element $Sn$. This element is with the outer shell $5s^{2}5p^2$ of valence electrons. Since $Sn$ is a metal we have that $Sn$ is in the state $\kappa^2<\frac13 h^2$.

Similar to the element $T_i$ for $Sn$ in the state of superconductivity we have that $h^2$ is needed to be reduced such that $\kappa^2\geq \frac{1}{(2+\sqrt{3})^2}h^2$. In this case we need a mechanism for the decreasing of $h$ to reach the state of superconductivity.

Again let the Cooper pair be such a mechanism. There are two $p$ valence electrons with the same spin and two $s$ valence electrons with different spins.

When the temperature $T$ is low enough one of the two $s$ electrons with the wave factor $e^{-ih_{Sn1}}$ is coupled with a $s$ electron with wave factor $e^{ih_{Sn1}}$ of another atom $Sn$.

This forms a Cooper $ss$-pair with the wave factor $e^{i(-h_{Sn1}+h_{Sn1})}$. This gives a reduction of energy, as explained in the above analysis. Because of the $ss$-pairing the energy parameter $h_{Sn}$ is proportional to the first ionization energy of $Sn$, as the $ss$-pair of the atom $Ti$. In this case the superconductivity is the Type I superconductivity.

On the other hand when the temperature $T$ is low enough the spin of one of the $p$ valence electron is changed, as similar to the case of $Nb$. Then one of the two $p$ electrons with the wave factor $e^{-ih_{Sn2}}$ is coupled with a $p$ electron with wave factor $e^{ih_{Sn2}}$ of another atom $Sn$.
This forms a Cooper $pp$-pair with the wave factor $e^{i(-h_{Sn2}+h_{Sn2})}$. This gives a reduction of energy.

Because the $pp$-pairing is obtained by the changing of the spin of a $p$ electron, as the case of $Nb$ we have  that the energy parameter $h_{Sn2}$ is proportional to the third ionization energy of $Sn$. In this case the superconductivity is also the Type I superconductivity because this is a $pp$-pairing.

Then the forming of the two Cooper pairs reduces the total energy $2|h_{Sn1}|+2|h_{Sn2}|$ of the two $s$ electrons and the two $p$ electrons for forming the two Cooper pairs.

Then as the case of $Mo$, we have the following formula of critical temperature $T_c$ of $Sn$:
\begin{equation}
k_BT_c=\kappa=D_{Sn}\Delta_{Sn}
\label{Sn}
\end{equation}
where $D_{Sn}=\frac{1}{2+\sqrt{3}}$ is mainly for Type I superconductivity ($D_{Sn}$ may be greater than $\frac{1}{2+\sqrt{3}}$), and
$\Delta_{Sn}=h=|h_{Sn1}|+|h_{Sn2}|$ is the energy gap of superconductivity of $Sn$.

Then from the existing table of ionization energies we have that the
first and third ionization energies of $Sn$ are approximately equal to $708.6$ kJ/mol and $2943.0$ kJ/mol respectively.

Let $|h_{Sn1}|\approx \xi 708.6$ kJ/mol and $|h_{Sn2}|\approx \xi 2943.0$ kJ/mol where $\xi$ is the proportional constant.
Then from (\ref{Sn}) we can compute the critical temperature $T_c$ of $Sn$:
\begin{equation}
T_c \approx 3.332K \quad \mbox{(Computed value of $T_c$ of $Sn$)}
\label{Sn5}
\end{equation}
This agrees with  the experimental value of $T_c\approx 3.722K$ of $Sn$ under ambient pressure.

As explained in the above we identify the superconductivity of $Sn$ as the Type I superconductivity.

\subsection{\bf Nontransition element $\bf{Pb}$}.
Let us then consider the element $Pb$. This element is with the outer shell $6s^{2}6p^2$ of valence electrons. Since $Pb$ is a metal we have that $Pb$ is in the state $\kappa^2<\frac13 h^2$.

Similar to the element $T_i$ for $Pb$ in the state of superconductivity we have that $h^2$ is needed to be reduced such that $\kappa^2\geq \frac{1}{(2+\sqrt{3})^2}h^2$. In this case we need a mechanism for the decreasing of $h$ to reach the state of superconductivity.

Again let the Cooper pair be such a mechanism. There are two $p$ valence electrons with the same spin and two $s$ valence electrons with different spins.

When the temperature $T$ is low enough one of the two $s$ electrons with the wave factor $e^{-ih_{Pb1}}$ is coupled with a $s$ electron with wave factor $e^{ih_{Pb1}}$ of another atom $Pb$.

This forms a Cooper $ss$-pair with the wave factor $e^{i(-h_{Pb1}+h_{Pb1})}$. This gives a reduction of energy, as explained in the above analysis. Because of the $ss$-pairing the energy parameter $h_{Pb}$ is proportional to the first ionization energy of $Pb$, as the $ss$-pair of the atom $Ti$. In this case the superconductivity is the Type I superconductivity.

On the other hand when the temperature $T$ is low enough the spin of one of the $p$ valence electron is changed, as similar to the case of $Nb$. Then one of the two $p$ electrons with the wave factor $e^{-ih_{Pb2}}$ is coupled with a $p$ electron with wave factor $e^{ih_{Pb2}}$ of another atom $Pb$.
This forms a Cooper $pp$-pair with the wave factor $e^{i(-h_{Pb2}+h_{Pb2})}$. This gives a reduction of energy.

Because the $pp$-pairing is obtained by the changing of the spin of a $p$ electron, as the case of $Nb$ we have  that the energy parameter $h_{Pb2}$ is proportional to the third ionization energy of $Pb$. In this case the superconductivity is also the Type I superconductivity because this is a $pp$-pairing.

Then the forming of the two Cooper pairs reduces the total  energy $2|h_{Pb1}|+2|h_{Pb2}|$ of the two $s$ electrons and the two $p$ electrons for forming the two Cooper pairs.

Then as the case of $Mo$, we have the following formula of critical temperature $T_c$ of $Pb$:
\begin{equation}
k_BT_c=\kappa=D_{Pb}\Delta_{Pb}
\label{Pb}
\end{equation}
where $D_{Pb}=\frac{1}{\sqrt{3}}$ (As the case of $Hg$, here $D_{Pb}$ is not equal to the usual value $\frac{1}{2+\sqrt{3}}$ for Type I superconductivity. This is because that $Pb$ under ambient pressure is with the highest critical temperature $T_c$ such that the critical temperature $T_c$ is decreasing when the pressure is increasing), and
$\Delta_{Pb}=h=|h_{Pb1}|+|h_{Pb2}|$ is the energy gap of superconductivity of $Pb$.

Then from the existing table of ionization energies we have that the
first and third ionization energies of $Pb$ are approximately equal to $715.6$ kJ/mol and $3081.5$ kJ/mol respectively.

Let $|h_{Pb1}|\approx \xi 715.6$ kJ/mol and $|h_{Pb2}|\approx \xi 3081.5$ kJ/mol where $\xi$ is the proportional constant.
Then from (\ref{Pb}) we can compute the critical temperature $T_c$ of $Pb$:
\begin{equation}
T_c \approx 7.463K \quad \mbox{(Computed value of $T_c$ of $Pb$)}
\label{Pb5}
\end{equation}
This agrees with  the experimental value of $T_c\approx 7.193K$ of $Pb$ under ambient pressure.

As explained in the above we identify the superconductivity of $Pb$ as the Type I superconductivity.

\subsection{\bf Transition element $\bf{La}$}

Let us then consider the element $La$. Since $La$ is a metal we have that $La$ is in the state $\kappa^2<\frac13 h^2$. Similar to the element $T_i$ for $La$ in the state of superconductivity we have that $h^2$ is needed to be reduced such that $\kappa^2\geq \frac{1}{(2+\sqrt{3})^2}h^2$. In this case we need a mechanism for the decreasing of $h$ to reach the state of superconductivity.

Again let the Cooper pair be such a mechanism.
Let us first consider the case that the crystal structure of $La$ is of the hexagon type (h.c.p). In this case this element $La$ is with the outer shells $5d^16s^{2}$ of valence electrons.

There are two $s$ valence electrons with different spins and one $d$ valence electron. Because of the special hexagon structure there are $5$ $La$ atoms clustered in a crystal of $La$ atoms of the hexagon type (This is similar to that one $Nb$ atom is (clustered) in the center of a crystal of $Nb$ atoms of the cube type (b.c.c). The numbe $5=6-1$ where $6$ is the of the h.c.p. structure). Let us call these $5$ $La$ atoms as a unit element.

When the temperature $T$ is low enough one of the two $s$ electrons of a $La$ atom of a unit element with the wave factor $e^{-ih_{La1}}$ is coupled with a $s$ electron with wave factor $e^{ih_{La1}}$ of a $La$ atom of another unit element.

This forms a Cooper $ss$-pair with the wave factor $e^{i(-h_{La1}+h_{La1})}$. This gives a reduction of energy, as explained in the above analysis. 

Because of the $ss$-pairing the energy parameter $h_{La1}$ is proportional to the first ionization energy of $La$, as the $ss$-pair of the atom $Ti$. In this case the superconductivity is the Type I superconductivity.

Then because of the clustering of the five $La$ atoms we have that these five $La$ atoms can be regarded as a larger atom. In this case there are ten Cooper $ss$-pairs formed from this larger atom.
Then the forming of the ten Cooper $ss$-pairs reduces the total  energy $10|h_{La1}|$ of the ten $s$ electrons of the clustering of $La$ atoms for forming the ten Cooper $ss$-pairs.

Then as the case of $Ti$, we have the following formula of critical temperature $T_c$ of $La$:
\begin{equation}
k_BT_c=\kappa=D_{La}\Delta_{La}
\label{La}
\end{equation}
where $D_{La}=\frac{1}{\sqrt{3}}$ (As the case of $Hg$), and
$\Delta_{La}=h=5|h_{La1}|$ is the energy gap of superconductivity of $La$.

Then from the existing table of ionization energies we have that the
first ionization energy of $La$ is approximately equal to $538.1$ kJ/mol.

Let $|h_{La1}|\approx \xi 538.1$ kJ/mol where $\xi$ is the proportional constant.
 Then from (\ref{La}) we can compute the critical temperature $T_c$ of $La$ (h.c.p):
\begin{equation}
T_c \approx 5.29K \,\, \mbox{(Computed value of $T_c$ of $La$ (h.c.p))}
\label{La5}
\end{equation}
This agrees with  the experimental value of $T_c\approx 5.2K$ of $La$ (h.c.p) under ambient pressure.

Let us then consider the case that the crystal structure of $La$ is of the face centered cubic type (f.c.c). It is known that under high pressure the crystal structure of $La$ is changed from h.c.p. to f.c.c.. In this case there is a $s\to d$ transition of electrons \cite{Sch}. In this case this element $La$ is with the outer shells $5d^26s^{1}$ of valence electrons.
Thus there are two $d$ valence electrons with different spins and one $s$ valence electron. We can regard  the f.c.c structure as a structure such that there are $3$ $La$ atoms clustered in a crystal of $La$ atoms of the f.c.c. type (This is similar to that one $Nb$ atom is (clustered) in the center of a crystal of $Nb$ atoms of the cube type (b.c.c)) (The numbe $3=4-1$ where $4$ is the of the f.c.c. structure). Let us call these $3$ $La$ atoms as a unit element.

When the temperature $T$ is low enough one of the two $d$ electrons of a $La$ atom of a unit element with the wave factor $e^{-ih_{La2}}$ is coupled with a $d$ electron with wave factor $e^{ih_{La2}}$ of a $La$ atom of another unit element.
This forms a Cooper $dd$-pair with the wave factor $e^{i(-h_{La2}+h_{La2})}$. This gives a reduction of energy, as explained in the above analysis.
 Because the $dd$-pairing is obtained by the $f\to d$ transition of electrons we have that the energy parameter $h_{La2}$ is proportional to the second ionization energy of $La$. In this case the superconductivity is the Type I superconductivity.
Then because of the clustering of the three $La$ atoms we have that these three $La$ atoms can be regarded as a larger atom. In this case there are six Cooper $dd$-pairs formed from this larger atom.
Then the forming of the six Cooper $dd$-pairs reduces the total  energy $6|h_{La1}|$ of the six $d$ electrons of the clustering of $La$ atoms for forming the six Cooper $dd$-pairs.
Then as the case of $Ti$, we have the following formula of critical temperature $T_c$ of $La$:
\begin{equation}
k_BT_c=\kappa=D_{La}\Delta_{La}
\label{La2}
\end{equation}
where $D_{La}=\frac{1}{\sqrt{3}}$ (As the case of $Hg$), and
$\Delta_{La}=h=3|h_{La2}|$ is the energy gap of superconductivity of $La$ (f.c.c.).
Then from the existing table of ionization energies we have that the
second ionization energy of $La$ is approximately equal to $1067$ kJ/mol.
Let $|h_{La2}|\approx \xi 1067$ kJ/mol where $\xi$ is the proportional constant.
Then from (\ref{La2}) we can compute the critical temperature $T_c$ of $La$(f.c.c.):
\begin{equation}
T_c \approx 6.293K \,\, \mbox{(Computed value of $T_c$ of $La$ (f.c.c.))}
\label{La6}
\end{equation}
This agrees with  the experimental value of $T_c\approx 6K$ of $La$ (f.c.c.) under ambient pressure (This value $T_c\approx 6K$ is predicted from the experimental value of $T_c$ under high pressure for the state of $La$ (f.c.c.)).
As explained in the above we identify the superconductivity of $La$ as the Type I superconductivity.

\subsection{\bf Nontransition element $\bf{Be}$}

Let us then consider the element $Be$. This element is with the electron configuration $1s^22s^2$ and the outer shell $2s^2$ of valence electrons. Since $Be$ is a metal we have that $Be$ is in the state $\kappa^2<\frac13 h^2$.
 Similar to the element $T_i$ for $Be$ in the state of superconductivity we have that $h^2$ is needed to be reduced such that $\kappa^2\geq \frac{1}{(2+\sqrt{3})^2}h^2$. In this case we need a mechanism for the decreasing of $h$ to reach the state of superconductivity.
Again let the Cooper pair be such a mechanism.
There are two $s$ valence electrons with different spins  in the same $Be$ atom.
In this case when the temperature $T$ is very low,
the two $s$ valence electrons in the same $Be$ atom are separated such that these two $s$ electrons (with wave factors $e^{\pm ih_{Be}}$) are coupled to the separated $s$ electrons (with wave factors $e^{\pm ih_{Be}}$) of other $Be$ atoms to form Cooper $ss$-pairs with the wave factor $e^{i(-h_{Be}+h_{Be})}$. This gives a reduction of energy, as explained in the above analysis.

Because such Cooper $ss$-pairs are from the $s$ valence electrons we have that the energy parameter
$h_{Be}$ is proportional to the first ionization energy of $Be$. In this case the superconductivity is the Type I superconductivity.
Then the forming of Cooper pair reduces the total energy $2|h_{Be}|$ of two $s$ electrons for forming the Cooper $ss$-pair.
Then as the case of $T_i$, let $h_1$ be the energy reduced from this total energy by the mechanism of Cooper pairing and
let $h_2=2|h_{Be}|-h_1$ be the energy after the reduction of energy.
Then by the same reason as that for the element $Ti$,
from the formula $h^2=3\kappa^2$ we have that $h_1=h=h_2$. Thus we have the following formula of critical temperature $T_c$ of $Be$:
\begin{equation}
k_BT_c=\kappa=D_{Be}\Delta_{Be}
\label{Be}
\end{equation}
where $\frac{1}{2+\sqrt{3}}\leq D_{Be}\leq\frac{1}{\sqrt{3}}$ is for Type I superconductivity and
$\Delta_{Be}=h=|h_{Be}|$ is the energy gap of superconductivity of $Be$.

Let us then determine the energy parameter $|h_{Be}|$.
Since the next shell to the shell of $s$ valence electrons is a full shell of $s$ electrons which is not a shell of valence elecrons it is more difficult to separate the $s$ valence electrons to form Cooper pairs than the case that the next shell to the shell of $s$ valence electrons is the shell of $d$ valence electrons. Thus the proportional constant for defining $|h_{Be}|$ is one-step smaller than $\xi$ and is the proportional constant $\xi_1=0.0214\xi$ for antiferromagnetic effect, as explained in the case of $Ir$.

Then from the existing table of ionization energies we have that the
first ionization energy of $Be$ is approximately equal to $899.5$ kJ/mol.
Thus we have
$|h_{Be}|\approx \xi_1 899.5$ kJ/mol where the proportional constant is the constant $\xi_1$.
 Then from (\ref{Be}) we can compute the critical temperature $T_c$ of $Be$:
\begin{equation}
T_c \approx 0.0175\sim 0.0378 K \quad \mbox{(Computed value of $T_c$ of $Be$)}
\label{Be5}
\end{equation}
where we use the constant $D_{Be}$ such that $\frac{1}{2+\sqrt{3}}\leq D_{Be}\leq\frac{1}{\sqrt{3}}$.
This agrees with  the experimental value of $T_c\approx 0.026K$ of $Be$ under ambient pressure.

As explained in the above we identify the superconductivity of $Be$ as the Type I superconductivity.

\subsection{\bf Nontransition element $\bf{Li}$}

Let us then consider the element $Li$. This element is with the electron configuration $1s^22s^1$ and the outer shell $2s^1$ consisting of one valence electron. Since $Li$ is a metal we have that $Li$ is in the state $\kappa^2<\frac13 h^2$.

Similar to the element $T_i$ for $Li$ in the state of superconductivity we have that $h^2$ is needed to be reduced such that $\kappa^2\geq \frac{1}{(2+\sqrt{3})^2}h^2$. In this case we need a mechanism for the decreasing of $h$ to reach the state of superconductivity.

Again let the Cooper pair be such a mechanism.
There are two $s$ nonvalence electrons with different spins  in the inner shell of the same $Li$ atom.

In this case when the temperature $T$ is very low,
the two $s$ nonvalence electrons in the same $Li$ atom are separated such that these two $s$ electrons (with wave factors $e^{\pm ih_{Li}}$) are coupled to the separated $s$ nonvalence electrons (with wave factors $e^{\pm ih_{Li}}$) of other $Li$ atoms to form Cooper $ss$-pairs with the wave factor $e^{i(-h_{Li}+h_{Li})}$. This gives a reduction of energy, as explained in the above analysis.

Because such Cooper $ss$-pairs are from the $s$ electrons we have that the energy parameter
$h_{Li}$ is proportional to the first ionization energy of $Li$. In this case the superconductivity is the Type I superconductivity.
Then the forming of Cooper pair reduces the total energy $2|h_{Li}|$ of two $s$ electrons for forming the Cooper $ss$-pair.
Then as the case of $T_i$, let $h_1$ be the energy reduced from this total energy by the mechanism of Cooper pairing and
let $h_2=2|h_{Li}|-h_1$ be the energy after the reduction of energy.
Then by the same reason as that for the element $Ti$,
from the formula $h^2=3\kappa^2$ we have that $h_1=h=h_2$. Thus we have the following formula of critical temperature $T_c$ of $Li$:
\begin{equation}
k_BT_c=\kappa=D_{Li}\Delta_{Li}
\label{Li}
\end{equation}
where $ D_{Li}=\frac{1}{\sqrt{3}}$ is for Type I superconductivity and
$\Delta_{Li}=h=|h_{Li}|$ is the energy gap of superconductivity of $Li$.
Let us then determine the energy parameter $|h_{Li}|$.
Since the shell of $s$ electrons for forming Cooper pair is not a shell of valence elecrons it is more difficult to separate these $s$ electrons to form Cooper pairs than the case of $s$ valence electrons.
 Thus the proportional constant for defining $|h_{Li}|$ is one-step smaller than $\xi_1$ which is for the case of $Be$. Thus the proportional constant is the constant $\xi_2=0.0214\xi_1$, as explained in the case of $Ir$.
Then from the existing table of ionization energies we have that the
first ionization energy of $Li$ is approximately equal to $520.2$ kJ/mol.
Thus we let
$|h_{Li}|\approx \xi_2 520.2$ kJ/mol where the proportional constant is the constant $\xi_2$.
 Then from (\ref{Li}) we can compute the critical temperature $T_c$ of $Li$:
\begin{equation}
T_c \approx 0.00046 K \quad \mbox{(Computed value of $T_c$ of $Li$)}
\label{Li5}
\end{equation}
This agrees with  the experimental value of $T_c\approx 0.0004K$ of $Li$ under ambient pressure.
As explained in the above we identify the superconductivity of $Li$ as the Type I superconductivity.

Let us then investigate some well known $A15$-type superconductors \cite{Mat,Mat2,Mat3}.

\subsection{\bf $\bf{A15}$-type superconductor $\bf{Nb_{3}Ge_{1-x}Al_{x}}$}

Let us
consider the $A15$-type material $Nb_{3}Ge_{1-x}Al_{x}$ ($0\leq x\leq 1$).
While $Nb_{3}Ge$ and $Nb_{3}Al$ are conventional superconductors let us find a way to use the theory of unconventional superconductivity to investigate this $A15$-type material $Nb_{3}Ge_{1-x}Al_{x}$.
From the $A15$ crystal structure of $Nb_{3}Ge$ there are six $Nb$ atoms in a unit cell of $Nb_{3}Ge$.
Then in order to treat the $Nb$ atoms as a quasi-2D conducting layer we consider only the $Nb$ atoms at one direction. From the $A15$ crystal structure of $Nb_{3}Ge$ we have that along a direction there are two $Nb$ atoms in a unit cell of $Nb_{3}Ge$. These two $Nb$ atoms then give the quasi-2D conducting layer as the honeycomb of six $B$ atoms of $MgB_2$.
Then from the $A15$ crystal structure of $Nb_{3}Ge_{1-x}Al_{x}$
we have that along a direction there are two $Nb$ atoms connecting to one $Ge_{1-x}Al_{x}$ atom in the center of a unit cell of $Nb_{3}Ge_{1-x}Al_{x}$.
These two $Nb$ atoms and the centered $Ge(Al)$ are clustered as a larger atom.

Let us find out the doping mechanism. We consider the following function:
\begin{equation}
f(x)=1537(1-x)+1816.7x
 \label{Nb3GeAl83}
\end{equation}
 where $1537$ kJ/mol and $1816.7$ kJ/mole are approximately the second ionization energies of $Ge$ and $Al$ respectively.
 This function gives the relation that the increasing of $x$ giving the increasing of $h$. Then we set up the following relation:
\begin{equation}
f(x_0)=1537(1-x)+1816.7x =1380f_0
\label{Nb3GeAl8}
\end{equation}
for some $x_0$ ($0\leq x_0< 1$) where $1537$ kJ/mol is approximately the second
ionization energy of $Nb$, where we let $f_0=\frac{1537}{1380}$.  

When this relation holds  (or approximately holds) we have that the channel for the high-$T_c$ superconductivity is open. In this case the Cooper pairs of $d$ valence electrons of $Nb$, the Cooper pairs of $s$ valence electrons of $Al$
and the bifurcation region of high-$T_c$ superconductivity have been formed.
We notice that $x_0\approx 0$.

Thus $Nb_3Ge_{1-x}Al_x$ comes into the range  of high-$T_c$ superconductivity when $x_0\leq x\leq x_1$ for some $x_1$ such that $x_0< x_1\leq 1$.
When (\ref{Nb3GeAl8}) holds the $d$ valence electrons of $Nb$ and the $p$ valence electrons of $Ge$ are in the basic state of second ionization energy, and that other states are to be reached from this state.
Further, from the  $A15$ crystal structure of $Nb_{3}Ge_{1-x}Al_x$ we have that
the $d$ valence electrons of $Nb$ and the $p$ valence electrons of $Ge$ are unified to occupy a sequence of states such that the $d$ valence electrons of $Nb$ occupy the higher states while the $p$ valence electrons of $Ge$ occupy the lower states.
Thus, as $MgB_2$, the $d$ valence electrons of $Nb$ are in the state of third ionization energy, and the $p$ valence electrons of $Ge$ are in the state of second ionization energy. Denote the energy parameter for the $p$ valence electrons of $Ge$ by $|h_{Ge2}|$. Then the $s$ valence electrons of $Ge$ are in the state of first ionization energy. Denote the energy parameter of the $s$ valence electrons of $Ge$ by $|h_{Ge1}|$.
Thus,
the maximum value of the energy parameters $|h_{Nb}|$ for the $d$ valence electrons of $Nb$ is proportional to the
third ionization energy of $Nb$, the maximum value of the energy parameters $|h_{Ge1}|$ and $|h_{Ge2}|$ are proportional to the first and second ionization energies of $Ge$ respectively.
Then, as $MgB_2$, we have the following formula of $T_c$ of $Nb_{3}Ge_{1-x}Al_{x}$:
\begin{equation}
k_BT_c=\kappa=\frac{1}{\sqrt{3}}\Delta_{Nb3Ge_{1-x}Al_{x}}
\label{Nb3GeAl}
\end{equation}
where $\Delta_{Nb3Ge_{1-x}Al_{x}}=h=4|h_{Nb}|+|h_{Ge1}|+|h_{Ge2}|$ is the energy gap of the  superconductivity of $Nb_{3}Ge_{1-x}Al_{x}$ where the coeffecient $4=2\cdot2$ with one of the $2$ is the number of Cooper $dd$-pairs of $Nb$ and one of the $2$ is the number of $Nb$ atoms of the cluster of atoms.
We have $|h_{Nb}|\approx \xi 2416$ kJ/mol, $|h_{Ge1}|\approx \xi 762.1$ kJ/mol and
 $|h_{Ge2}|\approx \xi 1537$ kJ/mol
 where $2416$ kJ/mol is approximately the third ionization energy of $Nb$ and $762.1$ kJ/mol is approximately the first ionization energy of $Ge$.
 Then from (\ref{Nb3GeAl}) we can compute the highest critical temperature $T_c$ of $Nb_{3}Ge_{1-x}Al_{x}$:
\begin{equation}
T_c \approx  23.52 K\,\, \mbox{(Computed $T_c$ of $Nb_{3}Ge_{1-x}Al_{x}$)}
\label{Nb3GeAl5}
\end{equation}

Now since the experimental value $T_c$ of $Nb_{3}Ge$ is $T_c\approx 23.2 K$, we have that the highest critical temperature $T_c$ of $Nb_{3}Ge_{1-x}Al_{x}$ is approximately at $x_0=0$;
and that the computed
$T_c$ of $Nb_{3}Ge$ agrees with the experimental value $T_c\approx 23.2 K$ of $Nb_{3}Ge$.
Then, since the highest critical temperature $T_c$ of $Nb_{3}Ge_{1-x}Al_{x}$ is approximately at $x_0=0$;
we have that the optimal doping is approximately at $x_0=0$. Thus we deduce that for the $A15$ material $Nb_{3}Ge_{1-x}Al_{x}$, the $T_c$ is decreasing as $x$ increasing. 
This agrees with  the experimental values $T_c\approx 23.2 K$, $T_c\approx 20.5 K$ and $T_c\approx 18.5 K$ for $Nb_{3}Ge$ (for $x=0$) , $Nb_{3}Ge_{1-x_1}Al_{x_1}$ (for some $x_1$ such that $0<x_1<1$) and $Nb_{3}Al$ (for $x=1$) respectively.

\subsection{\bf $\bf{A15}$-type superconductor $\bf{Nb_{3}Ga}$}

Then let us
consider the well known $A15$-type material $Nb_{3}Ga$. From the $A15$ crystal structure of $Nb_{3}Ga$ we have that along a direction there are two $Nb$ atoms connecting to one $Ga$ atom in the center of a unit cell of $Nb_{3}Ga$.
These two $Nb$ atoms and the centered $Ge$ are clustered as a larger atom.
Then the doping mechanism of superconductivity of $Nb_{3}Ga$ is similar to the above doping mechanism  for $Nb_{3}Ge$.
We notice that the second ionization energies of $Nb$ and $Ga$ are approximately of $1380$ kJ/mol and $1979$ kJ/mol respectively. Since these two second ionization energies are quite close to each other that the channel connecting the states of ionization energies of $Nb$ and $Ga$ is open, as similar to $Nb_{3}Ge$.
In this case the Cooper pairs of $d$ valence electrons of $Nb$,  the Cooper pairs of $s$ valence electrons of $Ga$,
and the bifurcation region of superconductivity can be formed ($Ga$ has only one $p$ valence electron and thus has no Cooper pair of $p$ valence electrons).
In this case the $d$ valence electrons of $Nb$ and the $p$ valence electrons of $Ga$ are in the basic state of second ionization energy, and that other states are to be reached from this state.
Further, from the  $A15$ crystal structure of $Nb_{3}Ga$ we have that
the $d$ valence electrons of $Nb$ and the $p$ valence electrons of $Ga$ are unified to occupy a sequence of states such that the $d$ valence electrons of $Nb$ in the $Nb$ plane occupy the higher states while the $p$ valence electrons of $Ga$ occupy the lower states.
Thus, as $MgB_2$, the $d$ valence electrons of $Nb$ are in the state of third ionization energy, and the $p$ valence electrons of $Ga$ are in the state of second ionization energy. Denote the energy parameter by $|h_{Ga2}|$. Then the $s$ valence electrons of $Ga$ are in the state of first ionization energy. Denote the energy parameter by $|h_{Ga1}|$.
Thus,
the maximum value of the energy parameters $|h_{Nb}|$ for the $d$ valence electrons of $Nb$ is proportional to the
third ionization energy of $Nb$, the maximum value of the energy parameters $|h_{Ga1}|$ are proportional to the first and second ionization energies of $Ge$ respectively.
Then, as $MgB_2$, we have the following formula of $T_c$ of $Nb_{3}Ga$:
\begin{equation}
k_BT_c=\kappa=\frac{1}{\sqrt{3}}\Delta_{Nb3Ga}
\label{Nb3Ga}
\end{equation}
where $\Delta_{Nb3Ga}=h=4|h_{Nb}|+|h_{Ga1}|$ is the energy gap of the  superconductivity of $Nb_{3}Ga$.
We have $|h_{Nb}|\approx \xi 2416$ kJ/mol, $|h_{Ga1}|\approx \xi 578.8$ kJ/mol and
 where $578.8$ kJ/mol is approximately the first ionization energy of $Ga$.
Then from (\ref{Nb3Ga}) we can compute the highest critical temperature $T_c$ of $Nb_{3}Ga$:
\begin{equation}
T_c \approx  20.14 K\,\, \mbox{(Computed $T_c$ of $Nb_{3}Ga$)}
\label{Nb3Ga5}
\end{equation}
This agrees with the experimental value $T_c\approx 20.3 K$ of $Nb_{3}Ga$.

\section{Superfluid and color superconductivity in neutron star}\label{color}

In the phase of quark-gloun-plasma (QGP) of netutron and protons, it is argued that there are pairing mechansium of flavor quarks and color quarks  for the states of superfuid and superconductivity of  flavor quarks and color quarks  \cite{Ald}.  In the literature the possible state of superconductivity of quarks is usually called the color superconductivity of  quarks since the existing theory of the strong interaction of quarks is the QCD theory \cite{Ald}.

 In this section let us use the above renormalization group analysis to investigate the states of superfluid and  superconductivity of  quarks  in the core of neutron star where the  neutrons and protons are in the phase of  QGP.

In the phase of QGP of neutrons and protons, since quarks are still confined by the quantum knots of protons and neutrons,  each neutron or proton can be  regarded as a unit cell of the condensed matter QGP of  netutrons and protons, just as the crystal  unit cell of  metals.
Then since there are three (more than two)  quarks and two quarks with up and down spins in a neutron or a proton which are as a unit cell of a crystal,  it is possible to form Cooper pairs of quarks beween two neutrons or two protons.

In the above section on renormalization group equations we have shown that QCD has two 
 nontrivial fixed points.  Since QCD has two nontrivial fixed points, there are two special phase lines {\bf e}.  Thus the region for the forming of the quasi-2D phase of superconductivity (or superfluid) of one nontrivial fixed point  overlaps with the region of 3D phase of superconductivity  (or superfluid) of another nontrivial fixed point  in the phase diagram of QGP. Thus the quasi-2D phase of superconductivity always  exists stably  in the case of  the existence of the 3D phase of superconductivity  (or superfluid). This also means that there are layers of quasi-2D dimensional  superconducting protons (or quasi-2D dimensional  superfluid of neutrons).
 
Since the  quasi-2D phase of superconductivity  (or superfluid) can  stably exist, this gives that the color superconductivity  (or superfluid) can exsit with no  upper limit on the critical temperature $T_c$. 
This agrees with the recent analysis of data of  cooling rate of the neutron star of   Cassiopeia A (Cas A)  that the critical temperature $T_c$ is very high that $T_c \geq 2\times 10^9 K$ \cite{Pag}\cite{Sht}.

 Further, since protons are of about $10$ percent in the neutron star of Cas A, protons may  form  layers in the neutron star of Cas A.   Then the  quarks of neutrons and of the protons are approximately with the same ionization energies (or the enegy of the nearly free of the quarks), we have that a channel opening between neutron and proton is possible at some energy of the nearly free of quarks. Thus 3D superfluid  of quarks of neutrons near the layers of protons and quisi-2D superconductivity of quarks of the layers of protons both exist and this gives the 
quasi-2D high-$T_c$ superconductivity of quarks of protons. 
Then since the enegy of the nearly free of the quarks is very high that the critical temperature $T_c$ is very high by the formula of $T_c$ as that in the QED case.
In the forming of pairing of quarks of protons,  energy is not released in  the neutrino form and the  neutron star is not cooled by the energy  released in  the neutrino form.

Then more and more quarks of neutrons form the 3D superfluid through the channel opening while the quisi-2D superconductivity of quarks of the layers of protons  exists.
In the pairing of quarks of neutrons in  forming the 3D superfluid, energy is released in the neutrino form and the  neutron star is cooled.
 This is as the period of rapidly cooling rate of the  neutron star of Cas A.
This gives the superfluid and color superconductivity in the core of the neutron star of Cas A.


\begin{thebibliography}{99}

\bibitem{Wit}
E. Witten,
Quantum field theory and the Jones polynomial,
Commun. Math. Phys. {\bf 121}, 351 (1989)



\bibitem{Di}
 P.A.M. Dirac, 
{\em Directions in Physics},
(John Wiley, 1978).

\bibitem{Zub}
C. Itzykson and J-B. Zuber,
{\em Quantum Field Theory},
(McGraw-Hill Inc., 1980).


\bibitem{Kau}
L. Kauffman,
{\em Knots and Physics},
(World Scientific 1993).


\bibitem{Baez}
J. Baez and J. Muniain, 
{\em Gauge Fields, Knots and Gravity},
(World Scientic 1994).

\bibitem{Chari}
V. Chari and A. Pressley,
{\em A Guide to Quantum Groups},
(Cambridge University Press 1994).

\bibitem{Koh}
T. Kohno,
Ann. Inst. Fourier (Grenoble) {\bf 37}  139-160 (1987).

\bibitem{Dri}
V. G. Drinfel'd.
Leningrad Math. J. {\bf 1}  1419-57 (1990).

\bibitem{Ng1}
S.K. Ng,  arXiv hep-th/0209143v6.

\bibitem{Fra}
P. Di Francesco, P. Mathieu and D. Senechal,
{\em Conformal Field Theory},
(Springer-Verlag 1997).

\bibitem{Fuc}
J. Fuchs,
{\em Affine Lie Algebras and Quantum Groups},
(Cambridge University Press 1992).

\bibitem{Kni}
V.G. Knizhnik and A.B. Zamolodchikov,
Current algebra and Wess-Zumino model in two dimensions,
Nucl. Phys. B {\bf 247}, 83 (1984).

\bibitem{Mur}
K. Murasugi,
Knot Theory and Its Applications,
(Birkhauser Verlag, 1997).

\bibitem{Lic}
W.B.R. Lickorish,
An Introduction to Knot Theory,
(Springer, 1997).

\bibitem{Ng}
S.K. Ng,  arXiv hep-th/0210024v4.


\bibitem{Wilson}
K.G. Wilson and M.E. Fisher, {\it Phys. Rev. Lett.}, {\bf 28},
p.248 (1972).

\bibitem{Huang}
K. Huang, Statistical Mechanics, John Wiley and Sons, 1987.

\bibitem{Yu}
 L.Yu and B. Hao, Advances in Statistical Physics (in Chinese),
Science Press, 1981.

\bibitem{Creswick}
R.J. Creswick, H.A. Farach and  C.P. JR. Poole, Introduction to
Renormalization Group Methods in Physics, John Wiley and Sons,
1992.

\bibitem{Onsager}
L. Onsager, {\it Phys. Rev.}, {\bf 65}, p.117 (1944).



\bibitem{Bar}
J. Bardeen, L.N. Cooper and J.R. Schrieffer, {\it Phys. Rev.}, {\bf 108}, 1175 (1957).


\bibitem{Abr}
A.A. Abrikosov, {\it Sov. Phys. JETP}, {\bf 5}, 1174 (1957).

\bibitem{Bog}
N.N. Bogoliubov, {\it Sov. Phys. JETP}, {\bf 7}, 41 (1959).

\bibitem{Hub}
J. Hubbard, {\it Proc. Roy. Soc. London Ser. A}, {\bf 276}, 238 (1963).


\bibitem{Gor}
L.P. Gor'kov, {\it Zh. Eksp, Teor. Fiz.}, {\bf 36}, 1364 (1986).

\bibitem{Tin}
M.Tinkham, Introduction to superconductivity, McGraw-Hill, 1975.




\bibitem{Lev7}
L.P. Levy, {\it et al.} {\it Phys. Rev. Lett.}, {\bf 64}, 2074 (1990).

\bibitem{Cha6}
V. Chandrasekhar, {\it et al.} {\it Phys. Rev. Lett.}, {\bf 67}, 3578 (1991).

\bibitem{Mai6}
D. Mailly, {\it et al.} {\it Phys. Rev. Lett.}, {\bf 70}, 2020 (1993).

\bibitem{Bar6}
H. Bary-Soroker, {\it et al.} {\it Phys. Rev. Lett.}, {\bf 101}, 057001 (2008).

\bibitem{Blu6}
H. Bluhm, {\it et al.} {\it Phys. Rev. Lett.}, {\bf 102}, 36802 (2009).




\bibitem{Gei}
A.K. Geim, {\it et al.}, {\it Nature}, {\bf 396}, 144 (1998).

\bibitem{Tho}
D.J. Thompson, {\it et al.}, {\it Phys. Rev. Lett.}, {\bf 75}, 529 (1995).

\bibitem{Kos}
P. Kostic, {\it et al.}, {\it Phys. Rev. B}, {\bf 53}, 791 (1996).

\bibitem{Sve}
P. Svedlindh, {\it et al.}, {\it Physica C}, {\bf 162}, 1365 (1989).

\bibitem{Bra5}
W. Braunisch, {\it et al.}, {\it Phys. Rev. B}, {\bf 48}, 4030 (1993).

\bibitem{Kho}
D. Khomskii, {\it et al.}, {\it J. Low. Temp. Phys.}, {\bf 95}, 205 (1994).

\bibitem{Sig}
M. Sigrist and T.M. Rice, {\it Rev. Mod. Phys.}, {\bf 67}, 503 (1995).

\bibitem{Rie}
S. Riedling, {\it et al.}, {\it Phys. Rev. B}, {\bf 49}, 13283 (1994).

\bibitem{Li}
M.S. Li, {\it Physics Report}, {\bf 376}, 133 (2003).




\bibitem{Wis3}
W.D. Wise, {\it et al.} {\it Nature Physics}, {\bf 4}, 696 (2008).

\bibitem{Deg}
O. Degtyareva, {\it et al.}, {\it Phys. Rev. Lett.}, {\bf 99}, 155505 (2007).

\bibitem{Ber3}
E. Berg, {\it et al.}, {\it Phys. Rev. Lett.}, {\bf 100}, 027003 (2008).

\bibitem{Seo}
K. Seo, {\it et al.}, {\it Phys. Rev. B}, {\bf 78}, 094510 (2008).


\bibitem{Sax}
S.S. Saxena, {\it et al.} {\it Nature}, {\bf 406}, 587 (2000).

\bibitem{Aok}
D. Aoki, {\it et al.}, {\it Nature}, {\bf 413}, 613 (2001).

\bibitem{Huy}
N.T. Huy, {\it et al.}, {\it Phys. Rev. Lett.}, {\bf 99}, 67006 (2007).

\bibitem{Nav}
A.H. Navidomsky, cond-mat/0412247.

\bibitem{Mac}
K. Machida and T. Ohmi, cond-mat/0008345.







\bibitem{Lia}
W.Y. Liang, {\it et al.}, {\it J. Phys.: Condens. Matter}, {\bf 10}, 11365 (1998).

 \bibitem{Bra}
W. Braunich, {\it et al.}, {\it Phys. Rev. Lett.}, {\bf 68}, 1908 (1992).

\bibitem{Tim}
T. Timusk and B. Statt, {\it Rep. Prog. Phys.}, {\bf 62}, 61 (1999).




\bibitem{And}
P.W. Anderson, {\it Science}, {\bf 235}, 1196 (1987).

\bibitem{And5}
P.W. Anderson, The Theory of Superconductivity in the
High-$T_c$ Cuprate Superconductors,
Princeton, 1997.

\bibitem{Mul}
K.A. Muller, {\it 
J. Phys.: Condens. Phys.}, {\bf 19}, 25102 (2007).

\bibitem{Bus}
A. Bussman-Holder,  {\it at al.}, {\it Phys. Rev.}, {\bf B55}, 11751 (1997).





\bibitem{Ber}
I.B. Bersuker,  The Jahn-Teller effect, Cambridge University Press (2006).





\bibitem{Chu}
M.K. Wu, {\it et al.},
{\it Phys. Rev. Lett.} {\bf 58}, 908 (1987).




\bibitem{Nag2}
J. Nagamatsu, {\it et al.}, {\it Nature}, {\bf 410}, 63 (2001).





\bibitem{Tra}
J.M. Tranquada, {\it et al.}
{\it Nature}, {\bf  375}, 561 (1995).

\bibitem{Fuj}
M. Fujita, {\it et al.}
{\it Phys. Rev. B}, {\bf  70}, 104517 (2004).

\bibitem{Ich}
N. Ichikawa, {\it et al.}
{\it Phys. Rev. Lett.}, {\bf 85}, 1738 (2000).

\bibitem{Eme3}
V.J. Emery, S.A. Kivelson, and O. Zachar,
{\it Phys. Rev. B}, {\bf 56}, 6120 (1997).

\bibitem{Zaa5}
J. Zaanen, {\it et al.} {\it Phil. Mag, B}, {\bf 81}, 1485 (2001).

\bibitem{Mac3}
K. Machida, {\it Physica C}, {\bf 158}, 192 (1989).






\bibitem{Bia}
A. Bianconi, {\it et al.}
{\it Physica C: Superconductivity}, 341-348 (No. 1-4) , 1719-1722 (2000).

\bibitem{Kiv}
S.A. Kivelson, {\it et al.}
{\it Rev. Mod. Phys.}, {\bf 75}, 1201 (2003).






\bibitem{Cha}
X.K. Chan, {\it et al.}, {\it Phys. Rev. Lett.}, {\bf 87}, 1570021 (2001).

\bibitem{Liu}
A.Y. Liu, {\it et al.}, {\it Phys. Rev. Lett.}, {\bf 87}, 0870051 (2001).

\bibitem{Cho}
H.J. Choi, {\it et al.}, {\it Nature}, {\bf 418}, 758 (2002).

\bibitem{Iav}
M. Iavarone, {\it et al.}, {\it Phys. Rev. Lett.}, {\bf 89}, 1870021 (2002).

\bibitem{Bou}
F. Bouquet, {\it et al.}, {\it Phys. Rev. Lett.}, {\bf 89}, 257001-1 (2002).

\bibitem{Sol}
A.V. Sologubenko, {\it et al.}, {\it Phys. Rev. B}, {\bf 66}, 014504 (2002).

\bibitem{Fel}
I. Felner, cond-mat/0102508.





\bibitem{Zhao}
Y.G. Zhao, {\it et al.}, cond-mat/0103077.

\bibitem{Tan}
K. Tanaka, {\it et al.}, {\it, Science}, {\bf 314}, 1910 (2006).


\bibitem{Val}
T. Valla, {\it et al.}, {\it, Science}, {\bf 314}, 1914 (2006).

\bibitem{Ma}
J-H. Ma, {\it et al.}, {\it Phys. Rev. Lett.}, {\bf 101}, 207002 (2008).

\bibitem{Man3}
N. Mannella, {\it et al.}, {\it, Nature}, {\bf 438}, 474 (2005).






\bibitem{Mei2}
C. Meingast, {\it et al.}, {\it Phys. Rev. Lett.}, {\bf 86}, 1606 (2001).

\bibitem{Ido}
M. Ido, {\it et al.}, {\it J. Low. T. Phys.}, {\bf 117}, 329 (1999).





\bibitem{Tra2}
J.M. Tranquada, {\it et al.}, {\it Phys. Rev. B}, {\bf 78}, 174529 (2008).




\bibitem{Han2}
T. Hanaguri, {\it et al.},
{\it Nature}, {\bf 430}, 1001 (2004).

\bibitem{Mce}
K. McElroy, {\it et al.}, {\it Phys. Rev. Lett.}, {\bf 94}, 197005 (2005).

\bibitem{Che}
H.D. Chen, {\it et al.}, {\it Phys. Rev. Lett.}, {\bf 89}, 137004 (2002).

\bibitem{Pol}
A. Polkovnikov, {\it et al.}, {\it Phys. Rev. B.}, {\bf 65}, 220509 (2002).

\bibitem{Zhu2}
J.X. Zhu, {\it et al.}, {\it Phys. Rev. Lett.}, {\bf 89}, 067003 (2002).

\bibitem{Che2}
Y. Chen, {\it et al.}, {\it Phys. Rev. B.}, {\bf 66}, 104501 (2002).

\bibitem{Pol2}
D. Podolsky, {\it et al.}, {\it Phys. Rev. B.}, {\bf 67}, 094514 (2002).

\bibitem{And3}
Y. Ando, {\it et al.}, {\it Phys. Rev. Lett.}, {\bf 75}, 4662 (1995).

\bibitem{And6}
Y. Ando, {\it et al.}, {\it Phys. Rev. Lett.}, {\bf 77}, 2065 (1996).

\bibitem{Boe}
G.S. Boebinger, {\it et al.}, {\it Phys. Rev. Lett.}, {\bf 77}, 5417 (1996).

\bibitem{Hil}
R.W. Hill, {\it et al.}, {\it Nature}, {\bf 414}, 711 (2001).

\bibitem{Cuk}
T. Cuk, {\it et al.}, {\it Phys. Rev. Lett.}, {\bf 100}, 217003 (2008).





\bibitem{Tok}
Y. Tokura, H. Takagi, S. Uchida, Nature {\bf 337} (1989) 345.


\bibitem{Smi}
M.G. Smith, A. Manthiram, J. Zhou, J.B. Goodenough, and J.T. Markert,  Nature {\bf 351} (1991) 549.


\bibitem{Nai}
M. Naito and M. Hepp,   {\it Jpn. J. Appl. Phys.},  {\bf 39}, L485 (2000).

\bibitem{Nai2}
M. Naito,  S. Karimoto and  A. Tsukada,  {\it Supercond. Sci. Technol.},  {\bf 15}, 1663 (2002).

\bibitem{Saw}
A. Sawa,  M. Kawasaki, H. Takagi, Y. Tokura,  {\it Phys. Rev. B}, {\bf 66}, 014531 (2002).

\bibitem{Zue}
Y. Zuev,  {\it et al.}, {\it Phys. Sat. Sol.}, {\bf b 286}, No.2,  412 (2003).

\bibitem{Tan9}
Y. Tanaka, M. Karppinen and H. Yamauchi,
  {\it Supercond. Sci. Technol.},  {\bf 22}, 065004 (2009).

\bibitem{Koj}
K.M. Kojima,  {\it et al.}, {\it Physica B}, {\bf 374-375},  207 (2006).





\bibitem{Put}
S.N. Putilin, {\it et al.}, {\it Nature}, {\bf 362}, 226 (1993).


\bibitem{Sch5}
A. Schilling, {\it et al.}, {\it Nature}, {\bf 363}, 56 (1993).

\bibitem{Gao}
L. Gao, {\it et al.}, {\it Phys. Rev. B}, {\bf 50}, 4260 (1994).

\bibitem{She5}
Z.Z. Sheng and A.M. Hermann, {\it Nature}, {\bf 332}, 55 (1988).

\bibitem{Tor}
C.C. Torardi, {\it et al.}, {\it Phys. Rev. B}, {\bf 38}, 225 (1988).

\bibitem{Kan}
T. Kaneko, {\it et al.}, {\it Physica C}, {\bf 178}, 377 (1991).

\bibitem{Par}
S.S.P. Parkin, {\it et al.}, {\it Phys. Rev. Lett.}, {\bf 60}, 2539 (1988).

\bibitem{Par5}
S.S.P. Parkin, {\it et al.}, {\it Phys. Rev. Lett.}, {\bf 61}, 250 (1988).

\bibitem{Sch}
J.S. Schilling, {\it Physica C}, {\bf 460-462}, 182 (2007).

\bibitem{Schi}
J.E Schiriber and C.A. Swenson {\it Phys. Rev.}, {\bf 123}, 1115 (1961).

\bibitem{Sch2}
S. Deemyad and J.S. Schilling, {\it Phys. Rev. Lett.}, {\bf 91}, 167001 (2003).

\bibitem{Mat}
B. Matthias B., T.H. Geballe and E. Corenzwit, {\it Rev. Mod. Phys.}, {\bf 35}, 1 (1963).

\bibitem{Mat2}
J.M. Vandenberg and B. Matthias, {\it Science}, Vol. 198, No. 4313, 194 (1977).

\bibitem{Mat3}
B. Matthias, {\it et al.}, {\it Phys. Rev.}, {\bf 95}, 1435 (1954).






\bibitem{Chak}
D.J. Chakrabari and D.E. Laughlin, {\it Bull. Alloy Phase Diagram}, {\bf 2}, 319 (1981).

\bibitem{Sub}
P.R. Subramanian and D.E. Laughlin, {\it Bull. Alloy Phase Diagram}, {\bf 9}, 316 (1988).



\bibitem{Vuc}
J.H.N. van Vucht, {\it et al.}, 
{\it Philips Res. Rept.}, {\bf 25}, 133 (1970).



\bibitem{Vin}
C.A. Vincent and B. Scrosati, {\it Modern Batteries}, (John Wiley 1997).



\bibitem{Tan2}
A. Taniguchi, {\it et al.}, {\it J. of Power Sources}, {\bf 100}, 117 (2001). 


\bibitem{Lir}
R. Li, J. Wu, S. Zhou and J. Qian,  {\it J. of Rare Earths}, {\bf 24}, 341 (2006).






\bibitem{Er}
J.K. Erbacher, Nickel-Metal Hydride Technology Evaluation, Final Report Wright Laboratories,WL-TR-96-2069, Fairborn, OH, 1996.




\bibitem{Lin}
D. Linden and D. Magnusen, Portable Sealed Nickle-Metal Hydride Batteries, in  Chapter 29, Hanbook of Batteries, 
MaGraw-HIll (2004).

\bibitem{Ter}
Hirohito Teraoka, Development of Low Self-dischargeble Nickel-metal Hydride Battery, SANYO Energy Twicell Co., Ltd. 






\bibitem{Ald}
M.G. Alford,  {\it et al.}, {\it Rev. Mod. Phys.}, {\bf 80}, 1455 (2008).

\bibitem{Pag}
D. Page,    {\it et al.}, {\it Phys. Rev. Lett.}, {\bf 106}, 081101 (2011).

\bibitem{Sht}
P.S. Shternin,  {\it et al.}, {\it  Not. R. Astron. Soc. Mon. Not. Roy. Astron. Soc. }, {\bf 412}, L108-L112 (2011).



\end{thebibliography}
\end{document}